\documentclass[reqno]{amsart}
\usepackage{comment,amssymb,latexsym,upref,enumerate,fouridx}
\usepackage{finalT}

\numberwithin{equation}{section}
\input Tdef.def
\input flfb.def

\begin{document}

\title[A priori estimates for relativistic liquid bodies]{A priori estimates for relativistic liquid bodies}

\author[T.A. Oliynyk]{Todd A. Oliynyk}
\address{School of Mathematical Sciences\\
Monash University, VIC 3800\\
Australia}
\email{todd.oliynyk@sci.monash.edu}

%\subjclass[2000]{83C25}

\begin{abstract}
\noindent We demonstrate that a sufficiently smooth solution of the
relativistic Euler equations that represents a dynamical compact liquid body, when expressed in Lagrangian coordinates,
determines a solution to a system of non-linear wave equations with acoustic boundary conditions. Using this wave
formulation, we prove that these solutions satisfy energy estimates without loss of derivatives. Importantly,
our wave formulation does not require the liquid to be irrotational, and the
energy estimates do not rely on divergence and curl type estimates employed in previous works.
\end{abstract}

\maketitle 
\sect{intro}{introduction}

\subsect{relativistic}{Relativistic Euler equations} On a 4-dimensional spacetime, the relativistic Euler equations are given by\footnote{With the exception of
Section \ref{linIBVP}, we use lower case Greek indices, i.e. $\mu,\nu,\gamma$, to
label spacetime coordinate indices which run from $0$ to $3$.}
\leqn{eulint1}{
\nabla_\mu T^{\mu\nu} = 0
}
where
\eqn{eulint2}{
T^{\mu \nu} = (\rho + p)v^\mu v^\nu + p g^{\mu \nu}
}
is the stress energy tensor,
\eqn{metric}{
g = g_{\mu \nu} dx^\mu dx^\nu
}
is a Lorentzian metric of signature $(-,+,+,+)$, $\nabla_\mu$ is the Levi-Civita connection of
$g_{\mu\nu}$, $v^\mu$ is the fluid four-velocity normalized
by\footnote{Following standard conventions, we lower and raise spacetime coordinate indices, i.e. $\mu,\nu,\gamma$,  using
the metric $g_{\mu\nu}$ and inverse metric $g^{\mu\nu}$, respectively.}
\eqn{vnorm}{
g_{\mu\nu}v^\mu v^\nu = -1,
}
$\rho$ is the proper energy density of the fluid, and $p$ is the pressure.
Projecting \eqref{eulint1} into the subspaces parallel and orthogonal to $v^\mu$ yields
the following well known form of the relativistic Euler equations
\lalin{eulint3}{
v^\mu \nabla_\mu \rho + (\rho+p)\nabla_\mu v^\mu & =  0, \label{eulint3.1}\\
(\rho + p)v^\mu \nabla_\mu v^\nu + h^{\mu\nu}\nabla_\mu p & = 0, \label{eulint3.2}
}
where
\leqn{eulint4}{
h_{\mu\nu} = g_{\mu\nu} + v_\mu v_\nu
}
is the induced positive  definite metric on the subspace orthogonal to $v^\mu$.
In this article, we will be concerned
with fluids with a barotropic equation of state of the form
\eqn{eos}{
\rho = \rho(p)
}
where $\rho$ satisfies
\lgath{eosA}{
\rho \in C^{\infty}([0,\infty),[\rho_0,\rho_1]), \quad  \rho(0) = \rho_0, \label{eosA.1}
\\
 \intertext{and}
 \rho'(p) >0, \quad 0\leq p < \infty, \label{eosA.2}
}
for some constants $0<\rho_0 < \rho_1$.

For fluid bodies with compact support, the timelike matter-vacuum boundary is defined by the vanishing
of the pressure. Due to the above restrictions on the equation of state, the type of
fluids considered in this article are \emph{liquids}, which are characterized by having a jump discontinuity in the proper energy
density at the matter-vacuum boundary. The main aim of this article is to derive a priori estimates for sufficiently smooth
solutions of the relativistic Euler equations that represent dynamical compact liquid bodies. The precise form
of the a priori estimates can be found in Theorem \ref{apthm},
which represent the main result of this article. The key analytic difficulties
in establishing the a priori estimates are due to the presence of the matter-vacuum boundary, which is free.

Our approach to establishing a priori estimates begins with showing
that sufficiently smooth solutions, which represent dynamical liquid bodies, of the
relativistic Euler equations
satisfy, when expressed in Lagrangian coordinates, a system of non-linear wave equations
with acoustic boundary conditions; see \eqref{cbvpA.1}-\eqref{cbvpA.9}
for the complete initial boundary value problem (IBVP). Although it is well known
that the Euler equations can be reduced to a non-linear scalar wave equation in
the special case of irrotational fluids,
see \cite{Christodoulou:2007,VisserMolina:2010} for details, our formulation is different
and able to handle the general case where rotation is present. 

The starting point for the derivation of our wave equation is the relativistic continuity equation given by
\begin{equation} \label{relcont} 
\nabla_\mu ( \nc v^\mu) = 0
\end{equation}
where $\nc = (\rho+p)/\zeta$ is the number density with $\zeta$ the
\textit{specific enthalpy} of the fluid. Introducing the one form
\begin{equation*}
\theta^0_\mu = -\zeta v_\mu,
\end{equation*}
we write the continuity equation \eqref{relcont} as
\begin{equation*}
\nabla^\mu\biggl(\frac{1}{f}\theta_\mu^0\biggr) = 0,
\end{equation*}
where $f$ is a function of the specific enthalpy $\zeta$ defined below by \eqref{fdefB}. Applying the covariant derivative
to $\nabla_\nu$ to this expression, we obtain, after commuting the covariant derivative $\nabla_\nu$ through and anti-symmetrizing 
$\nabla_\nu \theta^0_\mu$, a wave equation for $\theta^0_\nu$ that includes a term at the highest order which involves 
the fluid vorticity $2$-form $\Lambda = -d\theta^0$. What makes our wave formulation useful is that we can control the vorticity term in the wave equation by using
the well-known propagation equation $\mathrm{L}_v \Lambda =0$ satisfied by the vorticity.      
It is worth noting that in all previous works on a priori estimates for
the relativistic fluid bodies, the propagation equation  $\mathrm{L}_v \Lambda =0$ is also used in an essential way.

The other essential element of our wave formulation are the boundary conditions. These boundary conditions are non-linear
and of acoustic type, which means that they are Neumann boundary conditions with source terms that involves two time derivatives, in Lagrangian coordinates, along
the boundary, the same order as appears in the wave equation. The orgin of this non-standard boundary condition is the geometric identity
\begin{equation*}
\frac{\nabla_v \theta^0_\nu}{|\nabla_v \theta^0|_g} = - n_\nu,
\end{equation*}
where $n_\nu$ is the unit, outward pointing co-normal to the matter-vacuum boundary,    
which holds along the matter-vacuum boundary for solutions of relativistic Euler equations that satisfy the liquid boundary conditions and the Taylor sign condition.

 The importance of our wave formulation is
that it is well suited for deriving energy estimates \textit{without derivative loss} in the presence of a free
matter-vacuum boundary. This
is due, in part, to the wave structure of the equations, and in part, to the nature of the acoustic
boundary conditions. Indeed, in Section \ref{linIBVP}, we establish a local existence and uniqueness theory
for linear systems of wave equations with acoustic boundary conditions; see, in particular, Theorem \ref{linlocthmC}.
This linear theory provides the key technical result needed to establish
our a priori estimates. We anticipate that the linear theory developed in Section \ref{linIBVP} may be
of independent interest as it can be applied more generally to other systems of wave equations having
acoustic boundary conditions.

\subsect{previous}{Comparison with existing results} The only existing work that contains
a mathematical analysis of compact, relativistic liquid bodies is \cite{Trakhinin:2009}. There,
the local existence and uniqueness of solutions is established using the theory of symmetric
hyperbolic systems. However,
the energy estimates derived from the symmetric hyperbolic theory
involve a derivative loss that is repaired using a Nash-Moser iteration scheme.
This leads to a rather high requirement on the regularity of the initial data in order
to close the scheme. In contrast, the energy estimates established in this work
do not involve derivative loss and require less regularity on the initial data.

In the non-relativistic limit, the relativistic Euler equations reduce to the
compressible Euler equations. In this setting, there are more results available
with the first local existence and uniqueness result for compressible liquids due to Lindblad \cite{Lindblad:2005b}.
The estimates derived in \cite{Lindblad:2005b} also required the use of a Nash-Moser scheme
due to derivative loss. We note that the work of Trakhinin \cite{Trakhinin:2009}
applies in the non-relativistic setting and provides an alternate approach.

More recently, a local existence and uniqueness theory without derivative loss
for non-relativistic, compact, compressible liquid bodies has been developed in \cite{Coutand_et_al:2013}. The regularity requirements,
as measured by the amount
of regularity assumed on the initial data for the map that defines the Lagrangian coordinates,
for the energy estimates derived in this article
are comparable with those of \cite{Coutand_et_al:2013} with both needing $4.5$ derivatives bounded in an $L^2$ sense.
One key technical difference between the approach taken here compared to that taken in \cite{Coutand_et_al:2013}
and also \cite{Lindblad:2005b} is the energy estimates derived in this article do not rely on
the divergence and curl estimates developed in \cite{Coutand_et_al:2013,Lindblad:2005b}.

It is worth noting that a priori estimates for solutions of the relativistic Euler equations that
represent dynamical gaseous bodies have been established in the recent works \cite{Hadzic_et_al:2015,Jang_et_al:2016}.
There are also a number of other related results for the non-relativistic Euler
equations that involve either incompressible, or
gaseous fluid bodies with compact support; for example, see \cite{CoutandShkoller:2007,Coutand_et_al:2010,CoutandShkoller:2012,Lindblad:2005a,LindbladNordgren:2009,ShatahZeng:2008}
and references cited therein.

\subsect{future}{Future directions} In work that is currently in preparation, we use
the techniques developed in this article to establish the local existence and uniqueness of solutions
to the relativistic Euler equations that represent dynamical, compact liquid bodies. The key technical
step in going from the a priori estimates presented here to existence and uniqueness is to
view the relations \eqref{eiC.1}-\eqref{eiC.3} as constraints, and show that
solutions of IBVP \eqref{cbvpA.1}-\eqref{cbvpA.9} that satisfy these constraints initially continue
to satisfy these constraints to the future; that is, we show that the constraints \textit{propagate}. Since
it is relatively straightforward to establish the existence and uniqueness of solutions
for \eqref{cbvpA.1}-\eqref{cbvpA.9} using the linear theory developed in
Section \ref{linIBVP}, the key difficulties are reduced to showing that
the IBVP satisfied by
the constraints \eqref{eiC.1}-\eqref{eiC.3} has unique solutions. This uniqueness problem is
solved by applying standard results from hyperbolic theory. This work will be presented in
a separate article. We also are able to show that once the local existence and uniqueness problem is
settled for a fixed metric, this theory can
be used in conjunction with the techniques developed
in \cite{AnderssonOliynyk:2014,Andersson_et_al:2014}, suitably adapted, to establish
the local existence of solutions to the Einstein-Euler equations. This work will be presented
in a separate article.

\subsect{overview}{Overview of this paper} In Section \ref{fwe}, we review the Frauendiener-Walton formulation of the relativistic Euler
equations, which is the starting point for the derivation of our wave formulation. Here, we set out
the class of solutions for
which we establish a priori estimates. The derivation of our wave formulation, in the Eulerian picture, is carried out in
Section \ref{frame} with the key equation being \eqref{waveE}.  Lagrangian coordinates adapted to our problem are introduced
in Section \ref{frameLag} and the wave equation \eqref{waveE} is transformed into these coordinates with
the resulting wave equation given by \eqref{waveG}.
A time differentiated version of the wave equations \eqref{waveG} is also derived in this section, see \eqref{waveK}.
In the follwing section, Section \ref{boundary}, we show that the liquid boundary condition, which corresponds to the
vanishing of the pressure at the matter-vacuum interface, implies acoustic type boundary conditions for time differentiated wave equation \eqref{waveK}. The complete IBVP
given by \eqref{cbvpA.1}-\eqref{cbvpA.9}
is then presented in Section \ref{IVBP}. A linear existence and uniqueness theory, which includes energy estimates, for this type of IBVP is developed in
Section \ref{linIBVP}. These energy estimates are then used in the final section, Section \ref{nlinIBVP}, to obtain the desired a priori estimates
that are presented in Theorem \ref{apthm}, which represent the main result of this article. Finally, a number of useful calculus inequalities, elliptic estimates,
and determinant formulas are listed in Appendices \ref{calculus}, \ref{elliptic} and \ref{LA}, respectively.

\sect{fwe}{The Frauendiener-Walton formulation of the relativistic Euler equations}

\subsect{FWeqns}{Fraudendiner-Walton equations}

The derivation of our wave formulation of the Euler equations is based on the Frauendiener-Walton formulation
of the Euler equations \cite{Frauendiener:2003,Walton:2005}, which we quickly review. In the Frauendiener-Walton formulation of the Euler equations, the proper energy density $\rho$ and
the normalized fluid 4-velocity $v^\mu$ are combined into a single timelike vector field $w^\mu$ defined by
\begin{equation}
w^\mu = \frac{1}{\zeta}v^\mu \label{rhov.2}
\end{equation}
where, as above,
\leqn{eul3}{
\zeta =\frac{1}{\sqrt{w^2}} \qquad (w^2 := -w_\nu w^\nu > 0)
}
is the specific enthalpy of the fluid. As shown in  \cite{Frauendiener:2003,Walton:2005}, $w^\mu$
satisfies the symmetric hyperbolic equation
\leqn{eul1}{
A_{\mu\nu}{}^\gamma \nabla_\gamma w^\nu = 0
}
where
\eqn{eul2}{
A_{\mu\nu}{}^\gamma = \left(3 + \frac{1}{s^2}\right) \frac{w_\mu w_\nu}{w^2} w^\gamma + \delta^\gamma_\nu w_\mu
+\delta^\gamma_\mu w_\nu + w^\gamma g_{\mu\nu},
}
and the square of the sound speed $s^2$ is a function of $\zeta$.
An explicit formula for $s^2=s^2(\zeta)$ can be calculated as follows: first, the pressure $p=p(\zeta)$ is
determined by solving the initial value problem
\lalin{prel}{
\frac{d p}{d \zeta} &= \frac{1}{\zeta}\left(\rho(p) + p\right), \label{prel.1}\\
p(\zeta_0) &= p_0, \label{prel.2}
}
for any particular choice of $\zeta_0>0$ and  $p_0\geq 0$. To be definite, we set
\leqn{zpfix}{
p_0 = 0 \AND \zeta_0 = 1.
}
Solving \eqref{prel.1}-\eqref{prel.2} then yields, by standard ODE theory and \eqref{eosA.1}, a solution
\leqn{preg}{
p \in C^\infty((0,\infty)).
}
With $p(\zeta)$ determined, $s^2$ is then given by the formula
\eqn{s2form}{
s^2= \frac{1}{\rho'\bigl(p(\zeta)\bigr)}.
}
From this formula, it is then clear from \eqref{preg} and the assumption \eqref{eosA.2} on the equation
of state that $s^2$ satisfies
\eqn{s2reg}{
s^2 \in C^\infty((0,\infty)) \AND s^2(\zeta) >0, \quad 0<\zeta <\infty.
}
It is also clear that the proper energy density of the fluid can be recovered using the formula
\begin{equation} \label{rhov.1}
\rho =\rho\bigl(p(\zeta)).
\end{equation}

\subsect{FW2Eul}{Equivalence of $w^\mu$ and $(\rho,v^\mu)$.}
The formulas \eqref{rhov.2}, \eqref{eul3} and \eqref{rhov.1} define the explicit transformation that takes solutions $(\rho,v^\nu)$
of the Euler equations \eqref{eulint3.1}-\eqref{eulint3.2} to solutions  $w^\mu$  of \eqref{eul1} and vice versa. This
shows that the system \eqref{eul1} is completely equivalent to the standard formulation of the Euler equations
given by \eqref{eulint3.1}-\eqref{eulint3.2}.

\subsect{envort}{Vorticity}
By definition (e.g. see (6.22) and (6.25) from \cite{Gourgoulhon:2006}),  the fluid \textit{vorticity} is
the $2$-form, denote here by $\Lambda_{\mu\nu}$, constructed from the specific enthalpy $\zeta$ and
four-velocity $v_\mu$ of the fluid as follows
\eqn{vorticityA}{
\Lambda_{\mu\nu} = \nabla_{\mu}(\zeta v_\nu) - \nabla_\nu (\zeta v_\mu).
}
Consequently, using \eqref{rhov.2} and \eqref{eul3}, we see that
the fluid vorticity takes the form
\leqn{vorticityB}{
\Lambda_{\mu\nu} = -\nabla_{\mu}\biggl(\frac{1}{g(w,w)}w_\nu\biggr) +  \nabla_{\nu}\biggl(\frac{1}{g(w,w)}w_\mu\biggr)
}
when expressed in terms of the Frauendiener-Walton vector field $w^\mu$.

\subsect{solassump}{Solution assumptions}

Before proceeding, we fix the class
of relativistic solutions that will be of interest to us. As noted in the introduction,
our assumptions on the equation of state imply that
this class of solutions represent dynamical, compact liquid bodies.

\bigskip

\noindent \textbf{Assumptions:} We assume the following:
\begin{enumerate}
\item[(A.1)] The $(x^\mu)$, $\mu=0,1,2,3$, are (global) Cartesian coordinates on $\Rbb^4$, and $g=g_{\mu\nu}d x^\nu d x^\mu$ is
a smooth Lorentzian metric on $\Rbb^4$.
\item[(A.2)] $U$ is an open, bounded set in $\Rbb^4$ that is diffeomorphic to a timelike cylinder with base
\eqn{Omega0def}{
\Omega_0 = \{0\}\times \Omega
}
where $\Omega$ is a bounded, open set in $\Rbb^3$ with smooth boundary.
\item[(A.3)] The vector field $w=w^\mu \del{\mu}$ is timelike, has components $w^\mu \in C^{k}(\overline{U})$ for $k\geq 8$, and
satisfies the Frauendiener-Walton-Euler equations \eqref{eul1} on $U$. Here, and in the following,
we employ the notation
\eqn{deldefa}{
\del{\mu} = \frac{\del{}\;}{\del{}x^{\mu}}
}
for partial derivatives with respect to the coordinates $(x^\mu)$.
\item[(A.4)] The set $U$ is invariant under the flow of $w$. Letting
$\Fc_\tau(x^\mu) = (\Fc^\mu_t(x^\mu))$ denote the flow map of the vector field $w$, so that
$\Fc^\mu_{\tau}(x^\nu)$ is the unique
solution to the initial value problem
\leqn{flowB}{
\frac{d\;}{d\tau} \Fc^\mu_\tau(x^\nu)  = w^\mu(x^\nu), \qquad \Fc^\mu_0(x^\nu) = x^\mu,
}
it is clear from this assumption that
there exists a $T>0$ such that
\eqn{UTdef}{
U_T :=  \bigcup_{0\leq \tau \leq T} \Fc_\tau\bigl(\{0\}\times \Omega \bigr) \subset U.
}
\item[(A.5)] The vector field $w$ is tangent
to the timelike boundary
\eqn{BTdef}{
B_T := \bigcup_{0\leq \tau \leq T}\Fc_\tau \bigl(\{0\}\times \del{}\Omega \bigr)
}
of $U_T$, i.e.
\leqn{wtanB}{
w|_{B_T} \in T B_T.
}
This is not actually an independent assumption since it is a consequence of
the previous assumption. However, we state it here separately to emphasize
the condition \eqref{wtanB}.
\item[(A.6)] There
exists constants $0<c_{s}^{-}< c_{s}^+<1$ such that
\eqn{sbnd}{
0< c_s^-\leq s^2 \leq c_{s}^+ < 1 \hspace{0.4cm} \text{in $U$}.
}
\item[(A.7)] The pressure vanishes on the timelike boundary $B_T$, i.e.
\leqn{passumpA}{
p|_{B_T} = 0,
}
and satisfies, for some constant $c_p>0$, the \emph{Taylor sign condition}
\leqn{Taylor}{
0 < c_p \leq -\nabla_n p \hspace{0.4cm} \text{in $B_T$}
}
where $n= n^\nu\del{\nu}$ is the outward pointing unit normal to $B_T$.
\end{enumerate}

\begin{rem} \label{assumpArem} $\;$

\begin{enumerate}
\item[(i)] Since the vector field $w^\mu$ is timelike by assumption (A.3), there exists a constant $0<c_{w}$ such that
\eqn{w2bnd}{
 0 < c_w \leq  w^2 \hspace{0.4cm} \text{in $U$.}
}
\item[(ii)] Assumption (A.5) implies via \eqref{rhov.2} that
\leqn{vbc}{
v|_{B_T} \in T B_T.
}
Together, \eqref{passumpA} and \eqref{vbc} make up the standard representation
of the free boundary conditions for a fluid body.
\item[(iii)] Using \eqref{prel.1}, we have
\eqn{TaylorA}{
-\nabla_n p = \frac{\zeta^2}{2}\bigl(\rho(p)+p\bigr)\nabla_n w^2.
}
Evaluating this on the boundary yields
\eqn{TaylorB}{
-\nabla_n p|_{B_T} = \frac{\rho_0}{2} \nabla_n w^2
}
by \eqref{eosA.1}, \eqref{zpfix} and \eqref{passumpA}. From this it is then clear that
\leqn{TaylorC}{
0 < \frac{2 c_p}{\rho_0} \leq \nabla_n w^2 \hspace{0.3cm} \text{in $B_T$}
}
is equivalent to Taylor sign condition \eqref{Taylor}.
\end{enumerate}
\end{rem}

%Away from the boundary $B_T$, it is well-known that the local existence and uniqueness of
%solutions to the relativistic Euler equations
%follows from standard hyperbolic PDE theory, e.g. see \cite{?????}. In light of this,
%we only need to understand the solution in the neighborhood of the boundary $B_T$.
%Since the Taylor sign condition \eqref{Taylor} implies that
%$-\nabla_n p$ is bounded away from zero in a neighborhood of the boundary $B_T$,
%it follows that we lose no generality if we assume that

\sect{frame}{An Eulerian wave formulation of the Euler equations}

\subsect{fform}{A frame formulation of the relativistic Euler equations}
The starting point for the derivation of our wave formulation is
the  frame formulation of the relativistic Euler equations from \cite{Oliynyk:PRD_2012}.
Following \cite{Oliynyk:PRD_2012}, we introduce a frame\footnote{With the exception
of Section \ref{linIBVP}, lower case Latin indices (i.e. $i,j,k$) will denote frame indices that run from $0$ to $3$.}
\eqn{eidef}{
e_i = e^\mu_i\del{\mu}
}
where
\begin{enumerate}[(i)]
\item
\leqn{e0def}{
e_0 := w
}
and $w=w^\mu\del{\mu}$ is a solution of the Frauendiener-Walton-Euler equations \eqref{eul1} that
satisfies the Assumptions (A.1)-(A.7) from Section \ref{solassump},
\item and the remaining frame fields\footnote{Upper case Latin indices
(i.e. $I, J, K$) will always run from $1$ to $3$ and denote the spatial frame indices.} $\{e_I\}_{I=1}^3$ are determined by solving the Lie transport equations
\leqn{eIdefa}{
[e_0,e_I] = 0.
}
\end{enumerate}

\begin{rem} \label{eIsol}
Writing \eqref{eIdefa} as
\leqn{eIsolA}{
e^\mu_0 \del{\mu} e^\nu_I - (\del{\nu} e^\mu_0) e^\nu_I = 0,
}
it is clear, since $U_T$ is invariant under the flow of $e_0=w$ by assumption, that we can solve \eqref{eIsolA} for given initial
data
\leqn{eIsolB}{
e^\mu_I|_{\Omega_0} = f^\mu_I \in C^{k}(\overline{\Omega})
}
using the method of characteristics to get a solution
\eqn{eIsolC}{
e^\mu_I \in C^{k-1}(\overline{U}_T).
}
\end{rem}

Next, we let
\eqn{theta}{
\theta^i = \theta^i_\mu dx^\mu \qquad \bigl( (\theta^i_\mu):=(e^\mu_j)^{-1}\bigr)
}
denote the co-frame, and we recall that the connection coefficients $\omega_i{}^k{}_j$ are defined
via the relation
\eqn{omegadef1}{
\nabla_{e_i} e_j = \omega_i{}^k{}_j e_k.
}
We define the associated connection 1-forms $\omega^k{}_j$ in standard fashion by
\eqn{omegadef2}{
\omega^k{}_j = \omega_i{}^k{}_j\theta^i,
}
and we set\footnote{In this article, we will follow standard convention and lower and raise the frame indices (i.e. $i,j,k$) with the
frame and inverse frame metrics $\fm_{ij}$ and $\fm^{ij}$, respectively.}
\eqn{omegadef3}{
\omega_{kj} = \fm_{kl}\omega^l{}_j = \omega_{ikj}\theta^i
}
where
\eqn{fmet}{
\fm_{ij} := g(e_i,e_j) = g_{\mu\nu}e_i^\mu e_j^\nu,
}
is the frame metric
and
\eqn{omegadef4}{
\omega_{ikj} := \fm_{kl}\omega_i{}^l{}_j = g(\nabla_{e_i} e_j, e_k).
}
We also let $\fm^{ij}$ denote the inverse frame metric, i.e.
\eqn{fmetinv}{
(\fm^{ij})^{-1} = (\fm_{ij}),
}
and note that
\eqn{w2}{
w^2 = -\fm_{00} \AND \zeta =  \left(\frac{1}{-\fm_{00}}\right)^{\frac{1}{2}}
}
by  \eqref{eul3} and \eqref{e0def}. Due to the choice \eqref{zpfix},
the boundary condition \eqref{passumpA} is equivalent to
\leqn{fbcB}{
\fm_{00}|_{B_T} = -1.
}

For use below, we let
\eqn{Freg}{
F\in C^\infty\bigl((0,\infty),(0,\infty)\bigr)
}
denote the unique solution to the initial value problem
\lalin{Fdef}{
F'(\zeta) &= -\frac{F(\zeta)}{\zeta s^2(\zeta)}, \label{Fdef.1} \\
F(1) &= 1 \label{Fdef.2}.
}
We also choose initial data \eqref{eIsolB} so that
\leqn{eiA}{
\fm_{0J}|_{\Omega_0} = 0 \AND \det\bigl(\fm_{IJ}|_{\Omega_0}\bigr)= F\bigl((-\fm_{00}|_{\Omega_0})^{-\frac{1}{2}}\bigr)^2.
}
As we show in
Section \ref{IV}, this is always possible.

\subsect{orthogonality}{Constraint propagation and vorticity}
Since
\begin{enumerate}[(i)]
\item $e_0^\mu =w^\mu$, by assumption, satisfies the Frauendiener-Walton-Euler equations
\eqref{eul1} on $\overline{U}_T$,
\item and the frame vectors $e_I^\mu$, by construction, satisfy the Lie transport
equations \eqref{eIdefa} on $\overline{U}_T$ for initial data \eqref{eIsolB} satisfying the conditions
\eqref{eiA},
\end{enumerate}
we are in a position to apply
Proposition IV.3 of \cite{Oliynyk:PRD_2012} to conclude that the following equations
are satisfied on $\overline{U}_T$:
\lalin{eiC}{
[e_0,e_J] & = 0, \label{eiC.1} \\
\fm_{0J} &= 0, \label{eiC.2} \\
F\bigl( (-\fm_{00})^{-1/2}\bigr)^2 &= \det(\fm_{IJ}), \label{eiC.3}
\intertext{and}
e_0(\sigma_l{}^j{}_k) & = 0\label{eiC.4}
}
where
\leqn{eiD}{
\sigma_l{}^j{}_k := \theta^j([e_l,e_k]) = \theta^j_\lambda
\bigl( e^\sigma_l \del{\sigma} e^\lambda_k  - e^\sigma_k \del{\sigma} e^\lambda_l \bigr).
}

First, we observe that \eqref{eiC.2}-\eqref{eiC.3} show that the constraints \eqref{eiA}
on the initial data \eqref{eIsolB} propagate, while \eqref{eiC.1} and \eqref{eiD} show that $\sigma_i{}^j{}_k$
is \textit{anti-symmetric} in the $i,k$ indices and satisfies
\leqn{eiE}{
\sigma_0{}^j{}_k = \sigma_{k}{}^j{}_0 = 0.
}
Next, we note that \eqref{eiC.2} implies that the co-vector
\eqn{eiEa}{
\frac{g_{\mu\nu}}{\fm_{00}}e^\nu_0
}
satisfies
\eqn{eiEa}{
 \biggl(\frac{g_{\mu\nu}}{\fm_{00}}e^\nu_0\biggr) e^\mu_i  = \delta^0_i.
}
But $\theta^0_\mu$ also satisfies $\theta^0_\mu e^\mu_i = \delta^0_i$, and so, the equality
\leqn{hC}{
\theta^0_\mu = \frac{g_{\mu\nu}}{\fm_{00}}e^\nu_0
}
must hold.
Since $e_0^\mu = w^\mu$ and $\fm_{00}=g(e_0,e_0)$, by definition, we can write the above expression as
\eqn{hCa}{
\theta^0_\mu = \frac{1}{g(w,w)}w_\mu,
}
which together with \eqref{vorticityB} shows that the fluid vorticity $\Lambda$ is given
by the formula
\leqn{vorticityC}{
\Lambda = -d\theta^0.
}
Using the Cartan structure equations
\leqn{CartanA.1}{
d\theta^i = -\omega^i{}_j \wedge \theta^j,
}
we can write \eqref{vorticityC} as
\eqn{vorticityCa}{
\Lambda = (\omega_{i}{}^0{}_j-\omega_{j}{}^0{}_i)\theta^i\wedge\theta^j,
}
which shows that
\leqn{vorticityDa}{
\Lambda = 0 \quad \Longleftrightarrow  \quad \omega_{i}{}^0{}_j-\omega_{j}{}^0{}_i=0.
}
Writing the Cartan Structure equations \eqref{CartanA.1} in the alternative form
\eqn{CartanB}{
\theta^k([e_i,e_j]) = \omega_i{}^k_j - \omega_j{}^k{}_i,
}
we see from \eqref{eiD} that $\sigma_i{}^k{}_j$ may be expressed as
\leqn{CartanC}{
\sigma_i{}^k{}_j = \omega_i{}^k{}_j - \omega_j{}^k{}_i.
}
From this and \eqref{vorticityDa}, we conclude that
\eqn{vorticityD}{
\Lambda = 0 \quad \Longleftrightarrow  \quad \sigma_{i}{}^0{}_j=0,
}
or in other words, the solution to the relativistic Euler equations determined by $w^\mu$ is \textit{irrotational} in $U_T$
precisely when $\sigma_{i}{}^0{}_j$ vanishes in $U_T$.

\begin{rem} \label{vortrem}
The above discussion shows that $\sigma_{i}{}^0{}_j$ provides an alternate characterization of the
fluid vorticity and that the fluid is irrotational precisely when $\sigma_{i}{}^0{}_j$ vanishes.
We also observe that \eqref{eiC.4} implies the statement:
\eqn{vorticityE}{
\sigma_{i}{}^0{}_j|_{\Omega_0}=0 \quad \Longrightarrow \quad \sigma_{i}{}^0{}_j=0\quad \text{in $U_T$,}
}
which is equivalent to the well known fact that if the fluid is initially irrotational, then it remains
irrotational under evolution.

It is also worthwhile noting that $\sigma_{i}{}^0{}_j|_{\Omega_0}=0$ is
only a condition on the initial data for $e_0^\mu=w^\mu$ as can be seen via the formula \eqref{hC}. In
the case that $\sigma_{i}{}^0{}_j|_{\Omega_0}=0$ is satisfied initially, it is possible to choose the initial
data \eqref{eIsolB} for the other frame fields $e_I^\mu$ so that $\sigma_{i}{}^K{}_j|_{\Omega_0}=0$
is also satisfied\footnote{This amounts to choosing the $e_I$ to be commuting vector fields on $\Omega_0$, i.e. $e_I|_{\Omega_0}\in T\Omega_0$ and $[e_I,e_J]|_{\Omega_0}=0$.} in addition to the other constraints \eqref{eiA}. By \eqref{eiC.4}, the constraints
 $\sigma_{i}{}^K{}_j|_{\Omega_0}=0$ then propagate to yield $\sigma_{i}{}^k{}_j=0$ in $U_T$.
\end{rem}

\subsect{wderive}{The wave equation}

Using the Cartan structure equations \eqref{CartanA.1} and
\leqn{CartanA.2}{
d\fm_{ij} = \omega_{ij}+\omega_{ji},
}
it is not difficult, see the proof of Proposition IV.1 in \cite{Oliynyk:PRD_2012} for details, to show that
\eqref{eiC.1}-\eqref{eiC.4} imply that connection coefficients $\omega_i{}^j{}_k$  satisfy
\lalin{connB}{
\frac{1}{s^2 \fm_{00}} \omega_{k00} - \fm^{IJ}\omega_{kIJ} & = 0 \notag %\label{connB.1}
\intertext{and}
\omega_{k0J}+\omega_{kJ0} & = 0 \notag %\label{connB.2}
}
in $\overline{U}_T$.
We also note that, due to \eqref{eiC.2}, the inverse frame metric $\fm^{ij}$ satisfies
\lalin{tiA}{
\fm^{0J} &= 0 \label{tiA.1}\\
\bigl(\fm^{IJ}\bigr) &= \bigl(\fm_{IJ}\bigr)^{-1} %\label{tiA.2}
\intertext{and}
\fm^{00} &= \frac{1}{\fm_{00}}, %\label{tiA.3}
}
which, in turn, allows us to express \eqref{eiC.3} as
\leqn{tiAa}{
F\bigl( \sqrt{-\fm^{00}}\bigr)^2 = \frac{1}{\det(\fm^{IJ})}.
}

Appealing to the definition of the Hodge star operator and the linear independence of the co-frame $\theta^i$, we
know that the
1-form $* \bigl(\theta^1 \wedge \theta^2 \wedge \theta^3\bigr) $ is
non-vanishing and orthogonal to the $\theta^I$, and  consequently,
\leqn{tiB}{
\theta^0 = \ff * \bigl(\theta^1 \wedge \theta^2 \wedge \theta^3\bigr)
}
for some non-vanishing function $\ff$ by \eqref{tiA.1}. To determine $\ff$, we apply the Hodge star operator to \eqref{tiB} to get\footnote{Recall
that $**\lambda = (-1)^{p(4-p)+1}\lambda$ for p-forms $\lambda$.}
\leqn{tiC}{
\frac{1}{\ff} * \theta^0 = \theta^1 \wedge \theta^2 \wedge \theta^3.
}
Wedging this with $\theta^0$ gives
\eqn{tiD}{
\frac{1}{\ff}\theta^0 \wedge * \theta^0 = \theta^0\wedge \theta^1 \wedge \theta^2 \wedge \theta^3,
}
which we can write as\footnote{We are using the well known identity
$ \lambda \wedge *\beta = g(\lambda,\beta)\mu$, which holds for all one forms $\alpha$ and $\beta$
with $\mu$ the volume form. In terms of the co-frame $\theta^i$, $\mu$ is given by
$\mu = (-\det(\fm^{ij}))^{-1/2} \theta^0\wedge \theta^1 \wedge \theta^2 \wedge \theta^3$.}
\eqn{tiE}{
\frac{1}{\ff} \frac{\fm^{00}}{\sqrt{-\det(\fm^{ij})}} \theta^0\wedge \theta^1 \wedge \theta^2 \wedge \theta^3= \theta^0\wedge \theta^1 \wedge \theta^2 \wedge \theta^3.
}
From this, we conclude that
\leqn{fdefA}{
\ff = -\frac{-\fm^{00}}{\sqrt{-\det(\fm^{ij})}},
}
which, with the help of \eqref{tiA.1} and \eqref{tiAa}, implies
\leqn{fdefB}{
\ff=\ff(\zeta) := -\zeta F(\zeta)
}
where $\zeta = \sqrt{-\fm^{00}}$.

Applying the operator $*d$ to \eqref{tiC} yields
\eqn{tiF}{
*d\left(\frac{1}{\ff}*\theta^0\right) = *\bigl(d\theta^1 \wedge \theta^2 \wedge \theta^3
- \theta^1 \wedge d\theta^2 \wedge \theta^3
+\theta^1 \wedge \theta^2 \wedge d\theta^3\bigr).
}
Noting that the right hand side of this vanishes, since
\leqn{CartanD}{
d\theta^i = -\Half \sigma_{M}{}^i{}_L \theta^M\wedge \theta^L % \quad \Longleftrightarrow
%\quad \nabla_{\nu}\theta^i_\mu  - \nabla_{\mu}\theta^i_\nu = \sigma_{M}{}^i{}_L\theta^M_\mu\theta^L_\nu,
}
by \eqref{eiE} and \eqref{CartanA.1}, and
\eqn{tiH}{
\theta^I\wedge \theta^J \wedge \theta^K \wedge \theta^L =0
}
for any choice of $I,J,K,L \in \{1,2,3\}$, we see that\footnote{This is just the relativistic continuity equation \eqref{relcont}.}
\eqn{tiG}{
*d\left(\frac{1}{\ff}*\theta^0\right) = 0,
}
or equivalently, in terms of components,
\eqn{tiL}{
\nabla^\mu\left(\frac{1}{\ff} \theta^0_\mu\right) = 0.
}
Applying the covariant derivative $\nabla_\nu$ to this expression, yields, after commuting the covariant
derivatives,\footnote{Here, we use $\nabla_\mu \nabla_\nu \lambda_\gamma - \nabla_\nu \nabla_\mu \lambda_\gamma
= R_{\mu\nu\gamma}{}^\sigma \lambda_\sigma$ and $R_{\mu\gamma}=R_{\mu\nu\gamma}{}^\nu$.} the relation
\leqn{tiM}{
\nabla^\mu\left( \nabla_\nu\left(\frac{1}{\ff} \theta^0_\mu\right)\right) = \frac{1}{\ff}R_{\nu}{}^{\lambda}\theta^0_\lambda.
}

Continuing on, we compute
\eqn{waveA}{
\nabla_\nu \zeta = -\frac{g(\theta^0,\nabla_\nu \theta^0)}{\sqrt{-\fm^{00}}} \qquad (\zeta = \sqrt{-\fm^{00}}),
}
and using \eqref{fdefB} and \eqref{Fdef.1},
\eqn{waveB}{
\frac{d\,}{d\zeta} \ln(-\ff(\zeta)) = \frac{1}{\zeta} + \frac{1}{F(\zeta)}\frac{d F}{d\zeta} = \frac{1}{\zeta} - \frac{1}{\zeta s^2(\zeta)}.
}
From these two expressions, we find that
\leqn{waveC.1}{
\nabla_{\nu} \left(\frac{1}{\ff}\theta^0_\mu\right) = \frac{1}{\ff}\left[\nabla_\nu \theta^0_\mu
- \frac{d\,}{d\zeta} \ln(-\ff(\zeta))\nabla_\nu \zeta \theta^0_\mu  \right]
= \frac{1}{\ff}\left[\nabla_\nu \theta^0_\mu
- \left(1-\frac{1}{s^2}\right)\frac{g(\theta^0,\nabla_\nu \theta^0)}{\fm^{00}} \theta^0_\mu  \right].
}
Using \eqref{rhov.2} and \eqref{e0def}, it is clear from the definition \eqref{eulint4} that
 the projection tensor $h_{\mu\nu}$ can be written as
\eqn{hA}{
h_{\mu\nu} = g_{\mu\nu} - \frac{e_{0\mu} e_{0\nu}}{\fm_{00}},
}
or equivalently
\leqn{hD}{
h_{\mu\nu} = g_{\mu\nu} - \frac{1}{\fm^{00}}\theta^0_\mu \theta^0_\nu = g_{\mu\nu} - \fm_{00}\theta^0_\mu \theta^0_\nu
}
using \eqref{hC}.

Setting $i=0$ in \eqref{CartanD} and expressing the result in terms of covariant derivatives,\footnote{Note that $d\theta^i = -\Half \sigma_{M}{}^i{}_L \theta^M\wedge \theta^L  \quad \Longleftrightarrow
\quad \nabla_\nu\theta^i_\mu  - \nabla_\mu\theta^i_\nu = \sigma_{M}{}^i{}_L\theta^M_\mu\theta^L_\nu$.} we get that
\leqn{CartanE}{
\nabla_\nu \theta^0_\mu = \nabla_\mu \theta^0_\nu + \sigma_M{}^0{}_L\theta^M_\mu \theta^L_\nu.
}
Using \eqref{hD} and \eqref{CartanE}, we can write \eqref{waveC.1} as
\eqn{waveD}{
\nabla_{\nu} \left(\frac{1}{\ff}\theta^0_\mu\right) = \frac{1}{\ff}\biggl(h^\alpha_\mu +
\frac{1}{s^2}\frac{\theta^{0\alpha}\theta^0_\mu}{g(\theta^0,\theta^0)}\biggr)\bigl[\nabla_\alpha \theta^0_\nu
+ \sigma_M{}^0{}_L \theta^M_\alpha \theta^L_\nu
\bigr].
}
Substituting this into \eqref{tiM} yields the wave equation
\leqn{waveE}{
\nabla_{\alpha}\left(\frac{1}{\ff}a^{\alpha \beta}\bigl[\nabla_\beta \theta^0_\nu
+ \sigma_M{}^0{}_L \theta^M_\beta \theta^L_\nu
\bigr]\right)= \frac{1}{\ff}R_{\nu}{}^{\lambda}\theta^0_\lambda,
}
where
\eqn{adefA}{
a^{\alpha\beta} = h^{\alpha\beta} +
\frac{1}{s^2}\frac{\theta^{0\alpha}\theta^{0\beta}}{\fm^{00}},
}
which defines our \emph{Eulerian
wave formulation of the Euler equations}.
\begin{rem} \label{Eulwaverem}
$\;$
\begin{itemize}
\item[(i)] From \eqref{hC}, it is clear that we can write $a^{\alpha\beta}$ as
\leqn{adefB}{
a^{\alpha\beta} = h^{\alpha\beta} +
\frac{1}{s^2}\frac{e^{\alpha}_0 e^{\beta}_0}{\fm_{00}} = g^{\alpha\beta} -\left(1-\frac{1}{s^2}\right)\frac{e^{\alpha}_0 e^{\beta}_0}{\fm_{00}},
}
or equivalently, using \eqref{rhov.2} and \eqref{e0def},
\eqn{adefC}{
a^{\alpha\beta} = h^{\alpha\beta} - \frac{1}{s^2}v^\mu v^\nu,
}
which we recognize as the inverse of the \textit{acoustic metric}
\eqn{adefD}{
a_{\alpha\beta} = h_{\alpha\beta} - s^2 v_\mu v_\nu.
}
\item [(ii)]  Although
$\theta^0_\mu$ is
equivalent
by \eqref{hC} to the vector field $e_0=w$, which we know from the discussion in
Section \ref{fwe} is completely equivalent to the usual parameterization of a barotropic
 perfect fluid in terms of $\rho$ and $v^\nu$, the wave
 equation \eqref{waveE} does not, in general, represent a complete evolution equation for
  the fluid. This is due to the presence of the co-frame fields $\theta^M_\beta$, which
  satisfy their own evolution equations that can be derived from \eqref{eiC.1}. The one exception
  to this is when the fluid is irrotational, which was shown in Section \ref{orthogonality} to be
  equivalent to the condition $\sigma_i{}^0{}_{\!j}=0$. In this case, it is clear that \eqref{waveE} reduces
  to a wave equation involving only $\theta^0_\nu$, and hence provides a complete evolution
  equation for the fluid.
\end{itemize}
\end{rem}
Although \eqref{waveE} does not, in general, provide a complete evolution
equation for the perfect fluid, we show in the next section that this
defect can be remedied by transforming to Lagrangian coordinates. The introduction of
the Lagrangian coordinates is also essential to fix the free boundary and allow us to work
on a fixed domain. To prepare for the change to Lagrangian coordinates, we
express the covariant derivatives in \eqref{waveE} in terms of the partial derivatives and Christoffel
symbols to get
\leqn{waveF}{
\del{\alpha} Y^\alpha_\nu = -\Gamma_{\alpha\gamma}^\alpha Y^\gamma_\nu + \Gamma^\gamma_{\alpha\nu}Y^\alpha_\gamma + \frac{1}{\ff}R_{\nu}{}^{\lambda}\theta^0_\lambda
}
where
\leqn{Ydef}{
Y^\alpha_\nu = \frac{1}{\ff}a^{\alpha \beta}\Bigl[\del{\beta} \theta^0_\nu - \Gamma_{\beta\nu}^\gamma \theta^0_\gamma
+ \sigma_M{}^0{}_L \theta^M_\beta \theta^L_\nu\Bigr].
}

\sect{frameLag}{The Lagrangian wave formulation of the Euler equations}

\subsect{lag}{Lagrangian coordinates} We introduce Lagrangian coordinates $(\xb^\mu)$ adapted to the vector field
$e_0 = w$ via the formula
\leqn{phidefA}{
x^\mu = \phi^\mu(\xb^\lambda) := \Fc_{\xb^0}^\mu(0,\xb^\Lambda)
}
where $\Fc_\tau$ is the flow of $w$ defined by \eqref{flowB}. From the regularity assumption
(A.3) for the components of the vector field $w^\mu$ and the standard properties of
flows, it is not difficult to verify that $\phi$
defines a $C^k$ diffeomorphism
\eqn{phidiff}{
\phi \: :\: \Omega_T := [0,T] \times \Omega \longrightarrow U_T \: :\: (\xb^\lambda) \longmapsto (\phi^\mu(\xb^\lambda)),
}
and  that
\leqn{phiOmega}{
\phi|_{\Omega_0} = \mathrm{id}_{\Omega_0}.
}
Moreover, since $\phi$ is generated by the flow of $e_0$, the pullback
\eqn{wbdefAa}{
\eb_0 = \phi^* e_0
}
satisfies
\leqn{wbdefA}{
\eb_0 = \delb{0}\quad \Longleftrightarrow  \quad \eb^\mu_0 = \delta_0^\mu
}
where here and below, we use
\eqn{delbdef}{
\delb{\mu} = \frac{\partial\;}{\partial \xb^\mu}
}
to denote the partial derivative with respect to the Lagrangian coordinates $(\xb^\mu)$.

For use below, we let
\leqn{Jdef}{
J = (J^\mu_\nu) := (\delb{\nu}\phi^\mu)
}
denote the Jacobian matrix of the coordinate transformation \eqref{phidefA} and
\leqn{Jinv}{
\Jch = (\Jch^{\mu}_\nu) := J^{-1}
}
its inverse. Using this notation, the components of the pullback $\eb_0^\mu$ can be computed via the
formula
\eqn{wpback}{
\eb_0^\mu = \Jch^\mu_\nu e_0^\nu\circ\phi.
}
From this, it follows immediately that
\leqn{wbdefB}{
\delb{0}\phi^\mu = e^\mu_0\circ \phi
}
is equivalent to \eqref{wbdefA}.

\subsect{IV}{Initial conditions} Since the vector field $e_0^\mu=w^\mu$ completely determines the proper energy
density $\rho$ and the fluid four-velocity $v^\nu$,
\eqn{f0def}{
f^\mu_0 := e_0^\mu|_{\Omega_0}
}
contains all of the initial data for the fluid. We are then free to choose the other
initial data $f^\mu_I$, see \eqref{eIsolB}, as we like as long as the constraints
\eqref{eiA} are satisfied. For our purposes, we need to fix the $f^\mu_I$ in a specific fashion beyond
just satisfying the constraints \eqref{eiA}. However, before we discuss this, we first make some observations
starting with the pullback frame
\leqn{fwev25}{
\eb_i = \phi^* e_i \quad \Longleftrightarrow  \quad \eb_i^\mu = \Jch^\mu_\nu \et^\nu_i,
}
where
\leqn{etdef}{
\et^\nu_i = e^\nu_i \circ \phi,
}
or equivalently
\leqn{theta3C}{
\et^\mu_i = \eb_i^\lambda \delb{\lambda} \phi^\mu.
}

Since the $e_i$ satisfy \eqref{eiC.1} and the Lie bracket is natural\footnote{That is, the Lie bracket $[\cdot,\cdot]$ satisfies $\psi^*[X,Y]=[\psi^*X,\psi^*Y]$ for
all diffeomorphisms $\psi$ and vector fields $X,Y$.}
with respect to
pullbacks by diffeomorphisms, the $\eb_i$ must also satisfy $[\eb_0,\eb_I]=0$, or equivalently,
by \eqref{wbdefA},
\leqn{fwev27}{
\delb{0} \eb^\mu_j = 0.
}
In particular, this shows that the frame fields $\eb^\mu_I$ are independent of $\xb^0$.
We also note that
\lalin{etidata}{
\et^\mu_j|_{\{0\}\times \Omega} &= f^\mu_j, \label{etidata.1}
\intertext{and}
\bigl(J^\mu_\nu \bigr)|_{\{0\}\times \Omega} &= \begin{pmatrix} f^\mu_0 & \delta^\mu_\Sigma \end{pmatrix} \label{etidata.2}
}
by \eqref{phiOmega} and \eqref{wbdefB}. Using \eqref{etidata.2}, we compute
\eqn{Jchidata}{
\bigl(\Jch^\mu_\nu\bigr)|_{\{0\}\times \Omega}  = \begin{pmatrix}\begin{displaystyle} \frac{1}{f^0_0}\end{displaystyle} & 0\\
\begin{displaystyle} -\frac{1}{f^0_0}f^\Lambda_0\end{displaystyle} & \delta^\Lambda_\Sigma \end{pmatrix}.
}
From this, \eqref{etidata.1} and the fact that the $\eb^\mu_I$ are $\xb^0$-independent, we find
that
\leqn{ebform}{
(\eb_i^\mu) =  \begin{pmatrix} 1 & \begin{displaystyle} f^0_I \end{displaystyle}  \\
0 & \begin{displaystyle} -\frac{1}{f^0_0}f^\Lambda_0 f^0_I + f^\Lambda_I \end{displaystyle} \end{pmatrix}.
}

By assumption, $\eb_0=\delb{0}$ is tangent to the matter-vacuum boundary, which is given in
Lagrangian coordinates by
\eqn{Gammabdef}{
\Gamma_T := [0,T]\times \del{}\Omega.
}
Using the Gram-Schmidt algorithm, we can complete $\eb_0$ to
a basis
\leqn{basisTan}{
\{\eb_0,\zb_1=\zb_1^\mu\delb{\mu},\zb_2=\zb_2^\mu\delb{\mu}\}
}
of $\Gamma_T$
such that the frame is orthogonal
at $\xb^0=0$, and the components $\zb_\Af^\mu$, $\Af=1,2$, are $\xb^0$-independent.
Here, the orthogonality at $\xb^0=0$ is determined with respect to the pull-back metric
\eqn{gbdefA}{
\gb :=\phi^*g = \gb_{\alpha\beta}dx^\alpha dx^\beta
}
with components given by
\leqn{gbdef}{
\gb_{\alpha\beta} = J^\mu_\alpha J^\nu_\beta g_{\mu\nu}\circ\phi.
}
We can then extend the basis \eqref{basisTan} to
$\Omega_T$, while keeping the components  $\zb_\Af^\mu$, $\Af=1,2$, $\xb^0$-independent, and complete it to a full frame
\eqn{baisFull}{
\{\eb_0,\zb_1=\zb_1^\mu\delb{\mu},\zb_2=\zb_2^\mu\delb{\mu},\zb_3=\zb_3^\mu\delb{\mu}\}
}
such that: $\zb_3$ is outward pointing at the boundary $\Gamma_T$, the components $\zb_3^\mu$
are $\xb^0$-independent, and the frame is
orthogonal at $\xb^0=0$. Thus, in particular,
\leqn{OmegaSpanBa}{
\text{Span}\bigl\{\, \{\eb_0, \zb_\Af\}|_{T\Gamma_T} \, | \, \Af=1,2 \bigr\} = T\Gamma_T
}
and
\eqn{OmegaOrthAa}{
\gb(\eb_0,\zb_I)|_{\Omega_0} = 0.
}
We also observe via \eqref{phiOmega}, \eqref{etidata.2} and \eqref{gbdef} that
\eqn{fsetA}{
\gb(\zb_I,\zb_J)\bigl|_{\Omega_0} = \bigl(\zb^0_I f^\mu_0 + \zb^\Sigma_I \delta_\Sigma^\mu\bigr)
\bigl(\zb^0_J f^\nu_0 + \zb^\Lambda_J \delta_\Lambda^\nu\bigr)g_{\mu\nu}|_{\Omega_0}.
}

Next, we set
\lalin{fsetB}{
f^0_I &= \det\bigl(\gb(\zb_I,\zb_J)\bigr)^{-\frac{1}{6}}F\bigl((-g_{00})^{\frac{1}{2}}\bigr)\zb^0_I \Bigl|_{\Omega_0}
\label{fsetB.1}
\intertext{and}
f^\Lambda_I &= \det\bigl(\gb(\zb_I,\zb_J)\bigr)^{-\frac{1}{6}}F\bigl((-g_{00})^{\frac{1}{2}}\bigr)
\biggl(\zb^\Lambda_I + \frac{f^\Lambda_0}{f^0_0} \zb^0_I\biggr)\biggl|_{\Omega_0}. \label{fsetB.2}
}
From the above two expressions and \eqref{ebform}, we obtain
\eqn{fsetC}{
\eb^\mu_I = \det\bigl(\gb(\zb_I,\zb_J)\bigr)^{-\frac{1}{6}}F\bigl((-g_{00})^{\frac{1}{2}}\bigr)^{\frac{1}{3}}
\zb^\mu_I,
}
which, in turn, implies that
\alin{fsetD}{
g(e_0,e_J)|_{\Omega_0} &= \gb(\eb_0,\eb_J)|_{\Omega_0} = 0
\intertext{and}
\det(g(e_I,e_J))|_{\Omega_0} &= \det(\gb(\eb_I,\eb_J))|_{\Omega_0} = F\bigl((-g_{00})^{\frac{1}{2}}\bigr)^2|_{\Omega_0}.
}
This shows the constraints \eqref{eiA} are satisfied for the choice of
initial data $f^\mu_I$ given by \eqref{fsetB.1} and \eqref{fsetB.2}, and
that
\leqn{OmegaSpanB}{
\text{Span}\Bigl\{\, \eb_\af \bigl|_{\Gamma_T } \, | \, \af=0,1,2 \Bigr\} = T\Gamma_T
}
by virtue of \eqref{OmegaSpanBa}. It is also not difficult to see from the above construction
that we can always ensure that
\leqn{deteb}{
 c_{f} > \det(\eb) > \frac{1}{c_{f}} > 0 \hspace{0.4cm}\text{in $\Gamma_T$}
}
holds for some positive constant $c_{f}$.

In terms of the co-frame, we
have that
\leqn{fwev28}{
\thetab^i = \phi^*\thetab^i \quad \Longleftrightarrow  \quad \thetab^i_\mu = J_\mu^\nu \thetat_\nu^i
}
where
\leqn{etdefA}{
(\thetat_\nu^i) := (\theta_\nu^i\circ\phi) = (\et^\nu_i)^{-1} \AND (\thetab_\nu^i) = (\eb^\nu_i)^{-1}.
}
By definition, $\thetab^3_\mu \eb^\mu_j = \delta^3_j$,
and so, by \eqref{wbdefA} and \eqref{OmegaSpanB}, $\thetab^3_\mu$
satisfies $\theta^3_0 = 0$,
and, consequently,
\leqn{normalC}{
\nu_\mu := \thetab^3_\mu = \delta_\mu^\Sigma \theta^3_\Sigma
}
defines a $\xb^0$-independent, outward pointing (non-normalized) co-normal to the boundary $\Gamma_T$,
while $\nu_\Sigma$ defines an outward pointing co-normal to the boundary
$\del{}\Omega$.
Finally, letting
\eqn{fchdef}{
(\fch^j_\mu) = (f^j_\mu)^{-1}
}
denote the inverse of $f^j_\mu$, a short calculation shows that
\leqn{thetabform}{
\bigl(\thetab^{j}_\mu\bigr) = \begin{pmatrix} \delta^j_0 & \fch^j_\Lambda \end{pmatrix}.
}

\subsect{wave}{The wave formulation in Lagrangian coordinates}
The key to transforming the wave equation \eqref{waveF}, or equivalently \eqref{waveE}, into Lagrangian coordinates is the following well known
transformation formula for the divergence of a vector field:
\leqn{divXa}{
|\gb|^{-\frac{1}{2}}\delb{\mu}\bigl(|\gb|^{\frac{1}{2}}\Xb^\mu\bigr) = \phi^*\bigl(|g|^{-\frac{1}{2}}\del{\mu}\bigl(|g|^{\frac{1}{2}}X^\mu\bigr)\bigr)
}
where
\eqn{divXb}{
\Xb^\mu = \Jch^\mu_\nu X^\nu\circ \phi,
}
$|g|=-\det(g_{\mu\nu})$, and $|\gb|=-\det(\gb_{\mu\nu})$.
Setting
\eqn{divXc}{
X^\alpha = |g|^{-\frac{1}{2}}Y^\alpha_\nu
}
where $Y^\alpha_\nu$ is as defined previously by \eqref{Ydef},
a short calculation, using the chain rule, shows that
\leqn{divXd}{
|\gb|^{-\frac{1}{2}}\Xb^\alpha = \biggl(\frac{|\gb|}{|g|\circ\phi}\biggr)^{\frac{1}{2}}\biggl[
\frac{1}{\fft}\ab^{\alpha \beta}\delb{\beta} \thetat^0_\nu
+\frac{1}{\fft}  \ab^{\alpha \beta}\Bigl(- J_\beta^\omega\Gammat_{\omega\nu}^\gamma \thetat^0_\gamma
+ \sigmat_M{}^0{}_L \thetab^M_\beta \thetat^L_\nu\Bigr)\biggr]
}
where
\eqn{abdef}{
\ab^{\alpha \beta} = (\phi^* a)^{\alpha \gamma} = \Jch^\alpha_\mu \Jch^\beta_\nu a^{\mu\nu}\circ\phi
}
denotes the pull-back of the acoustic metric $a^{\mu\nu}$ by $\phi$, and
we have introduced the definitions
\leqn{thetatdef}{
\thetat^0_\nu = \theta^0_\nu \circ \phi,
}
\leqn{sigmatdef}{
\sigmat_{M}{}^0{}_L = \sigma_{M}{}^0{}_L\circ \phi,
}
\leqn{Gammatdef}{
\Gammat^\mu_{\omega \sigma} = \Gamma^\mu_{\omega \sigma}\circ\phi,
}
\leqn{fftdef}{
\fft = \ff \circ \phi
}
and
\leqn{atdef}{
\at^{\mu\nu} = a^{\mu\nu}\circ \phi.
}

Next, we define
\leqn{gtdef}{
\gt_{\mu\nu} = g_{\mu\nu}\circ \phi \AND \gt^{\mu\nu} = g^{\mu\nu} \circ \phi,
}
and note that
\eqn{detgbA}{
|\gb| = -\det\bigl(J^\alpha_\nu \gt_{\alpha\beta} J^\alpha_\nu\bigr) = \det(J)^2 |\gt| = \det(J)^2 |g|\circ \phi,
}
or equivalently
\leqn{detgbB}{
\biggl(\frac{|\gb|}{|g|\circ\phi}\biggr)^{\frac{1}{2}} = \det(J).
}
Additionally, we define
\leqn{gammatdef}{
\fmt_{ij} = \fm_{ij}\circ \phi \AND \fmt^{ij} = \fm^{ij}\circ \phi,
}
and observe that
\leqn{fftdefa}{
\fft = \ff\bigl((-\fmt_{00})^{-1/2}\bigr),
}
where
\leqn{fmt00}{
\fmt_{00} = \gt_{\alpha\beta}\et^\alpha_0\et^\beta_0 =\frac{1}{\gt^{\alpha\beta}\thetat^0_\alpha \thetat^0_\beta},
}
and
\leqn{atexp}{
\at^{\mu\nu} = \gt^{\mu\nu} - \frac{1}{\fmt_{00}}\left(1-\frac{1}{\st^2}\right)\et^\mu_0\et^\nu_0
= \gt^{\mu\nu} - \fmt_{00}\left(1-\frac{1}{\st^2}\right)\gt^{\mu\alpha}\gt^{\nu\beta}\thetat^0_\alpha
\thetat^0_\beta,
}
where
\leqn{stdef}{
\st^2 = s^2\bigl((-\fmt_{00})^{-1/2}\bigr).
}

Using \eqref{divXa}, \eqref{divXd}, and \eqref{detgbB}, it follows that the wave equation \eqref{waveF}
when expressed in Lagrangian coordinates becomes
\leqn{waveG}{
\delb{\alpha}\bigl(\At^{\alpha\beta}\delb{\beta}\thetat^0_\nu + \Lt^\alpha_\nu\bigr) = \Ft_\nu
}
where
\lalin{waveGvars}{
\At^{\alpha\beta} &= -\frac{\det(J)}{\fft}\ab^{\alpha\beta}, \label{waveGvars.1}\\
\Lt^\alpha_\nu &= -\frac{\det(J)}{\fft}  \ab^{\alpha \beta}\Bigl(- J_\beta^\omega\Gammat_{\omega\nu}^\gamma \thetat^0_\gamma
+ \sigmat_M{}^0{}_L \thetab^M_\beta \thetat^L_\nu\Bigr), \label{waveGvars.2} \\
\Ft_\nu &= -\det(J)\biggl(-\Gammat_{\alpha\gamma}^\alpha \Yt^\gamma_\nu +
\Gammat^\gamma_{\alpha\nu}\Yt^\alpha_\gamma + \frac{1}{\ff}\Rt_{\nu}{}^{\lambda}\thetat^0_\lambda\biggr), \label{waveGvars.3} \\
\Yt^\alpha_\nu & = \frac{1}{\fft}\at^{\alpha \beta}\Bigl[\Jch_\beta^\gamma\delb{\gamma} \thetat^0_\nu -
\Gammat_{\beta\nu}^\gamma \thetat^0_\gamma
+ \sigmat_M{}^0{}_L \thetat^M_\beta \thetat^L_\nu\Bigr], \label{waveGvars.4}
\intertext{and}
\Rt_{\nu}{}^{\lambda} &= R_{\nu}{}^{\lambda} \circ \phi. \label{waveGvars.5}
}
From the chain rule, \eqref{eiC.4} and \eqref{wbdefB} it is straightforward
to verify that $\sigmat_i{}^k{}_j$ satisfies the evolution equation
\eqn{tvarsC}{
\delb{0}\sigmat_i{}^k{}_j = 0.
}
Using this and \eqref{phiOmega}, it follows immediately
that
\leqn{tvarsD}{
\sigmat_i{}^k{}_j(\xb^0,\xb^\Sigma) = \sigma_i{}^k{}_j(0,\xb^\Sigma).
}
We also observe that
\leqn{phiEv}{
\delb{0}\phi^\mu = \gammat_{00}\gt^{\mu\nu}\thetat^0_\nu
}
follows from \eqref{hC} and \eqref{wbdefB},
while
\leqn{thetatIEv}{
\thetat^I_\mu = \Jch^\nu_\mu \thetab^I_\nu
}
is a consequence of \eqref{fwev28}.

The two evolution equations \eqref{waveG} and \eqref{phiEv} along with
the definitions \eqref{Jdef}, \eqref{Jinv}, \eqref{Gammatdef},
\eqref{gtdef},  \eqref{fftdefa}-\eqref{stdef}, \eqref{waveGvars.1}-\eqref{tvarsD}, and
\eqref{thetatIEv} constitute our \emph{Lagrangian wave formulation of the Euler equations}.

\begin{rem} \label{completerem}
In the simpler setting where the pressure never vanishes and there
is no free boundary, it already follows from known results that our Lagrangian wave formulation is complete in the sense
that the evolution equations \eqref{waveG} and \eqref{phiEv} for the pair $(\phi^\mu,\thetat^0_\mu)$ are well-posed, and in particular, solutions of \eqref{waveG} and \eqref{phiEv}
satisfy energy estimates. Indeed, the evolution equations \eqref{waveG} and
\eqref{phiEv} fit within the class of wave equations considered in \cite{Koch:1993}. Although boundary conditions,
which do not include the boundary conditions consider in this article,
were imposed in \cite{Koch:1993},
it is clear that results of \cite{Koch:1993} continue to hold in the absence of a boundary, and
hence apply to the system \eqref{waveG} and \eqref{phiEv} immediately when no boundary is present.
\end{rem}

\subsect{tdiffwave}{The time differentiated wave equation}
Due to the free boundary, it is not enough to consider the evolution equations \eqref{waveG} and \eqref{phiEv}
alone.
Instead, we must supplement them with what amounts to a $\xb^0$-differentiated version of the wave equation \eqref{waveG}.
With this in mind, we define
\eqn{psidefA}{
\psi_\nu = \bigl(\nabla_{e_0} \theta^0_\nu\bigr) \circ \phi,
}
and note using \eqref{wbdefA} and the chain rule
that this can be written as
\eqn{psidefB}{
\psi_\nu = \delb{0}\thetat^0_\nu + \beta^\lambda_\nu \thetat^0_\lambda
}
where we have set
\eqn{psidefC}{
\beta^\lambda_\nu := - \et^\gamma_0 \Gammat^\lambda_{\gamma\nu}
\oset{\eqref{hC}}{=} -\gammat_{00}\gt^{\gamma \sigma}\thetat^0_\sigma
\Gammat^\lambda_{\gamma\nu}.
}
Differentiating \eqref{waveG} with respect to $\xb^0$,
a short calculation shows that $\psi_\nu$ satisfies the wave equation
\leqn{waveH}{
\delb{\alpha}\bigl(\At^{\alpha\beta}\delb{\beta}\psi_\nu + L^\alpha_\nu\bigr) = F_\nu
}
where
\alin{waveHvars}{
L^\alpha_\nu &= \delb{0}\At^{\alpha\beta}\delb{\beta}\thetat^0_\nu+\delb{0}\Lt^\alpha_\nu
+\beta^\omega_\nu \Lt^\alpha_\omega - \At^{\alpha\beta}\delb{\beta}\beta^\omega_\nu \thetat^0_\omega,
\intertext{and}
F_\nu &= \delb{0}\Ft_\nu + \beta^\omega_\nu \Ft_\omega + \delb{\alpha}\beta^\omega_\nu\bigl(\At^{\alpha\beta}
\delb{\beta}\thetat^0_\omega + \Lt^\alpha_\omega\bigr).
}

Next, we define a positive definite, symmetric 2-tensor by
\leqn{mupdef}{
m^{\alpha\beta} := g^{\alpha\beta} - \frac{1}{\fm_{00}} e^\alpha_0 e^\beta_0
\overset{\eqref{hC}}{=} g^{\alpha\beta} - \fm_{00}g^{\mu\alpha}g^{\nu\beta}\theta^0_\mu \theta^0_\nu,
}
and we let
\eqn{mdndefA}{
(m_{\alpha\beta}) = (m^{\alpha\beta})^{-1}
}
denote its inverse. A short calculation then shows that
\eqn{mdndefB}{
m_{\alpha\beta} = g_{\mu\nu} - \frac{1}{\fm^{00}} \theta^0_\mu \theta^0_\nu.
}
Setting
\eqn{mtdef}{
\mt^{\mu\nu} = m^{\mu\nu}\circ \phi \AND \mt_{\mu\nu} = m_{\mu\nu}\circ \phi,
}
we then have that
\gath{mtdefC}{
\mt_{\alpha \beta} = \gt_{\alpha\beta} -\frac{2}{\gammat^{00}}\thetat^0_\alpha\theta^0_\beta, %\label{mtdefC.1}
\\
\mt^{\alpha\beta} = \gt^{\alpha\beta} - \frac{2}{\gammat_{00}}\et_0^\alpha \et_0^\beta
=  \gt^{\alpha\beta} - 2\fmt_{00}\gt^{\mu\alpha}\gt^{\nu\beta}\thetat^0_\alpha \thetat^0_\beta, %\label{mtdefC.2}
\intertext{and}
\mt_{\alpha \beta} \mt^{\beta \gamma} = \delta_\alpha^\gamma. %\label{mtdefC.3}
}
Continuing on, we define\footnote{Here and
below, we use following standard notation for norms: $|X|_g=(g_{\mu\nu}X^\mu X^\nu)^{\frac{1}{2}}$,
$|\chi|_g = (g^{\mu\nu}\chi_{\mu}\chi_{\nu}))^{\frac{1}{2}}$,
$|Y|_m = (m_{\mu\nu}Y^\mu Y^\nu)^{\frac{1}{2}}$ and $|\xi|_m= (m^{\mu\nu}\xi_\mu\xi_\nu)^{\frac{1}{2}}$ where $X^\mu,\chi_\mu$ are spacelike and $Y^\mu,\xi_\mu$ are arbitrary, and similar notation when
$g$ and $m$ are replaced with $\gt$ and $\mt$. }
\leqn{mudef}{
\mu = |\psi|_{\mt} %:= \bigl(\mt^{\alpha\beta}\psi_\alpha \psi_\beta\bigr)^{\frac{1}{2}}
}
and
\leqn{psihdef}{
\psih_\nu = \frac{1}{\mu} \psi_\nu.
}
\begin{rem} \label{murem}
The variable $\psih_\nu$ is well defined wherever $\mu=|\nabla_{e_0}\theta^0|_m\circ \phi$ is positive.
As we show in Section \ref{acoustic} below, see  \eqref{theta3rep1}
and \eqref{e3g00}, the Taylor sign condition implies
that $\mu$ is bounded away from zero on the boundary $\Gamma_T$. From the smoothness of the solution
it follows immediately that $\mu$ is bounded away from zero in a neighborhood of the boundary. Away
from the boundary, we can, via the finite propagation speed, obtain energy estimates using the theory of symmetric hyperbolic systems. Because of this, it is enough to consider the problem in a neighborhood of the boundary, and
consequently, we lose no generality by assuming that
\leqn{muassump}{
\mu \geq c_\mu > 0 \hspace{0.4cm}\text{in $\overline{\Omega}_T$}
}
for some positive constant $c_\mu$.
\end{rem}
Differentiating $\psih_\nu$, we find that
\leqn{dpsih}{
\delb{\mu}\psih_\nu = \frac{1}{\mu}\pi^\lambda_\nu \delb{\mu}\psi_\lambda - \frac{1}{2} \delb{\mu}\mt^{\alpha\beta}
\psih_{\alpha}\psih_\beta \psih_\nu,
}
where
\leqn{pidef}{
\pi^\lambda_\mu = \delta^\lambda_\nu - \psih^\lambda\psih_\nu   \qquad (\psih^\lambda := \mt^{\lambda\gamma}\psih_\gamma)
}
is the projection operator uniquely defined by the properties
\leqn{piprop}{
\pi^\lambda_\nu \psih_\lambda = 0 \AND \pi^\lambda_\nu \pi_\lambda^\mu = \pi^\mu_\nu,
}
while differentiating $\mu$ yields
\leqn{dmu}{
\delb{\gamma}\mu = \psih^\alpha \delb{\gamma}\psi_\alpha + \frac{\mu}{2}\delb{\gamma}\mt^{\alpha\beta}\psih_\alpha\psih_\beta.
}
For use below, we define the following variants of the projection operator $\pi^\mu_\nu$:
\eqn{piupdn}{
\pi^{\mu\nu} = \mt^{\mu\lambda}\pi_\lambda^\nu \AND \pi_{\mu\nu} = \mt_{\mu\lambda} \pi^\lambda_\nu.
}

From \eqref{dpsih}, we see that
\eqn{Adpsih}{
\At^{\alpha\beta}\delb{\beta} \psih_\nu =
\frac{1}{\mu}\pi^\omega_\nu \At^{\alpha\beta}\delb{\beta}\psi_\omega - \frac{1}{2} \At^{\alpha\beta}
\delb{\beta}\mt^{\sigma\gamma}\psih_\sigma\psih_\gamma \psih_\nu,
}
and hence that
\eqn{Adpsihdiv}{
\delb{\alpha}\bigl(\mu\At^{\alpha\beta}\delb{\beta} \psih_\nu\bigr) =
\delb{\alpha}\biggl(\pi^\omega_\nu \At^{\alpha\beta}\delb{\beta}\psi_\omega - \frac{\mu}{2} \At^{\alpha\beta}
\delb{\beta}\mt^{\sigma\gamma}\psih_\sigma\psih_\gamma \psih_\nu\biggr).
}
Using this and \eqref{waveH}, we see that $\psih_\nu$ satisfies the
wave equation
\leqn{waveI}{
\delb{\alpha}\biggl(\mu \At^{\alpha\beta}\delb{\beta}\psih_\nu +
\frac{\mu}{2} \At^{\alpha\beta}\delb{\beta}\mt^{\sigma\gamma}\psih_\sigma\psih_\gamma\psih_\nu
+ \pi^\omega_\nu L_\omega^\nu\biggr) = \pi^\omega_\nu F_\omega + \delb{\alpha}\pi^\omega_\nu\bigl(\At^{\alpha\beta}\delb{\beta}\psi_\omega
+L^\alpha_\omega\bigr).
}
A similar computation starting from the identity
\eqn{Amu}{
\At^{\alpha\beta}\delb{\beta}\mu = \frac{\psi^\gamma}{\mu} \At^{\alpha\beta}
\delb{\beta}\psi_\gamma + \frac{\mu}{2} \delb{\beta}\mt^{\sigma\gamma}\At^{\alpha\beta} \psih_\sigma\psih_\gamma,
}
which holds by \eqref{dmu},
shows that $\mu$ satisfies the wave equation
\leqn{waveJ}{
\delb{\alpha}\biggl(\At^{\alpha\beta}\del{\beta}\mu -\frac{\mu}{2}\At^{\alpha\beta}\delb{\beta}\mt^{\sigma\gamma}\psih_\sigma\psih_\gamma
+\psih^\gamma L^\alpha_\gamma \biggr) = \psih^\gamma F_\gamma+ \delb{\alpha}\psih^{\gamma}\bigl(
\At^{\alpha\beta}\delb{\beta}\psi_\gamma+L^\alpha_\gamma\bigr).
}
Setting
\leqn{Psidef}{
\Psi = \begin{pmatrix}\psih_\nu \\ \mu \end{pmatrix},
}
we combine \eqref{waveI} and \eqref{waveJ} into the single equation
\leqn{waveK}{
\delb{\alpha}\bigl(\Ac^{\alpha\beta}\delb{\beta}\Psi+\Lc^\alpha\bigr) = \Fc
}
where
\alin{waveKvars}{
\Ac^{\alpha\beta} &= \begin{pmatrix}\mu^2 \At^{\alpha\beta}\mt^{\mu\nu} & 0 \\ 0 & \At^{\alpha\beta} \end{pmatrix},
%\label{waveKvars.1}
\\
\Lc^{\alpha} &= \begin{pmatrix}\frac{\mu^2}{2} \At^{\alpha\beta}\delb{\beta}\mt^{\sigma\gamma}\psih_\sigma\psih_\gamma\psih^\mu
+ \mu\pi^{\omega\mu} L_\omega^\nu & 0 \\
0 & -\frac{\mu}{2}\At^{\alpha\beta}\delb{\beta}\mt^{\sigma\gamma}\psih_\sigma\psih_\gamma
+\psih^\gamma L^\alpha_\gamma \end{pmatrix}, %\label{waveKvars.2}
\intertext{and}
\Fc &= \begin{pmatrix} \fc_\nu \\ \fc \end{pmatrix}% \label{waveKvars.3}
}
with
\alin{waveKvarsA}{
\fc_\nu &= \mu\mt^{\mu\nu}\Bigl(\pi^\omega_\nu F_\omega + \delb{\alpha}\pi^\omega_\nu\bigl(\At^{\alpha\beta}\delb{\beta}\psi_\omega
+L^\alpha_\omega\bigr)\Bigr)  \\
&\hspace{2.0cm} + \delb{\alpha}\bigl(\mu\mt^{\mu\nu}\bigr)\Bigl(\mu \At^{\alpha\beta}\delb{\beta}\psih_\nu +
\frac{\mu}{2} \At^{\alpha\beta}\delb{\beta}\mt^{\sigma\gamma}\psih_\sigma\psih_\gamma\psih_\nu
+ \pi^\omega_\nu L_\omega^\nu\Bigr), \\
\fc &= \psih^\gamma F_\gamma+ \delb{\alpha}\psih^{\gamma}\bigl(
\At^{\alpha\beta}\delb{\beta}\psi_\gamma+L^\alpha_\gamma\bigr).
}

\sect{boundary}{Boundary conditions}

By virtue of \eqref{wbdefA}, the boundary condition \eqref{wtanB} is automatically
incorporated into the definition of the Lagrangian coordinates \eqref{phidefA}.
 Therefore, from this perspective the only non-trivial boundary condition is \eqref{passumpA}, or equivalently
\eqref{fbcB}, which in Lagrangian coordinates becomes
\leqn{dirbc}{
\fmt_{00}|_{\Gamma_T} = -1.%\quad \oset{\eqref{tiA.3}}{\Longleftrightarrow}
%\quad \fmt^{00}|_{\Gamma_T} = -1.
}
Written this way, there does not seem to be enough boundary conditions to derive energy
estimates using our wave formulation. However, as we show below
this single boundary condition does, in fact, imply a sufficient set of boundary conditions.
We also note that  \eqref{fbcB} implies, by
\eqref{Fdef.2}, \eqref{fdefB} and \eqref{fftdefa},
that
\leqn{dirbcf}{
\fft|_{\Gamma_T} = -1.
}

\subsect{neumann}{Neumann boundary conditions} The first set of boundary conditions that
we derive from \eqref{dirbc}, or equivalently \eqref{fbcB} in the Eulerian picture, are of Neumann type, albeit
degenerate.
The derivation of these boundary conditions begins with the identity
\lalin{bcA}{
\theta^3_\alpha g^{\alpha\beta}\nabla_{\beta} \theta^0_\nu &=
\theta^3_\alpha g^{\alpha\beta}\nabla_{\nu}\theta_\beta^0 - \sigma_{M}{}^0{}_L \theta^3_\alpha g^{\alpha\beta}\theta_\beta^I\theta^J_\nu
&& \text{(by \eqref{CartanE})} \notag \\
&= g(\theta^3,\nabla_{\nu}\theta^0) - \sigma_{M}{}^0{}_L g(\theta^3,\theta^M)\theta^L_\nu  \notag\\
&= -g(\nabla_{\nu}\theta^3,\theta^0) - \sigma_{M}{}^0{}_L g(\theta^3,\theta^M)\theta^L_\nu && \text{(by \eqref{tiA.1})} \notag \\
&= -g^{\alpha\beta}\nabla_{\nu}\theta^3_\alpha \theta^0_\beta - \sigma_{M}{}^0{}_{L}g(\theta^3,\theta^M)\theta^L_\nu \notag\\
&=-g^{\alpha\beta}\bigl[\nabla_{\alpha}\theta^3_\nu-\sigma_{M}{}^3{}_L\theta^M_\nu\theta^L_\alpha\bigr]\theta^0_\beta
- \sigma_{M}{}^0{}_{L}g(\theta^3,\theta^M)\theta^L_\nu   && \text{(by setting $i=3$ in \eqref{CartanD})} \notag \\
&=-\theta_\beta^0 g^{\alpha\beta}\nabla_{\alpha}\theta^3_\nu - \sigma_{M}{}^0{}_{L}g(\theta^3,\theta^M)\theta^L_\nu
&& \text{(by \eqref{tiA.1})} \notag \\
&= -\frac{1}{\fm_{00}}e^\alpha_0\nabla_{\alpha}\theta^3_\nu - \sigma_{M}{}^0{}_{L}g^{\alpha\beta}\theta^3_\alpha \theta^M_\beta\theta^L_\nu
&& \text{(by \eqref{hC})} \label{bcA.1}.
}
Noting that
\eqn{bcB}{
\theta^3_\alpha g^{\alpha\beta} = \theta^3_\alpha a^{\alpha\beta}
}
follows from  \eqref{adefB}, we can write \eqref{bcA.1} as
\leqn{bcC}{
\theta^3_\alpha \bigl(a^{\alpha\beta}[\nabla_{\beta} \theta^0_\nu + \sigma_{M}{}^0{}_{L}\theta^M_\beta\theta^L_\nu]\bigr) = -\frac{1}{\fm_{00}}e^\alpha_0\nabla_{\alpha}\theta^3_\nu.
}
Transforming this expression into Lagrangian coordinates  gives
\leqn{bcEa}{
\nml_\alpha\bigl( \At^{\alpha\beta} \delb{\beta} \thetat^0_\nu + \Lt^\alpha_\nu\bigr)
\bigl|_{\Gamma_T}
= \frac{\det(J)}{\fft\fmt_{00}}
\bigl(\delb{0}\thetat^3_\nu +\beta_\nu^\gamma\thetat^3_\gamma\bigr)
}
where $\nml_\alpha$ is the outward co-normal to $\Gamma_T$ defined by \eqref{normalC}.

Next, we calculate
\alin{detgbC}{
|g|\circ\phi &= -\det(\gt_{\mu\nu}) \\
&= \det(\thetat^i_\mu \fmt_{ij}\theta^j_\nu) \\
&= \det(\thetat)^2 \det(\fmt_{ij}) \\
&= \det(\thetab\Jch)^2\bigl[\fmt_{00}\ff\bigl((-\fmt_{00})^{-1/2}\bigr)\bigr]^2
&& \text{(by \eqref{fdefA} and \eqref{fwev28})}\\
&= \left(\frac{\det(\thetab)}{\det(J)}\right)^2\bigl[\fmt_{00}\ff\bigl((-\fmt_{00})^{-1/2}\bigr)\bigr]^2,
}
which after taking the square root, gives
\leqn{detgbD}{
\det(J) = \frac{\det(\thetab)}{(|g|\circ \phi)^{\frac{1}{2}}}\fmt_{00}\fft,
}
and allows us to write \eqref{bcEa} as
\leqn{bcE}{
\thetab^3_\alpha\bigl( \At^{\alpha\beta} \delb{\beta} \thetat^0_\nu + \Lt_\nu^\alpha\bigr)\bigl|_{\Gamma_T}= \frac{\det(\thetab)}{|\gt|^{\frac{1}{2}}}
\bigl(\delb{0}\thetat^3_\nu + \beta^\gamma_\nu \thetat^3_\gamma\bigr).
}
We now examine the structure of the righthand side of \eqref{bcE}. By definition, $(\thetat^\mu_j)$ is the inverse of $(\et^j_\mu)$,
and so, we must have that
\eqn{bcF}{
\thetat^3_\mu = \frac{\cof(\et)^\mu_3}{\det(\et)}.
}
But,
\alin{detet}{
\det(\et) &= \det(J\eb) && \text{(by \eqref{fwev25})} \\
&= \det(J)\det(\eb) \\
&= \frac{\det(\thetab)\det(\eb)\fmt_{00}\fft}{\sqrt{|\gt|}} && \text{(by \eqref{detgbD})} \\
&= \frac{\fmt_{00}\fft}{|\gt|^{\frac{1}{2}}}  && \text{(since $\thetab = \eb^{-1}$),}
}
and so
\eqn{thetat3Aa}{
\thetat^3_\mu =  \frac{|\gt|^{\frac{1}{2}}}{\fmt_{00}\fft} \cof(\et)^\mu_3.
}
By definition of the cofactor matrix,\footnote{Recall for  $a= (a^\mu_\nu)\in \Mbb{n}$ that
\eqn{cofdef}{
\cof(a)^\mu_\nu = \frac{1}{(n-1)!} \ep^\mu{}_{\mu_1 \mu_2 \ldots \mu_{n-1}}\ep_j{}^{\nu_1 \nu_2 \ldots \nu_{n-1}} a^{\mu_1}_{\nu_1}a^{\mu_2}_{\nu_2}\cdots
a^{\mu_{n-1}}_{\nu_{n-1}}
}
where $\ep_{\mu_1 \mu_2 \ldots \mu_n}$ is the completely antisymmetric symbol and
all indices are raised using the  Kronecker delta $\delta^{\mu\nu}$.
}
\eqn{bcK}{
\cof(\et)^\mu_3 = - \ep_{\mu\alpha\beta\gamma}\et^\alpha_0\et^\beta_1\et^\gamma_2,
}
which gives
\eqn{thetat3A}{
\thetat^3_\mu =  -\frac{|\gt|^{\frac{1}{2}}}{\fmt_{00}\fft}\ep_{\mu\alpha\beta\gamma}\et^\alpha_0\et^\beta_1\et^\gamma_2.
}
Evaluating this at the boundary, we see, with the help of
\eqref{dirbcf},  that
\leqn{thetat3B}{
\thetat^3_\mu |_{\Gamma_T} =  -|\gt|^{1/2}\ep_{\mu\alpha\beta\gamma}\et^\alpha_0\et^\beta_1\et^\gamma_2.
}
Differentiating \eqref{thetat3B} with respect to $\xb^0$, we find, using
\eqref{theta3C} and \eqref{fwev27}, that
\leqn{theta3D}{
-\delb{0} \thetat^3_\mu  |_{\Gamma_T} =  |\gt|^{1/2}\bigl(\ep_{\mu\nu\alpha\gamma} \et^\alpha_1 \et^\gamma_2 \eb^\omega_0
+  \ep_{\mu\alpha \nu \gamma} \et^\alpha_0 \et^\gamma_2 \eb^\omega_1
+ \ep_{\mu\alpha \gamma \nu} \et^\alpha_0 \et^\gamma_1 \eb^\omega_2\bigr)\delb{\omega}\et_0^\nu +  \frac{1}{2|\gt|^{1/2}}\frac{\partial |\gt|}{\partial \phi^\kappa} \et^\kappa_0 \ep_{\mu\alpha\beta\gamma}\et^\alpha_0\et^\beta_1\et^\gamma_2,
}
while a short calculation, using \eqref{hC}, \eqref{OmegaSpanB} and \eqref{dirbc}, shows that
\leqn{theta3Da}{
\eb_{\af}(\et_0^\nu)\bigl|_{\Gamma_T} = -\gt^{\nu\omega}\eb_{\af}(\thetat^0_\omega)-
\eb_{\af}(\gt^{\nu\omega})\thetat^0_\omega = -\mt^{\nu\omega}\eb_{\af}(\thetat^0_\omega)-
\eb_{\af}(\gt^{\nu\omega})\thetat^0_\omega+ \eb_{\af}(\gt^{\sigma\omega})\et^\nu_0 \thetat^0_\sigma \thetat^0_\omega .
}

Taken together, \eqref{bcE}, \eqref{thetat3B}, \eqref{theta3D} and \eqref{theta3Da} imply that
\leqn{thetatBC}{
\nml_\alpha\bigl(\At^{\alpha\beta} \delb{\beta} \thetat^0_\nu + \Lt_\nu^\alpha\bigr)\bigl|_{\Gamma_T}=
\St_{\nu \mu}{}^\omega \mt^{\mu \gamma}\delb{\omega}\thetat^0_\gamma + \Zt_\nu
}
where
\eqn{StdefA}{
\St_{\nu\mu}{}^\omega = \det(\thetab)
\bigl(\ep_{\nu\mu\alpha\gamma} \et^\alpha_1 \et^\gamma_2 \eb^\omega_0
+  \ep_{\nu\alpha \mu \gamma} \et^\alpha_0 \et^\gamma_2 \eb^\omega_1
+ \ep_{\nu\alpha \gamma \mu} \et^\alpha_0 \et^\gamma_1 \eb^\omega_2\bigr),
}
and
\eqn{Ztdef}{
\Zt_{\nu} = \St_{\nu\mu}{}^\omega\delb{\omega}(\gt^{\mu\lambda})\thetat^0_\lambda-  \St_{\nu\mu}{}^\omega\delb{\omega}
(\gt^{\sigma\lambda})\et^\nu_0 \thetat^0_\sigma \thetat^0_\lambda
- \det(\thetab)\biggl(\frac{1}{2|\gt|}\frac{\partial |\gt|}{\partial \phi^\kappa} \et^\kappa_0
\delta^\lambda_\nu + \beta^\lambda_\nu\biggr)\ep_{\lambda\alpha\beta\gamma}\et^\alpha_0\et^\beta_1\et^\gamma_2.
}
Setting
\leqn{StdefB}{
\St^{\mu\nu\omega} := \mt^{\mu\alpha}\St_{\alpha\beta}{}^\omega \mt^{\beta\nu},
}
we observe that $\St^{\mu\nu\omega}$ satisfies
\leqn{Stprop}{
\St^{\mu\nu\omega} = -\St^{\nu\mu\omega} \AND \nu_\omega \St^{\mu\nu\omega}= 0.
}

\begin{rem} \label{neurem}
Equation \eqref{thetatBC} is the fundamental Neumann boundary condition satisfied by our system.
It is important to realize that this Neumann boundary condition is degenerate
in the sense that it does not yield coercive elliptic estimates of the type
\eqref{coerc}, and as such, is not directly useful for deriving energy estimates.
However, by considering the time differentiated version
of these boundary conditions, we show below that this degeneracy can be removed, although it should be
noted that the resulting boundary conditions are of acoustic type.
\end{rem}

To proceed, we differentiate \eqref{thetatBC} with respect to $\xb^0$ to find, after a short calculation, that
$\psi_\nu$ satisfies
\leqn{psiBC}{
\nml_\alpha\bigl(\At^{\alpha\beta} \delb{\beta} \psi_\nu + L_\nu^\alpha\bigr)\bigl|_{\Gamma_T}=
\St_{\nu \mu}{}^\omega \mt^{\mu \gamma}\delb{\omega}\psi_\gamma + Z_\nu
}
where
\eqn{Zdef}{
Z_\nu = -\St_{\nu \mu}{}^\omega \mt^{\mu \gamma}\delb{\omega}\beta_\gamma^\lambda\thetat^0_\lambda
+ \delb{0}\bigl(\St_{\nu \mu}{}^\omega \mt^{\mu \gamma}\bigr)\delb{\omega}\thetat^0_\gamma
+\delb{0}\Zt_{\nu}.
}
From \eqref{dpsih}, \eqref{piprop}, \eqref{dmu}, and \eqref{StdefB}, we see also that
\lgath{Sdecomp}{
\pi_\nu^\delta \St_{\delta \mu}{}^\omega \mt^{\mu \gamma} \delb{\omega}\psi_\gamma
= \mu \pi_{\nu\delta}\St^{\delta \mu \omega}\pi_\mu^\tau \delb{\omega}\psih_\tau
+ \pi_{\nu\delta}\St^{\delta \mu \omega}\psih_\mu \delb{\omega} \mu -
\frac{\mu}{2} \pi_{\nu\delta}\St^{\delta \mu \omega}\psih_\mu \delb{\omega}\mt^{\sigma\lambda}\psih_\sigma \psih_\lambda,
\label{SdecompA}
\intertext{and}
\psih^{\nu} \St_{\nu\mu}{}^\omega\mt^{\mu\gamma}\delb{\omega}\psi_\gamma
=\mu \psih_\nu \St^{\nu \mu \omega}\pi_\mu^\tau \delb{\omega}\psih_\tau. \label{SdecompB}
}
A straightforward calculation using \eqref{psiBC}-\eqref{SdecompB} then shows that
\lalin{psihBC}{
\nml_\alpha\biggl(\mu^2 \At^{\alpha\beta} \delb{\beta} \psih_\nu
+\frac{\mu^2}{2}\At^{\alpha\beta}\delb{\beta}\mt^{\sigma\gamma}\psih_\sigma\psih_\gamma\psih_\nu+
\mu\pi^\omega_\nu L_\omega^\alpha\biggr)\bigl|_{\Gamma_T} = &
\mu^2 \pi_{\nu\delta}\St^{\delta \mu \omega}\pi_\mu^\tau \delb{\omega}\psih_\tau
 \notag \\
  + \mu\pi_{\nu\delta}\St^{\delta \mu \omega}\psih_\mu \delb{\omega} \mu  -
\frac{\mu^2}{2} \pi_{\nu\delta}&\St^{\delta \mu \omega}\psih_\mu \delb{\omega}\mt^{\sigma\lambda}\psih_\sigma \psih_\lambda
+\mu\pi^\omega_\nu Z_\omega, \label{psihBC.1}
}
and
\leqn{muBC}{
\nml_\alpha\biggl(\At^{\alpha\beta} \delb{\beta} \mu - \frac{\mu}{2}\At^{\alpha\beta}\delb{\beta}\mt^{\sigma\gamma}\psih_\sigma\psih_\gamma +
\psih^\nu L^\alpha_\nu\biggr)\bigl|_{\Gamma_T}=
 \mu \psih_\nu \St^{\nu \mu \omega}\pi_\mu^\tau \delb{\omega}\psih_\tau +
\psih^\omega Z_\omega.
}

Employing the definition \eqref{Psidef}, we collect \eqref{psihBC.1}-\eqref{muBC} into the single
boundary condition
\leqn{PsiBC}{
\nml_\alpha\bigl(\Ac^{\alpha\beta}\delb{\beta}\Psi +\Lc^\alpha \bigr)\bigl|_{\Gamma_T} = \Sc^{\omega}\delb{\omega}\Psi + \Gc
}
where
\alin{PsiBCvars}{
\Sc^\omega &= \begin{pmatrix} \mu^2 \pi^\mu_\alpha\St^{\alpha \beta \omega}\pi_\beta^\nu &
\mu\pi^{\mu}_{\alpha}\St^{\alpha \beta \omega}\psih_\beta \\
 \mu \psih_\alpha \St^{\alpha \beta \omega}\pi_\beta^\nu & 0 \end{pmatrix}, %\label{PsiBCvars.1}
 \intertext{and}
\Gc &= \begin{pmatrix}\frac{\mu^2}{2} \pi^{\mu}_\alpha\St^{\alpha \beta \omega}\psih_\beta \delb{\omega}\mt^{\sigma\lambda}\psih_\sigma \psih_\lambda
+\mu\pi^{\mu\omega} Z_\omega \\
\psih^\omega Z_\omega \end{pmatrix}. %\label{PsiBCvars.2}
}
Additionally, we note that
\leqn{Scprops}{
(\Sc^{\omega})^{\tr}=-\Sc^{\omega} \AND \nml_\alpha S^\omega = 0
}
by \eqref{Stprop}.

\subsect{acoustic}{Acoustic boundary conditions} The boundary conditions \eqref{PsiBC}
for the time differentiated system are still degenerate in the sense of Remark \ref{neurem}.
The first step in removing this degeneracy is to observe that at the boundary
\alin{om00}{
\omega_0{}^0{}_i|_{B_T} &= -\omega_{00i} && \text{(by \eqref{fbcB} and \eqref{eiC.2}) }\\
&= -\omega_{i00} && \text{(by \eqref{eiE} and \eqref{CartanC})} \\
&= -\Half e_i(\gamma_{00}) && \text{(by \eqref{CartanA.2})}\\
&= -\Half e_3(\gamma_{00})\delta^3_i && \text{(by \eqref{fbcB} and \eqref{OmegaSpanB}).}
}
From this and the connection formula
\eqn{cotheta}{
\nabla_{e_0}\theta^0_\mu = -\omega_0{}^0{}_i\theta^i_\mu,
}
we obtain
\leqn{theta3rep1}{
\theta^3_\mu |_{B_T} = \frac{2}{e_3(\fm_{00})} \nabla_{e_0}\theta^0_\mu,
}
which is well defined since
\leqn{e3g00}{
e_3(\fm_{00})|_{B_T} \leq c < 0
}
for some constant $c$ by \eqref{TaylorC}. Clearly, \eqref{theta3rep1} implies that
\eqn{theta3rep2}{
\theta^3_\mu |_{B_T} = -|\theta^3|_g \frac{\nabla_{e_0}\theta^0_\mu}{|\nabla_{e_0}\theta^0|_g},
}
or equivalently, since $|\nabla_{e_0}\theta^0|_g|_{B_T} = |\nabla_{e_0}\theta^0|_m$ by \eqref{fbcB}
and \eqref{mupdef},
\eqn{theta3rep3}{
\theta^3_\mu |_{B_T} = -|\theta^3|_g \frac{\nabla_{e_0}\theta^0_\mu}{|\nabla_{e_0}\theta^0|_m}.
}
Using this, we can write \eqref{bcC} as
\leqn{acousticA}{
\theta^3_\alpha \bigl(a^{\alpha\beta}[\nabla_{\beta} \theta^0_\nu + \sigma_{M}{}^0{}_{L}\theta^M_\beta\theta^L_\nu]\bigr)|_{B_T} = \frac{1}{\fm_{00}}\nabla_{e_0}\biggl(|\theta^3|_g \frac{\nabla_{e_0}\theta^0_\mu}{|\nabla_{e_0}\theta^0|_m}\biggr).
}
We also observe that
\lalin{dtmuA}{
e_0(|\nabla_{e_0}\theta^0|^2_m)|_{B_T}
&= \nabla_{e_0}\biggl(g^{\mu\nu}\nabla_{e_0}\theta^0_\mu \nabla_{e_0}\theta^0_\mu - \frac{\bigl(g^{\mu\nu}\theta_\mu^0\nabla_{e_0}\theta^0_\nu\bigr)^2}{\gamma^{00}}\biggr) \notag &&\text{(by \eqref{mupdef})} \\
&= \nabla_{e_0}\biggl(g^{\mu\nu}\nabla_{e_0}\theta^0_\mu \nabla_{e_0}\theta^0_\mu - \frac{\bigl(e_0(\gamma^{00})\bigr)^2}{4\gamma^{00}}\biggr) \notag \\
& = 2 g^{\mu\nu}\nabla_{e_0}\theta^0_\mu \nabla_{e_0}\nabla_{e_0}\theta^0_\mu
&& \text{(by \eqref{fbcB})}\notag \\
&=2m^{\mu\nu}\nabla_{e_0}\theta^0_\mu \nabla_{e_0}\nabla_{e_0}\theta^0_\mu
+ \frac{2}{\fm^{00}}g^{\mu\nu}\theta_\mu^0\nabla_{e_0}\theta_\nu^0 g^{\alpha\beta} \theta_\alpha^0\nabla_{e_0}
\nabla_{e_0}\theta_\beta^0  \notag \\
&= 2m^{\mu\nu}\nabla_{e_0}\theta^0_\mu \nabla_{e_0}\nabla_{e_0}\theta^0_\mu
+ \frac{1}{\fm^{00}} e_{0}(\fm^{00})g^{\alpha\beta} \theta_\alpha^0\nabla_{e_0}
\nabla_{e_0}\theta_\beta^0 \notag \\
&= 2m^{\mu\nu}\nabla_{e_0}\theta^0_\mu \nabla_{e_0}\nabla_{e_0}\theta^0_\mu &&\text{(by \eqref{fbcB}).} \label{dtmuA.1}
}
When expressed in Lagrangian coordinates,  \eqref{acousticA} and \eqref{dtmuA.1} become
\leqn{acousticB}{
\thetab^3_\alpha \bigl(\At^{\alpha\beta}\delb{\beta} \thetat^0_\nu + \Lt^\alpha_\nu\bigr)|_{\Gamma_T} =
-\frac{\det(J)}{\fft\fmt_{00}}\Dc_0\bigl(|\thetat^3|_{\gt} \psih_\nu\bigr)
\oset{\text{\eqref{detgbD}}}= -\frac{\det(\thetab)}{|\gt|^{1/2}}\Bigl(|\thetat^3|_{\gt}\Dc_0 \psih_\nu + \delb{0}|\thetat^3|_{\gt} \psih_\nu\Bigr),
}
and
\leqn{dtmuB}{
\delb{0}\mu^2|_{\Gamma_T} = 2\mt^{\mu\nu}\psi_\mu \Dc_0\psi_\nu,
}
respectively,
where $\Dc_0$ is defined by
\eqn{Dcdef}{
\Dc_0 \zeta_\nu = \delb{0}\zeta_\nu + \beta_\nu^\lambda \zeta_\lambda.
}
From \eqref{pidef} and \eqref{dtmuB}, we then obtain
\eqn{dtpsihA}{
\Dc_0 \psih_\nu|_{\Gamma_T} = \Dc_0\biggl(\frac{1}{\mu}\psi_\nu\biggr)\Bigl|_{\Gamma_T}
= \frac{1}{\mu} \pi_\nu^\omega \Dc_0 \psi_\omega,
}
which in turn, implies that
\leqn{dtpsihB}{
\psih^\nu \Dc_0\psih_\nu |_{\Gamma_T} = 0.
}

Next, differentiating \eqref{acousticB} with respect to $\xb^0$, we see that
\eqn{acousticBb}{
\nml_\alpha \bigl(\At^{\alpha\beta}\delb{\beta} \psi_\nu + L^\alpha_\nu\bigr)|_{\Gamma_T}
= \Dc_0\bigl(\alpha\Dc_0 \psih_\nu + \lambda \psih_\nu\bigr)
}
where
\eqn{acbcdef}{
\alpha = -\frac{\det(\thetab)|\thetat^3|_{\gt}}{|\gt|^{1/2}}, \quad \lambda = -\frac{\det(\thetab)\delb{0}|\thetat^3|_{\gt}}{|\gt|^{1/2}},
}
and it is understood that $|\thetat^3|_{\gt}$ is calculated using the right hand side
of \eqref{thetat3B}.
From this and \eqref{dpsih}, we then get that
\leqn{acousticC}{
\nml_\alpha \biggl(\mu^2 \At^{\alpha\beta}\delb{\beta} \psih_\nu + \frac{\mu^2}{2}\At^{\alpha\beta}\delb{\beta}\mt^{\sigma\gamma}\psih_\sigma\psih_\gamma\psih_\nu+
\mu\pi^\omega_\nu L_\omega^\alpha
\biggr)\Bigl|_{\Gamma_T}
= \mu\pi^\omega_\nu\Dc_0\bigl(\alpha \Dc_0 \psih_\omega + \lambda\psih_\omega\Bigr).
}
But for any $\kappa \in \Rbb$,
\alin{acousticD}{
\pi^\omega_\nu\Dc_0\Bigl(\alpha\Dc_0 \psih_\omega + \lambda \psih_\omega\Bigr)
&= \alpha\pi^\omega_\nu\Dc_0 \Dc_0 \psih_\nu +
(\delb{0}\alpha+\lambda)\pi^\omega_\nu\Dc_0\psih_\omega &&\text{(by \eqref{piprop})} \\
&= \alpha\pi^\omega_\nu\Dc_0 \Dc_0 \psih_\nu +
\biggl((\delb{0}\alpha+\lambda)\delta^\omega_\nu + \frac{\kappa}{\mu}\psih^\omega \psih_\nu \biggr)\Dc_0\psih_\omega
&& \text{(by \eqref{dtpsihB})},
}
and so, we see that
\leqn{acousticEa}{
\nml_\alpha \biggl(\mu^2 \At^{\alpha\beta}\delb{\beta} \psih_\nu + \frac{\mu^2}{2}\At^{\alpha\beta}\delb{\beta}\mt^{\sigma\gamma}\psih_\sigma\psih_\gamma\psih_\nu+
\mu\pi^\omega_\nu L_\omega^\alpha
\biggr)\Bigl|_{\Gamma_T}
= \mu\alpha\pi^\omega_\nu\Dc_0 \Dc_0 \psih_\nu +
\bigl(\mu(\delb{0}\alpha+\lambda)\delta^\omega_\nu + \kappa\psih^\omega \psih_\nu \bigr)\Dc_0\psih_\omega
}
follows from \eqref{acousticC}.

To proceed, we calculate
\alin{acousticF}{
\Dc_0\Dc_0 \psi^\omega &= \delb{0}^2\psih_\omega + 2\beta^\lambda_\omega \Dc_0\psih_\lambda
+\delb{0}\beta^\lambda_\omega \psih_\lambda - \beta^\lambda_\omega\beta^\gamma_\lambda \psih_\gamma \\
& = \delb{0}^2\psih_\omega + 2\beta^\lambda_\omega\bigl(\pi^\delta_\lambda + \psih^\delta
\psih_\lambda\bigr) \Dc_0\psih_\lambda
+\delb{0}\beta^\lambda_\omega \psih_\lambda - \beta^\lambda_\omega\beta^\gamma_\lambda \psih_\gamma  \\
&= \delb{0}^2\psih_\omega + 2\beta^\lambda_\omega\pi^\delta_\lambda\Dc_0\psih_\delta
+\delb{0}\beta^\lambda_\omega \psih_\lambda - \beta^\lambda_\omega\beta^\gamma_\lambda \psih_\gamma
&& \text{(by\eqref{dtpsihB})}.
}
Using this, we can express \eqref{acousticEa} as
\leqn{acousticE}{
\nml_\alpha \biggl(\mu^2 \mt^{\mu\nu}\At^{\alpha\beta}\delb{\beta} \psih_\nu + \frac{\mu^2}{2}\At^{\alpha\beta}\delb{\beta}\mt^{\sigma\gamma}\psih_\sigma\psih_\gamma\psih^\mu+
\mu\pi^{\mu\omega} L_\omega^\alpha
\biggr)\Bigl|_{\Gamma_T}
= q^{\mu\nu}\delb{0}^2\psih_\nu + p^{\mu\nu}\delb{0}\psih_\nu + r^\mu
}
where
\lalin{acousticEvars}{
q^{\mu\nu}  &= \mu\alpha\pi^{\mu\nu}, \label{acousticEvars.1}\\
p^{\mu\nu} &= 2 \mu\alpha\pi^{\mu\omega}\beta^\lambda_\omega\pi^\nu_\lambda + \mu(\delb{0}\alpha+\lambda)\mt^{\mu\nu} + \kappa\psih^\mu \psih^\nu, \label{acousticEvars.2}
\intertext{and}
r^\mu &= \mt^{\mu\nu}\Bigl[\mu\alpha\pi^\omega_\nu\delb{0}\beta^\lambda_\omega \psih_\lambda - \mu\alpha\pi^\omega_\nu\beta^\lambda_\omega\beta^\gamma_\lambda \psih_\gamma
+\bigl(\mu(\delb{0}\alpha+\lambda)\delta^\omega_\nu + \kappa\psih^\omega \psih_\nu \bigr)\beta^\sigma_\omega
\psih_\sigma + 2 \mu\alpha\pi^\omega_\nu\beta^\lambda_\omega\pi^\delta_\lambda\beta_\delta^\sigma \psih_\sigma\Bigr].
\label{acousticEvars.3}
}

\begin{rem} \label{acousticrem}
The boundary conditions \eqref{acousticE} for the wave equation \eqref{waveI} are of acoustic type.\footnote{Acoustic boundary
conditions were first defined in \cite{BealeRosencrans:1974,MorseIngard:1968} and further
analyzed in \cite{Beale:1976}. See also \cite{Gal_et_al:2003} for
more recent work and relations to Wentzell boundary conditions.} Although it is not
obvious at the moment, these boundary conditions can be used to remove the
degeneracy from the Neumann boundary conditions \eqref{PsiBC}, and provide effective
boundary conditions for the wave equation \eqref{waveK}. For details, see Lemmas \ref{coerclem} and \ref{pqlem}.
\end{rem}

\sect{IVBP}{The complete initial boundary value problem}

From the above calculations and results, in particular, \eqref{hC}, \eqref{waveG}, \eqref{mudef}, \eqref{psihdef},
 \eqref{Psidef}, \eqref{waveK}, \eqref{thetatBC},
\eqref{PsiBC}, and \eqref{acousticE}, it follows from a straightforward calculation that,
for any choice of $\kappa,\ep,\delta\in \Rbb$, the triple $\{\phi=(\phi^\mu),\thetat^0=(\theta^0_\mu\circ\phi),\Psi=(\psih_\mu,\mu)^{\tr}\}$
derived from a solution of the Frauendiener-Walton-Euler equations satisfying the assumptions (A.1)-(A.7)
from Section \ref{fwe}, and \eqref{muassump}
solves the (overdetermined) IBVP:
\lalin{cbvpA}{
\delb{\alpha}\bigl(B^{\alpha\beta}\delb{\beta} \thetat^0 + M^\alpha\bigr) &= H
\hspace{3.15cm}\text{in $\Omega_T$,} \label{cbvpA.1}\\
\nml_\alpha \bigl(B^{\alpha\beta}\delb{\beta} \thetat^0 + M^\alpha\bigr) &= K
\hspace{3.15cm} \text{in $\Gamma_T$,}  \label{cbvpA.2} \\
\delb{\alpha}(\Bc^{\alpha\beta}\delb{\beta}\Psi + \Mcal^\alpha) &= \Hc
\hspace{3.15cm}\text{in $\Omega_T$,}  \label{cbvpA.3}\\
\nml_{\alpha}(\Bc^{\alpha\beta}\delb{\beta}\Psi + \Mcal^\alpha) &= Q\delb{0}^2\Psi + P \delb{0}\Psi + \Kc
\hspace{0.5cm} \text{in $\Gamma_T$,} \label{cbvpA.4}\\
\delb{0}\phi^\mu - \fmt_{00}\gt^{\mu\nu}\thetat_\nu^0 &= 0 \hspace{3.30cm}\text{in $\Omega_T$,} \label{cbvpA.5}\\
\delb{0}\thetat_\mu^0 + \beta^\lambda_\nu\thetat^0_\lambda &= \psi_\mu \hspace{3.05cm}\text{in $\Omega_T$,} \label{cbvpA.6}\\
(\phi^0,\phi^\Sigma)&= (0,\text{id}_{\Omega_0}) \hspace{2.2cm}\text{in $\Omega_0$,} \label{cbvpA.7} \\
\thetat^0_\mu &= \theta^0_\mu|_{\Omega_0}
\hspace{2.65cm}\text{in $\Omega_0$,} \label{cbvpA.8} \\
(\Psi,\delb{0}\Psi) &= (\Psi_0,\Psi_1) \hspace{2.2cm}\text{in $\Omega_0$} \label{cbvpA.9}
}
where the initial data $(\Psi_0,\Psi_1)$ is given by
\alin{cbvpBa}{
\Psi_0 &= \biggl(\frac{\nabla_{e_0}\theta^0_\mu}{|\nabla_{e_0}\theta^0|_m},
|\nabla_{e_0}\theta^0|_m\biggr)^{\tr}\Bigl|_{\Omega_0}, %\label{cbvpB.0}
\\
\Psi_1 &= \biggl(\frac{e_0\bigl(\nabla_{e_0}\theta^0_\mu\bigr)|\nabla_{e_0}\theta^0|_m-
\nabla_{e_0}\theta^0_\mu e_0\bigl(|\nabla_{e_0}\theta^0|_m\bigr)}{|\nabla_{e_0}\theta^0|_m^2},e_0
\bigl(|\nabla_{e_0}\theta^0|_m\bigr)\biggr)^{\tr}\Bigl|_{\Omega_0},  %\label{cbvpB.1}
}
the system coefficients are defined by
\lalin{cbvpB}{
B^{\alpha\beta} &= \bigl((\mt^{\mu\nu}+\ep\pi^{\mu\nu})\At^{\alpha\beta}+2\St^{\mu\nu[\alpha}\nml^{\beta]}\bigr)
\qquad \biggr(\nml^\alpha := \frac{\tilde{m}{}^{\alpha\beta}\nml_\beta}{\tilde{m}{}^{\mu\nu}\nml_\mu\nml_\nu}\biggl),  \label{cbvpBa.1} \\
M^{\alpha} &= \Bigl(\bigl(\mt^{\mu\nu}+\ep\pi^{\mu\nu}\bigr)\Lt^\alpha_\nu\Bigr), \label{cbvpBa.2} \\
H &=\Bigl( 2\delb{\alpha}\bigl(\St^{\mu\nu[\alpha}\nml^{\beta]}\bigr)\delb{\beta}\thetat^0_\nu + \bigl(\mt^{\mu\nu}+\ep\pi^{\mu\nu}\bigr)\Ft_\nu + \delb{\alpha}\bigl(\mt^{\mu\nu}+\ep\pi^{\mu\nu}\bigr)\bigl(\At^{\alpha\beta}
\delb{\beta}\thetat^0_\nu + \Lt^\alpha_\nu\bigr)  \Bigr), \label{cbvpBa.3} \\
K&= \Bigl( \ep \alpha\pi^{\mu\nu}\bigl(\delb{0}\psih_\nu + \beta^\lambda_\nu \psih_\lambda\bigr)+
\mt^{\mu\nu}\Zt_\nu \Bigr), \label{cbvpBa.4} \\
\Pbb  &= \begin{pmatrix}\delta^{\mu}_{\nu} & 0 \\ 0 & 0 \end{pmatrix}, \label{cbvpB.2}\\
\Bc^{\alpha\beta} &= (\id + \delta \Pbb)\Ac^{\alpha\beta} + 2\nu^{[\beta}S^{\alpha]}, \label{cbvpB.3}\\
\Mcal^{\alpha} &= (\id + \delta\Pbb)\Lc^\alpha, \label{cbvpB.4} \\
\Hc &= (\id + \delta \Pbb)\Fc + 2\delb{\alpha}\bigl(\nu^{[\beta}S^{\alpha]}\bigr)\delb{\beta}\Psi, \label{cbvpB.5}\\
Q &= \begin{pmatrix}\delta q^{\mu\nu} & 0 \\ 0 & 0 \end{pmatrix},  \label{cbvpB.6} \\
P &= \begin{pmatrix}\delta p^{\mu\nu} & 0 \\ 0 & 0 \end{pmatrix}, \label{cbvpB.7} \\
\Kc &= \begin{pmatrix}\frac{\mu^2}{2} \pi^{\mu}_\alpha\St^{\alpha \beta \omega}\psih_\beta \delb{\omega}\mt^{\sigma\lambda}\psih_\sigma \psih_\lambda
+\mu\pi^{\mu\omega} Z_\omega + \delta r^\mu \\
\psih^\omega Z_\omega \end{pmatrix}, \label{cbvpB.8}
}
and the other quantities are as previously defined:

\smallskip

%\noindent \underline{Eulerian variables}
%\lalin{defrecEul}{
%e_0^\mu &= w^\mu, \label{defrecEul.1} \\
%(\theta^j_\mu) &= (e^\mu_j)^{-1}, \label{defrecEul.2}\\
%\fm_{ij} &= g(e_i,e_j), \label{defrecEul.3}\\
%m^{\alpha\beta} &= g^{\alpha\beta}-\frac{1}{\fm_{00}}e^{\alpha}_0 e^{\beta}_0, \quad \bigl((m_{\alpha\beta}) = %(m^{\alpha\beta})^{-1} \bigr), \label{defrecEul.4}\\
%a^{\mu\nu}&= g^{\mu\nu}-\frac{1}{\fm_{00}}
%\biggl(1-\frac{1}{s^2}\biggr)e^\mu_0 e^\nu_0, \label{defrecEul.5}
%}

\noindent \underline{Lagrangian variables}
\lalin{defrecLagA}{
J &= (J^\mu_\nu), \label{defrecLag.2}\\
J^\mu_\nu &= \delb{\nu}\phi^\mu, \notag %\label{defrecLag.1}
\\
(\Jch^\mu_\nu) &= J^{-1}, \notag %\label{defrecLag.3}
\\
\et^\mu_i &= e^\mu_i \circ \phi,  \notag %\label{defrecLag.4}
\\
(\thetat^i_\mu) &= (\et_i^\mu)^{-1},\notag %  \label{defrecLag.5}
\\
\gt_{\alpha\beta} & = g_{\alpha\beta}\circ \psi, \notag % \label{defrecLag.6}
\\
\Gammat^\gamma_{\alpha\beta} &= \Gamma^\gamma_{\alpha\beta}\circ\phi, \notag %\label{defrecLag.7}
\\
\Rt_{\alpha\beta} &= R_{\alpha\beta}\circ\phi  \notag %\label{defrecLag.8}
\\
\fmt_{ij} &:= \fm_{ij}\circ \phi = \gt_{\mu\nu}\et^\mu_i \et^\nu_j, \notag %\label{defrecLag.9}
\\
\fmt^{ij} &:= \fm^{ij}\circ \phi = \gt^{\mu\nu}\thetat_\mu^i \thetat_\nu^j, \notag % \label{defrecLag.9a}
}
\lalin{defrecLagB}{
\mt^{\alpha\beta}& := m^{\alpha\beta}\circ\phi = \gt^{\alpha\beta}-\frac{1}{\fmt_{00}}\et^{\alpha}_0 \et^{\beta}_0 , \label{defrecLag.12} \\
\at^{\mu\nu} & := a^{\mu\nu}\circ\phi = \gt^{\mu\nu}-\frac{1}{\fmt_{00}}
\biggl(1-\frac{1}{\tilde{s}^2}\biggr)\et^\mu_0 \et^\nu_0, \label{defrecLag.13}\\
\ab^{\alpha\beta} &= \Jch^{\alpha}_\mu\Jch^\beta_\nu\at^{\mu\nu}, \label{defrecLag.14}\\
\tilde{s}^2 &:= s^2\circ \phi = s^2\bigl(\sqrt{-\gammat^{00}}\bigr)  , \notag %\label{defrecLag.10}
\\
\fft &:= \ff\circ \phi = \ff(\sqrt{-\gammat^{00}}), \notag %\label{defrecLag.11}
\\
\psi_\nu &= \bigl(\nabla_{e_0}\theta^0_\nu)\circ \phi, \notag %\label{defrecLag.15}
\\
\mu &= \enorm{\psi}_{\mt}, \notag %\label{defrecLag.16}
\\
\psih_\nu &= \frac{1}{\mu}\psi_\nu,  \notag %\label{defrecLag.17}
\\
\psih^\mu &= \mt^{\mu\nu}\psih_\nu,  \notag %\label{defrecLag.18}
}

\smallskip

\noindent \underline{Projection operator}
\lalin{defrecproj}{
\pi^\lambda_\nu &= \delta^\lambda_\nu - \psih^\nu\psih_\nu, \label{defrecproj.1}\\
\pi^{\mu\nu} &= \mt^{\mu\lambda}\pi_\lambda^\nu, \label{defrecproj.2}
}

\smallskip

\noindent \underline{Bulk coefficients}
\lalin{defrecA}{
\At^{\alpha\beta} &= -\frac{\det(J)}{\fft}\ab^{\alpha\beta},\label{defrec.1} \\
\Ac^{\alpha\beta} &= \begin{pmatrix}\mu^2 \At^{\alpha\beta}\mt^{\mu\nu} & 0 \\ 0 & \At^{\alpha\beta} \end{pmatrix},
\label{defrec.8}\\
\beta^\lambda_\nu &= -\et^\gamma_0\Gammat^\lambda_{\gamma\nu}, \label{defrec.5}
}
\lalin{defrecB}{
\Lt^\alpha_\nu &= -\frac{\det(J)}{\fft}\ab^{\alpha \beta}\Bigl(- J_\beta^\omega\Gammat_{\omega\nu}^\gamma \thetat^0_\gamma
+ \sigmat_M{}^0{}_L \thetab^M_\beta \thetat^L_\nu\Bigr), \notag % \label{defrec.2}
\\
\Ft_\nu &= -\det(J)\biggl(-\Gammat_{\alpha\gamma}^\alpha \Yt^\gamma_\nu +
\Gammat^\gamma_{\alpha\nu}\Yt^\alpha_\gamma + \frac{1}{\ff}\Rt_{\nu}{}^{\lambda}\thetat^0_\lambda\biggr),
\notag %\label{defrec.3}
\\
\Yt^\alpha_\nu & = \frac{1}{\fft}\at^{\alpha \beta}\Bigl[\Jch_\beta^\gamma\delb{\gamma} \thetat^0_\nu -
\Gammat_{\beta\nu}^\gamma \thetat^0_\gamma
+ \sigmat_M{}^0{}_L \thetat^M_\beta \thetat^L_\nu\Bigr], \notag %\label{defrec.4}
\\
L^\alpha_\nu &= \delb{0}\At^{\alpha\beta}\delb{\beta}\thetat^0_\nu+\delb{0}\Lt^\alpha_\nu
+\beta^\omega_\nu \Lt^\alpha_\omega - \At^{\alpha\beta}\delb{\beta}\beta^\omega_\nu \thetat^0_\omega,
\notag %\label{defrec.6}
\\
F_\nu &= \delb{0}\Ft_\nu + \beta^\omega_\nu \Ft_\omega + \delb{\alpha}\beta^\omega_\nu\bigl(\At^{\alpha\beta}
\delb{\beta}\thetat^0_\omega + \Lt^\alpha_\omega\bigr), \notag %\label{defrec.7}
\\
\Lc^{\alpha} &= \begin{pmatrix}\frac{\mu^2}{2} \At^{\alpha\beta}\delb{\beta}\mt^{\sigma\gamma}\psih_\sigma\psih_\gamma\psih^\mu
+ \mu\pi^{\omega\mu} L_\omega^\nu & 0 \\
0 & -\frac{\mu}{2}\At^{\alpha\beta}\delb{\beta}\mt^{\sigma\gamma}\psih_\sigma\psih_\gamma
+\psih^\gamma L^\alpha_\gamma \end{pmatrix}, \notag %\label{defrec.9}
\\
\fc_\nu &= \mu\mt^{\mu\nu}\Bigl(\pi^\omega_\nu F_\omega + \delb{\alpha}\pi^\omega_\nu\bigl(\At^{\alpha\beta}\delb{\beta}\psi_\omega
+L^\alpha_\omega\bigr)\Bigr) \notag \\
&\hspace{2.0cm} + \delb{\alpha}\bigl(\mu\mt^{\mu\nu}\bigr)\Bigl(\mu \At^{\alpha\beta}\delb{\beta}\psih_\nu +
\frac{\mu}{2} \At^{\alpha\beta}\delb{\beta}\mt^{\sigma\gamma}\psih_\sigma\psih_\gamma\psih_\nu
+ \pi^\omega_\nu L_\omega^\nu\Bigr), \notag %\label{defrec.10}
\\
\fc &= \psih^\gamma F_\gamma+ \delb{\alpha}\psih^{\gamma}\bigl(
\At^{\alpha\beta}\delb{\beta}\psi_\gamma+L^\alpha_\gamma\bigr), \notag % \label{defrec.11}
\\
\Fc &= \begin{pmatrix} \fc_\nu \\ \fc \end{pmatrix}, \notag %\label{defrec.12}
}

\smallskip

\noindent \underline{Boundary coefficients}

\lalin{defrecC}{
\St_{\nu\mu}{}^\omega &= \det(\thetab) \bigl(\ep_{\nu\mu\alpha\gamma} \et^\alpha_1 \et^\gamma_2 \eb^\omega_0
+  \ep_{\nu\alpha \mu \gamma} \et^\alpha_0 \et^\gamma_2 \eb^\omega_1
+ \ep_{\nu\alpha \gamma \mu} \et^\alpha_0 \et^\gamma_1 \eb^\omega_2\bigr), \label{defrec.13}\\
\St^{\mu\nu\omega} &= \mt^{\mu\alpha}\St_{\alpha\beta}{}^\omega \mt^{\beta\nu}, \label{defrec.14}\\
\Sc^\omega &= \begin{pmatrix} \mu^2 \pi^\mu_\alpha\St^{\alpha \beta \omega}\pi_\beta^\nu &
\mu\pi^{\mu}_{\alpha}\St^{\alpha \beta \omega}\psih_\beta  \\
 \mu \psih_\alpha \St^{\alpha \beta \omega}\pi_\beta^\nu & 0 \end{pmatrix}, \label{defrec.17}\\
\Zt_{\nu} &=  \St_{\nu\mu}{}^\omega\delb{\omega}(\gt^{\mu\lambda})\thetat^0_\lambda- \St_{\nu\mu}{}^\omega\delb{\omega}
(\gt^{\sigma\lambda})\et^\nu_0 \thetat^0_\sigma \thetat^0_\lambda
- \det(\thetab)\biggl(\frac{1}{2|\gt|}\frac{\partial |\gt|}{\partial \phi^\kappa} \et^\kappa_0
\delta^\lambda_\nu + \beta^\lambda_\nu\biggr)\ep_{\lambda\alpha\beta\gamma}\et^\alpha_0\et^\beta_1\et^\gamma_2,
\notag %\label{defrec.15}
\\
Z_\nu &= -\St_{\nu \mu}{}^\omega \mt^{\mu \gamma}\delb{\omega}\beta_\gamma^\lambda\thetat^0_\lambda
+ \delb{0}\bigl(\St_{\nu \mu}{}^\omega \mt^{\mu \gamma}\bigr)\delb{\omega}\thetat^0_\gamma
+\delb{0}\Zt_{\nu},  \notag %\label{defrec.16}
}
\lalin{defrecD}{
\alpha &= -\frac{\det(\thetab)|\thetat^3|_{\gt}}{|\gt|^{1/2}}, \label{defrec.18}\\
\lambda &= -\frac{\det(\thetab)\delb{0}|\thetat^3|_{\gt}}{|\gt|^{1/2}}, \label{defrec.19}\\
q^{\mu\nu}  &= \mu\alpha\pi^{\mu\nu}, \label{defrec.20}\\
p^{\mu\nu} &= 2 \mu\alpha\pi^{\mu\omega}\beta^\lambda_\omega\pi^\nu_\lambda + \mu(\delb{0}\alpha+\lambda)\mt^{\mu\nu} + \kappa\psih^\mu \psih^\nu,  \label{defrec.21}
\\
r^\mu &= \mt^{\mu\nu}\Bigl[\mu\alpha\pi^\omega_\nu\delb{0}\beta^\lambda_\omega \psih_\lambda - \mu\alpha\pi^\omega_\nu\beta^\lambda_\omega\beta^\gamma_\lambda \psih_\gamma
+\bigl(\mu(\delb{0}\alpha+\lambda)\delta^\omega_\nu + \kappa\psih^\omega \psih_\nu \bigr)\beta^\sigma_\omega
\psih_\sigma + 2 \mu\alpha\pi^\omega_\nu\beta^\lambda_\omega\pi^\delta_\lambda\beta_\delta^\sigma \psih_\sigma\Bigr],
\notag  %\label{defrec.22}
}
where in \eqref{defrec.18} and \eqref{defrec.19}, $|\thetat^3|_{\gt}$ is computed using the formula
\eqref{conprev.9} for $\thetat^3_\mu$ on the boundary $\Gamma_T$,

\bigskip

\noindent \underline{Propagated bulk constraints}

\lalin{conprev}{
\et^\mu_0 &= \fmt_{00}\gt^{\mu\nu}\thetat^0_\nu, \label{conprev.1}\\
\fmt_{0J} & = \fmt^{0J} = 0, \label{conprev.2} \\
\fmt_{00} & = \frac{1}{\fmt^{00}}, \label{conprev.3} \\
(\fmt_{IJ}) & = (\fmt^{IJ})^{-1}, \notag %\label{conprev.4}
\\
F\bigl(\sqrt{-\fmt^{00}}\bigr)^2 &= \frac{1}{\det(\fmt^{IJ})}, \notag % \label{conprev.5}
}
where $F$ is defined by \eqref{Fdef.1}-\eqref{Fdef.2}, and
%\eqn{conprevA}{
%\ft = -\sqrt{-\fmt^{00}}F\bigl(\sqrt{-\fmt^{00}}\bigl)^2.
%}

\smallskip

\noindent \underline{Boundary conditions and propagated boundary constraints}

\lalin{defrecBndry}{
\gammat_{00}|_{\Gammab_T} &=-1, \label{conprev.6} \\
\psih_\mu|_{\Gammab_T} &= -|\fmt^{33}|^{-1/2}\thetat_\mu^3, \label{conprev.7}\\
\ft|_{\Gammab_T} &=-1, \label{conprev.8}\\
\thetat^3_\mu |_{\Gamma_T} &=  -|\gt|^{1/2}\ep_{\mu\alpha\beta\gamma}\et^\alpha_0\et^\beta_1\et^\gamma_2.
\label{conprev.9}
}

\bigskip

For the purposes of deriving energy
estimates, we view \eqref{cbvpA.3}-\eqref{cbvpA.6} as the primary evolution equations, while we treat \eqref{cbvpA.1}-\eqref{cbvpA.2}
as an elliptic constraint for $\thetat^0_\mu$ by using \eqref{cbvpA.6}
to express the time derivatives
of $\thetat^0_\mu$ in terms of the variables $\{\thetat^0,\phi,\psih\}$.
We further remark that in the following, it turns out to be convenient to ``forget'' that
$\psih_\nu$ satisfies $|\psih_\nu|=1$. This necessitates
redefining $\pi^\mu_\nu$, see \eqref{defrec.21}, as
\leqn{pirem1}{
\pi^\mu_\nu  =  \delta^\mu_\nu - \frac{\psih^\mu\psih_\nu}{|\psih|_{\mt}^2}
}
so that it remains a projection operator that agrees with the previous definition \eqref{pidef}
for $\psih_\mu$ satisfying $|\psih|_{\mt}=1$. We also redefine $p^{\mu\nu}$
by
\leqn{pirem2}{
p^{\mu\nu} = 2 \mu\alpha\pi^{\mu\omega}\beta^\lambda_\omega\pi^\nu_\lambda + \mu(\delb{0}\alpha+\lambda)
\mt^{\mu\nu} + \kappa\frac{\psih^\mu \psih^\nu}{|\psih|^2_{\mt}}
}
in order to agree with the previous definition \eqref{acousticEvars.2}.

\subsect{coefsmoothstruc}{Coefficient smoothness and structure}
We make the following observations about the smoothness of the coefficients as functions of their
variables which are a straightforward consequence of the above formulas, the relations (see \eqref{theta3C},
\eqref{ebform}, and \eqref{normalC})
\lgath{funrel}{
\et_i^\mu = \eb_i^\omega \delb{\omega}\phi^\mu, \label{funrel.1} \\
(\eb_i^\mu) =  \begin{pmatrix} 1 & \begin{displaystyle} f^0_I \end{displaystyle}  \\
0 & \begin{displaystyle} -\frac{1}{f^0_0}f^\Lambda_0 f^0_I + f^\Lambda_I \end{displaystyle} \end{pmatrix},
\label{funrel.2}  \\
\nu_\omega \eb_i^\omega |_{\Gamma_T} = 0, \quad i=0,1,2,  \label{funrel.3}
%\nu_\omega \St^{\mu\nu\omega}|_{\Gamma_T} = 0,  \label{funrel.4}
}
and the time-independence of the initial data (see \eqref{etidata.1} and \eqref{tvarsD})
\eqn{tind}{
f^\lambda_I=f^\lambda_I(\xb^\Lambda) \AND \sigma_{i}{}^k{}_{j}=\sigma_{i}{}^k{}_j(\xb^\Lambda).
}

\begin{itemize}
\item[(i)]
From the definitions \eqref{defrecLag.2}, \eqref{defrec.1} and \eqref{defrec.8},
it is not difficult to verify that
the maps
\eqn{AtAc1}{
\At^{\alpha\beta} = \At^{\alpha\beta}_{\nu\mu}\bigl(\phi,J\bigr), \AND
\Ac^{\alpha\beta} = \Ac^{\alpha\beta}\bigl(\phi,J,\mu\bigr)
}
are smooth for $((\phi,J),\mu)\in\Uc\times\Rbb_{> 0}$
where
\eqn{Ucdef}{
\Uc = \Bigl\{\, (\phi,J) \in \Rbb^4\times \Mbb{4} \, \Bigl|\, \det(J) > 0, \;\; -g_{\mu\nu}(\phi)J_0^\mu
J^\nu_0>0, \;\ s^2\Bigl( \bigl(-g_{\mu\nu}(\phi)J_0^\mu
J^\nu_0\bigr)^{-1/2}\Bigr) > 0 \, \Bigr\},
}
and moreover, that
\leqn{AtAc2}{
\At^{\alpha\beta} = \At^{\beta\alpha} \AND \Ac^{\alpha\beta} = \bigl(\Ac^{\beta\alpha}\bigr)^{\tr},
}
respectively, and
for any bounded open subsets $\widetilde{\Uc} \subset \Uc$ and $\Ic\subset\Rbb_{>0}$, there exists a constants
$c_{\At}^0,c_{\Ac}^0 > 0$ such that
\alin{AtAc3}{
\ipe{\Psi}{\Ac^{00}(\phi,J,\mu)\Psi} &\leq -c_{\Ac}^0 |\Psi|^2 && \forall \; \bigl((\phi,J),\mu,\Psi\bigr)\in \widetilde{\Uc} \times \Ic \times \Rbb^5, \\
\At^{00}(\phi,J) &\leq -c_{\At}^0 && \forall \; (\phi,J)\in \widetilde{\Uc} . \\
%\xi_\Sigma\xi_\Lambda\ipe{\Psi}{\Ac^{\Sigma\Lambda}(\phi,J,\mu)\Psi} &\geq c_{\Ac} |\xi|^2|\Psi|^2 && \forall \;
%(\phi,J,\mu,\Psi,\xi)\in \Vc\times \Ic \times \Rbb^5\times \Rbb^3,
%\intertext{and}
%\xi_\Sigma \xi_\Lambda \At^{\Sigma \Lambda}(\phi,J) &\geq c_{\At} |\xi|^2 && \forall \; (\phi,\Phi,\xi,)\in \Vc \times \Rbb^3,
}
Here and below,  we use $\ipe{\cdot}{\cdot}$ and $|\cdot|$ to denote the Euclidean inner product and norm, respectively.
\item[(ii)] From \eqref{AtAc2}, the antisymmetry conditions  from \eqref{Stprop} and \eqref{Scprops}, and the obvious symmetry $\pi^{\mu\nu}=\pi^{\nu\mu}$,
 it is clear from the definitions \eqref{cbvpBa.1} and \eqref{cbvpB.3} that $B^{\alpha\beta}$ and $\Bc^{\alpha\beta}$ satisfy the symmetry conditions
\eqn{BBc}{
(B^{\alpha\beta})^{\tr} = B^{\beta\alpha} \AND \Bc^{\alpha\beta} = \bigl(\Bc^{\beta\alpha}\bigr)^{\tr},
}
respectively.
\item[(iii)]
Setting
\alin{sigmafdef}{
%\sigmat^0 = (\sigmat_I{}^0_J(\xb^\Lambda)), \quad f  = (f^\mu_j(\xb^\Lambda)), \quad
J_0 = (J_0^\mu),
}
%and using
%\eqn{deldef}{
% \Db(\cdot) = \bigr(\delb{1}(\cdot),\delb{2}(\cdot),\delb{3}(\cdot)\bigr) \AND \delb{}(\cdot) = \bigl(\delb{0}(\cdot) , %\Db(\cdot)\bigr)
%}
%to denote the spatial and spacetime gradients, respectively,
and letting
\eqn{delbdef}{
\delb{}(\cdot)=\bigl(\delb{\lambda}(\cdot)\bigr),\quad \lambda=0,1,2,3, \AND \Db(\cdot)=\bigl(\delb{\Lambda}(\cdot)\bigr), \quad \Lambda =1,2,3,
}
denote the spacetime and spatial gradients, respectively,
we see from the definitions \eqref{cbvpBa.1}-\eqref{cbvpBa.3} and \eqref{cbvpB.3}-\eqref{cbvpB.5}
that the maps
\lalin{Fsmooth}{
B^{\alpha\beta} &= B^{\alpha\beta}(\xb,\phi,J,\Psi\bigr), \label{Fsmooth.0a}
\\
\Bc^{\alpha\beta} &= \Bc^{\alpha\beta}(\xb,\phi,J,\Psi\bigr), \label{Fsmooth.0b}
\\
M^\alpha &= M^\alpha\bigl(\xb,\phi,J,\Psi\bigr), \label{Fsmooth.1}
\\
\Mcal^\alpha &= \Mcal^\alpha\bigl(\xb,\phi,J,\delb{} J_0,\Psi\bigr), \label{Fsmooth.2}
\\
H &= H\bigl(\xb,\phi,J,\delb{}J,\Psi,\delb{}\Psi\bigr) \label{Fsmooth.3}
\intertext{and}
\Hc &= \Hc(\xb,\phi,J,\delb{}J,\delb{}\delb{0}J_0,\Psi,\delb{}\Psi\bigr),\label{Fsmooth.4}
}
are well defined and smooth in all variables
provided that\footnote{Here, we are using the notation $\Rbb^4_{\times} = \Rbb^4\setminus \{0\}$.}
$((\phi,J),\Psi)\in \Uc\times (\Rbb\times \Rbb^4_{\times})$. It is also worth noting that the dependence
of the above maps on the spatial coordinates $\xb=(\xb^\Lambda)$ arises via the time-independent
functions $f_I^\mu$ and
$\sigmat_{i}{}^k{}_j$.
\item[(iv)] Letting
\eqn{delslbdef}{
\delslb(\cdot) \AND \Dslb(\cdot)
}
denote the derivatives tangent to the spacetime boundary $\Gamma_T$
and the spatial boundary $\del{}\Omega$, respectively, it is also not difficult to see from
\eqref{cbvpBa.4}, \eqref{cbvpB.6}-\eqref{cbvpB.8}, \eqref{defrec.14}, \eqref{defrec.17},
and the formulas \eqref{conprev.9} and \eqref{funrel.1}-\eqref{funrel.3} that
\lalin{bcsmooth}{
%\alpha & = \alpha\bigl(f,\phi,\delslb\phi\bigr),\\
%\lambda &= \lambda \bigl(f,\phi,\delslb\phi,\delslb J_0\bigr),\\
\St^\alpha &= \bigl(\St^{\mu\nu\alpha}\bigl(\xb,\phi,\delslb\phi\bigr)\bigr),  \label{bcsmooth.1}
\\
\Sc^\alpha &= \Sc^\alpha\bigl(\xb,\phi,\delslb \phi,\Psi\bigr),  \label{bcsmooth.2}
\\
K &= K\bigl(\xb,\phi,\delslb \phi,\Psi,\delb{0}\Psi\bigr),\label{bcsmooth.3}
\\
\Kc &= \lc^\mu_1(\xb,\phi,\delslb\phi,\Psi)\delb{\mu} J_0 + \kc_1\bigl(\xb,\phi,\delslb\phi,\delb{0}J_0,\Psi\bigr),\label{bcsmooth.4}
\\
P &= \lc^\mu_2(\xb,\phi,\delslb\phi,\Psi)\delb{\mu} J_0 + \kc_2\bigl(\xb,\phi,\delslb\phi,\delb{0}J_0,\Psi\bigr), \label{bcsmooth.5}
\intertext{and}
Q &= Q\bigl(\xb,\phi,\delslb\phi,\Psi\bigr) \label{bcsmooth.6}
}
where
\eqn{bcsmoothA}{
\nu_\mu\lc^\mu=0, \quad \lc^0=0,
}
and the maps are smooth
in all their variables.
\end{itemize}

\subsect{timerel}{Fundamental relations under time differentiation}

From \eqref{conprev.1}, \eqref{conprev.3}, \eqref{funrel.1} and \eqref{funrel.2}, we have that
\eqn{timrelA}{
\delb{0}\phi^\mu = \frac{\gt^{\mu\nu}}{\fmt^{00}}\thetat^0_\nu,
}
or schematically,
\leqn{timerelB}{
J_0 = \Jc(\phi,\thetat^0) \qquad \bigl(J_0=(\delb{}\phi^\mu)\bigr)
}
for a smooth map $\Jc$ satisfying $\Jc(\phi,0)=0$. Differentiating this with respect to $\xb^0$, we find, with the help
of the evolution equation \eqref{cbvpA.6}, that
\leqn{timerelC}{
\delb{0}J_0 = \Jf(\phi,\thetat^0,\Psi)
}
for a smooth map $\Jf$ satisfying $\Jf(\phi,0)=0$. Using \eqref{timerelB} and \eqref{timerelC}, we can alternatively express the
functional dependence of the maps \eqref{Fsmooth.1}-\eqref{Fsmooth.4},
\eqref{bcsmooth.4} and \eqref{bcsmooth.5} as
\footnote{Here, we are abusing notation and not changing the name of the maps under the change
of variables. This will cause no difficulties since the precise form of the functions are not required beyond
knowing that they are smooth in their variables on the appropriate domains.}
\lalin{cfsmoothA}{
M^\alpha &= M^\alpha\bigl(\xb,\phi,\thetat^0,\Psi\bigr), \label{cfsmoothA.1}\\
\Mcal^\alpha &= \Mcal^\alpha\bigl(\xb,\phi,\delb{}\phi,\thetat^0,\delb{}\thetat^0,\Psi\bigr),
 \label{cfsmoothA.2} \\
H &= H\bigl(\xb,\phi,\bar{D}\phi,\bar{D}^2\phi, \thetat^0,\delb{}\thetat^0,\Psi,\delb{}\Psi\bigr),
 \label{cfsmoothA.3}\\
\Hc &= \Hc(\xb,\phi,\bar{D}\phi,\bar{D}^2\phi,\thetat^0,\delb{}\thetat^0,\Psi,\delb{}\Psi\bigr),
 \label{cfsmoothA.4}\\
\Kc &= \lc^\mu_1(\xb,\phi,\delslb\phi,\Psi)\delb{\mu}\bigl(\Jc(\phi,\thetat^0)\bigr) +
\kc_1\bigl(\xb,\phi,\delslb\phi,\thetat^0,\Psi\bigr),  \label{cfsmoothA.5}
\intertext{and}
P &= \lc^\mu_2(\xb,\phi,\delslb\phi,\Psi)\delb{\mu}\bigl(\Jc(\phi,\thetat^0)\bigr) +
\kc_2\bigl(\xb,\phi,\delslb\phi,\thetat^0,\Psi\bigr). \label{cfsmoothA.6}
}

Differentiating $\Bc^{\alpha\beta}$, $\Mcal$, $Q$, $\kc$, and $\lc^\mu$ with respect to $\xb^0$, we also observe with the help
of \eqref{cbvpA.6} and \eqref{timerelB} that
\lalin{cfdsmoothA}{
\delb{0}\lc^\mu_a &= \dot{\lc}^\mu_a\big(\xb,\phi,\delslb\phi,\thetat^0,\delslb\thetat^0,\Psi,\delb{0}\Psi\bigr), \quad a=1,2, \label{cfdsmoothA.1} \\
\delb{0}\kc_a &= \dot{\kc}_a\bigl(\xb,\phi,\delslb\phi,\thetat^0,\delslb\thetat^0,\Psi,\delb{0}\Psi\bigr),\quad a=1,2, \label{cfdsmoothA.2} \\
\delb{0}Q  &= \dot{Q}\bigl(\xb,\phi,\delslb\phi,\thetat^0,\delslb\thetat^0,\Psi,\delb{0}\Psi\bigr), \label{cfdsmoothA.3} \\
\delb{0} B^{\alpha\beta} &=  \dot{B}^{\alpha\beta}(\xb,\phi,\delb{}\phi,\thetat^0,\Psi,\delb{0}\Psi), \label{cfdsmoothA.4} \\
\delb{0} \Bc^{\alpha\beta} &= \dot{\!\Bc}^{\alpha\beta}(\xb,\phi,\delb{}\phi,\thetat^0,\Psi,\delb{0}\Psi) , \label{cfdsmoothA.5}
\intertext{and}
\delb{0} \Mcal^\alpha &= \dot{\!\!\Mcal}^\alpha (\xb,\phi,\delb{}\phi,\thetat^0,\delb{}\thetat^0,\Psi,\delb{}\Psi),
\label{cfdsmoothA.6}
}
where the maps $\dot{\kc}_a$, $\dot{\lc}^\mu_a$, $\dot{Q}$, $\dot{B}^{\alpha\beta}$, $\dot{\!\Bc}^{\alpha\beta}$
 and $\;\dot{\!\!\Mcal}^{\alpha}$ are smooth in their arguments.

%\begin{rem}
%The particular manner in which the above maps depend on the variables %$(\phi,\delb{}\phi,\thetat^0,\delb{}\thetat^0,\Psi,\delb{}\Psi)$ is crucial, because, as we
%shall see in Section \ref{nlimIBVP}, it allows us to use the linear estimates from Theorem \ref{linlocthmAaa}
%to derive a prior estimates for solutions of the \eqref{cbvpA.1}-\eqref{cbvpA.9}.
%\end{rem} 
\sect{linIBVP}{Linear wave equations}
With our wave formulation complete, we now turn in this section to the problem
of establishing the existence, uniqueness, and energy estimates for solutions to linear equations
that include equations of the form \eqref{cbvpA.3}-\eqref{cbvpA.4}.

\subsect{prelim}{Preliminaries} Before proceeding, we first introduce some notation and fix our conventions that will
be used throughout this section. Unlike the previous sections, we work here in arbitrary dimensions.

\subsubsect{notation}{Notation} We use $(x^{\mu})_{\mu=0}^n$  to denote Cartesian coordinates
on $\Rbb^{n+1}$, and we use  $x^0$ and $t$, interchangeably, to denote the time coordinate,
and $(x^i)_{i=1}^n$ to denote spatial coordinates. We also use $x=(x^1,\ldots,x^n)$
%and $\xv = (x^0,\ldots,x^n)$
to denote spatial points.
%and spacetime points, respectively.

As before, partial derivatives are denoted by
\eqn{partial}{
\del{\mu} = \frac{\partial \;}{\partial x^\mu},
}
and
we use $Du(x) = (\del{1}u(x),\ldots,\del{n}u(x))$ and
$\del{}u(t,x) =  (\del{0}u(t,x),Du(t,x))$ to denote the spatial and spacetime gradients, respectively. For time derivatives, we often employ the notation
\leqn{ft}{
u_\ell = \del{t}^\ell u,
}
and use
\leqn{fvect}{
\uv_r = (u_0, u_1, \ldots, u_r )^{\tr}
}
to denote the collection of partial derivatives of $u$ with respect to $t$.

The notation \eqref{fvect} will be employed more generally for vectors with components $u_\ell$, $0\leq \ell \leq r$,
that are not necessarily of the type \eqref{ft}. It will be clear from the context whether $u_\ell$
represents some general vector component, or is of the type \eqref{ft}. We will find it convenient to use the notation
\leqn{dtfvect}{
\del{t}^\ell \uv_r = (u_\ell,u_{\ell+1},\ldots,u_{\ell+r})^{\tr}
}
on occasion, where again the $u_{\ell+j}$ can be thought of as either the time derivative $\del{t}^{\ell+j}u$ of some function $u$, or just a
component of a vector.

\subsubsect{funct}{Function spaces}

$\;$

\bigskip

\noindent\textit{Spatial function spaces}: In the following, we let $\Omega$ denote
a bounded, open set in $\Rbb^n$ with a $C^\infty$ boundary, and employ similar notation
when $\Omega$ is replaced by a smooth, closed $n$-dimensional manifold. Given a finite dimensional vector space $V$, we let
$W^{s,p}(\Omega,V)$, $s\in \Rbb$ $(s\in \Zbb_{\geq 0})$, $1< p < \infty$ $(1\leq p\leq \infty)$,
denote the space of $V$-valued maps on $\Omega$ with fractional (integral) Sobolev regularity $W^{s,p}$. Particular
cases of interest for $V$ will be $V=\Rbb^N$ and $V=\Mbb{N}$, where, here, we use
$\Mbb{N}$ to denote the set of $N$-by-$N$ matrices. In the special case of $V=\Rbb$,
we employ the more compact notation $W^{s,p}(\Omega)=W^{s,p}(\Omega,\Rbb)$.

When $p=2$, we use the standard
notation $H^s(\Omega,V)=W^{s,2}(\Omega,V)$, and on $L^2(\Omega,\Rbb^N)$, we denote the inner product by
\eqn{ipdef}{
\ip{u}{v}_{\Omega} = \int_{\Omega} \ipe{u(x)}{v(x)}\, d^n x \qquad u,v\in L^2(\Omega,\Rbb^N)
}
where, as previously,
\eqn{eipdef}{
\ipe{\xi}{\zeta} = \xi^{\tr}\zeta,\qquad \xi,\zeta \in \Rbb^N,
}
is the Euclidean inner product on $\Rbb^N$. We will also have occasion to use the standard inner-product
on $H^s(\Omega)$, which we denote by $\ip{\cdot}{\cdot}_{H^s(\Omega)}$.

Given $s=k/2$, for $k\in \Zbb_{\geq 0}$, we define the spaces
\leqn{XdefA}{
X^{s,r}(\Omega,V) = \prod_{\ell=0}^r H^{s-\frac{\ell}{2}}(\Omega,V)
}
for $0\leq r \leq 2s$, and set
\leqn{XTdefBa}{
X^s(\Omega,V) = X^{s,2s}(\Omega,V).
}
Using the vector notation \eqref{fvect}, we can write the norms for the spaces \eqref{XdefA} and \eqref{XTdefBa} as
\eqn{Xnorm}{
\norm{\uv_r}^2_{X^{s,r}} = \sum_{\ell=0}^r \norm{u_\ell}_{H^{s-\frac{\ell}{2}}(\Omega)}^2 \AND \norm{\uv_{2s}}_{X^{s}} = \norm{\uv_{2s}}_{X^{s,2s}},
}
respectively.

\bigskip

\noindent\textit{Spacetime function spaces:} Given $T>0$, and $s=k/2$ for $k\in \Zbb_{\geq 0}$, we define the spaces
\leqn{XTdefA}{
X^{s,r}_T(\Omega,V) = \bigcap_{\ell=0}^{r} W^{\ell,\infty}\bigl([0,T],H^{s-\frac{\ell}{2}}(\Omega,V)\bigr),
%C^{\ell}\bigl([0,T),\Hc^{m_{s-\ell},s-\ell}(\Tbb^n)\bigr),
}
for $0\leq r \leq 2s$,
\leqn{XTdefB}{
X^s_T(\Omega,V) = X^{s,2s}_T(\Omega,V),
}
%\leqn{XcTdefA}{
%\Xc^{s,2s-2}_T(\Omega,V) =  \bigcap_{\ell=0}^{2s-2} W^{\ell,\infty}\bigl([0,T],H^{s-\frac{m_\ell}{2}}(\Omega,V)\bigr),
%}
and
\leqn{XcTdefB}{
\Xc^{s}_T(\Omega,V) = \bigcap_{\ell=0}^{2s-1} W^{\ell,\infty}\bigl([0,T],H^{s-\frac{m_\ell}{2}}(\Omega,V)\bigr)
}
where
\leqn{melldef}{
m_\ell = \begin{cases} \ell & \text{if $0\leq \ell \leq 2s-2$} \\
2s &\text{if $\ell = 2s-1$}
\end{cases}.
}

\begin{rem} \label{Xcrem}
$\;$
\begin{itemize}
\item[(i)]
The spaces $\Xc^{s}_T(\Omega,V)$ consist of time-dependent
functions with spatial regularity of $H^s(\Omega)$ that lose $1/2$ a derivative of spatial regularity
for each time derivative taken until the
spatial regularity $H^1(\Omega)$ is reached after which the last time derivative
reduces the spatial regularity by $1$ leaving highest time derivative in
$L^2(\Omega)$. This type of regularity is necessary for
 the class of wave equations with acoustic boundary conditions that we consider in
 the following sections; see, in particular, Section \ref{simple} for a simplified example
 that illustrates the need for these spaces.
\item[(ii)]
The $s$ in the definition of $m_\ell$ given by \eqref{melldef} refers to the space
defined by \eqref{XcTdefB}. Thus if $s$ is replaced by $s+1$ in \eqref{XcTdefB},
then $s$ must also be replaced by $s+1$
in the definition of $m_\ell$, that is
\eqn{melldefsp1}{
m_\ell = \begin{cases} \ell & \text{if $0\leq \ell \leq 2s$} \\
2s+2 & \text{if $\ell = 2s+1$}
\end{cases}.
}
\end{itemize}
\end{rem}

We define the following \emph{energy norms}:
\gath{XTnormA}{
\norm{u}_{E^{s,r}}^2 = \sum_{\ell=0}^r \norm{\del{t}^\ell u}^2_{H^{s-\frac{\ell}{2}}(\Omega)}, \quad\norm{u}^2_{E^s} =  \norm{u}_{E^{s,2s}}^2, \AND
%\norm{u}_{\Ec^{s,2s-2}}^2 = \sum_{\ell=0}^{2s-2} \norm{\del{t}^\ell u}^2_{H^{s-\frac{m_\ell}{2}}(\Omega)} \AND
\norm{u}_{\Ec^{s}}^2  = \sum_{\ell=0}^{2s-1} \norm{\del{t}^\ell u}^2_{H^{s-\frac{m_\ell}{2}}(\Omega)}.
}
In terms of these energy norms, we can write the norms of the spaces \eqref{XTdefA}-\eqref{XcTdefB} as
\gath{XTnormB}{
\norm{u}_{X^{s,r}_T} = \sup_{0\leq t \leq T} \norm{u(t)}_{E^{s,r}}, \quad
\norm{u}_{X^s_T} = \sup_{0\leq t \leq T} \norm{u(t)}_{E^s},
%\norm{u}_{\Xc^{s,2s-2}_T}  = \sup_{0\leq t \leq T} \norm{u(t)}_{\Ec^{s,2s-2}}
\AND
\norm{u}_{\Xc^{s}_T} =  \sup_{0\leq t \leq T} \norm{u(t)}_{\Ec^s},
}
respectively. We also define the
subspace
\alin{XVdef}{
C\Xc^s_T(\Omega,V) &= \bigcap_{\ell=0}^{2s-1} C^{\ell}\bigl([0,T],H^{s-\frac{m_\ell}{2}}(\Omega,V)\bigr).
}

We conclude with the following elementary, but useful, relations
\alin{eleminq}{
\norm{u}_{E^{s_1,r_1}} &\leq \norm{u}_{E^{s_2,r_2}}, &&s_1\leq s_2, \; r_1\leq r_2, \\
\norm{u}_{\Ec^{s_1}}&\leq \norm{u}_{\Ec^{s_2}}, && s_1\leq s_2,\\
\norm{u}_{\Ec^{s}}^2 &= \norm{u}_{E^{s,2s-2}}^2 + \norm{\del{t}^{2s-1}u}_{L^2(\Omega)}^2,\\
\norm{u}_{E^{s,r}}^2& = \norm{u}_{H^{s}(\Omega)}^2 + \norm{\del{t}u}^2_{E^{s-\frac{1}{2},r-1}},\\
\norm{\del{t}u}_{E^{s,r}} &\leq \norm{u}_{E^{s+\frac{1}{2},r+1}}, \\
\norm{D u}_{E^{s,r}}  &\leq \norm{u}_{E^{s+1,r}}, \\
\norm{\del{t}u}_{\Ec^s} &\leq \norm{u}_{\Ec^{s+\frac{1}{2}}},
\intertext{and}
\norm{Du}_{\Ec^s} &\leq \norm{u}_{E^{s+1,2s-1}} \leq \norm{u}_{\Ec^{s+1}}.
}

\subsubsect{cost}{Estimates and constants}

We employ that standard notation
\eqn{lesssimA}{
a \lesssim b
}
for inequalities of the form
\eqn{lesssimB}{
a \leq C b
}
in situations where the precise value or dependence on
other quantities of the constant $C$ is not required.  On the other hand,  when the dependence of the constant
on other inequalities needs to be specified, for example if the constant depends on the norms $\norm{u}_{L^\infty(\Tbb^n)}$ and $\norm{v}_{L^\infty(\Omega)}$, we use the notation
\eqn{lesssimC}{
C = C(\norm{u}_{L^\infty(\Tbb^n)},\norm{v}_{L^\infty(\Omega)}).
}
Constants of this type will always be non-negative, non-decreasing, continuous functions of their arguments.

\subsect{modlin}{A model class of linear wave equations}
Rather than directly considering linear wave equations that include equations of
the form \eqref{cbvpA.3}-\eqref{cbvpA.4},  we instead consider a related model class
of equations for which it is easier to establish an existence and uniqueness result.
The desired existence and uniqueness result will then follow from this one.

The model class that we consider are wave equations of the form:
\lalin{linB}{
\del{\alpha}\bigl(b^{\alpha\beta}\del{\beta}v+L^\alpha\bigr) +\lambda c v &= F \hspace{2.85cm}
\text{in $\Omega_T=[0,T]\times \Omega$,}
\label{linB.1}\\
\nu_\alpha\bigl(b^{\alpha\beta}\del{\beta} v + L^\alpha) & = q\del{t}^2 v + P \del{t} v + G
 \hspace{0.5cm} \text{in $\Gamma_T=[0,T]\times\del{}\Omega$, } \label{linB.2}\\
(v,\del{t}v) &= (\vt_0,\vt_1) \hspace{2.0cm} \text{in $\Omega_0$,}\label{linB.3}\\
\Pbb_{q}\del{t}v & = \wt_1 \hspace{2.70cm} \text{in $\Gamma_0$}\label{linB.4}
}
where
\begin{enumerate}[(i)]
\item $\lambda \in \Rbb$,
\item $\Omega\subset \Rbb^n$ is open and bounded with smooth boundary,
\item $\nu_\alpha = \delta^i_\alpha \nu_i$ where $\nu_i$ is the outward pointing unit normal to $\del{}\Omega$,
\item $v=v(t,x)$ is a $\Rbb^N$-valued map,
\item $L=(L^\alpha) \in W^{1,2}\bigl([0,T],L^2(\Omega,\Rbb^N)\bigr)$ and
$F \in L^2(\Omega_T,\Rbb^N)$,
\item  the matrix valued maps $c\in L^{\infty}\bigl([0,T],L^n(\Omega,\Mbb{N})\bigr)$,
$P \in W^{1,\infty}(\Omega_T,\Mbb{N})$, and $q \in W^{1,\infty}(\Omega_T,\Mbb{N})$
satisfy $\del{t}q \in W^{1,\infty}(\Omega_T,\Mbb{N})$,\footnote{Given $A,B\in \Mbb{N}$, we define $A\leq B$ $\Longleftrightarrow$ $\ipe{\xi}{A\xi} \leq
\ipe{\xi}{B\xi}$ for all $\xi\in \Rbb^N$.}
\lgath{qPcdef}{
q^{\tr} = q, \quad q \leq 0, \quad -\frac{1}{\gamma} q \leq q^2 \leq -\gamma q,
\quad \text{rank}(q)= N_q \hspace{0.4cm} \text{in $\Gamma_T$}, \label{qPdef}
}
and
\leqn{cdef}{
 c\leq -\sigma\hspace{0.4cm}\text{in $\Omega_T$}
 }
for some positive constants $\sigma, \gamma>0$ and $N_q \in \{0,1,\ldots,N\}$,
\item $\Pbb_q$ is the projection onto $\text{ran}(q)$,
\item  the matrix valued maps $b^{\alpha\beta} \in W^{1,\infty}\bigl([0,T],L^\infty(\Omega,\Mbb{N})\bigr)$
satisfy
\lgath{bsym}{
(b^{\alpha\beta})^{\tr} = b^{\beta\alpha} \label{bsym.1}
\intertext{and}
 b^{00} \leq -\kappa_0, \label{bsym.2}
}
in $\Omega_T$ for some positive constant $\kappa_0>0$,
\item there exists constants $\kappa_1 >0$ and $\mu \geq 0$ such that
\leqn{coerc}{
\ip{\del{i}v}{b^{ij}(t)\del{j}v}_{\Omega} \geq \kappa_1 \norm{v}^2_{H^1(\Omega)} -\mu \norm{v}^2_{L^2(\Omega)}
}
for all $v\in H^1(\Omega)$ and $t\in [0,T]$,
\item and the vector valued map $G$ is of the form
\leqn{GformA}{
G = k^i \del{i}\theta + g
}
were $k^i \in W^{1,\infty}\bigl([0,T],W^{1,\infty}(\Omega,\Mbb{N})\bigr)$,
$\theta\in W^{1,\infty}\bigl([0,T],H^1(\Omega,\Rbb^N)\bigr)$,
$g\in H^{1}(\Omega_T,\Rbb^N)$ and
\eqn{GformB}{
\nu_i k^i =0.
}
\end{enumerate}

\begin{rem} \label{coercrem}$\;$

\begin{enumerate}[(i)]
\item The coercive condition \eqref{coerc} is known to be equivalent to the matrix $b^{ij}$
being strongly elliptic at each point in $\overline{\Omega}$ and satisfying the strong complementing
condition for every point on the boundary $\del{}\Omega$. For a proof of this equivalence, see
Theorem 3 in Section 6 of \cite{SimpsonSpector:1987}.
\item Letting $\Pbb_q^\perp$ denote the
projection onto the orthogonal complement of $\text{ran}(q)$, we can decompose $\Rbb^N$ as
\eqn{RNdecomp}{
\Rbb^N = \Pbb_q\Rbb^N\!\oplus \Rbb_q^\perp\Rbb^N.
}
It then follows from the assumptions \eqref{qPdef} on the matrix $q$ that
\eqn{qdecompB}{
q = \Pbb_q q \Pbb_q
}
and
\leqn{qdecompC}{
-\frac{1}{\gamma} \ipe{\xi}{q\xi} \leq \ipe{q\xi}{q\xi} \leq  -\gamma \ipe{\xi}{q\xi}, \qquad \forall \, \xi \in \Rbb^N,
}
which imply that the map
\leqn{qdecompD}{
q|_{\Pbb_q\Rbb^N} \: : \: \Pbb_q\Rbb^N \longrightarrow \Pbb_q\Rbb^N
}
is invertible and satisfies
\leqn{qdecompE}{
|(q|_{\Pbb_q\Rbb^N})^{-1}\xi| \leq \frac{1}{\gamma}|\xi|, \qquad \forall\, \xi \in \Pbb_q\Rbb^N.
}
It also follows directly from \eqref{qdecompC} that the inequality
\leqn{qnorm}{
\frac{1}{\sqrt{\gamma}}\ip{w}{(-q(t))w}_{\del{}\Omega}^{\frac{1}{2}}
\leq \norm{q(t)w}_{L^2(\del{}\Omega)} \leq \sqrt{\gamma}\ip{w}{(-q(t))w}_{\del{}\Omega}^{\frac{1}{2}}
}
holds on the boundary for each $t\in [0,T]$ and all $w\in L^2(\del{}\Omega)$.
\end{enumerate}
\end{rem}

\subsect{linweak}{Weak solutions}

\begin{Def} \label{weaksoldef}
A pair $(v,w)\in H^1(\Omega_T,\Rbb^N)\times L^2(\Gamma_T,\Rbb^N)$ is called a \emph{weak solution} of \eqref{linB.1}-\eqref{linB.4} if  $(v,w)$ define maps $v \, : \, [0,T]\rightarrow H^1(\Omega,\Rbb^N)$, $\del{t}v \, :
\, [0,T]\rightarrow L^2(\Omega,\Rbb^N)$ and $w \,:\, [0,T] \longrightarrow L^2(\Gamma,\Rbb^N)$ that
satisfy\footnote{We use the standard notation
$\rightharpoonup$ to denote weak convergence.}
\eqn{weaksoldef1}{
(u(t),\del{t}v(t)) \rightharpoonup (\vt_0,\vt_1) \quad\text{in $H^1(\Omega,\Rbb^N)\times L^2(\Omega,\Rbb^N)$},  \AND w(t) \rightharpoonup \wt_1 \quad\text{in $L^2(\Gamma,\Rbb^N)$ }
}
as $t\searrow 0$,
\eqn{weaksoldef2}{
w \in \text{ran}(q)\hspace{0.4cm} \text{in $\Gamma_T$,}
}
and
\leqn{weaksoldef2a}{
\ip{qw}{\phi}_{\Gamma_T} = - \ip{\del{t}q v}{\phi}_{\Gamma_T}-\ip{q v}{\del{t}\phi}_{\Gamma_T}
}
and
\leqn{weaksoldef3}{
\ip{b^{\alpha\beta}\del{\beta}v + L^\alpha}{\del{\alpha} \phi}_{\Omega_T} +
\ip{(\del{t}q-P)\del{t}v}{\phi}_{\Gamma_T}-\ip{g}{\phi}_{\Gamma_T}-\ip{k^i\del{i}\theta}{\phi}_{\Gamma_T}+\ip{qw}{\del{t}\phi}_{\Gamma_T}= \ip{\lambda c v -F}{\phi}_{\Omega_T}
}
for all $\phi \in C^1_0\bigl([0,T],C^1(\overline{\Omega},\Rbb^N)\bigr)$.
\end{Def}

\begin{rem} \label{weakrem}$\;$

\begin{itemize}
\item[(i)]
As in \cite{Koch:1993}, the boundary
terms $\ip{g}{\phi}_{\Gamma_T}$ and $\ip{(\del{t}q-P)\del{t}v}{\phi}_{\Gamma_T}$ are defined via the
expressions\footnote{For sufficiently differentiable vector valued and
matrix valued maps $\{v,\phi\}$ and $S$, respectively, such that $\phi|_t=\phi|_{t=T}=0$, the identities
\alin{weakremfoot1}{
\ip{v}{\phi}_{\Gamma_T} &= \int_{\Omega_T}\del{\alpha}\bigl[ \nu^\alpha \ipe{v}{\phi}\bigr]\, d^{n+1} x =
\ip{\nu(v)+\del{i}\nu^i v}{\phi}_{\Omega_T} + \ip{v}{\nu(\phi)}_{\Omega_T}, \\
0&=\int_{\Omega_T}\del{\beta}\bigl[\delta^\beta_0\ipe{S\nu(v)}{\phi}\bigr]\, d^{n+1} x =
\ip{\del{t}S\nu(v)+S\del{t}\nu(v)}{\phi}_{\Omega_T}+\ip{S\nu(v)}{\del{t}\phi}_{\Omega_T}
}
follow from the divergence theorem. The second identity together with one more application of
the divergence theorem then yields
\eqn{weakremfoot2}{
\ip{S\del{t}v}{\phi}_{\Gamma_T} = \int_{\Omega_T}\del{\beta}\bigl[\nu^\beta\ipe{S\del{t}v}{\phi}\bigr]\, d^{n+1} x = \ip{\nu(S)\del{t}v-\del{t}S\nu(v) +
\del{\alpha}\nu^{\alpha}
S\del{t}v}{\phi}_{\Omega_T}+ \ip{S\del{t}v}{\nu(\phi)}_{\Omega_T}-
\ip{S\nu(v)}{\del{t}\phi}_{\Omega_T}.
}
}
\leqn{weakrem2a}{
\ip{g}{\phi}_{\Gamma_T} = \ip{\nu(g)+\del{i}\nu^i g}{\phi}_{\Omega_T} + \ip{g}{\nu(\phi)}_{\Omega_T},
}
and
\lalin{weakrem2}{
&\ip{(\del{t}q-P)\del{t}v}{\phi}_{\Gamma_T} = \ip{\nu(\del{t}q-P)\del{t}v-\del{t}(\del{t}q-P)\nu(v)}{\phi}_{\Omega_T} \notag \\
 &\hspace{0.2cm} +
\ip{\del{\alpha}\nu^{\alpha}
(\del{t}q-P)\del{t}v}{\phi}_{\Omega_T}
+ \ip{(\del{t}q-P)\del{t}v}{\nu(\phi)}_{\Omega_T}-
\ip{(\del{t}q-P)\nu(v)}{\del{t}\phi}_{\Omega_T}, \label{weakrem2.1}
}
respectively,
where $\nu(\cdot) = \nu^\alpha\del{\alpha}(\cdot)$, $\nu^\alpha = \delta^{\alpha i}\nu_i$, and
$\nu_i$ is any smooth extension to $\Omega$ of the outward pointing unit normal to $\del{}\Omega$.
\item[(ii)] The condition \eqref{weaksoldef2a}
implies that $v$ \textit{weakly} satisfies
\eqn{weakrem1}{
q w = q \del{t} v \hspace{0.4cm} \text{in $\Gamma_T$,}
}
where here, we are again defining the
boundary terms on the right hand side of \eqref{weaksoldef2a} using the same type of formula
as \eqref{weakrem2a}.
\item[(iii)] The interpretation of the boundary term $\ip{k^i\del{i}\theta}{\phi}_{\Gamma_T}$ is
more involved. First, we let
\eqn{weakrem3}{
\perp^i_j = \delta^i_j - \nu^i\nu_j
}
denote the orthogonal projection onto the the subspace orthogonal to the normal vector $\nu_i$,
and we set
\eqn{weakrem4}{
\Dsl_i = \perp_i^j\del{j},
}
which defines a complete collection of derivatives tangent to $\del{}\Omega$. This allows us to write
\eqn{weakrem5}{
k^i\del{i} = k^i\Dsl_{i}
}
since
$k^i = \perp^i_j k^j$ by assumption.
Letting $\Deltasl$ denote the Laplacian arising from the restriction of the Euclidean metric on $\Rbb^n$
to $\del{}\Omega$, we then obtain
\lalin{weakrem6}{
\bigl|\ip{\phi}{k^i\del{i}\theta}_{\del{}\Omega}\bigr| &= \bigl|\ip{\phi}{k^i \Dsl_i\theta}_{\del{}\Omega}\bigr| \notag \\
&= \bigl|\ip{(k^i)^{\tr}\phi}{\Dsl_i\theta}_{\del{}\Omega}\bigr|\notag \\
&=\bigl|\ip{\Dsl_i\bigl((k^i)^{\tr}\phi\bigr)}{\theta}_{\del{}\Omega} \bigr|\notag \\
&=\bigl|\ip{(1-\Deltasl)^{\frac{1}{4}}(1-\Deltasl)^{-\frac{1}{4}}\Dsl_i\bigl((k^i)^{\tr}\phi\bigr)}{\theta}_{\del{}\Omega}
\bigr|\notag\\
&=\bigl|\ip{(1-\Deltasl)^{-\frac{1}{4}}\Dsl_i\bigl((k^i)^{\tr}\phi\bigr)}{(1-\Deltasl)^{\frac{1}{4}}\theta}_{\del{}\Omega}
\bigr|\notag\\
&\lesssim \norm{(k^i)^{\tr}\phi}_{H^{\frac{1}{2}}(\del{}\Omega)}
\norm{\theta}_{H^{\frac{1}{2}}(\del{}\Omega)}  \notag\\
& \lesssim \norm{(k^i)^{\tr}\phi}_{H^1(\Omega)}\norm{\theta}_{H^1(\Omega)}  \label{weakrem6.1}
}
where, in deriving the last four lines, we used the self-adjointness and well-known mapping properties of
the pseudo-differential operator $(1-\Deltasl)^s$, the Cauchy-Schwartz inequality, and the Trace Theorem, see Theorem \ref{trace}. Using \eqref{weakrem6.1}, we see with the
help of \eqref{weakremA3} that
\leqn{weakrem7}{
\bigl|\ip{\phi}{k^i\del{i}\theta}_{\del{}\Omega}\bigl| \lesssim \norm{k}_{L^\infty(\Omega)\cap W^{1,n}(\Omega)} \norm{\phi}_{H^1(\Omega)}\norm{\theta}_{H^1(\Omega)},
}
which we can use to obtain the estimate
\eqn{weakrem8}{
\bigl|\ip{\phi}{k^i\del{i}\theta}_{\Gamma_T}\bigr|
%\leq  \norm{k}_{W^{1,\infty}(\Omega_T)} \int_0^T
%\norm{\phi(t)}_{H^1(\Omega)}\norm{\theta(t)}_{H^{1}(\Omega)}\, dt
\leq  \norm{k}_{L^\infty([0,T],L^\infty(\Omega)\cap W^{1,n}(\Omega))}  \norm{\phi}_{H^{1}(\Omega_T)}
\norm{\theta}_{H^1(\Omega_T)}.
}
This estimate implies the continuity of the bilinear map
\eqn{weakrem9}{
H^1(\Omega_T,\Rbb^n)\times H^1(\Omega_T,\Rbb^n) \ni (\phi,\theta) \longmapsto \ip{\phi}{k^i\del{i}\theta}_{\Gamma_T}
\in \Rbb,
}
and gives meaning to the boundary term $\ip{k^i\del{i}\theta}{\phi}_{\Gamma_T}$.
\end{itemize}
\end{rem}

When $q=0$, $k^i=0$, and $P\leq 0$, the existence and uniqueness of weak solutions to the
IBVP \eqref{linB.1}-\eqref{linB.3} (in this case, $w=0$ and \eqref{linB.4} is redundant)
is a consequence of Theorem 2.2 from \cite{Koch:1993}. Using similar arguments, we
establish the following generalization.

\begin{thm} \label{weakthm}
Suppose $\vt_0\in H^1(\Omega,\Rbb^N)$, $\vt_1\in L^2(\Omega,\Rbb^N)$, $\wt_1 \in L^2(\del{}\Omega,\Rbb^N)$
 such that $\wt_1\in \text{\emph{ran}}(q|_{t=0})$,
the assumptions (i)-(x) from Section \ref{modlin} are fulfilled, \eqn{weakthm2a}{
P-\Half \del{t}q-\chi q \leq 0 \hspace{0.5cm}\text{in $\Gamma_T$,}
}
and let
\eqn{weakthm2ab}{
\kf = \norm{k}_{W^{1,\infty}([0,T],L^\infty(\Omega)\cap W^{1,n}(\Omega))}.
}
Then
there exists a unique weak solution $(v,w)$ to the IBVP \eqref{linB.1}-\eqref{linB.4},
and this solution satisfies the energy estimate
\lalin{weakthm1}{
E(t) &\leq E(0) + C\int_0^t\Bigl(1+\norm{\del{t}b(\tau)}_{L^\infty(\Omega)}+\norm{\chi(\tau)}_{L^\infty(\del{}\Omega)}
+ |\lambda| \norm{c(\tau)}_{L^n(\Omega)}\Bigr)
E(\tau) \notag \\
&\hspace{1.0cm}+ \norm{\del{t}L(\tau)}_{L^2(\Omega)}^2  + \norm{F(\tau)}_{L^2(\Omega)}^2+
\norm{g(\tau)}_{H^1(\Omega)}^2+ \norm{\del{t}g(\tau)}_{L^2(\Omega)}^2+\norm{\del{t}\theta(\tau)}_{H^1(\Omega)}^2 \, d\tau
\label{weakthm1.1}
}
for $0\leq t\leq T$ where $C=C\bigl(\kappa_0,\kappa_1,\mu,\gamma,\kf\bigr)$,
\lalin{weakthm2}{
&E(t) = \frac{1}{2}\ip{\del{i}v(t)}{b^{ij}(t)\del{j}v(t)}_{\Omega}-\frac{1}{2}\ip{\del{t}v(t)}{b^{00}(t)\del{t}v(t)}_{\Omega}
- \frac{1}{2}\ip{w(t)}{q(t)w(t)}_{\del{}\Omega} \notag \\
&\hspace{1.0cm} +\ip{\del{i}v(t)}{L^i(t)}_{\Omega}-\ip{v(t)}{k^i\del{i}\theta(t)}_{\del{}\Omega}+ \frac{\mu}{2}\norm{v(t)}_{L^2(\Omega)}^2 +
\frac{2}{\kappa_1}\norm{\vec{L}(t)}^2_{L^2(\Omega)}+\frac{3\kf^2}{\kappa_1}\norm{\theta(t)}^2_{H^1(\Omega)}, \label{weakthm2.1}
}
$\vec{L}=(L^i)$, $b=(b^{\alpha\beta})$,
and $E(t)$ satisfies
\eqn{weakthm3}{
E(t) \geq \min\biggl\{\frac{\kappa_0}{2},\frac{\kappa_1}{8},\frac{1}{2}\biggr\}\norm{(v(t),w(t))}^2_E
}
with $\norm{(v(t),w(t))}^2_E$ given by
\leqn{weakthm4}{
\norm{(v(t),w(t))}^2_E = \norm{v(t)}_{H^1(\Omega)}^2 + \norm{\del{t}v}^2_{L^2(\Omega)} + \ip{w(t)}{(-q(t))w(t)}_{\del{}\Omega} .
}
Moreover, $(v,w) \in \bigcap_{j=0}^1 C^j\bigl([0,T],H^{1-j}(\Omega,\Rbb^N)\bigr)
\times C^0\bigl([0,T],L^2(\del{}\Omega,\Rbb^N)\bigr)$.
\end{thm}
\begin{proof}
%We will only consider the existence of the weak solution in any detail. The rest of the proof follows from adapting the
%arguments used in the proof of Theorem 2.2. from \cite{Koch:1993}.
\noindent\underline{Existence (restricted initial data)}:
Following the proof of Theorem 2.2 from \cite{Koch:1993}, we employ a Galerkin
method to establish the existence of weak solutions to the IBVP \eqref{linB.1}-\eqref{linB.4}.
However, before proceeding, we note that we can absorb $g$ into the coefficients $L^\alpha$ and $F$
under the following transformation
\eqn{weakthm5}{
g\longmapsto 0, \quad L^\alpha \longmapsto L^\alpha-\nu^\alpha g \AND F \longmapsto F-\del{\alpha}(\nu^\alpha g).
}
Thus, without loss of generality, we may assume that $g=0$.

To start the Galerkin scheme, we consider the following alternate weak formulation of the IBVP
obtained by testing \eqref{linB.1} with $\phi\in C^\infty(\overline{\Omega},\Rbb^N)$ to
get
\eqn{weakthm6a}{
\ip{\phi}{\del{t}\bigl(b^{0\beta}\del{\beta}v+L^0\bigr)}_{\Omega}
+\ip{\phi}{\del{i}\bigl(b^{i\beta}\del{b}v+L^i\bigr)}_{\Omega}
+ \lambda\ip{\phi}{cv}_{\Omega} = \ip{\phi}{F}_{\Omega},
}
which we can express after integrating by parts as
\eqn{weakthm6b}{
\ip{\phi}{\del{t}\bigl(b^{0\beta}\del{\beta}v+L^0\bigr)}_{\Omega}
-\ip{\del{i}\phi}{b^{i\beta}\del{\beta}v+L^i}_{\Omega}
+\ip{\phi}{\nu_\alpha(b^{\alpha\beta}\del{\beta}v+L^\alpha)}_{\del{}\Omega}
+ \lambda\ip{\phi}{cv}_{\Omega} = \ip{\phi}{F}_{\Omega}.
}
Applying the boundary condition \eqref{linB.2} to the above equation, we obtain
\leqn{weakthm6c}{
\ip{\phi}{\del{t}\bigl(b^{0\beta}\del{\beta}v+L^0\bigr)}_{\Omega}
-\ip{\del{i}\phi}{b^{i\beta}\del{\beta}v+L^i}_{\Omega}
+\ip{\phi}{q\del{t}^2v+P\del{t}v+k^i\del{i}\theta}_{\del{}\Omega}
+ \lambda\ip{\phi}{cv}_{\Omega} = \ip{\phi}{F}_{\Omega}.
}
We now look for solutions of \eqref{weakthm6c} where
\leqn{weakthm7}{
\phi \in \text{Span}\{\phi_1,\phi_2,\ldots,\phi_M\}
}
and $\{\phi_I\}_{I=1}^\infty \subset H^2(\Omega,\Rbb^N)$ is a basis for $L^2(\Omega,\Rbb^N)$ satisfying:
\lalin{weakthm8}{
&\quad \ip{\phi_I}{\phi_J}_{\Omega} = \delta_{IJ}, \label{weakthm8.1}\\
&\quad\ip{\phi_I}{\phi_J}_{H^1(\Omega)} = 0, \hspace{1.0cm} I\neq J, \label{weakthm8.2}\\
&\norm{\pi_M(u)}_{L^2(\Omega)} \leq \norm{u}_{L^2(\Omega)}, \label{weakthm8.3}\\
&\norm{\pi_M(u)}_{H^1(\Omega)} \leq \norm{u}_{H^1(\Omega)} \label{weakthm8.4}
\intertext{and}
\pi_M u \rightarrow u  \quad &\text{as $M\rightarrow \infty$ in $L^2(\Omega,\Rbb^N)$ and $H^1(\Omega,\Rbb^N)$}, \label{weakthm8.5}
}
where $\pi_M \: :\: L^2(\Omega,\Rbb^N)\longrightarrow L^2(\Omega,\Rbb^N)$ is the self-adjoint projection operator
defined by
\eqn{parsum}{
\pi_M(u) = \sum_{I=1}^M \ip{\phi_I}{u}_{\Omega} \phi_I.
}
The existence of such a basis is established in \cite{Shkoller:2012}; see, in particular, Problem 5.9.

We consider solutions to \eqref{weakthm6c} of the form
\leqn{weakthm9}{
v_M(t,x)=\sum_{I=1}^M V_I(t)\phi_I(x).
}
Substituting $v_M$ into \eqref{weakthm6c} and setting $\phi=\phi_J$, we see that coefficients $V_I(t)$ must
satisfy a linear differential equation of the form
\leqn{weakthm10}{
\sum_{J=1}^M \bigl(A_{IJ}(t)V_J''(t) + B_{IJ}(t)V_I'(t)+C_{IJ}V_J(t)\Bigr) + D_{I}(t)= 0, \quad 0\leq t \leq T,
}
where
\eqn{weakthm11}{
A_{IJ}(t) = \ip{\phi_I}{b^{00}(t)\phi_J}_{\Omega}+\ip{\phi_I}{q(t)\phi_J}_{\del{}\Omega},
}
and the coefficients $A_{IJ}(t)$ $B_{IJ}(t)$, $C_{IJ}(t)$ and $D_{I}(t)$ are bounded on $[0,T]$.
Moreover, it follows from the assumptions \eqref{qPdef}, \eqref{bsym.1} and \eqref{bsym.2} that
$A_{IJ}(t)$ is invertible with bounded inverse on  $[0,T]$. To see this, we first note that
\leqn{weakthm12}{
A_{IJ}(t) = \ip{\phi_I}{b^{00}(t)\phi_J}_{\Omega}+\ip{\phi_I}{q(t)\phi_J}_{\del{}\Omega}
= \ip{b^{00}(t)\phi_I}{\phi_J}_{\Omega}+\ip{q(t)\phi_I}{\phi_J} = A_{JI}(t),
}
or in other words, $A_{IJ}(t)$ is a symmetric matrix for all $t\in [0,T]$. Next, we compute
\lalin{weakthm13}{
\sum_{I,J=1}^M c_I A_{IJ}(t)c_J &= \sum_{I,J=1}^M\Bigl(
\ip{c_I\phi_I}{b^{00}(t)c_J\phi_J}_{\Omega}+\ip{c_I\phi_I}{q(t)c_J\phi_J}_{\del{}\Omega}\Bigr) \notag \\
&=\biggl\langle\sum_{I=1}^M c_I \phi_I\biggl| b^{00}(t)\sum_{J=1}^M c_J\phi_J\biggr\rangle_{\Omega}+
\biggl\langle\sum_{I=1}^M c_I \phi_I\biggr|q(t)\sum_{J=1}^M c_J\phi_J\biggr\rangle_{\del{}\Omega} \notag \\
&\leq -\kappa_0  \biggl\langle\sum_{I=1}^M c_I \phi_I\biggl|\sum_{J=1}^M c_J\phi_J\biggr\rangle_{\Omega}
\notag \\
&=-\kappa_0 \sum_{I=1}^M c_I^2\label{weakthm13.1}
}
where in deriving the second to last line we used the assumptions $b^{00}\leq \kappa_0$ and
 $q\leq 0$, and in deriving the final line we used the orthonormality property \eqref{weakthm8.1}.
Since $\kappa_0 < 0$ by assumption, \eqref{weakthm12} and \eqref{weakthm13.1} imply that $A_{IJ}(t)$
is invertible for $t\in[0,T]$.

We are now in a position to apply standard existence and uniqueness theorems for linear differential equations
to conclude that there exists a unique solution $V_I \in W^{2,\infty}([0,T])$, $1\leq I\leq M$, to \eqref{weakthm10} satisfying the
initial conditions
\leqn{weakthm14}{
V_I(0)=\ip{\phi_I}{\vt_0}_{\Omega} \AND V_{I}'(0)=\ip{\phi_I}{\vt_1}_{\Omega}, \quad I=1,2,\ldots,M,
}
where $\vt_0$, $\vt_1$ is the initial data \eqref{linB.3}, which, for the moment, we assume satisfies
\eqn{weakthm14a}{
(\vt_0,\vt_1)\in H^1(\Omega,\Rbb^N)\times H^1(\Omega,\Rbb^N).
}
By construction, the solution \eqref{weakthm9} solves \eqref{weakthm6c} for any $\phi$ satisfying
\eqref{weakthm7}, and hence, in particular, for
\eqn{weakthm15}{
\phi(t,x) = \del{t}v_M(t,x)=\sum_{I=1}^M V_I'(t)\phi_I.
}
Substituting this into \eqref{weakthm6c} with $v=v_M$ yields
\lalin{weakthm16}{
&\ip{\del{t}v_M}{\del{t}\bigl(b^{0\beta}\del{\beta}v_M+L^0\bigr)}_{\Omega}
-\ip{\del{i}\del{t}v_M}{b^{i\beta}\del{\beta}v_M+L^i}_{\Omega}
\notag \\
&\hspace{1.0cm}+\ip{\del{t}v_M}{q\del{t}^2v_M+P\del{t}v_M+k^i\del{i}\theta}_{\del{}\Omega}
+ \lambda\ip{\del{t}v_M}{cv}_{\Omega} - \ip{\del{t}v_M}{F}_{\Omega}=0. \label{weakthm16.1}
}
Next, from the symmetry \eqref{bsym.1} of the $b^{\alpha\beta}$, it is not difficult
to see that
\alin{weakthm17}{
\ip{\del{i}\del{t}v_M}{b^{i\beta}\del{\beta}v_M+L^i}_{\Omega} & =
\ip{\del{t}v_M}{\del{t}\bigl(b^{0i}\del{i}v_M+L^0\bigr)}_{\Omega}
-\ip{\del{\alpha}v_M}{\del{t}L^\alpha}_{\Omega}-\ip{\del{t}v_M}{\del{t}b^{0i}\del{i}v_M}_{\Omega} \\
&-\Half\ip{\del{i}v_M}{\del{t}b^{i j}\del{j}v_M}_{\Omega}
+\del{t}\Bigl[\Half\ip{\del{i}v_M}{b^{ij}\del{j}v_M}_{\Omega}+\ip{\del{i}{v_M}}{L^i}_{\Omega}\Bigr].
}
Substituting this into \eqref{weakthm16.1}, we obtain with the help
the symmetry conditions $(b^{00})^{\tr}=b^{00}$ and $q^{\tr}=q$ that
\lalin{weakthm18}{
&\del{t}\Bigl[\Half\Bigl(\ip{\del{t}v_M}{b^{00}\del{t}v_M}_{\Omega}-
\ip{\del{i}v_M}{b^{ij}\del{j}v_M}_{\Omega}\Bigr)-\ip{\del{i}{v_M}}{L^i}_{\Omega}
 +\Half\ip{\del{t}v_M}{q\del{t}v_M}_{\del{}\Omega}\notag \\
&\hspace{2.2cm} +\ip{v_M}{k^i\del{i}\theta}_{\del{}\Omega}\Bigr]
-\Half\ip{\del{t}v_M}{\del{t}b^{00}\del{t}v_M}_{\Omega}
+\Half\ip{\del{i}v_M}{\del{t}b^{ij}\del{j}v_M}_{\Omega} \notag\\
&\hspace{3.0cm}+\ip{\del{t}v_M}{\del{t}b^{0j}\del{j}v_M}_{\Omega}
+\ip{\del{i}v_M}{\del{t}L^i}_{\Omega}+2\ip{\del{t}v_M}{\del{t}L^0}_{\Omega} \notag\\
&\hspace{2.0cm}+\ip{\del{t}v_M}{\bigl(P-\Half\del{t}q\bigr)\del{t}v_M}_{\del{}\Omega}
-\ip{v_M}{\del{t}(k^i\del{i}\theta)}_{\del{}\Omega}+\ip{\del{t}v_M}{\lambda c v_M-F}_{\Omega}=0.
\label{weakthm18.1}
}

Introducing the energy
\lalin{weakthm19}{
&E_M(t) = -\frac{1}{2}\ip{\del{t}v_M(t)}{b^{00}(t)\del{t}v_M(t)}_{\Omega}+
\frac{1}{2}\ip{\del{i}v_M(t)}{b^{ij}(t)\del{j}v_M(t)}_{\Omega} \notag \\
&\hspace{1.3cm}+\ip{\del{i}{v_M(t)}}{L^i(t)}_{\Omega} -\frac{1}{2}\ip{\del{t}v_M(t)}{q(t)\del{t}v_M(t)}_{\del{}\Omega}
 -\ip{v_M(t)}{k^i\del{i}\theta(t)}_{\del{}\Omega} \notag\\
 &\hspace{4.6cm}+\frac{\mu}{2}\norm{v_M(t)}_{L^2(\Omega)}^2+\frac{2}{\kappa_1}\norm{\vec{L}(t)}_{L^2(\Omega)}^2
 +\frac{3\kf^2}{\kappa_1}\norm{\theta(t)}^2_{H^1(\Omega)}, \label{weakthm19.1}
}
where
\eqn{weakthm20}{
\vec{L}=(L^i) \AND \kf = \norm{k}_{W^{1,\infty}([0,T],L^\infty(\Omega)\cap W^{1,n}(\Omega))},
}
we find after differentiating $E_M(t)$ that
\lalin{weakthm20}{
&E_M'(t)= \del{t}\Biggl[-\frac{1}{2}\ip{\del{t}v_M(t)}{b^{00}(t)\del{t}v_M(t)}_{\Omega}+
\frac{1}{2}\ip{\del{i}v_M(t)}{b^{ij}(t)\del{j}v_M(t)}_{\Omega}\Bigr)\notag \\
&\hspace{1.5cm} +\ip{\del{i}v_M(t)}{L^i(t)}_{\Omega}
  -\frac{1}{2}\ip{\del{t}v_M(t)}{q(t)\del{t}v_M(t)}_{\del{}\Omega}
+\ip{v_M(t)}{k^i(t)\del{i}\theta(t)}_{\del{}\Omega}\Biggr] \notag \\
&\hspace{3.0cm}+
\mu\ip{v_M(t)}{\del{t}v_M(t)}_{\Omega}
+\frac{4}{\kappa_1}\ip{\vec{L}(t)}{\del{t}\vec{L}(t)}_{\Omega}+\frac{6\kf^2}{\kappa_1}\ip{\theta(t)}{\del{t}\theta(t)}_{H^1(\Omega)}.
\notag %\label{weakthm20.1}
}
From \eqref{weakthm18.1}, the assumption $-\Half\del{t}q+P-\chi q\leq 0$,
and the estimate \eqref{weakrem7}, it then follows
that $E_M'$ satisfies the inequality
\lalin{weakthm21}{
&E'_M \leq C(\mu,\kappa_1,\kf)\Bigl(\Bigl[1+|\lambda|\norm{c}_{L^n(\Omega)}+\norm{\del{t}b}_{L^\infty(\Omega)}
+\norm{\chi}_{L^\infty(\del{}\Omega)}
\Bigr]\Bigl(\norm{v_M}_{H^1(\Omega)}^2
 \notag \\
&\hspace{1.7cm} +\norm{\del{t}v_M}^2_{L^2(\Omega)}
-\ip{\del{t}v_M}{q\del{t}v_M}_{\del{}\Omega}\Bigr)+\Bigl(\norm{v_M}_{H^1(\Omega)}  +\norm{\del{t}v_1}_{L^2(\Omega)}
\notag \\
&\hspace{2.0cm}+\norm{\vec{L}}_{L^2(\Omega)}
 +\norm{\theta}_{H^1(\Omega)}
\Bigr)\Bigr(\norm{\theta}_{H^1(\Omega)}
 +\norm{\del{t}\theta}_{H^1(\Omega)} +\norm{\del{t}L}_{L^2(\Omega)}+\norm{F}_{L^2(\Omega)}\Bigr) . \label{weakthm21.1}
}

For the energy inequality \eqref{weakthm21.1} to be useful, we first must verify that $E_M$ controls
$\norm{(v_M,\del{t}v_M)}_{E}$
where $\norm{(\cdot,\cdot)}_E$ is the energy norm defined by \eqref{weakthm4}.  Using Cauchy-Schwartz and Young's
(i.e. $ab\leq \frac{1}{2\ep}a^2+\frac{\ep}{2}b^2$, $a,b\geq 0$ and $\ep>0$) inequalities together with the
assumptions \eqref{bsym.2} and \eqref{coerc}, it is clear from \eqref{weakthm19.1}
that the energy $E_M$ is bounded below by
\eqn{weakthm23}{
E_M \geq \frac{\kappa_1}{4}\norm{v_M}_{H^1(\Omega)}^2+\frac{\kappa_0}{2}\norm{\del{t}v_M}_{L^2(\Omega)}
-\frac{1}{2}\ip{\del{t}v_M}{q\del{t}v_M}_{L^2(\del{}\Omega)} +\frac{1}{\kappa_1}\norm{\vec{L}}^2_{L^2(\Omega)}
-\ip{v_M}{k^i\del{i}\theta}_{\del{}\Omega}+\frac{3\kf^2}{\kappa_1}\norm{\theta}^2_{H^1(\Omega)}.
}
From this estimate, we see with the help of Young's inequality and \eqref{weakrem7} that
\leqn{weakthm28}{
E_M \geq \min\left\{\frac{\kappa_1}{8},\frac{\kappa_0}{2},\frac{1}{2}\right\}\norm{(v_M,\del{t}v_M)}^2_E
+\frac{1}{\kappa_1}\norm{\vec{L}}_{L^2(\Omega)}^2 + \frac{\kf^2}{\kappa_1}\norm{\theta}^2_{H^1(\Omega)}.
}
Using this lower bound, it follows immediately from
\eqref{weakthm21.1} that $E_M$ satisfies
\lalin{weakthm20}{
&E_M' \leq C(\mu,\kappa_0,\kappa_1,\kf)\Bigl(\Bigl[1+|\lambda|\norm{c}_{L^n(\Omega)}+\norm{\del{t}b}_{L^\infty(\Omega)}
+\norm{\chi}_{L^\infty(\del{}\Omega)}
\Bigr]E_M'\notag \\
&\hspace{6.0cm} +\norm{\del{t}L}_{L^2(\Omega)}^2+
\norm{\del{t}\theta}_{H^1(\Omega)}^2+\norm{F}_{L^2(\Omega)}^2\Bigr). \notag %\label{weakthm29.1}
}
Integrating this in time yields the estimate
\lalin{weakthm30}{
&E_M(t) \leq C(\mu,\kappa_0,\kappa_1,\kf)\biggl(E_M(0)
+\int_0^t \Bigl(1+|\lambda|\norm{c(\tau)}_{L^n(\Omega)}+\norm{\del{t}b(\tau)}_{L^\infty(\Omega)}
\notag \\
&\hspace{2.0cm} +\norm{\chi(\tau)}_{L^\infty(\del{}\Omega)} \Bigr)E_M(\tau)
 +\norm{\del{t}L(\tau)}_{L^2(\Omega)}^2+
\norm{\del{t}\theta(\tau)}_{H^1(\Omega)}^2+\norm{F(\tau)}_{L^2(\Omega)}^2\, d\tau\biggr), \label{weakthm30.1}
}
which in turn, shows by Gronwall's inequality that
\leqn{weakthm31}{
E_M(t) \leq C\biggl(E_M(0)
+\int_0^t \norm{\del{t}L(\tau)}_{L^2(\Omega)}^2+
\norm{\del{t}\theta(\tau)}_{H^1(\Omega)}^2+\norm{F(\tau)}_{L^2(\Omega)}^2\, d\tau\biggr)
}
where
\eqn{weakthm32}{
C=C\bigl(\mu,\kappa_0,\kappa_1,\kf,\norm{\del{t}b}_{L^\infty(\Omega_T)},\norm{\chi}_{L^\infty(\Gamma_T)},
\norm{c}_{L^\infty([0,T],L^n(\Omega))} \bigr).
}

The next step in constructing the solutions is to
to let $M\rightarrow \infty$ and use the inequality \eqref{weakthm31} to
provide a uniform bound on the energy norm $\norm{(v_M(t),\del{t}v_M(t))}$ for $t\in [0,T]$. However, before
we can do this, we must first establish a uniform bound on $E_M(0)$. We begin by expressing
$v_M(0)$ and $\del{t}v_M(0)$ as
\eqn{weakthm46a}{
v_M(0) = \pi_M(\vt_0) \AND \del{t}v_M(0) = \pi_M(\vt_1),
}
which is clear from \eqref{weakthm9} and \eqref{weakthm14}.
Using properties \eqref{weakthm8.3} and \eqref{weakthm8.4} of the projection operator $\pi_M$,
we obtain the bound
\leqn{weakthm46b}{
\norm{v_M(0)}_{H^1(\Omega)}^2+\norm{\del{t}v_M(0)}_{L^2(\Omega)}^2 \leq
\norm{\vt_0}_{H^1(\Omega)}^2+\norm{\vt_1}_{L^2(\Omega)}^2.
}
Additionally, we observe that
\lalin{weakthm46}{
-\ip{\del{t}v_M(0)}{q(0)\del{t}v_M(0)}_{\del{}\Omega} &\lesssim \norm{q(0)\del{t}v_M(0)}_{L^2(\del{}\Omega)}^2
&& \text{(by \eqref{qnorm})}\notag \\
&\lesssim \norm{q(0)}_{W^{1,\infty}(\Omega)}^2\norm{\del{t}v_M(0)}_{H^1(\Omega)}^2
&& \text{(by Theorem \ref{trace})} \notag \\
& \leq \norm{\vt_1}_{H^1(\Omega)}^2 &&
\text{(by \eqref{weakthm8.4})}. \label{weakthm47}
}
Taken together, the inequalities \eqref{weakthm46b} and \eqref{weakthm47} imply
that
\eqn{weakthm47}{
E_M(0) \lesssim 1, \quad M \geq 1,
}
from which the uniform bound
\leqn{weakthm48}{
\norm{v_M(t)}_{H^1(\Omega)}+ \norm{\del{t}v_M(t)}_{L^2(\Omega)}
+\norm{q(t)\del{t}v_M(t)|_{\del{}\Omega}}_{L^2(\del{}\Omega)}  \leq C, \quad (M,t)\in \Nbb\times [0,T],
}
follows by \eqref{qnorm}, \eqref{weakthm28} and \eqref{weakthm31}.  By the Banach-Alaoglu Theorem,
there exists $v\in L^p\bigl([0,T],H^1(\Omega,\Rbb^N)\bigr)$, $\dot{v}\in L^p\bigl([0,T],L^2(\Omega,\Rbb^N)\bigr)$,
$\wh \in L^p\bigl([0,T],L^2(\del{}\Omega,\Rbb^N)\bigr)$, and a subsequence of $v_M$, again denoted $v_M$, such
that
\lalin{weakthm49}{
&v_M \rightharpoonup v \hspace{1.4cm} \text{in $L^p\bigl([0,T],H^1(\Omega,\Rbb^N)\bigl)$}, \label{weakthm49.1}\\
&\del{t}v_M \rightharpoonup \dot{v}\hspace{1.1cm} \text{in $L^p\bigl([0,T],L^2(\Omega,\Rbb^N)\bigr)$}, \label{weakthm49.2}
\intertext{and}
&q\del{t}v_M|_{\del{}\Omega} \rightharpoonup \wh \hspace{0.4cm} \text{in $L^p\bigl([0,T],L^2(\del{}\Omega,\Rbb^N)\bigr)$}
\label{weakthm49.3}
}
for any $p\in (1,\infty)$. Moreover, since the constant $C>0$ in the inequality \eqref{weakthm48} is
independent of $p\in (1,\infty)$, we can use the property $\lim_{p\rightarrow \infty}\norm{f(t)}_{L^p([0,T])}
= \norm{f(t)}_{L^\infty([0,T])}$ of $L^p$ norms to conclude that
\leqn{weakthm49a}{
(v,\dot{v},\wh) \in L^\infty([0,T],H^1\bigl(\Omega,\Rbb^N)\bigr)\times L^\infty\bigl([0,T],L^2(\Omega,\Rbb^N)\bigr)
\times L^\infty([0,T],L^2\bigl(\del{}\Omega,\Rbb^N)\bigr).
}

Since $q\del{t}v_M$ satisfies $\Pbb_{q}q\del{t}v_M = q\del{t}v_M$ it is not difficult to see that
the limit $\wh$ must also satisfy $\Pbb_{q}\wh = \wh$. This together with the invertibility of
the map $q|_{\Pbb_q}$, see \eqref{qdecompD}, allows us to write \eqref{weakthm49.3} as
\leqn{weakthm50}{
q\del{t}v_M|_{\del{}\Omega} \rightharpoonup q w \hspace{0.4cm} \text{in $L^\infty\bigl([0,T],L^2(\del{}\Omega,\Rbb^N)\bigr)$}
}
where
\leqn{weakthm51}{
w := (q|_{\Pbb_{q}})^{-1} \wh \in L^\infty\bigl([0,T],L^2(\del{}\Omega,\Rbb^N)\bigr).
}
Furthermore, since $\del{t}v_M$ satisfies
\lalin{weakthm52}{
\ip{\del{t}v_M}{\phi}_{\Omega_T} &= -\ip{v_M}{\del{t}\phi}_{\Omega_T} \notag
\intertext{and}
\ip{q\del{t}v_M}{\phi}_{\Gamma_T} &= -\ip{\del{t}q v_M}{\phi}_{\Gamma_T} - \ip{q v_M}{\del{t}\phi}_{\Gamma_T}
\label{weakthm52.2}
}
for all $\phi \in C^1_0([0,T],C^1(\Omega,\Rbb^N))$,
it follows from  \eqref{weakthm49.1}, \eqref{weakthm49.2} and \eqref{weakthm50} that\footnote{Here, we are interpreting
the boundary terms on the righthand side of \eqref{weakthm52.2} via the formula \eqref{weakrem2a}.}
\leqn{weakthm52a}{
\ip{\dot{v}}{\phi}_{\Omega_T} = -\ip{v}{\del{t}\phi}_{\Omega_T} \AND
\ip{q w}{\phi}_{\Gamma_T} = -\ip{\del{t}q v}{\phi}_{\Gamma_T} - \ip{q v}{\del{t}\phi}_{\Gamma_T}.
}
Thus, we the formulas
\lalin{weakthm53}{
\dot{v} &= \del{t}v  \label{weakthm53.1}
\intertext{and}
qw&=q\del{t}v \notag
}
hold in a weak sense.

Next, in \eqref{weakthm6c} with $v=v_M$, we set
\leqn{weakthm54}{
\phi(t,x)=\sum_{I=1}^M c_I(t)\phi_I(x), \quad c_I \in C^1_0([0,T]),
}
and integrate in time from
$t=0$ to $t=T$ and by parts to obtain
\eqn{weakthm55}{
\ip{b^{\alpha\beta}\del{\beta}v_M + L^\alpha}{\del{\alpha} \phi}_{\Omega_T} +
\ip{(\del{t}q-P)\del{t}v_M}{\phi}_{\Gamma_T}-\ip{k^i\del{i}\theta}{\phi}_{\Gamma_T}+
\ip{q\del{t}v_M}{\del{t}\phi}_{\Gamma_T}= \ip{\lambda c v_M -F}{\phi}_{\Omega_T},
}
and hence, that
\leqn{weakthm56}{
\ip{b^{\alpha\beta}\del{\beta}v + L^\alpha}{\del{\alpha} \phi}_{\Omega_T} +
\ip{(\del{t}q-P)\del{t}v}{\phi}_{\Gamma_T}-\ip{k^i\del{i}\theta}{\phi}_{\Gamma_T}+
\ip{qw}{\del{t}\phi}_{\Gamma_T}= \ip{\lambda c v -F}{\phi}_{\Omega_T}.
}
by \eqref{weakthm49.1}, \eqref{weakthm49.2}, \eqref{weakthm50} and \eqref{weakthm53.1}. Furthermore, since
$M$ can be chosen arbitrarily large and we can approximate maps $C^1_0([0,T],C^1(\overline{\Omega},\Rbb^N)$
by maps of the form \eqref{weakthm54}, it follows that \eqref{weakthm56} continues to hold for
all $\phi  \in C^1_0([0,T],C^1(\overline{\Omega},\Rbb^N)$.

Thus far, we have shown that the pair $(v,w)$ defined by \eqref{weakthm49.1} and
\eqref{weakthm51} satisfy \eqref{weakthm56} for all $\phi  \in C^1_0([0,T],C^1(\overline{\Omega},\Rbb^N)$. The next step is to verify
that the pair $(v,w)$ satisfies the initial conditions $v|_{t=0}=\vt_0$, $\del{t}v|_{t=0}=\vt_1$ and $w|_{t=0}= \Pbb_{q(0)}\vt_1|_{\del{}\Omega}$.
To see this, we choose a $\phi \in C^1(\overline{\Omega}_T,\Rbb^N)$ that satisfies $\phi|_{t=T}=0$. Since
\eqn{vMlimreg}{
v\in C^{0}\bigl([0,T],L^2(\Omega,\Rbb^N)\bigr) \subset W^{1,\infty}\bigl([0,T],L^2(\Omega,\Rbb^N)\bigr)
\AND \del{t}v \in L^\infty\bigl([0,T],L^2(\Omega,\Rbb^N)\bigr)
}
by \eqref{weakthm49a} and \eqref{weakthm53.1}, we see from integrating
\eqn{weaktm57}{
\del{t}\ipe{\phi}{v} = \ipe{\del{t}\phi}{v} + \ipe{\phi}{\del{t}v}
}
over $\Omega_T$ that
\leqn{weakthm58}{
\ip{\phi(0)}{v(0)}_{\Omega} = \ip{\del{t}\phi}{v}_{\Omega_T} + \ip{\phi}{\del{t}v}_{\Omega_T}.
}
By a similar calculation, we also have that
\eqn{weakthm59}{
\ip{\phi(0)}{\pi_M(\vt_0)}_{\Omega} = \ip{\del{t}\phi}{v_M}_{\Omega_T} + \ip{\phi}{\del{t}v_M}_{\Omega_T},
}
which in turn, implies, after letting $M\rightarrow \infty$, that
\leqn{weakthm60}{
\ip{\phi(0)}{\vt_0}_{\Omega} = \ip{\del{t}\phi}{v}_{\Omega_T} + \ip{\phi}{\del{t}v}_{\Omega_T},
}
by \eqref{weakthm8.5}, \eqref{weakthm49.1}, \eqref{weakthm49.2} and \eqref{weakthm53.1}. Since \eqref{weakthm58} and \eqref{weakthm60} hold
for all $\phi \in C^1(\overline{\Omega}_T,\Rbb^N)$ satisfying $\phi|_{t=T}=0$, we conclude that
\leqn{weakthm61}{
v(0) = \vt_0.
}

To proceed, we set
\eqn{weakthm62}{
\phi(t,x)=\sum_{I=1}^M c_I(t)\phi_I(x),
}
where $c_I\in C^1([0,T])$ and  $c_I(T)=0$, in \eqref{weakthm6c} with $v=v_M$, which we know by construction solves \eqref{weakthm6c}. Integrating the resulting expression over $t$ from $t=0$ to $t=T$, we
find, with the help of \eqref{weakthm14}, that
\alin{weakthm63}{
&\ip{q(0)\pi_M \vt_1}{\phi(0)}_{\del{}\Omega}  + \ip{b^{00}(0)\pi_M \vt_1 + b^{0i}(0)\del{i}\pi_M \vt_0+L^0(0)}{\phi(0)}_{\Omega} \\
&\hspace{1.0cm} +\ip{(P(0)-\del{t}q(0))\nu(\pi_M\vt_0)}{\phi(0)}_{\Omega}  + \ip{b^{\alpha\beta}\del{\beta}v_M+L^\alpha}{\del{\alpha}\phi}_{\Omega_T} +\ip{q\del{t}v_M}{\del{t}\phi}_{\Gamma_T}\\
&\hspace{4.0cm}
+ \ip{(\del{t}q-P)\del{t}v_M}{\phi}_{\Gamma_T} -
\ip{k^i\del{i}\theta}{\phi}_{\Gamma_T}+\ip{\lambda c v_M-F}{\phi}_{\Omega_T} = 0,
}
where we recall that the ``boundary term'' $\ip{(\del{t}q-P)\del{t}v_M}{\phi}_{\Gamma_T}$ is defined via the formula \eqref{weakrem2.1}. Letting $M\rightarrow \infty$ in the above
expression, it follows from \eqref{weakthm8.5}, \eqref{weakthm49.1}, \eqref{weakthm49.2}, \eqref{weakthm50} and \eqref{weakthm53.1} that
\lalin{weakthm64}{
&\ip{q(0)\vt_1}{\phi(0)}_{\del{}\Omega}  + \ip{b^{00}(0)\vt_1 + b^{0i}(0)\del{i}\vt_0+L^0(0)}{\phi(0)}_{\Omega} \notag \\
&\hspace{0.5cm} +\ip{(P(0)-\del{t}q(0))\nu(\vt_0)}{\phi(0)}_{\Omega}  + \ip{b^{\alpha\beta}\del{\beta}v+L^\alpha}{\del{\alpha}\phi}_{\Omega_T} +\ip{q w}{\del{t}\phi}_{\Gamma_T} \notag\\
&\hspace{3.5cm}
+ \ip{(\del{t}q-P)\del{t}v}{\phi}_{\Gamma_T} -
\ip{k^i\del{i}\theta}{\phi}_{\Gamma_T}+\ip{\lambda c v-F}{\phi}_{\Omega_T} = 0, \label{weakthm64.1}
}
which, by approximation, holds for all $\phi \in C^1(\overline{\Omega}_T,\Rbb^N)$ satisfying $\phi_{t=T}=0$. Choosing $\phi$ so that
\eqn{weakthm65}{
\phi(t,x) = \eta(t)\psi(x),
}
where $\eta \in C^1_0([0,T])$ and $\psi\in C^1(\overline{\Omega})$, we see from \eqref{weakthm64.1}, or alternatively \eqref{weakthm56}, that
\lalin{weakthm66}{
&\int_0^T \eta'\Bigl[\ip{b^{0\beta}\del{\beta}v+L^0+(P-\del{t}q)\nu(v)}{\psi}_{\Omega}+\ip{qw}{\psi}_{\Omega}\Bigr]\, dt = \int_0^T \eta \Bigl[-\ip{b^{i\beta}\del{\beta}v+L^i}{\del{i}\psi} \notag \\
&\hspace{4.5cm} + \ip{k^i\del{i}\theta}{\psi}_{\del{}\Omega}
+\ip{\nu(P-\del{t}q)\del{t}v-\del{t}(P-\del{t}q)\nu(v)}{\psi}_{\Omega} \notag \\
&\hspace{3.3cm}+\ip{\del{\alpha}v^\alpha(P-\del{t}q)\del{t}v}{\psi}_{\Omega}+ \ip{(P-\del{t}q)\del{t}v}{\nu(\psi)}_{\Omega}+\ip{F-\lambda c v}{\psi}_{\Omega}\Bigr]\, dt. \label{weakthm66.1}
}
From this we conclude that there exists a map
$\Vc \in W^{1,\infty}([0,T],H^1(\Omega)) \subset C^{0}([0,T],H^1(\Omega))$
such that
\lalin{weakthm68}{
\ip{\Vc}{\psi}_{H^1(\Omega)} &= \ip{b^{0\beta}\del{\beta}v+L^0+(P-\del{t}q)\nu(v)}{\psi}_{\Omega}+\ip{qw}{\psi}_{\del{}\Omega} \label{weakthm68.1}
\intertext{and}
\ip{\del{t}\Vc}{\psi}_{H^1(\Omega)} &= \ip{b^{i\beta}\del{\beta}v+L^i}{\del{i}\psi}_{\Omega} - \ip{k^i\del{i}\theta}{\psi}_{\del{}\Omega}
-\ip{\nu(P-\del{t}q)\del{t}v-\del{t}(P-\del{t}q)\nu(v)}{\psi}_{\Omega}   \notag \\
&\qquad -\ip{\del{\alpha}v^\alpha(P-\del{t}q)\del{t}v}{\psi}_{\Omega}- \ip{(P-\del{t}q)\del{t}v}{\nu(\psi)}_{\Omega}-\ip{F-\lambda c v}{\psi}_{\Omega}
\label{weakthm68.2}
}
for all $\psi \in H^1(\Omega,\Rbb^N)$.

Fixing $\phi \in C^1(\overline{\Omega}_T,\Rbb^N)$ satisfying $\phi|_{t=T}=0$, we see from integrating the identity
\eqn{weakthm69}{
\del{t}\ip{\Vc}{\phi}_{H^1(\Omega)} = \ip{\del{t}\Vc}{\phi}_{H^1(\Omega)} + \ip{\Vc}{\del{t}\phi}_{H^1(\Omega)}
}
over $t$ from $t=0$ to $t=T$ and using \eqref{weakthm68.1}-\eqref{weakthm68.2} that
\lalin{weakthm70}{
\ip{\Vc(0)}{\phi(0)}_{H^1(\Omega)} +&  \ip{b^{\alpha\beta}\del{\beta}v+L^\alpha}{\del{\alpha}\phi}_{\Omega_T} +\ip{q w}{\del{t}\phi}_{\Gamma_T}
 \notag\\
&+ \ip{(\del{t}q-P)\del{t}v}{\phi}_{\Gamma_T} -
\ip{k^i\del{i}\theta}{\phi}_{\Gamma_T}+\ip{\lambda c v-F}{\phi}_{\Omega_T} = 0. \label{weakthm70.1}
}
Comparing \eqref{weakthm64.1} and \eqref{weakthm70.1}, we see that
\leqn{weakthm71}{
\lim_{t\searrow 0}\ip{\Vc(t)}{\psi}_{H^1(\Omega)} = \ip{q(0)\vt_1}{\psi}_{\del{}\Omega}  + \ip{b^{00}(0)\vt_1 + b^{0i}(0)\del{i}\vt_0+L^0(0)+(P(0)-\del{t}q(0))\nu(\vt_0)}{\psi}_{\Omega}
}
for all $\psi \in H^1(\Omega,\Rbb^N)$.
From the formulas \eqref{weakthm68.1} and \eqref{weakthm71}, and the regularity properties
\gath{weakthm72}{
v\in C^{0}\bigl([0,T],L^2(\Omega,\Rbb^N)\bigr)\cap L^\infty\bigl([0,T],H^1(\Omega,\Rbb^N)\bigr), \quad \del{t}v\in L^\infty([0,T],L^2(\Omega,\Rbb^N)) \notag \\
w\in L^\infty\bigl([0,T],L^2(\del{}\Omega,\Rbb^N)\bigr) \AND \Vc \in C^{0}\bigl([0,T],H^1(\Omega,\Rbb^N)\bigr),
}
and \eqref{weakthm61}, it is not difficult to verify that
\eqn{weakthm73}{
\bigl\langle\lim_{t\searrow 0} b^{00}(t)(\del{t}v(t)-\vt_1) \bigl|\psi\bigr\rangle_{\Omega}= 0, \quad \forall \, \psi \in C^1_0(\Omega,\Rbb^N),
}
from which we conclude that the weak limit $\lim_{t\searrow 0} b^{00}(t)(\del{t}v(t)-\vt_1)$ exists in $L^2(\Omega,\Rbb^N)$ and is zero. We conclude via
the invertibility of $b^{00}$ that
\leqn{weakthm74}{
\del{t}v(t) \rightharpoonup  \vt_1 \quad \text{as $t\searrow 0$ in $L^2(\Omega,\Rbb^N)$.}
}
Using this fact, it follows from \eqref{weakthm64.1} and \eqref{weakthm70.1} and similar reasoning that
\leqn{weakthm75}{
\bigl\langle\lim_{t\searrow 0} q(t)(w(t)-\vt_1|_{\del{}\Omega}) \bigl|\psi\bigr\rangle_{\del{}\Omega}= 0, \quad \forall \, \psi \in C^1(\Omega,\Rbb^N).
}
From this we see that the weak limit $\lim_{t\searrow 0} q(t)(w(t)-\vt_1|_{\del{}\Omega})$ exists in $L^2(\del{}\Omega,\Rbb^N)$ and is zero, and hence, that
\leqn{weakthm76}{
w(t)\rightharpoonup \Pbb_{q(0)}\vt_1|_{\del{}\Omega} \quad \text{as $t\searrow 0$ in $L^2(\del{}\Omega,\Rbb^N)$}.
}

Collectively, \eqref{weakthm52a}, \eqref{weakthm56}, \eqref{weakthm61}, \eqref{weakthm74} and \eqref{weakthm76} show that the pair
\eqn{weakthm77}{
(v,w) \in \cap_{\ell=0}^1 W^{\ell,\infty}([0,T],H^{1-\ell}(\Omega))\times L^\infty([0,T],L^2(\del{}\Omega))
}
define a weak solution of the IBVP \eqref{linB.1}-\eqref{linB.4} for the restricted class of initial data
\eqn{weakthm78}{
(v,\del{t}v)|_{\Omega_0} = (\vt_0,\vt_1) \in H^1(\Omega)\times H^1(\Omega)  \AND w|_{\Gamma_0} = \Pbb_{q(0)}\vt_1|_{\del{}\Omega} \in L^2(\del{}\Omega).
}
\begin{comment}
We
can now argue using a straightforward modification of the argument from the bottom of pg 21 in
\cite{Koch:1993} to show $(v,w)$ are weakly continuous at $t=0$, i.e. \eqref{weaksoldef1} is
satisfied. This establishes the existence of a weak solution to the IBVP \eqref{linB.1}-\eqref{linB.4}
for initial data $(\vt_0,\vt_1,\wt)$ satisfying \eqref{weakthm14a} and $\wt=\Pbb_{q(0)} \vt_1|_{\del{}\Omega}$.
\end{comment}

\bigskip

\noindent \underline{Uniqueness and regularity}:
By adapting the arguments used by Koch in the proof of Theorem 2.2 from \cite{Koch:1993}, it is
not difficult to show that weak solutions $(v,w)$, in the sense of Definition \ref{weaksoldef}
and not necessarily the ones constructed from the Galerkin scheme above,
to \eqref{linB.1}-\eqref{linB.4} are
unique, and satisfy the additional regularity
\eqn{weakthm79}{
(v,w) \in \bigcap_{j=0}^1 C^j\bigl([0,T],H^{1-j}(\Omega,\Rbb^N)\bigr)\times C^0\bigl([0,T],L^2(\del{}\Omega,\Rbb^N)\bigr),
}
and the energy estimate \eqref{weakthm1.1},
which, not surprisingly, is exactly the same as the
energy estimate \eqref{weakthm30.1} satisfied by the Galerkin approximations.

\bigskip

\noindent\underline{Existence (general initial data)}:
The final step is to approximate general initial data
\leqn{weakthm80}{
(\vt_0,\vt_1,\wt) \in H^1(\Omega,\Rbb^N)\times L^2(\Omega,\Rbb^N)\times L^2(\del{}\Omega,\Rbb^N),
}
where $\Pbb_{q(0)}\wt = \wt$, by the more regular initial data
\leqn{weakthm59}{
(\vt_0,\vt_1,\wt=\Pbb_{q(0)}\vt_1|_{\del{}\Omega}) \in H^1(\Omega,\Rbb^N)\times H^1(\Omega,\Rbb^N)\times L^2(\del{}\Omega,\Rbb^N)
}
for which we can prove existence via the Galerkin scheme. Weak solutions generated by
initial data of the type \eqref{weakthm58} are then obtained by taking weak limits of sequences of solutions generated from
initial data satisfying \eqref{weakthm59}.
\end{proof}

\begin{rem} \label{weakremA}
$\;$

\begin{itemize}
\item[(i)] It is clear from the proof of Theorem \ref{weakthm} that it continues to hold for $G$ of the form
\eqn{weakrem1a}{
G = k^i_{1}\del{i}\theta_{1}+k^i_{2}\del{i}\theta_{2}+g,
}
where $\nu_i k^i_{a}=0$, $a,b=1,2$,
provided that we make the replacements
\gath{weakrem1b}{
\theta \longmapsto \vec{\theta} = (\theta_{1},\theta_{1}) \AND
k=(k^i) \longmapsto \vec{k}=(k_{1}^i,k_{2}^i)
}
along with the obvious substitutions, e.g. $\norm{\theta}_{H^1(\Omega)}^2$ $\mapsto$ $\norm{\vec{\theta}}_{H^1(\Omega)}^2$ $=$
$\norm{\theta_1}^2_{H^1(\Omega)}$$+$$\norm{\theta_2}^2_{H^1(\Omega)}$, in the energy estimate \eqref{weakthm1.1} and
energy norm \eqref{weakthm2.1}.
\item[(ii)]
It is not difficult to show using an iteration method that the existence and uniqueness statement and the energy estimate
from Theorem \ref{weakthm} continues to hold if we allow the coefficients $L^\alpha$, $F$,
$G$ and $g$ to depend linearly on $v$ provided that they satisfy estimates that preserve the form of the energy estimate.
For example, we could have
\eqn{weakrem1}{
F=\bar{F} + \hat{F}(v)
}
where $\bar{F}\in L^2(\Omega_T,\Rbb^N)$ and
\eqn{weakrem2}{
\norm{\hat{F}(v(t))}_{L^2(\Omega)}^2 \leq  \hat{f}(t)E(t), \quad 0\leq t\leq T.
}
In this case, we would just replace $\norm{F(t)}^2_{L^2(\Omega)}$ in \eqref{weakthm1.1}
with $\norm{\bar{F}(t)}_{L^2(\Omega)}^2 + \hat{f}(t)E(t)$.
\end{itemize}
\end{rem}

\subsect{simple}{Higher time derivatives: an instructive example}

Theorem \ref{weakthm} establishes the existence and uniqueness of (weak) solutions to the model problem. However,
for the linear results to be applicable to non-linear problems, the existence and uniqueness result
of Theorem \ref{weakthm} must be generalized to include solutions with more regularity.
To give insight into our proof strategy for establishing this type of generalization, we first consider a simplified example
that contains the essential ideas, and highlights the major steps involved in
the proof. We note that this type of strategy was employed in \cite{AnderssonOliynyk:2014} to establish
the existence and uniqueness of solutions to wave equations with jump discontinuities.

\bigskip

\subsubsect{tss}{The simplified problem} The simplified IBVP that we consider is
\lalin{ss}{
-\del{t}^2 v + \delta^{ij}\del{i}\del{j}v -v &= f  \hspace{2.0cm} \text{in $\Omega_T$}, \label{ss.1} \\
\nu_i \delta^{ij}\del{j}v &= q \del{t}^2 v + p\del{t}v   \hspace{0.4cm}\text{in $\Gamma_T$}, \label{ss.2} \\
(v,\del{t}v) &= (\vt_0,\vt_1)  \hspace{1.1cm} \text{in $\Omega_0$}, \label{ss.3}
}
where we assume that $q$ and $P$ are constant matrices, i.e. $\del{\alpha}q=\del{\alpha}p = 0$, $q$
verifies the conditions \eqref{qPdef}, $p$ satisfies $p -\chi q \leq 0$ for some constant $\chi$,
and $f\in X_T^{\frac{3}{2}}(\Omega)$.

\subsubsect{tftds}{The formally time differentiated system} The first step in obtaining solutions
to \eqref{ss.1}-\eqref{ss.3} with more regularity is to \textit{formally} differentiate
the equations \eqref{ss.1}-\eqref{ss.2} in time. We do this three times to obtain
\lalin{ftds}{
\delta^{ij}\del{i}\del{j}v_0 -v_0 &= v_2+f_0  \hspace{1.425cm} \text{in $\Omega_T$}, \label{ftds.1} \\
\nu_i \delta^{ij}\del{j}v_0 &= q v_2 + pv_1   \hspace{1.075cm}\text{in $\Gamma_T$}, \label{ftds.2} \\
\delta^{ij}\del{i}\del{j}v_1 -v_1 &= \del{t}v_2+f_1  \hspace{1.125cm} \text{in $\Omega_T$}, \label{ftds.3} \\
\nu_i \delta^{ij}\del{j}v_1 &= q \del{t}v_2 + p v_2   \hspace{0.775cm}\text{in $\Gamma_T$}, \label{ftds.4} \\
-\del{t}^2 v_2 + \delta^{ij}\del{i}\del{j}v_2 -v_2 &= f_2  \hspace{2.175cm} \text{in $\Omega_T$}, \label{ftds.5} \\
\nu_i \delta^{ij}\del{j}v_2 &= q \del{t}^2 v_2 + p\del{t}v_2   \hspace{0.4cm}\text{in $\Gamma_T$}, \label{ftds.6}
}
where
\eqn{tdssA}{
v_\ell = \del{t}^\ell v \AND f_\ell = \del{t}^\ell f.
}
Of course, we could consider taking even more time derivatives provided that $f$ was sufficiently regular, but three
time derivatives will more than suffice to illustrate the method.

\subsubsect{thes}{The elliptic-hyperbolic IBVP} The next step is to view \eqref{ftds.1}-\eqref{ftds.6}
as the following elliptic-hyperbolic IBVP:
\lalin{tdss}{
\delta^{ij}\del{i}\del{j}v_0 -v_0 &= v_2+f_0  \hspace{1.35cm} \text{in $\Omega_T$}, \label{tdss.1} \\
\nu_i \delta^{ij}\del{j}v_0 &= q v_2 + pv_1   \hspace{1.0cm}\text{in $\Gamma_T$}, \label{tdss.2} \\
\delta^{ij}\del{i}\del{j}v_1 -v_1 &= \del{t}v_2+f_1  \hspace{1.05cm} \text{in $\Omega_T$}, \label{tdss.3} \\
\nu_i \delta^{ij}\del{j}v_1 &= q w_3 + p v_2   \hspace{0.95cm}\text{in $\Gamma_T$}, \label{tdss.4} \\
-\del{t}^2 v_2 + \delta^{ij}\del{i}\del{j}v_2 -v_2 &= f_2  \hspace{2.15cm} \text{in $\Omega_T$}, \label{tdss.5} \\
\nu_i \delta^{ij}\del{j}v_2 &= q \del{t}^2 v_2 + p\del{t}v_2   \hspace{0.4cm}\text{in $\Gamma_T$}, \label{tdss.6} \\
(v_2,\del{t}v_2) &= (\vt_2,\vt_3) \hspace{1.4cm}\text{in $\Omega_0$}, \label{tdss.7}\\
\del{t}v_2 &= \wt_3 \hspace{2.1cm}\text{in $\Gamma_0$}, \label{tdss.8}
}
where
\leqn{tdssidata}{
(\vt_2,\vt_3,\wt_3)\in H^1(\Omega)\times L^2(\Omega) \times L^2(\del{}\Omega),
}
$\wt_3$ satisfies $\Pbb_{q(0)} \wt_3 = \wt_3$, the variables $\{u_0,u_1,u_2\}$ are now treated as independent, and
$f_\ell=\del{t}^\ell f$, $\ell=0,1,2$, as defined above.
By Theorem \ref{weakthm}, we know that there exists a unique weak solution
\leqn{tdssA}{
(v_2,w_3) \in \bigcap_{j=0}^1 C^j\bigl([0,T]H^{1-j}(\Omega,\Rbb^N)\bigr)\times C^0\bigl([0,T],L^2(\del{}\Omega,\Rbb^N)\bigr)
}
of \eqref{tdss.5}-\eqref{tdss.8} that is generated by the initial data \eqref{tdssidata} and satisfies the energy estimate
\eqn{tdssB}{
\norm{(v_2(t),w_3(t))}_{E}^2 \lesssim \norm{(v_2(0),w_3(0))}_E^2 + \int_{0}^t
\norm{(v_2(\tau),w_3(\tau))}_{E}^2 + \norm{\del{t}^2f(\tau)}^2_{L^2(\Omega)} \, d\tau, \quad 0\leq t \leq T.
}

Next, we see that \eqref{tdss.3}-\eqref{tdss.4} can be solved for $v_1$ by appealing to Theorem \ref{ellipthmA}.
It clear from the regularity of the solution $(v_2,w_3)$ given by \eqref{tdssA} and Theorem \ref{ellipthmB} that
\leqn{tdssC}{
v_1 \in C^0\bigl([0,T],H^{\frac{3}{2}}(\Omega,\Rbb^N)\bigr).
}
It is important to observe that the $3/2$-derivative spatial regularity is a consequence of the fact that we
have $L^2(\del{}\Omega)$ control of the boundary term $q(t) w_3(t)$ and no better.
Solving \eqref{tdss.1}-\eqref{tdss.2} in
the same manner, we again see by Theorems \eqref{ellipthmA} and \eqref{ellipthmB}, and the regularity of the solutions
\eqref{tdssA} and \eqref{tdssC} that there exists a solution $v_0$ of \eqref{tdss.1}-\eqref{tdss.2}
satisfying
\leqn{tdssD}{
v_0 \in C^0\bigl([0,T],H^{2}(\Omega,\Rbb^N)\bigr).
}
Again, it is important to observe that $v_0$ gains only a $1/2$-derivative in spatial regularity over $v_1$. This is due
to the boundary control $v_2(t)|_{\del{}\Omega}\in H^{1/2}(\Omega)$ for $v_2$, which follows from the Trace Theorem, see Theorem \ref{trace}, and
the interior regularity $v_2(t)\in H^1(\Omega)$.

\subsubsect{sss}{Solutions to the simplified system} Collectively, \eqref{tdssA}, \eqref{tdssC} and \eqref{tdssD} define a solution of the the elliptic-hyperbolic
system \eqref{tdss.1}-\eqref{tdss.8}. However, it not clear if $v_0$, in any sense, solves the original system
\eqref{ss.1}-\eqref{ss.3}.  In order to establish this, we need to verify that the relations
\eqn{sssA}{
\del{t}v_0 = v_1 \AND \del{t}v_1 = v_2
}
hold. Since $v_0$ and $v_1$ only depend continuously on time, it is not yet clear in what sense we can differentiate
them in time. To address this issue, we smooth, in time, using a mollifier $J_\omega$ defined by
\leqn{sssB}{
J_\omega u(t,x) = \frac{1}{\omega}\int_\Rbb \rho\biggl(\frac{t-\tau}{\omega}\biggr)u(\tau,x)\,d\tau
}
where $\rho \in C^\infty_0(\Rbb)$ satisfies $\int_\Rbb \rho(\tau)\,d\tau =1$. Restricting our attention
to $t\in (T_1,T_2)$, where $0<T_1<T_2<T$, and choosing $\omega>0$ small enough so that $J_\omega u(t,x)$, for
$t\in (T_1,T_2)$, only depends on $u(t,x)$, for $t\in (0,T)$, we apply $J_\omega$
to \eqref{tdss.1}-\eqref{tdss.4} to get
\lalin{sssC}{
\delta^{ij}\del{i}\del{j}J_\omega v_0(t) -J_\omega v_0(t) &= J_\omega v_2(t)+ J_\omega f_0(t)  \hspace{0.8cm} \text{in
$\Omega$}, \label{sssC.1} \\
\nu_i \delta^{ij}\del{j}J_\omega v_0(t) &= q J_\omega v_2(t) + p J_\omega v_1(t)
\hspace{0.45cm}\text{in $\del{}\Omega$}, \label{sssC.2} \\
\delta^{ij}\del{i}\del{j}J_\omega v_1(t) - J_\omega v_1(t) &= \del{t} J_\omega v_2(t)+
J_\omega f_1(t)  \hspace{0.5cm} \text{in $\Omega$}, \label{sssC.3} \\
\nu_i \delta^{ij}\del{j} J_\omega v_1(t) &= q J_\omega w_3(t) + p J_\omega v_2(t)   \hspace{0.4cm}
\text{in $\del{}\Omega$}, \label{sssC.4}
}
for $t\in (T_1,T_2)$.

Next, since $(v_2,w_3)$ is a weak solution of \eqref{tdss.6}-\eqref{tdss.8}, it satisfies
\eqn{sssD}{
-\ip{\del{t}v_2}{\del{t}\phi}_{\Omega_T} - \ip{p\del{t}v_2}{\phi}_{\Gamma_T}
+ \ip{q w_3}{\del{t}\phi}_{\Omega_T} = -\ip{v_2}{\phi}_{\Omega_T}-\ip{f_2}{\phi}_{\Omega_T}.
}
Setting $\phi=J_\omega(\psi\eta)$, where $\eta \in C^1(\overline{\Omega},\Rbb^N)$ and $\psi \in C_0^1([0,T])$,
we see that
\leqn{sssE}{
-\ip{\del{t}^2 J_\omega v_2}{\psi\eta}_{\Omega_T} - \ip{p\del{t} J_\omega v_2}{\psi\eta}_{\Gamma_T}
- \ip{q \del{t}^2 J_\omega u_2}{\psi\eta}_{\Gamma_T} = -\ip{J_\omega v_2}{\eta\psi}_{\Omega_T}-
\ip{J_\omega f_2}{\psi\eta}_{\Omega_T},
}
where in deriving this we have used the property \eqref{weaksoldef2a} of weak solutions, which here,
implies that
\leqn{sssF}{
\ip{q w_3}{\del{t}(J_\omega(\psi\eta))}_{\Gamma_T} = -\ip{q v_2}{\del{t}^2(J_\omega(\psi\eta))}_{\Gamma_T}
= -\ip{q \del{t}^2 J_\omega v_2}{\psi\eta}_{\Gamma_T}.
}
Since all of the time-dependent quantities appearing in \eqref{sssE} and \eqref{sssF} are continuous
in time, and the test function $\psi=\psi(t)$ is arbitrary, it
follows that
\eqn{sssG}{
-\ip{\del{t}^2 J_\omega v_2(t)}{\eta}_{\Omega} - \ip{p\del{t} J_\omega v_2(t)}{\eta}_{\del{}\Omega}
- \ip{q \del{t}^2 J_\omega u_2(t)}{\eta}_{\del{}\Omega} = -\ip{J_\omega v_2(t)}{\eta}_{\Omega}-
\ip{J_\omega f_2(t)}{\eta}_{\Omega}
}
and
\eqn{sssFa}{
\ip{q \del{t} J_\omega w_3(t)}{\eta}_{\Omega} = \ip{q \del{t}^2 J_\omega v_2(t)}{\eta}_{\Omega}
}
for $t\in (T_1,T_2)$ and $\eta\in C^1(\overline{\Omega})$. This shows that we can view  $J_\omega v_2(t)$, $T_1<t<T_2$, as a weak solution of
\lalin{sssH}{
\delta^{ij}\del{i}\del{j}J_\omega v_2(t) - J_\omega v_2(t) &= \del{t}^2 J_\omega v_2(t)+
J_\omega f_1(t)  \hspace{0.4cm} \text{in $\Omega$}, \label{sssH.1} \\
\nu_i \delta^{ij}\del{j} J_\omega v_2(t) &= q \del{t}^2 J_\omega v_2 + p J_\omega v_2(t)   \hspace{0.45cm}
\text{in $\del{}\Omega$,} \label{sssH.2}
}
that satisfies
\leqn{sssI}{
q \del{t} J_\omega w_3(t) = q \del{t}^2 J_\omega v_2(t) \hspace{0.4cm} \text{ in $\del{}\Omega$}.
}
We further observe that the equality
\leqn{sssCa}{
J_\omega w_3(t) = \del{t} J_\omega v_2(t) \quad \text{in $\del{}\Omega$,}
}
for $t\in (T_1,T_2)$,follows from \eqref{weaksoldef2a}, which is satisfied by the weak solution $(v_2,w_3)$. This is seen by setting
$\phi=J_\omega(\psi\eta)$ to get $\ip{q J_\omega w_3}{\psi\eta}_{\Gamma_T} = -\ip{q J_\omega v_2}{\psi\eta}_{\Gamma_T}$.
Since $\psi \in C^1_0([0,T])$ is arbitrary it follows that $\ip{q J_\omega w_3(t)}{\eta}_{\del{}\Omega}
= -\ip{q J_\omega v_2(t)}{\eta}_{\del{}\Omega}$, for $t\in (T_1,T_2)$ and all $\eta \in C^1(\overline{\Omega})$,
which proves the assertion.

Differentiating \eqref{sssC.1}-\eqref{sssC.4} in time, which is now possible due to the smoothing,
we see with the help of \eqref{sssCa} and \eqref{sssI} that
$\del{t}J_\omega v_0$ and $\del{t}J_\omega v_1$ satisfy
\lalin{sssJ}{
\delta^{ij}\del{i}\del{j} \del{t}J_\omega  v_0(t) -J_\omega \del{t}v_0(t) &= J_\omega \del{t} v_2(t)+ J_\omega f_1(t)
\hspace{1.05cm} \text{in
$\Omega$}, \label{sssJ.1} \\
\nu_i \delta^{ij}\del{j} \del{t}J_\omega  v_0(t) &= q J_\omega  w_3(t) + p \del{t} J_\omega v_1(t)
\hspace{0.4cm}\text{in $\del{}\Omega$}, \label{sssJ.2} \\
\delta^{ij}\del{i}\del{j}\del{t}J_\omega  v_1(t) - \del{t} J_\omega v_1(t) &= \del{t}^2 J_\omega v_2(t)+
J_\omega f_2(t)  \hspace{1.0cm} \text{in $\Omega$}, \label{sssJ.3} \\
\nu_i \delta^{ij}\del{j} \del{t} J_\omega v_1(t) &= q \del{t}^2 J_\omega v_2(t) + p \del{t} J_\omega v_2(t)   \hspace{0.375cm}
\text{in $\del{}\Omega$}. \label{sssJ.4}
}
for $t\in (T_1,T_2)$.
It then follows from \eqref{sssC.3}-\eqref{sssC.4}, \eqref{sssH.1}-\eqref{sssH.2}, and \eqref{sssJ.1}-\eqref{sssJ.4}
that the differences $J_\omega v_1-\del{t}J_\omega  v_0$ and $J_\omega v_2 - \del{t} J_\omega v_1$
satisfy the homogenous elliptic system
\alin{sssK}{
\delta^{ij}\del{i}\del{j}\bigl(J_\omega v_1(t)-\del{t}J_\omega  v_0(t)\big)
 -\bigl(\del{t}J_\omega v_1(t)- J_\omega \del{t}v_0(t)\bigr) &= 0  \hspace{0.4cm} \text{in
$\Omega$},  \\
\nu_i \delta^{ij}\del{j}\bigl(J_\omega v_1(t)- \del{t}J_\omega  v_0(t)\bigr) &= 0
\hspace{0.4cm}\text{in $\del{}\Omega$}, \\
\delta^{ij}\del{i}\del{j}\bigl(J_\omega v_2(t) - \del{t} J_\omega v_1(t)\bigr)
 - \bigl(J_\omega v_2(t) - \del{t} J_\omega v_1(t)\bigr) &= 0\hspace{0.4cm} \text{in $\Omega$}, \\
\nu_i \delta^{ij}\del{j} \bigl(J_\omega v_2(t) - \del{t} J_\omega v_1(t)\bigr) &= 0   \hspace{0.4cm}
\text{in $\del{}\Omega$}
}
for $t\in (T_1,T_2)$. By uniqueness of solutions to this system, see
Theorem \ref{ellipthmA}, we conclude that
\eqn{sssL}{
J_\omega v_1 =  \del{t} J_\omega v_0 \AND J_\omega v_2 =\del{t} J_\omega v_1 \hspace{0.4cm} \text{ in $(T_1,T_2)\times\Omega$}.
}
Testing these relations with $\phi \in C^1_0([0,T],C^1(\overline{\Omega},\Rbb^N))$, we obtain, for $\omega$ small enough,
\eqn{sssLa}{
\ip{J_\omega\phi}{v_1}_{\Omega_T}=  -\ip{\del{t}\phi}{J_\omega v_0}_{\Omega_T} \AND
\ip{J_\omega\phi}{v_2} =-\ip{\del{t}\phi}{J_\omega v_1}_{\Omega_T}.
}
Sending $\omega \searrow 0$ gives
\eqn{sssLb}{
\ip{\phi}{v_1}_{\Omega_T} =  -\ip{\del{t}\phi}{v_0}_{\Omega_T} \AND
\ip{\phi}{v_2}_{\Omega_T} =-\ip{\del{t}\phi}{v_1}_{\Omega_T}.
}
From this we conclude that
\eqn{sssc}{
 v_1 =  \del{t} v_0 \AND v_2 = \del{t} v_1 \hspace{0.4cm} \text{ in $\Omega_T$}.
}
Thus we have established the following: $v_0 \in C\Xc^{2}_T(\Omega,\Rbb^N)$, $v_0$
solves wave equation \eqref{ss.1}-\eqref{ss.3},
and the pair $(\del{t}^2 v_0,w_3)$ define a weak solution of the twice time differentiated
version of \eqref{ss.1}-\eqref{ss.3}.

%We further note that in the above construction, the initial data was chosen for the highest
%time derivatives, i.e. $(v_2,\del{t}v_2,w_3)|_{t=0} = (\vt_2,\vt_3,\wt_3)\in H^1(\Omega)\times L^2(\Omega)
%\times L^2(\del{\Omega})$. The initial data in the traditional sense is then determined in terms of
%the triple $(\vt_2,\vt_3,\wt_3)$ by setting
%\leqn{sssIDA}{
%(\vt_0,\vt_1) := (v_0,\del{t}v_0)|_{t=0},
%}
%where $v_0$ is the solution to \eqref{ss.1}-\eqref{ss.3} constructed using the method above. By construction,
%we see that
%\leqn{sssIDB}{
%v_0 \in
%}

\subsect{lintdif}{The rescaled system}
We now proceed with the general existence proof for solutions to the model problem with higher regularity following essentially the
same steps taken as for the simplified example from the previous section. All of the technical complications come from
the coupling that arises between the equations
satisfied by time derivatives due to the presence of variable coefficients. Because of this coupling,
the resulting elliptic systems
will not be solvable without the introduction of a smallness condition for certain coefficients. To facilitate
this, we introduce a small parameter $\ep>0$ into the IBVP. Rather than considering the system
\eqref{linB.1}-\eqref{linB.3} directly, we instead consider the system
\lalin{linC}{
\del{\alpha}\bigl(B^{\alpha\beta}\del{\beta}v+M^\alpha\bigr) +\lambda c v &= F \hspace{3.05cm} \text{in $\Omega_{T}$,}
\label{linC.1}\\
\nu_\alpha\bigl(B^{\alpha\beta}\del{\beta} v + M^\alpha) & = \ep q\del{t}^2 v + \ep P \del{t} v + G
 \hspace{0.4cm} \text{in $\Gamma_{T}$, } \label{linC.2}\\
(v,\del{t}v) &= (\vt_0,\vt_1) \hspace{2.2cm} \text{in $\Omega_0$,}\label{linC.3}
%\del{t}v &= \wt_1 \hspace{3.15cm} \text{in $\Gamma_0$}\label{linC.4}
}
where
\lalin{BMdef}{
B^{\alpha\beta} &= \ep^2 \delta^\alpha_0\delta^\beta_0 b^{00} + \ep\delta^\alpha_0\delta^\beta_j b^{0j}
+ \ep \delta^\alpha_i\delta^\beta_0 b^{i0} + \delta^\alpha_i\delta^\beta_j b^{ij}, \notag % \label{BMdef.1}
\\
M^\alpha &= \ep \delta^\alpha_0 L^0 + \delta^\alpha_i L^i. \notag % \label{BMdef.2}
}
Here, we do not assume that $G$ is of the form \eqref{GformA}.

\begin{rem} \label{rescalerem}
The form of the system \eqref{linC.1}-\eqref{linC.2} is consistent with rescaling the spatial coordinates $x=(x^i)$
in \eqref{linB.1}-\eqref{linB.2} according to
$x \longmapsto \ep x$. This fact is exploited in the proof of Theorem \ref{linlocthmC} where it is shown
that the problem of existence and uniqueness for the unscaled system \eqref{linB.1}-\eqref{linB.3} can be reduced to establishing existence and
uniqueness for the rescaled system \eqref{linC.1}-\eqref{linC.3} with $\ep$ chosen suitably small.
\end{rem}

\subsect{eest}{Elliptic estimates}

Formally differentiating \eqref{linC.1}-\eqref{linC.2} $\ell$-times with respect to $t$ for $\ell=0,1,\ldots,2s-1$,
we see, employing the notation \eqref{ft}, that the $v_\ell$, $0\leq \ell \leq 2s-1$,
satisfy
\lalin{linD}{
\del{i}\bigl(b^{i j}_0\del{\beta}v_\ell + \ep\dc^i_\ell v_\ell + \Lc^i_\ell\bigr) +
\ep \ac^i_\ell\del{i}v_\ell + \lambda \cc_\ell v_\ell &= \Fc_\ell \hspace{1.55cm} \text{in $\Omega$,}
\label{linD.1}\\
\nu_i\bigl(b^{i j}_0\del{\beta}v_\ell + \ep \dc^i_\ell v_\ell + \Lc^i_\ell\bigr) & = \ep\hc_\ell v_\ell
+\Gc_\ell
 \hspace{0.4cm} \text{in $\del{}\Omega$ } \label{linD.2}
}
where
\lalin{linE}{
\dc^i_\ell &= \ell b_1^{i0}, \label{linE.1} \\
\ac^i_\ell &= (1+\ell)b^{0i}_1, \label{linE.2} \\
\cc_{\ell} &= c + \frac{\ell(\ell+1)}{2\lambda}b^{00}_2, \label{linE.3} \\
\Lc^i_\ell &= \ep b^{i0}_0 v_{\ell+1} + L^i_\ell+ \ep b^{i0}_\ell v_1 +  \ep \sum_{r=1}^{\ell-2}
\binom{\ell}{r} b^{i0}_{\ell-r}v_{r+1} +b^{ij}_\ell\del{j}v_0  + \sum_{r=1}^{\ell-1}\binom{\ell}{r} b^{ij}_{\ell-r}\del{j} v_r , \label{linE.4} \\
\Fc_\ell &= \ep \Bigl(b^{0 j}_0\del{j} v_{\ell+1} + \ep(\ell+1)b^{00}_1 v_{\ell+1}
+\ep b^{00}_0 v_{\ell+2}\Bigr) + F_\ell - \lambda\biggl(c_\ell v_0 +\sum_{r=1}^{\ell-1} \binom{\ell}{r} c_{\ell-r}v_r \biggr)
\notag \\
&\hspace{0.9cm} -L^0_{\ell+1} -\ep\biggl(b^{0j}_{\ell+1}\del{j}v_0 + b^{0j}_\ell \del{j}v_1
+ \sum_{r=1}^{\ell-1}\binom{\ell}{r} b^{0j}_{\ell+1-r}\del{j}v_r + \sum^{\ell-2}_{r=1}
\binom{\ell}{r} b^{0j}_{\ell-r}\del{j}v_{r+1}\biggr) \notag \\
&\hspace{3.0cm}-\ep^2\biggl(b^{00}_{\ell+1}v_1 + b^{00}_\ell v_2
+ \sum_{r=1}^{\ell-2}\binom{\ell}{r} b^{00}_{\ell+1-r}v_{r+1} + \sum^{\ell-3}_{r=1}
\binom{\ell}{r} b^{00}_{\ell-r}\del{j}v_{r+2}\biggr),
\label{linE.5} \\
\hc_\ell &= \frac{\ell(\ell-1)}{2}q_2 + \ell P_1, \label{linE.6}
\intertext{and}
\Gc_\ell &= \ep q_0 v_{\ell+2} +\ep\bigl(\ell q_1 + P_0\bigr)v_{\ell+1} + G_\ell
\notag \\
&\hspace{1.0cm} + \ep^2 q_\ell v_2 + \ep \sum_{r=1}^{\ell-3}\binom{\ell}{r} q_{\ell-r} v_{r+2} + \ep P_\ell v_1
+ \ep\sum_{r=1}^{\ell-2} \binom{\ell}{r}P_{\ell-r}v_{r+1}. \label{linE.7}
}
Just as in the simplified example from Section \ref{simple}, we interpret
\eqref{linD.1}-\eqref{linD.2} as
a system of elliptic equations, which will allow us to use elliptic estimates
to bound the
$v_\ell$, $0\leq \ell \leq 2s-1$. Here, we use
$v_{2s+1}$ as shorthand notation for a pair $(v_{2s+1},w_{2s+1})\in L^2(\Omega)\times L^2(\del{}\Omega)$,
where we abuse notation and denote $qw_{2s+1}$ as $qv_{2s+1}|_{\del{}\Omega}$.

To prepare for the elliptic estimates, we first
estimate the coefficients appearing in \eqref{linD.1}-\eqref{linD.2}, which
we collect in the following lemma. In the following, we freely
employ the vector notation \eqref{fvect}, e.g. $\Pv_{\ell}=(P_0,P_1,\ldots,P_\ell)$.
\begin{lem}\label{elem}
Suppose $s>n/2$, $\st \in [0,s]$, $0\leq\ep\leq 1$, and $(s,\st)=(k/2,\kt/2)$ for $k,\kt\in \Zbb$.
Then the following estimates hold:
\begin{enumerate}[(i)]
\item
\alin{elem0}{
\norm{\Gc_\ell}_{H^{\st-\frac{\ell}{2}}(\Omega)}
&\lesssim \ep \Bigl(\norm{q_0}_{H^s(\Omega)}\norm{v_{\ell+2}}_{H^{\st+1-\frac{\ell+2}{2}}(\Omega)}
+ \bigl( \norm{q_{1}}_{H^s(\Omega)}+\norm{P_0}_{H^s(\Omega)}\bigr)
\norm{v_{\ell+1}}_{H^{\st+1-\frac{\ell+1}{2}}(\Omega)}\Bigr)\\
&\hspace{0.4cm} + \norm{G_\ell}_{H^{\st-\frac{\ell}{2}}(\Omega)}+ \bigl(\norm{\Pv_\ell}_{X^{s,\ell}}
+ \norm{\qv_\ell}_{X^{s,\ell}} + \norm{q_1}_{H^s(\Omega)}\bigr)\norm{\vv_{\ell-1}}_{X^{\st+1,\ell-1}},
}
for $0\leq \ell \leq 2\st-2$,
and
\alin{elem2a}{
&\norm{\Gc_{2\st-1}|_{\del{}\Omega}}_{L^{2}(\del{}\Omega)} \lesssim  \ep\Bigl(\norm{q_0 v_{2\st+1}}_{L^2(\del{}\Omega)}
+ \bigl( \norm{q_{1}}_{H^s(\Omega)}+\norm{P_0}_{H^s(\Omega)}\bigr)
\norm{v_{2\st}}_{H^1(\Omega)}\Bigr)
+ \norm{G_{2\st-1}|_{\del{}\Omega}}_{L^2(\del{}\Omega)} \\
& \hspace{1.0cm}+ \Bigl(\norm{\Pv_{2s-2}}_{X^{s,2s-2}}+\norm{P_{2s-1}}_{L^2(\del{}\Omega)}
+ \norm{\del{t}\qv_{2s-2}}_{X^{s,2s-2}} \Bigr)\norm{\vv_{2\st-2}}_{X^{\st+1,2\st-2}},
}

\item
\alin{elem3a}{
\norm{\Lc^i_\ell}_{H^{\st-\frac{\ell}{2}}(\Omega)}
\lesssim \ep\norm{b_0}_{H^s(\Omega)}\norm{v_{\ell+1}}_{H^{\st+1-\frac{\ell+1}{2}}(\Omega)}+ \norm{L^i_\ell}_{H^{\st-\frac{\ell}{2}}(\Omega)} + \norm{\bv_{\ell}}_{X^{s,\ell}}
\norm{\vv_{\ell-1}}_{X^{\st+1,\ell-1}},
}
for $0\leq \ell \leq 2\st-2$, and

\alin{elem4}{
\norm{\Lc^i_{2\st-1}}_{L^2(\Omega)}
\lesssim \ep\norm{b_0}_{H^s(\Omega)}\norm{v_{2s}}_{H^1(\Omega)}+\norm{L^i_{2\st-1}}_{L^2(\Omega)} + \norm{\bv_{2s-1}}_{X^{s,2s-1}}
\norm{\vv_{2\st-2}}_{X^{\st+1,2\st-2}},
}
\item and
\alin{elem5a}{
\norm{\Fc_\ell}_{H^{\st-1-\frac{\ell}{2}}(\Omega)}
&\lesssim \ep\bigl(\norm{b_0}_{H^s(\Omega)} + \norm{b_1}_{H^{s-\frac{1}{2}}(\Omega)}\bigr)
\bigl(\norm{v_{\ell+1}}_{H^{\st+1-\frac{\ell+1}{2}}(\Omega)}+\norm{v_{\ell+2}}_{H^{\st+1-\frac{\ell+2}{2}}(\Omega)}
\bigr)\\
& + \norm{F_\ell}_{H^{\st-1-\frac{\ell}{2}}(\Omega)}
+ \norm{L^0_{\ell+1}}_{H^{\st-1-\frac{\ell}{2}}(\Omega)}
+\bigl(\norm{\bv_{\ell+1}}_{X^{s,\ell+1}}+\norm{\cv_{\ell}}_{X^{s,\ell}}\bigr)
\norm{\vv_{\ell-1}}_{X^{\st+1,\ell-1}}
}
for $0\leq \ell \leq 2\st-2$,
and
\alin{elem5}{
&\norm{\Fc_{2\st-1}}_{L^2(\Omega)}
\lesssim \ep\bigl(\norm{b_0}_{H^s(\Omega)} + \norm{b_1}_{H^{s-\frac{1}{2}}(\Omega)}\bigr)
\bigl(\norm{v_{2\st}}_{H^{1}(\Omega)}+\norm{v_{2\st+1}}_{L^2(\Omega)}
\bigr)+ \norm{F_{2\st-1}}_{L^2(\Omega)} \\
&\hspace{2.0cm}
+ \norm{L^0_{2\st}}_{L^2(\Omega)}+\bigl(\norm{\bv_{2s}}_{X^{s,2s}}+\norm{\cv_{2s-1}}_{X^{s,2s-1}}\bigr)
\norm{\vv_{2\st-2}}_{X^{\st+1,2\st-2}}.
}
\end{enumerate}
\end{lem}
\begin{proof}
From the assumptions $s=k/2 > n/2$ and $\st=\kt/2\in [0,s]$, it follows directly
from the fractional multiplication inequality, Theorem \ref{calcpropB},
that the estimates
\alin{elem2}{
\norm{q_0 v_{\ell+2}}_{H^{\st-\frac{\ell}{2}}(\Omega)}
&\lesssim \norm{q_0}_{H^s(\Omega)}\norm{v_{\ell+2}}_{H^{\st+1-\frac{\ell+2}{2}}(\Omega)}, \\
\norm{\bigl(\ep\ell q_1 + P_0\bigr)v_{\ell+1}}_{H^{\st-\frac{\ell}{2}}(\Omega)}
&\lesssim \bigl( \norm{q_{1}}_{H^s(\Omega)}+\norm{P_0}_{H^s(\Omega)}\bigr)
\norm{v_{\ell+1}}_{H^{\st+1-\frac{\ell+1}{2}}(\Omega)}, \\
\norm{q_\ell v_2}_{H^{\st-\frac{\ell}{2}}(\Omega)} &\lesssim
\norm{q_\ell}_{H^{s-\frac{\ell}{2}}(\Omega)} \norm{v_2}_{H^{\st+1-1}(\Omega)}, \\
\norm{q_{\ell-r}v_{r+2}}_{H^{\st-\frac{\ell}{2}}(\Omega)}
&\lesssim \norm{q_{\ell-r}}_{H^{s-\frac{\ell-r}{2}}(\Omega)}\norm{v_{r+2}}_{H^{\st+1-\frac{r+2}{2}}(\Omega)},
&& 1\leq r \leq \ell-3,\\
\norm{P_\ell v_1}_{H^{\st-\frac{\ell}{2}}(\Omega)} &\lesssim \norm{P_\ell}_{H^{s-\frac{\ell}{2}}(\Omega)}
\norm{v_1}_{H^{\st+1-\frac{1}{2}}(\Omega)},
\intertext{and}
\norm{P_{\ell-1} v_{r+1}}_{H^{\st-\frac{\ell}{2}}(\Omega)} &\lesssim
\norm{P_{\ell-r}}_{H^{s-\frac{\ell-r}{2}}(\Omega)}
\norm{v_{r+1}}_{H^{\st+1-\frac{r+1}{2}}(\Omega)}, && 1\leq r\leq \ell-2,
}
hold for $0\leq \ell \leq 2 \st-1$. Using these estimates, it is clear from \eqref{linE.7}
that $\Gc_\ell$ can be estimated by
\alin{elem3}{
\norm{\Gc_\ell}_{H^{\st-\frac{\ell}{2}}(\Omega)}
&\lesssim \ep \Bigl(\norm{q_0}_{H^s(\Omega)}\norm{v_{\ell+2}}_{H^{\st+1-\frac{\ell+2}{2}}(\Omega)}
+ \bigl( \norm{q_{1}}_{H^s(\Omega)}+\norm{P_0}_{H^s(\Omega)}\bigr)
\norm{v_{\ell+1}}_{H^{\st+1-\frac{\ell+1}{2}}(\Omega)}\Bigr)\\
& \hspace{0.6cm} + \norm{G_\ell}_{H^{\st-\frac{\ell}{2}}(\Omega)} + \bigl(\norm{\Pv_\ell}_{X^{s,\ell}}
+ \norm{\qv_\ell}_{X^{s,\ell}} + \norm{q_1}_{H^s(\Omega)}\bigr)\norm{\vv_{\ell-1}}_{X^{\st+1,\ell-1}}
}
for $0\leq \ell \leq 2\st -2$, which proves the first estimate. The remaining estimates
can be established in a similar fashion.
\end{proof}

\begin{prop} \label{lineexist}
Suppose $s > n/2+1$, $s=k/2$ for $k\in \Zbb_{\geq 0}$, $m\in \{0,1,2,\ldots,2s-1\}$, and
$(v_{m+1},v_{m+2})$ satisfy
\gath{lineexist1}{
(v_{m+1},v_{m+2}) \in
H^{s+1-\frac{m+1}{2}}(\Omega,\Rbb^N)\times H^{s+1-\frac{m+2}{2}}(\Omega,\Rbb^N), \quad \text{if $0\leq m \leq 2s-2$}
\intertext{and}
(v_{2s},v_{2s+1},q_0 v_{2s+1}|_{\del{}\Omega}) \in H^1(\Omega,\Rbb^N)\times L^2(\Omega,\Rbb^N)\times L^2(\del{}\Omega,\Rbb^N) \quad \text{if $m=2s-1$.}
}
Then there exist constants
\gath{lineexist1ab}{
\lambda^*=\lambda^*(\sigma,\mu) \geq 1 \AND
\delta^*= \delta^*\Bigl(\kappa_1,\norm{b^{i0}_1}_{H^{s-\frac{1}{2}}(\Omega)},
\norm{q_2}_{H^{s-1}(\Omega)},\norm{P_1}_{H^{s-\frac{1}{2}}(\Omega)}\Bigr) \geq 1,
} such that for each
$(\lambda,\ep) \in  [\lambda^*,\infty)\times [0,\frac{1}{\delta^*}]$ there exists
a unique solution
\eqn{lineexist2}{
\vv_m=(v_0,v_1,\ldots,v_m) \in X^{s+1,m}(\Omega,\Rbb^N)
}
to the sequence of equations \eqref{linD.1}-\eqref{linD.2} for $0\leq \ell \leq m$. Moreover, this solution satisfies the estimate
\eqn{lineexist3}{
\norm{\vv_m}_{X^{s+1,m}} \leq C\bigl(\norm{\Fv_m}_{X^{s-1,m}}+ \norm{\Lv^0_{m+1}}_{X^{s-1,m}}+\norm{\Gv_m}_{X^{s,m}}+\norm{\Lv^i_{m}}_{X^{s,m}}
+ \ep V_m\bigr]\bigr)
}
for $0\leq m \leq 2s-2$ and
\alin{lineexist3a}{
&\norm{\vv_{2s-1}}_{X^{s+1,2s-1}} \leq C\bigl(\norm{\Fv_{2s-2}}_{X^{s-1,2s-2}}+\norm{F_{2s-1}}_{L^2(\Omega)}+ \norm{\Lv^0_{2s-1}}_{X^{s-1,2s-2}}+ \norm{\Lv^0_{2s}}_{L^2(\Omega)}\\
&\hspace{1.0cm} +\norm{\Gv_{2s-2}}_{X^{s,2s-2}}+ \norm{G_{2s-1}|_{\del{}\Omega}}_{L^2(\del{}\Omega)}
+\norm{\Lv^i_{2s-1}}_{X^{s,2s-1}}
+ \ep V_{2s-1}\bigr]\bigr)
}
if $m=2s-1$,
where $\{\Lv^0_{m+1},\Lv^i_m,\Fv_m,\Gv_m\}$ are defined using the vector notation \eqref{fvect},
\eqn{lineexist4}{
V_m = \begin{cases} \norm{v_{m+1}}_{H^{s+1-\frac{m+1}{2}}(\Omega)} + \norm{v_{m+2}}_{H^{s+1-\frac{m+2}{2}}(\Omega)} & \text{if $0\leq m\leq 2s-2$} \\
% \norm{v_{m+1}}_{H^{\frac{3}{2}}(\Omega)} + \norm{v_{m+2}}_{H^1(\Omega)} & \text{if $m =2s-2$} \\
  \norm{v_{m+1}}_{H^{1}(\Omega)} + \norm{v_{m+2}}_{L^2(\Omega)}+\norm{q_0 v_{m+2}}_{L^{2}(\del{}\Omega)} & \text{if $m=2s-1$}
\end{cases}
}
and
\eqn{lineexist5}{
C = C\bigl(\kappa_1,\mu,\sigma,\lambda,\rho\Bigr)
}
with
\eqn{lineexists5aa}{
\rho=\Bigl(\norm{\bv_{2s}}_{X^{s,2s}},\norm{\cv_{2s-1}}_{X^{s,2s-1}},
\norm{\Pv_{2s-2}}_{X^{s,2s-2}},\norm{\qv_{2s-2}}_{X^{s,2s-2}},\norm{\del{t}\qv_{2s-2}}_{X^{s,2s-2}},
\norm{P_{2s-1}}_{L^2(\del{}\Omega)}\Bigr).
}
\end{prop}
\begin{proof}
We use proof by induction.

\smallskip

\noindent \underline{Base case}: From Theorems \ref{ellipthmA} and \ref{ellipthmB}, Lemma \ref{elem}, and the Sobolev
inequalities  (Theorems \ref{ISobolev} and \ref{FSobolev}),
we see that there exists constants\footnote{We fix here $\delta^*$ to be the maximum of the constants computed by setting $a^i=\ac^i_\ell$, $d^i = \dc^i_\ell$, and
$h=h_\ell$ in Theorem \ref{ellipthmA} for $\ell=0,1,\ldots 2s-1$.}
\eqn{lineexist6}{
\delta^*= \delta^*\Bigl(\kappa_1,\norm{b^{i0}_1}_{H^{s-\frac{1}{2}}(\Omega)},
\norm{q_2}_{H^{s-1}(\Omega)},\norm{P_1}_{H^{s-\frac{1}{2}}(\Omega)}\Bigr) \geq 1 \AND \lambda^*=\lambda^*(\sigma,\mu) \geq 1
}
such that \eqref{linD.1}-\eqref{linD.2} has a unique solution $\vv_0=v_0\in X^{s+1,0}(\Omega,\Rbb^N)=H^{s+1}(\Omega,\Rbb^N)$
for $\ell=0$ and $(\lambda,\ep)\in [\lambda^*,\infty)\times [0,\frac{1}{\delta^*}]$ that satisfies the bound
\eqn{lineexist7}{
\norm{v_0}_{X^{s+1,0}} \leq C_0\bigl(\norm{\Fv_0}_{X^{s-1,0}}+ \norm{\Lv^0_{1}}_{X^{s-1,0}}+\norm{\Gv_0}_{X^{s,0}}+ \norm{\Lv^i_{0}}_{X^{s,0}}
+ \ep V_0\bigr]\bigr)
}
where
\leqn{lineexist8}{
C_0 = C_0\bigl(\kappa_1,\mu,\sigma,\lambda,\rho\bigr)
}
with $\rho$ and $V_0$ as defined above in the statement of the proposition.

\bigskip

\noindent \underline{Induction hypothesis}: With the base case covered, we fix $m\in \{0,1,\ldots,2s-3\}$ and assume that the system \eqref{linD.1}-\eqref{linD.2}, $0\leq \ell\leq m$, has for $(\lambda,\ep)\in [\lambda^*,\infty)\times [0,\frac{1}{\delta}]$ with $\delta\geq \delta^*$
a unique solution $\vv_{m} \in X^{s+1,m}(\Omega,\Rbb^N)$ satisfying the bound
\leqn{lineexist10}{
\norm{\vv_m}_{X^{s+1,m}} \leq C_m\bigl(\norm{\Fv_m}_{X^{s-1,m}}+ \norm{\Lv^0_{m+1}}_{X^{s-1,m}}+\norm{\Gv_m}_{X^{s,m}}+ \norm{\Lv^i_{m}}_{X^{s,m}}
+ \ep V_m\bigr]\bigr)
}
for some constant $C_m$ of the form \eqref{lineexist8} and $V_m$ is as defined in the statement of
the proposition.

\bigskip

\noindent \underline{Induction step}: Appealing again to Theorems \ref{ellipthmA} and \ref{ellipthmB}, we see using
 Lemma \ref{elem}, and the Sobolev inequalities that for $(\lambda,\ep)\in [\lambda^*,\infty)\times [0,\frac{1}{\delta^*}]$
the BVP \eqref{linD.1}-\eqref{linD.2} with $\ell=m+1$ has a unique solution
$v_{m+1}\in H^{s+1-\frac{m+1}{2}}(\Omega,\Rbb^N)$ satisfying
\lalin{lineexist11}{
&\norm{v_{m+1}}_{H^{s+1-\frac{m+1}{2}}(\Omega)} \leq  c_{m+1}\bigl(\norm{\Fv_{m+1}}_{X^{s-1,{m+1}}}+ \norm{\Lv^0_{m+2}}_{X^{s-1,m+1}} \notag \\
&\hspace{3.0cm} +\norm{\Gv_{m+1}}_{X^{s,{m+1}}}+ \norm{\Lv^i_{m+1}}_{X^{s,{m+1}}}
+ \norm{v_{m}}_{X^{s+1,m}} + \ep V_{m+1}\bigr) \label{lineexist11.1}
}
where $c_{m+1}$ is a constant of the form \eqref{lineexist8}. Fixing $\delta \geq 2 C_m c_{m+1}$,
\eqref{lineexist10} and \eqref{lineexist11.1} imply that
the estimate
\alin{lineexist12}{
&\norm{v_{m+1}}_{H^{s+1-\frac{m+1}{2}}(\Omega)}  \leq C_{m+1}\bigl(\norm{\Fv_{m+1}}_{X^{s-1,{m+1}}}+ \norm{\Lv^0_{m+2}}_{X^{s-1,m+1}} \\
&\hspace{5.0cm} +\norm{\Gv_{m+1}}_{X^{s,{m+1}}}+\norm{\Lv^i_{m+1}}_{X^{s,{m+1}}}
 + \ep V_{m+1}\bigr)
}
holds for all $(\lambda,\ep) \in [\lambda^*,\infty)\times [0,\frac{1}{\delta}]$
where $C_{m+1}$ is again a constant of the form \eqref{lineexist8}.
Combining this estimate with \eqref{lineexist10} yields the desired estimate
\lalin{lineexist13}{
\norm{\vv_{m+1}}_{X^{s+1,m+1}} &\leq C_{m+1}\bigl(\norm{\Fv_{m+1}}_{X^{s-1,{m+1}}}+ \norm{\Lv^0_{m+2}}_{X^{s-1,m+1}} \notag\\
&\hspace{1.0cm}+ \norm{\Gv_{m+1}}_{X^{s,m+1}}+ \norm{\Lv^i_{m+1}}_{X^{s,m+1}}
+ \ep V_{m+1}\bigr),\label{lineeexist13.1}
}
which holds for $m=0,1,\ldots,2s-3$.

\noindent \underline{Final estimate}: For $m=2s-1$, the final estimate
\alin{lineexist3a}{
&\norm{\vv_{2s-1}}_{X^{s+1,2s-1}} \leq C\bigl(\norm{\Fv_{2s-2}}_{X^{s-1,2s-2}}+\norm{F_{2s-1}}_{L^2(\Omega)}+ \norm{\Lv^0_{2s-1}}_{X^{s-1,2s-2}}+ \norm{\Lv^0_{2s}}_{L^2(\Omega)}\\
&\hspace{0.3cm} +\norm{\Gv_{2s-2}}_{X^{s,2s-2}}+ \norm{G_{2s-1}|_{\del{}\Omega}}_{L^2(\del{}\Omega)}+\norm{\Lv^i_{2s-1}}_{X^{s,2s-1}}
+ \ep V_{2s-1}\bigr]\bigr)
}
follows using similar arguments as above from the estimate \eqref{lineeexist13.1}, Lemma \ref{elem},
and the elliptic theory from Appendix \ref{elliptic}.
\end{proof}

\subsect{linlocmodel}{Existence and uniqueness for the model problem}
Following the analysis of the simplified problem contained in Section \ref{simple}, the first step to
obtaining existence and uniqueness for the rescaled model problem  \eqref{linC.1}-\eqref{linC.3}
is to consider the following elliptic-hyperbolic IBVP, which is obtained from formally
differentiating the
rescaled model problem in time:
\lalin{linF}{
\del{i}\bigl(b^{i j}\del{j}v_\ell + \ep\dc^i_\ell v_\ell + \Lc^i_\ell\bigr) +
\ep \ac^i_\ell\del{i}v_\ell + \lambda \cc_\ell v_\ell &= \Fc_\ell \hspace{1.655cm} \text{in $\Omega$,}
\label{linF.1}\\
\nu_i\bigl(b^{i j}\del{j}v_\ell + \ep\dc^i_\ell v_\ell + \Lc^i_\ell\bigr) & = \ep\hc_\ell v
+\Gc_\ell
 \hspace{0.675cm} \text{in $\del{}\Omega$, } \label{linF.2}\\
\del{\alpha}\bigl(B^{\alpha\beta}\del{\beta}v_{2s}+\Mcal^\alpha_{2s}\bigr) +\lambda c v_{2s} &= \Fc_{2s} \hspace{1.525cm} \text{in $\Omega_T$,}
\label{linF.3}\\
\nu_\alpha\bigl(B^{\alpha\beta}\del{\beta} v_{2s} + \Mcal^\alpha_{2s})- \ep^2q\del{t}^2 v_{2s} -\ep\Pc \del{t} v_{2s} & =  \Gc_{2s}
\hspace{1.6cm} \text{in $\Gamma_T$, } \label{linF.4}\\
%\del{t} \Mcal_{2s-1}^\alpha &= \Nc_{2s-1} \hspace{1.1cm}  \text{in $\Omega_T$,}    \label{linF.5}\\
(v_{2s},\del{t}v_{2s}) &= (\vt_{2s},\vt_{2s+1}) \hspace{0.4cm} \text{in $\Omega_0$,}\label{linF.5}\\
\Pbb_{q}\del{t}v_{2s} & = \wt_{2s+1} \hspace{1.225cm} \text{in $\Gamma_0$}\label{linF.6}
%\Mcal_{2s-1} & = \widetilde{\Mcal}_{2s-1} \hspace{1.05cm} \text{in $\{0\}\times \Omega$,}\label{linF.7}
}
where $0\leq \ell \leq 2s-1$, $\Pbb_q$ , as defined previously, is the projection onto $\text{ran}(q)$,
\alin{Mcalt}{
\Mcal_{2s}^\alpha &= 2s B^{\alpha 0}_1 v_{2s}
 + \sum_{r=0}^{2s-2}\binom{2s}{r}B^{\alpha 0}_{2s-r}v_{r+1} +
\sum_{r=0}^{2s-1}\binom{2s}{r}B^{\alpha i}_{2s-r}\del{i}v_{r}
+ M_{2s}, \\
\Pc &=  P+ \ep 2s\del{t}q, \\
\Fc_{2s} &= F_{2s} -\lambda\sum_{r=0}^{2s-1}\binom{2s}{r}c_{2s-r}v_{r},\\
\Gc_{2s} &= G_{2s} +\ep\biggl(P_{2s}v_1+\sum_{r=1}^{2s-2}\binom{2s}{r}P_{2s-r}v_{r+1}
 + 2s P_1 v_{2s}\biggr)\\
 &\hspace{3.0cm} + \ep^2 \biggl(
 \sum_{r=0}^{2s-3}\binom{2s-1}{r} q_{2s-r}v_{r+2}+  \frac{2s(2s-1)}{2}q_2 v_{2s}\biggr),
}
and we are using the notation
$v_{2s+1}=\del{t} v_{2s}$ and $q v_{2s+1}|_{\Gamma_T} = qw_{2s+1}$ while, otherwise,
treating the $v_\ell$, $0\leq \ell \leq 2s-1$, as independent variables. Here, we assume that the coefficients, $B^{\alpha\beta}$,
$P$, $q$, etc., are time dependent maps, and in the following, we employ the
notation \eqref{ft} for time derivatives of these coefficients, and \eqref{fvect} and \eqref{dtfvect}
for collections of time derivatives.

\begin{lem} \label{linloclem}
Suppose $s>n/2+1$, $0\leq\ep\leq 1$, and $s=k/2$ for $k\in \Zbb_{\geq 0}$.
Then the following estimates hold:
\alin{linloclem1}{
\norm{\Mcal_{2s}^i}_{L^2(\Omega)} &\lesssim  \norm{L^i_{2s}}_{L^2(\Omega)}+ \norm{b}_{E^{s}}\norm{\vv_{2s}}_{X^{s+1,2s}},\\
\norm{\del{t}\Mcal_{2s}}_{L^2(\Omega)}& \lesssim \norm{L_{2s+1}}_{L^2(\Omega)} +
\Bigl(\norm{b}_{E^{s}}+\norm{\del{t}b}_{E^{s}}\Bigr)
\Bigl(\norm{\vv_{2s}}_{X^{s+1,2s}}+\norm{\del{t}v_{2s}}_{L^2(\Omega)}\Bigr) , \\
\norm{\Fc_{2s}}_{L^2(\Omega)} &\lesssim  \norm{F_{2s}}_{L^2(\Omega)}
+ \lambda\norm{c}_{E^{s}}\norm{\vv_{2s-1}}_{X^{s+1,2s-1}}
\intertext{and}
\norm{\Gc_{2s}-(G_{2s}+\ep P_{2s}v_1)}_{H^1(\Omega)}&\lesssim \Bigl(\norm{\del{t}^2\qv_{2s-2}}_{X^{s,2s-2}}
+\norm{\del{t}\Pv_{2s-2}}_{X^{s,2s-2}}\Bigr)\norm{\vv_{2s}}_{X^{s+1,2s}}.
}
\end{lem}
\begin{proof}
The proof of this lemma follows from similar arguments used to prove the estimates
from Lemma \ref{elem}. We omit the details.
\end{proof}

With the preliminary estimates out of the way, we are now ready to prove
the existence and uniqueness of solutions to the elliptic-hyperbolic IBVP \eqref{linF.1}-\eqref{linF.6}.

\begin{rem} To simplify the statement of the following existence and uniqueness theorems, we will
not explicitly state the spaces in which the coefficients, i.e. $F$, $L^\alpha$, $b^{\alpha\beta}$, etc.,
lie. This will be clear from the various norms on the coefficients in that the coefficients will need
to lie in the appropriate spaces defined in Section \ref{funct} so that all of the norms that appears make sense.
\end{rem}

\begin{thm} \label{linlocthmA}
Suppose $s>n/2+1$, $s=k/2$ for $k\in \Zbb$, $\ep_0>0$, the coefficients $\{b^{\alpha\beta},c,q\}$ satisfy \eqref{qPdef}-\eqref{coerc} for constants $\kappa_0,\kappa_1,\sigma >0$
and $\mu \geq 0$, there exists a function $\chi=\chi(t,x)$ such that
\eqn{linlocthmAa}{
\ep P + \ep\bigl(2s-\Threehalfs\bigr)\del{t}q - \chi q \leq 0 \hspace{0.4cm} \text{in $\Gamma_T$}
}
for $0<\ep\leq \ep_0$,
$\vt_{2s}\in H^1(\Omega,\Rbb^N)$, $\vt_{2s+1}\in L^2(\Omega,\Rbb^N)$,
$\wt_{2s+1}\in L^2(\del{}\Omega,\Rbb^N)$ satisfies $\Pbb_{q(0)} \wt_{2s+1} = \wt_{2s+1}$,
\eqn{linlocthmAaa}{
\del{t}^{2s}P= k^i_1\del{i}\theta_1 + g_1 \AND
\del{t}^{2s}G= k^i_2\del{i}\theta_2 + g_2
}
where $\nu_ik^i_a=0$ for $a=1,2$.
Then there exist constants $\lambda^*=\lambda^*(\sigma,\mu) \geq 1$ and
\gath{lineexist1a}{
\delta^*= \delta^*\Bigl(\kappa_1,\sup_{0\leq t\leq T}\norm{\del{t}b^{i0}(t)}_{H^{s-\frac{1}{2}}(\Omega)},
\sup_{0\leq t\leq T}\norm{\del{t}^2 q(t)}_{H^{s-1}(\Omega)},
\sup_{0\leq t\leq T}\norm{\del{t}P(t)}_{H^{s-\frac{1}{2}}(\Omega)}\Bigr) \geq 1,
}
such that for each $(\lambda,\ep) \in [\lambda^*,\infty)\times \bigl(0,\min\bigl\{\frac{1}{\delta^*},\ep_0\bigr\}\bigr]$ there exists
a unique solution
\eqn{linlocthmA2}{
(\vv_{2s-1},v_{2s},w_{2s+1}) \in C^0\bigl([0,T],X^{s+1,2s-1}(\Omega,\Rbb^N)\bigr) \times
\bigcap_{j=0}^1 C^j\bigl([0,T],H^{1-j}(\Omega,\Rbb^N)\bigr) \times C([0,T],L^2(\del{}\Omega,\Rbb^N)),
}
to \eqref{linF.1}-\eqref{linF.6}.
Moreover, by choosing $\ep_0>0$ small enough,
$(\vv_{2s-1},v_{2s},w_{2s+1})$ satisfies the energy estimate
\lalin{linlocthmA3}{
&\nnorm{(\vv_{2s-1}(t),v_{2s}(t),w_{2s+1}(t))}_{s+1}^2 \leq C\biggl[ \nnorm{(\vv_{2s-1}(0),v_{2s}(0),w_{2s+1}(0))}_{s+1}^2
+ \alpha_0 \notag \\
&\hspace{5.0cm} + \int_0^t \alpha_1(\tau)\Bigl( \nnorm{(\vv_{2s-1}(\tau),v_{2s}(\tau),w_{2s+1}(\tau))}_{s+1}^2
+ \alpha_2(\tau)\Bigr) \, d\tau \biggr] \notag % \label{linlocthmA3.1}
}
for $(t,\ep)\in [0,T]\times (0,\ep_0]$, where $\alpha_0 = \norm{F(0)}_{E^{s-1}}^2 + \norm{\del{t}^{2s-1}F(0)}_{L^2(\Omega)}^2
+ \alpha_2(0)$,
\eqn{linlocthmA4}{
\nnorm{(\vv_{2s-1}(t),v_{2s}(t),w_{2s+1}(t))}_{s+1}^2 = \norm{\vv_{2s-1}(t)}_{X^{s+1,2s-1}}^2 + \norm{(v_{2s}(t),w_{2s+1}(t))}_E^2,
}
\lalin{linlocthmA5}{
\alpha_1(t) & = 1+\norm{\chi(t)}_{H^{s}(\Omega)}^2+\norm{\del{t}b(t)}_{E^{s}}^2+ \norm{c(t)}_{E^{s}}^2+
\norm{\del{t}P(t)}_{E^{s,2s-2}}^2  \notag \\
&\hspace{2.3cm}+\norm{\del{t}^2q(t)}_{E^{s,2s-2}}^2
+\norm{g_1(t)}_{\Ec^1}^2+\norm{\theta_1(t)}_{H^1(\Omega)}^2+\norm{\theta_1(t)}_{H^1(\Omega)}^2, \label{linlocthmA5.1}\\
\alpha_2(t) &= \norm{\del{t}F(t)}_{E^{s-1}}^2+ \norm{\del{t}^{2s}F(t)}_{L^2(\Omega)}^2+ \norm{L(t)}_{E^{s}}^2+\norm{\del{t}L(t)}_{E^{s}}^2
+\norm{G(t)}_{E^{s,2s-2}}^2 \notag \\
&\hspace{2.2cm} +\norm{\del{t}G(t)}_{E^{s,2s-2}}^2
+ \norm{g_2(t)}_{\Ec^1(\Omega)}^2
+\norm{\theta_1(t)}_{H^1(\Omega)}^2+\norm{\theta_1(t)}_{H^1(\Omega)}^2 \label{linlocthmA5.2},
}
and
\eqn{linlocthmA6}{
C = C\bigl(\kappa_0,\kappa_1,\mu,\sigma,\gamma,\lambda,\rho\bigr)
}
with
\lalin{linlocthmA7a}{
\rho &= \norm{b}_{X_T^{s}}+ \norm{c}_{X_T^{s,2s-1}}+ \norm{P}_{X_T^{s,2s-2}}
+
\norm{\del{t}^{2s-1}P}_{L^\infty([0,T],L^2(\del{}\Omega))}\notag \\
&\hspace{3.6cm}
+\norm{q}_{L^\infty([0,T],H^s(\Omega))}+\norm{\del{t}q}_{X_T^{s,2s-2}}+
\norm{\vec{k}}_{W^{1,\infty}([0,T],H^s(\Omega))}
\label{linlocthmA7}
}
and $\vec{k}=(k^i_1,k^i_2)$.
\end{thm}
\begin{proof}
Given $v_{2s}(t)\in H^1(\Omega,\Rbb^N)$, $\del{t}v_{2s}(t)\in L^2(\Omega,\Rbb^N)$, and
$w_{2s+1}(t)\in L^2(\del{}\Omega,\Rbb^N)$ with $\text{ran}(w_{2s+1}(t))\subset \text{ran}(q(t))$, it follows from Proposition \ref{lineexist} that
there exists constants $\lambda^*=\lambda^*(\sigma,\mu) \geq 1$ and
\gath{linlocthmA8}{
\delta^*= \delta^*\Bigl(\kappa_1,\norm{\del{t}b^{i0}}_{L^\infty([0,T],H^{s-\frac{1}{2}}(\Omega))},
\norm{\del{t}^2 q}_{L^\infty([0,T],H^{s-1}(\Omega))},
\norm{\del{t}P}_{L^\infty([0,T],H^{s-\frac{1}{2}}(\Omega))}\Bigr) \geq 1,
}
such that for each $(\lambda,\ep) \in [\lambda^*,\infty)\times (0,\frac{1}{\delta^*}]$, there exists
a unique solution $\vv_{2s-1}(t) \in X^{s+1,2s-1}(\Omega,\Rbb^N)$
of \eqref{linF.1}-\eqref{linF.2} for  $\ell=0,\ldots, 2s-1$ that satisfies
\lalin{linlocthmA9}{
&\norm{\vv_{2s-1}(t)}_{X^{s+1,2s-1}} \leq
C\Bigl(\norm{F(t)}_{E^{s-1}}+\norm{\del{t}^{2s-1}F(t)}_{L^2(\Omega)}+ \norm{\del{t}L^0(t)}_{E^{s-1}}
 + \norm{\del{t}^{2s}L^0(t)}_{L^2(\Omega)}\notag\\
&\hspace{1.7cm} +\norm{G(t)}_{E^{s,2s-2}}+ \norm{\del{t}^{2s-1}G(t)|_{\del{}\Omega}}_{L^2(\del{}\Omega)}
+\norm{L^i(t)}_{E^{s,2s-1}}
+ \ep \norm{(v_{2s}(t),w(t)_{2s+1}}_E\Bigr) \label{linlocthmA9.1}
}
where $C=C(\kappa_1,\mu,\sigma,\gamma,\lambda,\rho)$ with $\rho$ given by \eqref{linlocthmA7}.
Here, we are interpreting $\vv_{2s-1}(t)$ as a non-local map that depends
linearly on $(v_{2s}(t),w_{2s+1}(t))$ and satisfies the estimates \eqref{linlocthmA9.1}.

By assumption, $\del{t}^{2s}P=k^i_1\del{i}\theta_1 + g_1$ and $\del{t}^{2s}G=k^i_2\del{i}\theta_i
+g_2$, which implies that
\eqn{linlocthmA9a}{
\del{t}^{2s}G+\ep\del{t}^{2s}P v_1= \ep k^i_1\del{i}(\theta_1 v_1)
+ k^i_2\del{i}\theta_2  -k^i_1\theta_1\del{i}v_1 + \ep g_1 v_1 +g_2.
}
From this, the estimate \eqref{linlocthmA9.1} and Lemma \ref{linloclem}, we are
able to conclude from Theorem \ref{weakthm} the existence of a unique weak solution
\leqn{linlocthmA11}{
(v_{2s},w_{2s+1}) \in \bigcap_{j=0}^1 C^j\bigl([0,T],H^{1-j}(\Omega,\Rbb^N)\bigr)
\times C\bigl([0,T],L^2(\del{}\Omega,\Rbb^N)\bigr)
}
of the IBVP defined by \eqref{linF.4}-\eqref{linF.6} that satisfies the energy estimate
\lalin{linlocthmA12}{
&\norm{(v_{2s}(t),w_{2s+1}(t))}_{E}^2 \leq C\biggl(
\norm{(v_{2s}(0),w_{2s+1}(0))}_{E}^2 + \alpha_0
 \notag \\
&\hspace{2.0cm} + \int_0^t \alpha_1(\tau)\Bigl(\nnorm{(\vv_{2s-1}(\tau),v_{2s}(\tau),w_{2s+1}(\tau))}_{s+1}^2
+\beta(\tau)\Bigr) d\tau\label{linlocthmA12.1}
}
where $C=C(\kappa_0,\kappa_1,\sigma,\gamma,\mu,\lambda,\rho)$, $\alpha_1$ is defined by
\eqref{linlocthmA5.1},
\lalin{linlocthmA13}{
&\beta(t) = \norm{\del{t}^{2s}F(t)}_{L^2(\Omega)}^2 + \norm{\del{t}^{2s}L(t)}_{L^2(\Omega)}^2+ \norm{\del{t}^{2s+1}L(t)}_{L^2(\Omega)}^2 \notag \\
&\hspace{2.5cm} +\norm{g(t)}_{H^1(\Omega)}^2+ \norm{\del{t}g(t)}_{L^2(\Omega)}^2
+\norm{\theta(t)}_{H^1(\Omega)}^2+\norm{\del{t}\theta(t)}_{H^1(\Omega)}^2 \notag %\label{linlocthmA13.1}
}
and
\eqn{linlocthmA14}{
\nnorm{(\vv_{2s-1}(t),v_{2s}(t),w_{2s+1}(t))}_{s+1}^2 = \norm{\vv_{2s-1}(t)}_{X^{s+1,2s-1}}^2 + \norm{(v_{2s}(t),w_{2s+1}(t))}_E^2.
}
Since the linear map $H^1(\Omega,\Rbb^N)\times L^2(\Omega,\Rbb^N)\times L^2(\del{}\Omega,\Rbb^N) \in(v_{2s},\del{t}v_{2s},w_{2s+1})
\longmapsto \vv_{2s-2} \in X^{s+1,2s-1}(\Omega,\Rbb^N)$ is bounded, it is clear from
\eqref{linlocthmA11} that
\eqn{linlocthmA14a}{
\vv_{2s-1} \in C^0\bigl([0,T],X^{s+1,2s-1}(\Omega,\Rbb^N)\bigr).
}

\begin{comment}
Differentiating $\norm{\del{t}^\ell F(t)}_{H^{s-1,\ell}}^2$ in time, we get
\eqn{linlocthmA16}{
\del{t}\norm{\del{t}^\ell F(t)}_{H^{s-1-\frac{\ell}{2}}(\Omega)}^2 =
2\ip{\del{t}^\ell F(t)}{\del{t}^{\ell+1}F(t)}_{H^{s-1-\frac{\ell}{2}}(\Omega)}.
}
Integrating in time gives
\lalin{linlocthmA16a}{
\norm{\del{t}^\ell F(t)}_{H^{s-1-\frac{\ell}{2}}(\Omega)}^2 &\leq
\norm{\del{t}^\ell F(0)}_{H^{s-1-\frac{\ell}{2}}(\Omega)}^2
+ 2\int_0^t\Bigl|\ip{\del{t}^\ell F(\tau)}{\del{t}^{\ell+1}F(\tau)}_{H^{s-1-\frac{\ell}{2}}(\Omega)}\Bigr|\, d\tau \notag \\
&\lesssim \norm{\del{t}^\ell F(0)}_{H^{s-1-\frac{\ell}{2}}(\Omega)}^2
+\int_{0}^t \norm{\del{t}^\ell F(\tau)}_{H^{s-1-\frac{\ell}{2}}(\Omega)}\norm{\del{t}^{\ell}\del{t}F(\tau)}_{H^{s-1-\frac{\ell}{2}}(\Omega)} \, d\tau \notag \\
&\lesssim \norm{F(0)}_{E^{s-1}}^2 + \int_{0}^t \norm{\del{t}^\ell F(\tau)}_{H^{s-1-\frac{\ell}{2}}(\Omega)}^2
+ \norm{\del{t}^{\ell}\del{t}F(\tau)}_{H^{s-1-\frac{\ell}{2}}(\Omega)}^2\, d\tau, \notag
}
for $0\leq \ell \leq 2s-2$, which in turn, implies that
\leqn{linlocthmA16b}{
\norm{F(t)}_{E^{s-1}}^2  \lesssim \norm{F(0)}_{E^{s-1}}^2 + \int_{0}^t \norm{F(\tau)}_{E^{s-1}}^2
+\norm{\del{t}F(\tau)}_{E^{s-1}}^2 \, d\tau.
}
\end{comment}

%In order to establish an energy inequality, we need an integral estimate for the terms appearing on
%the right hand side of \eqref{linlocthm9.1} except for the last term. To this end, we
Next, we observe that the integral estimate
\lalin{linlocthmA16c}{
&\norm{F(t)}_{E^{s-1}}^2+\norm{\del{t}^{2s-1}F(t)}_{L^2(\Omega)}^2+\norm{\del{t}L^0(t)}_{E^{s-1}}^2+\norm{\del{t}^{2s}L^0(t)}_{L^2(\Omega)}^2
+\norm{L^i(t)}_{E^{s,2s-1}}^2
+\norm{G(t)}_{E^{s,2s-2}}^2  \notag \\
&\hspace{0.5cm}\lesssim \norm{F(0)}_{E^{s-1}}^2+\norm{\del{t}^{2s-1}F(0)}_{L^2(\Omega)}^2+\norm{L(0)}_{E^{s}}^2+\norm{G(0)}_{E^{s,2s-2}}^2
 +\int_0^t \norm{\del{t}F(\tau)}_{E^{s-1}}^2 \notag \\
&\hspace{1.5cm}+\norm{\del{t}^{2s}F(\tau)}_{L^2(\Omega)}^2
+\norm{L(\tau)}_{E^{s}}^2
+\norm{\del{t}L(\tau)}_{E^{s}}^2+
\norm{G(\tau)}_{E^{s,2s-2}}^2+\norm{\del{t}G(\tau)}_{E^{s,2s-2}}^2\, d\tau \label{linlocthmA16c.1}
}
follows directly from Proposition \ref{stpropF}. We also observe that
\lalin{linlocthmA16d}{
\del{t}\norm{\del{t}^{2s-1}G}^2_{L^2(\del{}\Omega)} &= 2\ip{\del{t}^{2s-1}G}{\del{t}^{2s}G}_{\del{}\Omega} \notag \\
&= 2\ip{\del{t}^{2s-1}G}{k^i_2\del{i}\theta_2}_{\del{}\Omega} + 2\ip{\del{t}^{2s-1}G}{g_2}_{\del{}\Omega}\notag \\
&\lesssim \norm{k_2}_{L^\infty([0,T],W^{1,\infty}(\Omega))}\norm{\del{t}^{2s-1}G}_{H^{1}(\Omega)}
\bigl(\norm{g_2}_{H^1(\Omega)}+\norm{\theta_2}_{H^1(\Omega)}\bigr)
\notag
}
where in deriving the last inequality we used \eqref{weakrem7} and Theorem \ref{trace}. Integrating this
inequality in time gives
\lalin{linlocthmA16e}{
&\norm{\del{t}^{2s-1}G(t)}^2_{L^2(\del{}\Omega)} \lesssim
\norm{\del{t}^{2s-1}G(0)}^2_{L^2(\del{}\Omega)} \notag \\
&\hspace{2.0cm}+ \norm{k_2}_{L^\infty([0,T],W^{1,\infty}(\Omega))}\int_{0}^t \norm{\del{t}^{2s-1}G(\tau)}_{H^{1}(\Omega)}^2
+ \norm{g_2(\tau)}_{H^1(\Omega)}^2+\norm{\theta_2(\tau)}_{H^1(\Omega)}^2\, d\tau. \label{linlocthmA16e.1}
}

Combining the estimates \eqref{linlocthmA9.1}, \eqref{linlocthmA16c.1}, and \eqref{linlocthmA16e.1},
we find that
\leqn{linlocthmA17}{
\norm{\vv_{2s-1}(t)}_{X^{s+1,2s-1}}^2 \leq
C\biggl(\alpha_0 + \int_0^t \alpha_2(\tau)\, d\tau + \ep^2 \norm{(v_{2s-1}(t),w_{2s}(t))}_{E}^2\biggr)
%\notag \\
%&\hspace{3.4cm}
% +\norm{G(\tau)}_{E^{s,2s-2}}^2+ \norm{\del{t}G(\tau)}_{E^{s,2s-2}}^2 + \norm{\del{t}^{2s-1}G(\tau)}_{H^1(\Omega)}^2 \notag \\
%&\hspace{3.7cm}+ \norm{g(\tau)}_{H^1(\Omega)}^2+\norm{\theta(\tau)}_{H^1(\Omega)}^2\,d\tau
%+ \ep^2 \norm{(v_{2s-1}(t),w_{2s}(t))}_{E}^2\biggr). \label{linlocthmA17.1}
}
where $C=C(\kappa_1,\mu,\sigma,\gamma,\lambda,\rho)$ and $\alpha_2$ is given
by \eqref{linlocthmA5.2}. Since the constant $C$ in \eqref{linlocthmA17} is
independent of $\ep$, we see that by choosing $\ep_0>0$ small enough that
\leqn{linlocthmA18}{
\norm{\vv_{2s-1}(t)}_{X^{s+1,2s-1}}^2-\frac{1}{2}\norm{(v_{2s-1}(t),w_{2s}(t))}_{E}^2  \leq
C\biggl(\alpha_0 + \int_0^t \alpha_2(\tau)\, d\tau\biggr).
}
for $\ep \in (0,\ep_0]$. Adding the estimates \eqref{linlocthmA12.1} and \eqref{linlocthmA18} yields
the desired energy estimate and
completes the proof.
\end{proof}

\begin{rem}
From the proof of Theorem \ref{linlocthmA}, it is not difficult to see that the energy estimate continues to hold under the following
weaker assumptions on the coefficients $L^\mu$ provided that we make the following changes:
\begin{gather*}
\norm{L(t)}_{E^s}^2 + \norm{\del{t}L(t)}^2_{E^s} \longmapsto \norm{\del{t}^2L^0(t)}_{E^{s-1}}^2+
\norm{\del{t}^{2s+1}L^0(t)}_{L^2(\Omega)}^2+ \norm{\del{t}L^i(t)}_{E^{s,2s-1}}^2
\intertext{in $\alpha_2(t)$, and}
\alpha_0 \longmapsto \alpha_2(0) + \norm{F(0)}_{E^{s-1}}^2+ \norm{\del{t}^{2s-1} F(0)}_{L^2(\Omega)}^2+
\norm{\del{t}L^0(0)}_{E^{s-1}}^2+
\norm{\del{t}^{2s}L^0(0)}_{L^2(\Omega)}^2+ \norm{L^i(0)}_{E^{s,2s-1}}^2.
\end{gather*}
Although we will not make use of this observation in this article, we note it here because it may be of use in other application of
this linear theory to non-linear problems. 
\end{rem} 

In order to go from solutions of \eqref{linF.1}-\eqref{linF.6} to solutions of \eqref{linC.1}-\eqref{linC.3}, we
need to ensure that the initial data
\leqn{compatdef1}{
(v,\del{t}v)|_{\Omega_0} = (\vt_0,\vt_1) \in H^{s+1}(\Omega,\Rbb^N)\times H^{s+\frac{1}{2}}(\Omega,\Rbb^N)
}
satisfies the \emph{compatibility conditions} given by
\lalin{compatdef2}{
\vt_\ell &:= \del{t}^\ell v|_{\Omega_0} \in H^{s+1-\frac{\ell}{2}}(\Omega,\Rbb^N), \qquad 2\leq \ell \leq 2s, \label{compatdef2.1}\\
\vt_{2s+1} &:= \del{t}^{2s+1} v|_{\Omega_0} \in L^2(\Omega),\label{compatdef2.2}
%\intertext{and}
%\wt_{2s+1}&:=\Pbb_q\del{t}^{2s+1}v|_{\Gamma_0} \in L^2(\del{}\Omega,\Rbb^N). \label{compatdef2.3}
}
where the time derivatives $\del{t}^\ell v |_{\Omega_0}$, $\ell \geq 2$, are generated
from the initial data \eqref{compatdef1} by formally differentiating \eqref{linC.1} with respect to
$t$ and evaluating at $t=0$ and are assumed to satisfy $2s-1$ formal time derivatives of
the boundary condtion \eqref{linC.2} at $t=0$.

\begin{comment}
\begin{rem}
Due to the $\ep$-dependence in the evolution equations \eqref{linC.1}-\eqref{linC.3}, the
$\vt_\ell$, $\ell \geq 2$, and $\wt_{2s}$, which are generated from the initial data,
are themselves $\ep$-dependent. Allowing for the possibility that the initial data
$(\vt_0,\vt_1)$ is $\ep$-dependent, we will say that the compatibility conditions
are satisfied \emph{uniformly for $0< \ep \leq \ep_0$} if the the compatibility
conditions are satisfied for $0<\ep\leq \ep_0$ and the norms
$\norm{\vt_\ell}_{H^{s+1-\frac{m_\ell}{2}}(\Omega)}$ and $\norm{q(0)\wt_{2s}}_{L^2(\del{}\Omega)}$
are uniformly bounded.
\end{rem}
\end{comment}

\begin{cor} \label{linlocthmB}
Suppose $s>n/2+1$, $s=k/2$ for $k\in \Zbb$, $\ep_0>0$,$\wt_{2s+1}\in L^2(\del{}\Omega,\Rbb^N)$ satisfies $\Pbb_{q(0)} \wt_{2s+1} = \wt_{2s+1}$, the coefficients $\{b^{\alpha\beta},c,q\}$ satisfy \eqref{qPdef}-\eqref{coerc} for constants $\kappa_0,\kappa_1,\sigma >0$
and $\mu \geq 0$, there exists a function $\chi=\chi(t,x)$ such that
\eqn{linlocthmBa}{
\ep P + \ep\bigl(2s-\Threehalfs\bigr)\del{t}q - \chi q \leq 0 \hspace{0.4cm} \text{in $\Gamma_T$}
}
for $0<\ep \leq \ep_0$,
\eqn{linlocthmAaa}{
\del{t}^{2s}P= k^i_1\del{i}\theta_1 + g_1 \AND \del{t}^{2s}G=k_2^i\del{i}\theta_2+g_2
}
where $\nu_i k^i_a =0$, $a=1,2$,
and the initial data $(\vt_0,\vt_1)
\in H^{s+1}(\Omega,\Rbb^N)\times H^{s+\frac{1}{2}}(\Omega,\Rbb^N)$ verify the compatibility conditions \eqref{compatdef2.1}-\eqref{compatdef2.2} for $0<\ep \leq \ep_0$.
Then there exist constants $\lambda^*=\lambda^*(\sigma,\mu) \geq 1$ and
\gath{linlocthmB1a}{
\delta^*= \delta^*\Bigl(\kappa_1,\sup_{0\leq t\leq T}\norm{\del{t}b^{i0}(t)}_{H^{s-\frac{1}{2}}(\Omega)},
\sup_{0\leq t\leq T}\norm{\del{t}^2 q(t)}_{H^{s-1}(\Omega)},
\sup_{0\leq t\leq T}\norm{\del{t}P(t)}_{H^{s-\frac{1}{2}}(\Omega)}\Bigr) \geq 1,
}
such that for each $(\lambda,\ep) \in [\lambda^*,\infty)\times \bigl(0,\min\bigl\{\frac{1}{\delta^*},\ep_0\bigr\}\bigr]$ there exists
a unique solution
$v\in C\Xc^{s+1}_T(\Omega,\Rbb^N)$ to the IBVP \eqref{linC.1}-\eqref{linC.3}. Moreover, there exists a map $w_{2s+1} \in C\bigl([0,T],L^2(\Omega,\Rbb^N)\bigr)$ such that
\begin{itemize}
\item[(i)]
$(\del{t}^{2s}v,w_{2s+1})$ is
the unique weak solution of the linear wave equations obtained by differentiating \eqref{linC.1}-\eqref{linC.2} $2s$-times with respect to $t$  satisfying $(\del{t}^{2s}v(0),\del{t}^{2s+1}v(0),w_{2s+1}(0))=(\vt_{2s},\vt_{2s+1},\wt_{2s+1})$, and
\item[(ii)] for $\ep_0>$ chosen small enough, the pair $(v,w_{2s+1})$
satisfy the energy estimate
\lalin{linlocthmB2}{
&\nnorm{(v(t),w_{2s+1}(t))}_{s+1}^2 \leq C\biggl( \nnorm{(v(0),w_{2s+1}(0))}_{s+1}^2
+ \alpha_0  \notag \\
&\hspace{4.5cm} + \int_0^t \alpha_1(\tau)\Bigl( \nnorm{(v(\tau),w_{2s+1}(\tau))}_{s+1}^2
+ \alpha_2(\tau)\Bigr) \, d\tau \biggr) \label{linlocthmB2.1}
}
where $C = C\bigl(\kappa_0,\kappa_1,\sigma,\gamma,\mu,\lambda,\rho\bigr)$, $\rho$ and the $\alpha_i$ are as defined in Theorem \ref{linlocthmA}, and
\eqn{linlocthmB3}{
\nnorm{(v(t),w_{2s+1}(t))}_{s+1}^2 = \norm{v(t)}_{\Ec^{s+1}}^2 + \ip{w_{2s+1}(t)}{(-q(t))w_{2s+1}(t)}_{\del{}\Omega}.
}
\end{itemize}
\end{cor}
\begin{proof}
Given initial data
$(\vt_0,\vt_1)
\in H^{s+1}(\Omega,\Rbb^N)\times H^{s+\frac{1}{2}}(\Omega,\Rbb^N)$ satisfying the compatibility conditions \eqref{compatdef2.1}-\eqref{compatdef2.2}, we let
\eqn{linlocthmB4}{
(\vv_{2s-1},v_{2s},w_{2s+1}) \in C^0\bigl([0,T],X^{s+1,2s-1}(\Omega,\Rbb^N)\bigr) \times
\bigcap_{j=0}^1 C^j\bigl([0,T],H^{1-j}(\Omega,\Rbb^N)\bigr) \times C\bigl([0,T],L^2(\del{}\Omega,\Rbb^N)\bigr),
}
denote the unique solution to \eqref{linF.1}-\eqref{linF.6}, which we know exists
for $(\lambda,\ep)\in [\lambda^*,\infty)\times \bigl(0,\min\bigl\{\frac{1}{\delta^*},\ep_0\bigr\}\bigr]$
by Theorem \ref{linlocthmA}.
\begin{comment}
Following the same strategy as we used in the analysis of the simplified
system we considered in Section \ref{simple}, we apply the mollifier $J_\omega$, defined
previously by \eqref{sssB}, to the system \eqref{linF.1}-\eqref{linF.4} to get
\lalin{linlocthmB5}{
&\del{i}\bigl(b^{i j}\del{j}J_\omega v_\ell +[J_\omega,b^{ij}]\del{j}v_\ell +
+\ep(\dc^i_\ell J_\omega v_\ell+[J_\omega,\dc^i]v_\ell) + J_\omega\Lc^i_\ell\bigr)  \notag\\
&+
\ep \ac^i_\ell\del{i}v_\ell + \lambda \cc_\ell v_\ell = \Fc_\ell \hspace{1.655cm} \text{in $\Omega$,}
\label{linlocthmB5.1}\\
&\nu_i\bigl(b^{i j}\del{j}v_\ell + \ep\dc^i_\ell v_\ell + \Lc^i_\ell\bigr)  = \ep\hc_\ell v
+\Gc_\ell
 \hspace{0.675cm} \text{in $\del{}\Omega$, } \label{linlocthmB5.2}\\
&\del{\alpha}\bigl(B^{\alpha\beta}\del{\beta}v_{2s}+\Mcal^\alpha_{2s}\bigr) +\lambda c v_{2s} = \Fc_{2s} \hspace{1.525cm} \text{in $\Omega_T$,}
\label{linlocthmB5.3}\\
&\nu_\alpha\bigl(B^{\alpha\beta}\del{\beta} v_{2s} + \Mcal^\alpha_{2s})- \ep^2q\del{t}^2 v_{2s} -\ep\Pc \del{t} v_{2s}  =  \Gc_{2s}
\hspace{1.6cm} \text{in $\Gamma_T$, } \label{linlocthmB5.4}
}
\end{comment}
To proceed, we assume that the $v_{\ell}(t)$, $\ell=0,\ldots,2s-1$, are differentiable
in time and satisfy $\del{t}v_{\ell}(t)\in H^{s+1-\frac{\ell+1}{2}}(\Omega,\Rbb^N)$, $\ell=0,\ldots,2s-1$.
This assumption can be made rigorous by applying the mollifier $J_\omega$, see \eqref{sssB},
to the equations \eqref{linF.1}-\eqref{linF.4}, in analogy to what was done in the analysis of the simplified problem
in Section \ref{simple}. Commuting the mollifiers with the coefficients, which are now variable, it follows
that the smoothed variables $J_\omega v_\ell$ satisfy the same type of equations up to remainder terms
involving commutators of the type $[J_\omega,(\cdot)]$, where $(\cdot)$ represents the system coefficients.
Using standard properties of mollifiers, it is not difficult to see that these remainder terms are harmless
as far as this proof is concerned.

\begin{comment}
This
assumption can be justified by
replacing the derivative $\del{t}v_{\ell}$ by the difference quotient
$\Delta_h v_{\ell}(t,x) = h^{-1}\bigl(\Delta_h v_{\ell}(t+h,x)-\Delta_h v_\ell(t,x)\bigr)$
and sending $h\searrow 0$ at the end of the computation.
\end{comment}

Under the differentiability assumption, a straightforward calculation shows that the differences
$\vf_{\ell} := \del{t}v_{\ell-1}-v_{\ell}$, $\ell=1,\ldots,2s-1$, define a
solution to a collection of elliptic equations of the form
\lalin{linlocthmB5}{
\del{i}\bigl(b^{i j}\del{\beta}\vf_\ell + \ep\dc^i_\ell \vf_\ell + \Lf^i_\ell\bigr) +
\ep \ac^i_\ell\del{i}v_\ell + \lambda \cc_\ell v_\ell &= \Ff_\ell \hspace{1.75cm} \text{in $\Omega$,}
\label{linlocthmB5.1}\\
\nu_i\bigl(b^{i j}\del{\beta}v_\ell + \ep \dc^i_\ell \vf_\ell + \Lf^i_\ell\bigr) & = \ep\hc_\ell \vf_\ell
+\Gf_\ell
 \hspace{0.5cm} \text{in $\del{}\Omega$ \label{linlocthmB5.2} }
}
where the source terms $\Lf^i_\ell$, $\Ff_\ell$, and $\Gf_\ell$ are homogenous in the variables $\vf_{\ell}$
and satisfy homogeneous versions of the estimates from Lemma \ref{elem}; that is, estimates
that arise from making the replacements: $(\Gc_\ell,\Lc^i_\ell,\Fc_\ell) \mapsto (\Gf_\ell,\Lf_\ell^i,\Ff_\ell)$
and $(G_\ell,L^i_\ell,F_\ell)\mapsto (0,0,0)$ for $1\leq \ell \leq 2s-1$. The only point that one has to
be somewhat careful about is to remember that on the boundary the equality $q\del{t}v_{2s}-q w_{2s+1}=0$ holds since
it is true in a weak sense, see Remark \ref{weakrem}.(i).
Furthermore, it is not difficult to verify that the last term $\vf_{2s}$ defines a weak solution of an equation of the form \eqref{linlocthmB5.1}-\eqref{linlocthmB5.2} with $\ell=2s$, where again the source terms $\Lf^i_{2s}$, $\Ff_{2s}$, and $\Gf_{2s}$ are homogenous in the variables $\vf_{\ell}$ and satisfy the estimates
\eqn{inlocthmB6}{
\norm{\Gf_{2s}|_{\del{}\Omega}}_{L^2(\del{}\Omega)}+\norm{\Lf^{i}_{2s}}_{L^2(\Omega)}
+\norm{\Ff_\ell}_{L^2(\Omega)} \lesssim \sum_{\ell=1}^{2s}\norm{\vf_\ell}_{H^{s+1-\frac{\ell}{2}}(\Omega)}.
}

By Proposition \ref{lineexist} and Theorem \ref{ellipthmA}, we know that
solutions to \eqref{linlocthmB5.1}-\eqref{linlocthmB5.2} are unique for $(\lambda,\ep)\in [\lambda*,\infty)\times \bigl(0,\frac{1}{\delta^*}\bigr]$, and so, we conclude that the trivial solution, given by $\vf_{\ell}=0$ for $1\leq \ell \leq 2s$,
is the unique solution. From this, we see that
$\del{t}v_{\ell} = v_{\ell+1}$, $0\leq \ell \leq 2s-1$,
and hence, that
\leqn{linlocthmB8}{
v_\ell = \del{t}^\ell v, \qquad 0\leq \ell \leq 2s.
}
The proof now follows since it is clear from the properties of the
solution $(\vv_{2s-1},v_{2s},w_{2s+1})$, see Theorem \ref{linlocthmA}, and
\eqref{linlocthmB8} that
$v\in C\Xc^{s+1}(\Omega,\Rbb^N)$, $v$ solves the rescaled model problem \eqref{linC.1}-\eqref{linC.3},
$(\del{t}^{2s}v,w_{2s+1})$ is a weak solution of the linear wave equation
obtained by differentiating \eqref{linC.1}-\eqref{linC.2} $2s$-times with respect to $t$,
and $(v,w_{2s+1})$ satisfies the desired energy estimate.
\end{proof}

\begin{rem} \label{correm}
$\;$

\begin{itemize}
\item[(i)] It is clear from the proofs of Theorem \ref{linlocthmA} and Corollary \ref{linlocthmB} that
 they continue to hold for $P$ and $G$ satisfying
\eqn{weakrem1a}{
\del{t}^{2s}P = k^i_{1,1}\del{i}\theta_{1,1}+k^i_{1,2}\del{i}\theta_{1,2}
\AND \del{t}^{2s}G = k^i_{2,1}\del{i}\theta_{2,1}+k^i_{2,2}\del{i}\theta_{2,2},
}
where $\nu_ik^i_{a,b}=0$, $a,b=1,2$,
provided that we make the replacements
\eqn{weakrem1b}{
\theta_{a} \longmapsto (\theta_{a,1},\theta_{a,2}), \quad a=1,2, \AND
\vec{k} \longmapsto (k_{1,1}, k_{1,2},k_{2,1},k_{2,2}).
}
\item[(ii)]
There is a generalization of Corollary \ref{linlocthmB}, which
can be established using an iteration method that
allows for coefficients the $L^\alpha$, $F$,
$G$ and $g$ to depend linearly on $v$ provided that they satisfy estimates
that preserve the form of the energy estimate. For example, we could have
\eqn{weakrem1}{
F=\bar{F} + \hat{F}(v)
}
where $\bar{F}\in X_T^{s}(\Omega)$ and
\eqn{weakrem2}{
\norm{\hat{F}(v(t))}_{E^{s}}^2 \leq  \hat{f}(t)\norm{(v(t),w_{2s+1}(t))}_{s+1}^2, \quad 0\leq t\leq T.
}
In this case, we would just replace the term $\norm{F(t)}^2_{E^{s}}$ that appears in the energy estimate (i.e. in $\alpha_2(t)$) from Corollary \ref{linlocthmB}
with $\norm{\bar{F}(t)}_{E^{s}}^2 + \hat{f}(t)\norm{(v(t),w_{2s+1}(t))}_{s+1}^2$.
\end{itemize}
\end{rem}

\subsect{linloc}{Local existence and uniqueness}
We
are now prepared to use the existence and uniqueness results for
the model problem to establish an existence and uniqueness result
for linear wave equations that
include equations of the form \eqref{cbvpA.3}-\eqref{cbvpA.4}. The precise class of linear
wave equations that we consider are:
\lalin{linG}{
\del{\alpha}\bigl(b^{\alpha\beta}\del{\beta}u+\ell^\alpha\bigr) &= f \hspace{2.65cm} \text{in $\Omega_T$,}
\label{linG.1}\\
\nu_\alpha\bigl(b^{\alpha\beta}\del{\beta} u + \ell^\alpha) & = q\del{t}^2 u + p \del{t} u + g
 \hspace{0.4cm} \text{in $\Gamma_T$,} \label{linG.2}\\
(u,\del{t}u) &= (\ut_0,\ut_1) \hspace{1.75cm} \text{in $\Omega_0$}, \label{linG.3}
%(\del{t}^{2s}u|_{\text{ran}(q)}) &= \wt_{2s} \hspace{2.45cm} \text{in $\Gamma_0$.} \label{linG.4}
}
where the initial data
\leqn{ucompatdef1}{
(u,\del{t}u)|_{\Omega_{0}} = (\ut_0,\ut_1) \in H^{s+1}(\Omega,\Rbb^N)\times H^{s+\frac{1}{2}}(\Omega,\Rbb^N)
}
satisfies \emph{compatibility conditions} given by
\lalin{ucompatdef2}{
\ut_\ell &:= \del{t}^\ell u|_{\Omega_0} \in H^{s+1-\frac{\ell}{2}}(\Omega,\Rbb^N), \qquad 2\leq \ell \leq 2s, \label{ucompatdef2.1}\\
\ut_{2s+1} &:= \del{t}^{2s+1} u|_{\Omega_0}  \in L^2(\Omega,\Rbb^N), \label{ucompatdef2.2}
%\intertext{and}
%\wt_{2s} &:=\Pbb_q\del{t}^{2s}u|_{\Gamma_0} \in L^2(\del{}\Omega,\Rbb^N). \label{ucompatdef2.3}
}
where, as before, the higher time derivatives $\del{t}^\ell u |_{\Omega_0}$, $\ell \geq 2$, are generated
from the initial data \eqref{ucompatdef1} by formally differentiating \eqref{linG.1} with respect to
$t$ at $t=0$ and are assumed to satisfy $2s-1$ formal time derivatives of the boundary conditions
\eqref{linG.2} at $t=0$.

\begin{thm} \label{linlocthmC}
Suppose $s>n/2+1$, $s=k/2$ for $k\in \Zbb$, $\wt_{2s+1}\in L^2(\del{}\Omega)$ satisfies $\Pbb_{q(0)}\wt_{2s+1}=\wt_{2s+1}$,
the coefficients $\{b^{\alpha\beta},q\}$ satisfy \eqref{qPdef} and \eqref{bsym.1}-\eqref{coerc}
for constants $\kappa_0,\kappa_1,\gamma >0$
and $\mu \geq 0$,
\gath{linlocthmC0a}{
p + \biggl(2s-\frac{3}{2}\biggr)\del{t}q - \chi q \leq 0, \quad
\del{t}^{2s}p=k^i_1\del{i}\theta_1+h_1  \AND
\del{t}^{2s}g=k^i_2\del{i}\theta_2 + h_2
}
where $\nu_i k^i_a=0$, $a=1,2$,
and the initial data $(\vt_0,\vt_1)$ satisfies the compatibility conditions \eqref{ucompatdef2.1}-\eqref{ucompatdef2.2}.
Then there exists a unique solution
$u\in C\Xc^{s+1}_T(\Omega,\Rbb^N)$ to the IBVP \eqref{linG.1}-\eqref{linG.3}. Moreover, there exists a map $w_{2s+1} \in C([0,T],L^2(\Omega,\Rbb^N))$ such that
\begin{itemize}
\item[(i)]
$(\del{t}^{2s}u,w_{2s+1})$ is
 the unique weak solution of the linear wave equations obtained by differentiating \eqref{linG.1}-\eqref{linG.2} $2s$-times with respect to $t$ satisfying
 $(\del{t}^{2s}u(0),\del{t}^{2s+1}u(0),w_{2s+1}(0))=(\ut_{2s},\ut_{2s+1},\wt_{2s+1})$, and
\item[(ii)] the pair $(u,w_{2s+1})$
satisfy the energy estimate
\alin{linlocthmC2}{
\nnorm{(u(t),w_{2s+1}(t))}_{s+1}^2 \leq
C &\biggl(\nnorm{(u(0),w_{2s+1}(0))}_{s+1}^2 \\
 &+\alpha_0+
 \int_0^t  \alpha_1(\tau)\nnorm{(u(\tau),w_{2s+1}(\tau))}_{s+1}^2+\alpha_2(\tau)\,d\tau \biggr)
}
where $C = C\bigl(\kappa_0,\kappa_1,\mu,\gamma,\rho\bigr)$, $\alpha_0= \alpha_2(0)+ \norm{f(0)}^2_{E^{s-1}}+
\norm{\del{t}^{2s-1}f(0)}^2_{L^2(\Omega)}$,
\eqn{linlocthmC3}{
\nnorm{(u(t),w_{2s+1}(t))}_{s+1}^2 = \norm{u(t)}_{\Ec^{s+1}}^2 + \ip{w_{2s+1}(t)}{(-q)w_{2s+1}(t)}_{\del{}\Omega},
}
\alin{linlocthmC3a}{
\alpha_1(t) & = 1+\norm{\chi(t)}_{H^{s}(\Omega)}^2+\norm{\del{t}b(t)}_{E^{s}}^2
 + \norm{\del{t}p(t)}_{E^{s,2s-2}}^2 \\
 &\hspace{2.5cm}+\norm{\del{t}q(t)}_{\Ec^{s+\frac{1}{2}}}^2
+\norm{\vec{h}(t)}^2_{\Ec^1}+\norm{\vec{\theta}(t)}^2_{H^1(\Omega)}+\norm{\del{t}\vec{\theta}(t)}^2_{H^1(\Omega)},\\
\alpha_2(t) &= \norm{f(t)}_{E^{s-1}}^2+ \norm{\del{t}f(t)}_{E^{s-1}}^2+ \norm{\del{t}^{2s}f(t)}_{L^2(\Omega)}^2  + \norm{\ell(t)}_{E^{s}}^2 +\norm{\del{t}\ell(t)}_{E^{s}}^2 \\
&\qquad +\norm{g(t)}_{E^{s,2s-2}}^2
+\norm{\del{t}g(t)}_{E^{s,2s-2}}^2+ \norm{\vec{h}(t)}_{\Ec^1}^2
 +\norm{\vec{\theta}(t)}^2_{H^1(\Omega)}+\norm{\del{t}\vec{\theta}(t)}^2_{H^1(\Omega)}
\intertext{and}
\rho &= \norm{b}_{X_T^{s}}+ \norm{p}_{X_T^{s,2s-2}}+
\norm{\del{t}^{2s-1}p}_{L^\infty([0,T],L^2(\del{}\Omega))}
 \notag\\
&\hspace{2.5cm}+\norm{q}_{L^\infty([0,T],H^s(\Omega))}+\norm{\del{t}q}_{X_T^{s,2s-2}}+ \norm{\vec{k}}_{W^{1,\infty}([0,T],H^{s}(\Omega))}.
}
\end{itemize}
\end{thm}
\begin{proof}
$\;$

\smallskip

\noindent \underline{Reduction to the model problem}: First, a short computation shows that the IBVP \eqref{linG.1}-\eqref{linG.3} transform into
\lalin{linlocthmC4}{
\del{\alpha}\bigl(b^{\alpha\beta}\del{\beta}v+L^\alpha) + \lambda cv &= F
\hspace{4.3cm} \text{in $\Omega_{T}$,}\label{linlocthmC4.1}\\
\nml_\alpha \bigl(b^{\alpha\beta}\del{\beta}v+L^\alpha)    & = q\del{t}^2 v + P \del{t}v + G
\hspace{1.95cm} \text{in $\Gamma_{T}$,} \label{linlocthmC4.2}\\
(v,\del{t}v) = (\vt_0,\vt_1) &:= \bigl(e^{-\omega} \ut_0, e^{-\omega}(\ut_1 - \del{t}\omega \ut_0)\bigr)
\hspace{0.5cm} \text{in $\Omega_{0}$,} \label{linlocthmC4.3}
}
the under change of variables
%\eqn{linlocthmC6}{
$u=e^{\omega}v$,
%}
where
\gath{linlocthmC7}{
c =\frac{1}{\lambda}b^{\alpha\beta}\del{\alpha}\omega\del{\beta} \omega, \quad F = e^{-\omega}f-e^{-\omega}\ell^\alpha\del{\alpha}\omega v-\del{\alpha}\omega b^{\alpha\beta}\del{\beta}v,\\
L^\alpha = e^{-\omega}\ell^\alpha + b^{\alpha\beta}\del{\beta}\omega v, \quad P = p+ 2\del{t}\omega q ,
\intertext{and}
G = e^{-\omega}g + \bigl[(\del{t}\omega)^2+\del{t}^2\omega\bigr]q v +\del{t}\omega p v.
}
Setting
\eqn{linlocthmC8}{
\omega = \sqrt{\lambda}t, \qquad \lambda > 0,
}
gives
\eqn{linlocthmC9}{
c = b_{00} \leq -\kappa_0,
}
while
\eqn{linlocthmC10}{
P-\bigl(2s-\Threehalfs\bigr)\del{t}q - (2\sqrt{\lambda}+\chi)q \leq 0
}
follows from the assumption $p-\bigl(2s-\Threehalfs\bigr)\del{t}q - \chi q \leq 0$. Differentiating
$G$, we obtain
\lalin{linlocthmC10aa}{
\del{t}^{2s}G &= e^{-\omega}k^i_2\del{i}\theta_2   + k^i_1\del{i}(\theta_1\del{t}\omega v)
 +\Bigl\{ -k^i_1\theta_1\del{i}(\del{t}\omega v)
 + e^{-\omega}h_2+ h_1\del{t}\omega v  \notag \\
&\hspace{1.8cm} +[\del{t}^{2s},e^{-\omega}]g
+\del{t}^{2s}\bigl(\bigl[(\del{t}\omega)^2+\del{t}^2\omega\bigr]q v\bigr) +
\del{t}^{2s}\bigl(p\del{t}\omega  v\bigr)-\del{t}^{2s}p\del{t}\omega v \Bigr\} \label{linlocthmC10aa.1}
}
since $\del{t}^{2s}p=k^i_1\del{i}\theta_1+h_1$ and $\del{t}^{2s}g=k^i_2\del{i}\theta_2 + h_2$
by assumption. We note that the term $\del{t}^{2s}p$ does not appear in the term
$\del{t}^{2s}\bigl(p\del{t}\omega  v\bigr)-\del{t}^{2s}p\del{t}\omega v$
from \eqref{linlocthmC10aa.1}.
Similarly, differentiating $P$, we see that
\eqn{linlocthmC10ad}{
\del{t}^{2s}P = k^i\del{i}\theta_1 + h_1 + \del{t}^{2s}(2\del{t}\omega q).
}

Using the definitions \eqref{linlocthmA5.1}-\eqref{linlocthmA7} for $\alpha_1(t)$, $\alpha_2(t)$,
and $\rho$, where in those formulas $b^{\alpha\beta}$, $c$, $F$, $L^\alpha$, $P$, $G$, $\theta$ are as defined above, and
the following replacements are made
\alin{linlocthmC10aaa}{
g_1 & \longmapsto h_1 + \del{t}^{2s}(2\del{t}\omega q), \\
g_2 &\longmapsto \text{$\bigl\{\cdots\bigl\}$ on the r.h.s of \eqref{linlocthmC10aa.1}}, \\
k^i_2 &\longmapsto (e^{-\omega}k^i_2,k^i_1), \\
\theta_2 &\longmapsto (\theta_2,\theta_1\del{t}\omega v),
}
and setting
\alin{linlocthmC11}{
\alphat_1(t) & = 1+\norm{\chi(t)}_{H^{s}(\Omega)}^2+\norm{\del{t}b(t)}_{E^{s}}^2
 + \norm{\del{t}p(t)}_{E^{s,2s-2}}^2 \\
 &\hspace{2.5cm}+\norm{\del{t}q(t)}_{\Ec^{s+\frac{1}{2}}}^2
+\norm{\vec{h}(t)}^2_{\Ec^1}+\norm{\vec{\theta}(t)}^2_{H^1(\Omega)}+\norm{\del{t}\vec{\theta}(t)}^2_{H^1(\Omega)},\\
\alphat_2(t) &=  \norm{f(t)}_{E^{s-1}}^2+ \norm{\del{t}f(t)}_{E^{s-1}}^2+ \norm{\del{t}^{2s}f(t)}_{L^2(\Omega)}^2 +\norm{\ell(t)}_{E^{s}}^2+\norm{\del{t}\ell(t)}_{E^{s}}^2
\notag \\
&\hspace{1.5cm} +\norm{g(t)}_{E^{s,2s-2}}^2
+\norm{\del{t}g(t)}_{E^{s,2s-2}}^2 + \norm{\vec{h}(t)}_{\Ec^1}^2
 +\norm{\vec{\theta}(t)}^2_{H^1(\Omega)}+\norm{\del{t}\vec{\theta}(t)}^2_{H^1(\Omega)}
\intertext{and}
\rhot &= \norm{b}_{X_T^{s}}+ \norm{p}_{X_T^{s,2s-2}}+
\norm{\del{t}^{2s-1}p}_{L^\infty([0,T],L^2(\del{}\Omega))}
 \notag\\
&\hspace{2.5cm}+\norm{q}_{L^\infty([0,T],H^s(\Omega))}+\norm{\del{t}q}_{X_T^{s,2s-2}}+ \norm{\vec{k}}_{W^{1,\infty}([0,1],H^{s}(\Omega))},
}
it is not difficult to verify the inequalities
\eqn{linlocthmC12}{
\rho \lesssim \rhot, \quad
\alpha_1(t) \leq C(\rhot) \alphat_1(t), \AND
\alpha_2(t) \leq  C(\rhot)\bigl( \alphat_2(t) + \alphat_1(t)\norm{v(t)}_{\Ec^{s+1}}\bigr)
}
follows from the spacetime calculus inequalities from Appendix \ref{STineq}.

From the above calculations, it follows that it is enough to consider the IBVP \eqref{linlocthmC4.1}-\eqref{linlocthmC4.3}
where the coefficients $b^{ij}$, $P$, $q$, $c$ and $G$ satisfy the conditions: \eqref{qPdef}-\eqref{coerc},
\leqn{linlocthmC10a}{
P-\bigl(2s-\Threehalfs\bigr)\del{t}q - \chi q \leq 0,
}
and
\leqn{linlocthmC10b}{
\del{t}^{2s}G = k^i\del{i}\theta + g
}
where $\nu_i k^i=0$.

\bigskip

\noindent \underline{Reduction to the rescaled model problem}:
Letting
\eqn{Hbbdef}{
\Hbb^n = \{\, (x^1,\ldots,x^{n-1},x^n) \, | (x^1,\ldots, x^{n-1})\in \Rbb^{n-1}, \; x^n > 0 \, \}
}
denote the half space, we can always use a smooth coordinate transformation to map a portion
of the smooth boundary of $\Omega$ to the region $Q_{1}\cap \del{}\Hbb^n$, where
here, we are using $Q_{\delta}$ to denote an open $n$ cube of width $\delta>0$ centered at $x=0$. Consequently, nothing is
lost by just assuming that a portion of the boundary of $\Omega$ is given by $Q_1\cap \del{}\Hbb^n$
and that $Q_{1}^+\subset \Omega$, where
\eqn{Qplus}{
Q_{\delta}^+ := Q_{\delta}\cap \Hbb^n.
}
Next, we rescale
the spatial coordinates $x=(x^i)$ by $x\mapsto \ep x$, $0<\ep \leq 1$, which maps the boundary portion $Q_1^+$ to
$Q_{1/\ep}^+$ and transforms \eqref{linlocthmC4.1}-\eqref{linlocthmC4.3} to
\lalin{linlocthmC13}{
\del{\alpha}\bigl(B^{\alpha\beta}_\ep\del{\beta}v_\ep+ \ep M^\alpha_\ep) + \ep^2 \lambda c_\ep v_\ep &= \ep^2 F_\ep
\hspace{3.4cm} \text{in $[0,T]\times \Omega^\ep$,} \notag %\label{linlocthmC13.1}
\\
\nml_\alpha \bigl(B^{\alpha\beta}_\ep\del{\beta}v_\ep+ \ep M^\alpha_\ep)    & = \ep q_\ep\del{t}^2 v_\ep +
\ep P_\ep \del{t}v_\ep + \ep G_\ep
\hspace{0.4cm} \text{in $[0,T]\times \del{}\Omega^\ep$,} \notag %\label{linlocthmC13.2}
\\
(v_\ep,\del{t}v_\ep) &= (\vt_0^\ep,\vt_1^\ep)
\hspace{2.95cm} \text{in $\Omega^\ep$,}\notag % \label{linlocthmC13.3}
}
where $\Omega^\ep=\{\, x \in \Rbb^n \, | \, \ep x \in \Omega \,\}$,
\alin{linlocthmC14}{
B^{\alpha\beta}_\ep &= \ep^2 \delta^\alpha_0\delta^\alpha_0 b_{\ep}^{00} +
\ep \delta^\alpha_0\delta^\beta_j b^{0j}_\ep + \ep \delta^\alpha_i\delta^\beta_0 b^{0j}_\ep
+ b^{ij}_\ep, \\
M^{\alpha}_\ep &= \ep \delta^\alpha_0 L^0_\ep + \delta^\alpha_i L^i_\ep,
}
and all other quantities appearing with a subscript or superscript $\ep$ are obtained by scaling (i.e.
$v_\ep(t,x)=v(t,\ep x)$, $F_\ep(t,x)=F(t,\ep x)$, etc.).

\begin{lem} \label{scalemA}
Suppose $k>n/2$ and $f_\ep(x)=f(\ep x)$. Then
\eqn{scalem1}{
\norm{f_\ep}_{H^k(Q_1^+)} \lesssim \norm{f}_{H^k(\Omega)}, \quad 0<\ep \leq 1,
}
for all $f\in H^k(\Omega)$.
\end{lem}
\begin{proof}
First, it is straightforward to check that
\eqn{scalem1}{
\norm{D^k f_\ep}^2_{L^2(\Omega^\ep)} = \ep^{2k-n} \norm{f}^2_{H^k(\Omega)}, \quad 0<\ep \leq 1.
}
Since $Q^+_1\subset Q^+_{1/\ep} \subset \Omega^\ep$, and $k>n/2$ by assumption, we see
that
\leqn{scalem2}{
\norm{D^k f_\ep}_{L^2(Q^+_1)} \leq \norm{f}^2_{H^k(\Omega)}, \quad 0<\ep \leq 1.
}
Next
\leqn{scalem3}{
\norm{f_\ep}_{L^2(Q^+_1)} \leq \norm{f_\ep}_{L^\infty(Q^+_1)}\norm{1}_{L^2(Q^+_1)}
\lesssim \norm{f}_{L^\infty(Q^+_{\ep})} \leq \norm{f}_{L^\infty(\Omega)} \lesssim \norm{f}_{H^k(\Omega)},
}
where in deriving the last inequality we used Sobolev's inequality, see Theorem \ref{FSobolev}. The
proof now follows from the two inequalities \eqref{scalem2} and \eqref{scalem3} together with an
application of Ehrling's lemma, see Lemma \ref{Ehrling}.
\end{proof}
The above lemma shows that we have a bound on the rescaled coefficients given by
\lalin{linlocthmC15}{
\norm{b_\ep^{\alpha\beta}}_{H(Q_1^+)}&
+ \norm{\del{t}b_\ep^{0i}}_{H^{s-\frac{1}{2}}(Q_1^+)}
+ \norm{\del{t}P_\ep}_{H^{s-\frac{1}{2}}(Q_1^+)} + \norm{\del{t}^2 q_\ep}_{H^{s-1}(Q_1^+)} \notag \\
&\lesssim \norm{b^{\alpha\beta}}_{H^s(\Omega)}
+ \norm{\del{t}b^{0i}}_{H^{s-\frac{1}{2}}(\Omega)}
+ \norm{\del{t}P}_{H^{s-\frac{1}{2}}(\Omega)} + \norm{\del{t}^2 q}_{H^{s-1}(\Omega)},
\label{linlocthmC13.1}
}
which holds for all $(t,\ep)\in [0,T]\times (0,1]$. Moreover, it is clear
from \eqref{linlocthmC10a} and \eqref{linlocthmC10b} that
\eqn{linlocthmC16}{
P_\ep- \ep\bigl(2s-\Threehalfs\bigr)\del{t}q_\ep - \chi_\ep q_\ep \leq 0 \AND
\del{t}^{2s} \ep G_\ep = k^i_\ep\del{i}\theta_\ep + \ep g_\ep,
}
where $\nu_i k^i_\ep = 0$.

By Morrey's inequality and \eqref{linlocthmC13.1}, we also know that
\eqn{linlocthmC17}{
\norm{b_\ep^{\alpha\beta}}_{C^{0,\alpha}(\overline{Q_1^+})} \lesssim \norm{b_\ep^{\alpha\beta}}_{H(Q_1^+)}
\lesssim 1
}
for $\alpha \in (0,\min\{1,s-n/2\})$.
This estimate together with $b_\ep^{ij}(t,0)=b^{ij}(t,0)$ implies that
for any choice of constant $\beta>0$, there exists
a $\delta \in (0,1]$ such that
\eqn{linlocthmC18}{
\norm{b^{ij}_\ep(t,\cdot)-b^{ij}(t,0)}_{L^\infty(Q_\delta^+)} \leq \beta
}
for all $(t,\ep)\in [0,T]\times (0,1]$. Choosing $\beta>0$ small enough, we can, see Remark \ref{coercrem}.(i), guarantee that
$b^{ij}_\ep$ is strongly elliptic at each point $x\in Q^+_\delta$ and satisfies
the strong complementing condition on $\overline{Q^+_\delta}\cap \del{}\Hbb^n \subset \del{}\Omega^\ep$.

From the above considerations, it is not difficult to verify that, after suitably modifying the coefficients
outside of a region of the form $B_{\delta'}(0)\cap \overline{Q_\delta^+}$ for $\delta'$ small enough,  there
exists a bounded domain $\Omegah \subset Q_\delta^+ \subset \Omega^\ep$ with $C^\infty$ boundary satisfying
$\del{}\Omegah\cap \del{}\Hbb^n \subset \del{}\Omega^\ep$, and an IBVP
\lalin{linlocthmC19}{
\del{\alpha}\bigl(\Bh^{\alpha\beta}_\ep \del{\beta}\vh_\ep+ \ep\Mh_\ep^\alpha) + \lambda \ch_\ep v_\ep &= \ep^2\Fh_\ep
\hspace{3.4cm} \text{in $\Omegah_{T}$,}\label{linlocthmC19.1}\\
\nml_\alpha \bigl(\Bh^{\alpha\beta}_\ep \del{\beta}\vh_\ep+\ep\Mh^\alpha_\ep)    & = \ep \qh_\ep\del{t}^2 \vh_\ep
+ \ep\Ph_\ep \del{t}\vh_\ep + \ep\Gh_\ep
\hspace{0.4cm} \text{in $\Gammah_{T}$,} \label{linlocthmC19.2}\\
(\vh_\ep,\del{t}\vh_\ep) &= (\vch_0^\ep,\vch_1^\ep)
\hspace{3.0cm} \text{in $\Omegah_{0}$,} \label{linlocthmC19.3}
}
that satisfies the following:
\begin{itemize}
\item[(i)] the conditions \eqref{qPdef}-\eqref{coerc} are satisfied for constants $\gammah$, $\sigmah$, $\kappah_0$, $\kappa_1$,
and $\muh$ which are all $\ep$-independent with the exception of $\kappah_0$,
\item[(ii)] all of the systems coefficients (i.e $\Bh^{\alpha\beta}_\ep$, $\Fh_\ep$, etc.) agree with
the ``unhatted'' coefficients  (i.e. $B^{\alpha\beta}_\ep$, $F^\ep$, etc.) on $B_{\delta'}(0)\cap \overline{Q_\delta^+}\subset \overline{\Omegah}$,
\item[(iii)] the conditions
\eqn{linlocthmC19a}{
\Ph_\ep- \ep\bigl(2s-\Threehalfs\bigr)\del{t}\qh_\ep - \chih_\ep \qh_\ep \leq 0,
\AND
\del{t}^{2s} \ep \Gh_\ep = \kh^i_\ep\del{i}\thetah_\ep + \ep \gh_\ep
}
are satisfied,
\item[(iv)] the initial data $(\vch_0^\ep,\vch_1^\ep)$ satisfies the compatibility conditions and
agrees with the initial data $(\vt_0^\ep,\vt_1^\ep)$ on  $B_{\delta'}(0)\cap Q_\delta^+$,
\item[(v)] and the uniform bound
\alin{linlocthmC20}{
\frac{1}{\ep}\norm{\del{t}\Bh_\ep^{0j}}_{H^{s-\frac{1}{2}}(\Omegah)}&
+ \norm{\del{t}\Ph_\ep}_{H^{s-\frac{1}{2}}(\Omegah)} + \norm{\del{t}^2 \qh_\ep}_{H^{s-1}(\Omegah)} \notag \\
&\lesssim
\norm{\del{t}b^{0i}}_{H^{s-\frac{1}{2}}(\Omega)}
+ \norm{\del{t}P}_{H^{s-\frac{1}{2}}(\Omega)} + \norm{\del{t}^2 q}_{H^{s-1}(\Omega)}
}
holds for $\ep \in (0,1]$.
\end{itemize}

Choosing $\ep$ small enough and $\lambda$ big enough, we are then guaranteed by Corollary \ref{linlocthmB} the existence
of a solution $\vh \in C\Xc^{s+1}_T(\Omegah,\Rbb^N)$ to \eqref{linlocthmC19.1}-\eqref{linlocthmC19.3} that
satisfies the energy estimates \eqref{linlocthmB2.1}.
Appealing
to the finite propagation speed property for wave equations, this implies
that the solution satisfies the original problem \eqref{linlocthmC4.1}-\eqref{linlocthmC4.3} in a neighborhood of the boundary point $x=0\in \del{}\Omega$. Using the finite propagation speed property,
this is enough to construct a local solution in the neighborhood
of the whole boundary by patching the local solutions together. In the interior of $\Omega$, the existence and
uniqueness of solutions follows from standard hyperbolic theory. Patching the interior solution with the one that
is valid in the neighborhood of the boundary yields the desired solution, which it is not difficult to
see satisfies the stated energy
estimate that is inherited from the energy estimates satisfied by the local solutions. Uniqueness of the
 full solution follows from the energy estimate in the standard fashion.

\end{proof}

%--------------------- end comment ------------------------------------------------

\sect{nlinIBVP}{A priori estimates}
We now apply the energy estimates from Theorem \ref{linlocthmC} to obtain a priori estimates for the class of
solutions to the relativistic Euler equations
detailed in Section \ref{fwe}.
The following theorem, which constitutes the main result of this article, contains the
precise statement. However, before stating the theorem, we recall the following definitions and assumptions:
\begin{itemize}
\item[(i)] The vector field $(w^\mu) \in  C^{8}(\overline{U},\Rbb^4)$ is a solution of the Frauendiener-Walton formulation of the Euler equations satisfying assumptions (A.1)-(A.7) from
Section \ref{solassump}. We recall that $w^\mu$ uniquely determines a solution $(\rho,v^\mu)$ to the standard formulation of the Euler equations given by \eqref{eulint3.1}-\eqref{eulint3.2} via
the formulas \eqref{rhov.1} and \eqref{rhov.2}.
\item[(ii)] The vector fields $(e_I^\mu)\in C^7(\overline{U}_T,\Rbb^4)$, $I=1,2,3$, define the frame field completion of
\eqn{recallAa}{
e_0^\mu := w^\mu
}
obtained
by solving the Lie transport equations \eqref{eIsolA}-\eqref{eIsolB} with initial data $e^\mu_I|_{\Omega_0}=f^\mu_I$ chosen
as in Section \ref{IV}.
\item[(iii)] $\theta^i_\mu$ denotes the frame dual to $e_i^\mu$, or in other words, $(\theta^i_\mu) := (e^\mu_i)^{-1}$.
\item[(iv)] The map $\phi=(\phi^\mu)\in C^{8}(\overline{\Omega}_T,\Rbb^4)$, see Section \ref{lag}, defines the change of coordinates from the Eulerian coordinates $(x^\mu)$
to the Lagrangian coordinates adapted to the vector field $e_0^\mu=w^\mu$ according to the formula
\eqn{recallA}{
x^\mu = \phi^\mu(\xb^\lambda).
}
\item[(v)] The field\footnote{From
the assumption $(e_0^\mu) \in C^{8}(\overline{U},\Rbb^4)$, which, we note implies
that $\phi=(\phi^\nu) \in C^{8}(\overline{\Omega}_T,\Rbb^4)$ by solving the ODE \eqref{wbdefB},
it follows from the definition \eqref{etdef} that $(\et^\mu_0) \in C^{8}(\overline{\Omega}_T,\Rbb^4)$.
Similarly, from the fact that $e_I=(e^\mu_I) \in C^7(\overline{U}_T,\Rbb^4)$, we also have
that $(\et^\mu_I) \in C^7(\overline{\Omega}_T,\Rbb^4)$. At first glance, these two
statements seem to imply via definition \eqref{etdefA}
that $\theta^0=(\thetat^0_\mu) \in C^{7}(\overline{\Omega}_T,\Rbb^4)$. However,
due to the relation \eqref{hC}, the definition $\thetat^0_\mu = \theta^0_\mu \circ \phi$,
and the smoothness of the metric $g_{\mu\nu}$, we, in fact, have that
$\thetat^0 \in  C^{8}(\overline{\Omega}_T,\Rbb^4)$.}
$\thetat^0=(\thetat^0_\mu) \in C^{8}(\overline{\Omega}_T,\Rbb^4)$ represents the coframe field components $\theta^0_\mu$  evaluated in the Lagrangian coordinates, that
is
\eqn{recallB}{
\thetat^0_\mu := \theta^0_\mu \circ \phi.
}
\item[(vi)] The field $\psi = (\psi_\nu) \in C^{7}(\overline{\Omega}_T,\Rbb^4)$ denotes the components of $\nabla_{e_0}\theta^0$ evaluated in Lagrangian coordinates so that
\eqn{recallC}{
\psi_\nu := (\nabla_{e_0}\theta^0_\nu)\circ\phi.
}
From $\psi_\nu$, we construct the fields
\eqn{recallD}{
\mu := |\psi|_{\mt} = \sqrt{\mt^{\alpha\beta}\psi_\alpha\psi_\beta}
}
and
\eqn{recallE}{
\psih_\nu := \frac{1}{\mu}\psi_\nu,
}
where $\mt^{\alpha\beta}$ is the positive definite metric defined by
\eqn{recallF}{
\mt^{\alpha\beta} := \gt^{\alpha\beta} - \frac{2}{\gt^{\lambda\sigma}\thetat^0_\lambda\thetat^0_\sigma} \gt^{\alpha\mu}\gt^{\beta\nu}\thetat^0_\mu\thetat^0_\nu  \qquad (\gt^{\alpha\beta}:=g^{\alpha\beta}\circ\phi).
}
The fields  $\psih_\nu$ and $\mu$ are collectively denoted by the vector
\eqn{recallH}{
\Psi := (\psih_\nu,\mu)^{\text{tr}} \in C^{7}(\overline{\Omega}_T,\Rbb^5).
}
To ensure that $\psih_\nu$ is well defined, we further assume that
\leqn{recallG}{
\mu \geq c_\mu > 0 \quad \text{in $\overline{\Omega}_T$}
}
for some positive constant $c_\mu$. As discussed in Remark \ref{murem}, we lose no generality in assuming that \eqref{recallG} holds in addition to the assumptions (A.1)-(A.7) for the
solution $w^\mu$ to the Frauendiener-Walton formulation of the Euler equations.
\item[(vii)] The symmetric matrix $Q$ is defined by
\eqn{recallI}{
Q = \begin{pmatrix} \delta \mu \alpha \pi^{\alpha\beta} & 0 \\ 0 & 0  \end{pmatrix},
}
where $\delta$ is a positive constant to be determined,  $\pi^{\alpha\beta}$ is the projection operator given by
\eqn{recallJ}{
\pi^{\alpha\beta} = \mt^{\alpha\beta} - \psih^\alpha \psih^b \qquad (\psih^\alpha = \mt^{\alpha\nu}\psih_\nu),
}
and $\alpha$, a negative function, is given by
\eqn{recallK}{
\alpha = - \frac{\det(\thetab)|\thetat^3|_{\gt}}{|\gt|^{1/2}}.
}
In the formula for $\alpha$, we recall that $\thetab = (\thetab^i_\mu)$ is the time-independent co-frame defined by \eqref{thetabform}, and $\thetat^3$ is computed using the formula \eqref{conprev.9}.
\end{itemize}

\begin{thm} \label{apthm}
Suppose that $(w^\mu) \in C^{8}(\overline{U},\Rbb^4)$ is a solution of the Frauendiener-Walton formulation of the Euler equations satisfying assumptions (A.1)-(A.7)
from Section \ref{fwe} and the inequality \eqref{recallG}, the fields
\eqn{apthm1}{
(\phi^\mu,\thetat^0_\mu,\Psi) \in C^{8}(\overline{\Omega}_T,\Rbb^4)\times
C^{8}(\overline{\Omega}_T,\Rbb^4) \times C^{7}(\overline{\Omega}_T,\Rbb^5),
}
are constructed from the solution $w^\mu$ according $(i)$-$(vi)$ above, and let
\alin{apthm2}{
\Rc(\xb^0) = \norm{\phi(\xb^0)}_{E^{s+\frac{3}{2},2s-1}}^2+ \norm{\phi(\xb^0)}_{\Ec^{s+1}}^2+
\norm{\thetat^0(\xb^0)}_{\Ec^{s+1}}^2+ \norm{\Psi(\xb^0)}_{\Ec^{s+1}}^2 +
\ip{\delb{0}^{2s+1}\Psi(\xb^0)}{(-Q(\xb^0))\delb{0}^{2s+1}\Psi(\xb^0)}_{\del{}\Omega},
}
where $s=3$.
Then there exists a $T_*=T_*\bigl(\Rc(0)\bigr) \in (0,T)$ such that
\eqn{apthm3}{
\norm{\thetat^0(\xb^0)}_{H^{s+\frac{3}{2}}(\Omega)}^2 + \Rc(\xb^0) \leq C(\Rc(0))
}
for $0\leq \xb^0 \leq T_*$.
\end{thm}
\begin{proof}
We begin by setting $s=3$, and we let
\leqn{apthm4}{
(\phi^\mu,\thetat^0_\mu,\Psi) \in C^{8}(\overline{\Omega}_T,\Rbb^4)\times
C^{8}(\overline{\Omega}_T,\Rbb^4) \times C^{7}(\overline{\Omega}_T,\Rbb^5)
}
be as given in the statement of the theorem. By the assumptions (A.1)-(A.7) from Section \ref{fwe}, the assumption
\eqref{recallG}, and the formulas \eqref{theta3rep1} and \eqref{e3g00}, we see that
\lalin{chdef}{
\ch(t) = &\norm{(\gammat_{00})^{-1}}_{L^\infty(\overline{\Omega}_{t})} +
\norm{\det(J)^{-1}}_{L^\infty(\overline{\Omega}_{t})} + \norm{\mu^{-1}}_{L^\infty(\overline{\Omega}_{t})} \notag \\
&\qquad +\norm{|\thetat^3|_{\gt}^{-1}}_{L^\infty(\Gamma_{t})}
 +\norm{(1-\tilde{s}^2)^{-1}}_{L^\infty(\overline{\Omega}_{t})}
+ \norm{(\tilde{s}^2)^{-1}}_{L^\infty(\overline{\Omega}_{t})}. \label{chdef.1}
}
is bounded for $t\in[0,T]$. We also recall that $\tilde{s}^2$, $-\gammat_{00}$, $\det(J)$, and $|\thetat^3|_{\gt}^2$ are all positive on $\overline{\Omega}_{T}$. For use below, we define
\leqn{cchdef}{
\cch(t) =  \norm{\phi}_{L^\infty(\overline{\Omega}_{t})}
+ \norm{\delb{}\phi}_{L^\infty(\overline{\Omega}_{t})} +
\norm{\delb{}\thetat^0}_{L^\infty(\overline{\Omega}_t)}+
\norm{\Psi}_{L^\infty(\overline{\Omega}_t)}+ \norm{\delb{0}\Psi}_{L^\infty(\overline{\Omega}_t)}.
}

From the analysis carried out in Sections 3 to 5, we know that the triplet \eqref{apthm4} defines a solution
of the IBVP \eqref{cbvpA.1}-\eqref{cbvpA.9} for any of the freely specifiable constants
$\delta, \ep, \kappa \in \Rbb$. In order to apply the energy estimates from
Theorem \ref{linlocthmC} and the elliptic estimates from Theorem \ref{ellipthmB}
to the solution \eqref{apthm4}, we first need
to show the free parameters $\delta,\ep,\kappa$ can be chosen so that
$B^{\Sigma\Lambda}$ and $\Bc^{\Sigma\Lambda}$ satisfy coercive estimates of
the form \eqref{coerc}, and $P$ and $Q$ satisfy the conditions $P+\bigl(2s-\Threehalfs\bigr)\delb{0}Q
-\rc Q$ and $c_Q^{-1}Q\leq Q^2\leq c_Q Q$ for some constants $\rc \in \Rbb$ and $c_Q > 0$. That this is possible
is the content of the following three lemmas.

\begin{lem} \label{coerclem}
There exists constants $\delta=\delta(\ch(t),\cch(t))>0$, $c_{\Bc}=c_\Bc(\ch(t),\cch(t)) > 0$, independent of $\ep,\kappa\in \Rbb$, such that
\eqn{coerclem1}{
%\ip{\delb{\Sigma}\chi}{\At^{\Sigma\Lambda}(\xb^0)\delb{\Lambda}\chi}_{\Omega} \geq c_{\At} \norm{\Db{}\chi}_{L^2(\Omega)}^2,
%\AND
\ip{\delb{\Sigma}\Xi}{\Bc^{\Sigma\Lambda}(\xb^0) \delb{\Lambda}\Xi}_{\Omega} \geq c_{\Bc}
\norm{\Db{}\Xi}_{L^2(\Omega)}^2
}
for all $\xb^0\in [0,t]$, $0\leq t\leq T$, and  $\Xi \in C^1(\overline{\Omega},\Rbb^5)$.
\end{lem}
\begin{proof}
%Since $\phi=(\phi^\mu) \in C^2(\overline{\Omega}_T,\Rbb^4)$, and $(\phi(\xb^\lambda),J(\xb^\lambda))$, where %$J(\xb^\lambda)=(\delb{\nu}\phi^\mu(\xb^\lambda))$, lies
%in an open, bounded subset $\widetilde{\Uc} \subset \Uc$ for all $(\xb^\lambda) \in \overline{\Omega}_T$, it follows from
%smoothness of the metric components $g_{\mu\nu}$ and
From \eqref{chdef.1}, \eqref{defrecLag.12}-\eqref{defrecLag.14},
 and \eqref{conprev.1}, it is clear that there exists constants $c_{\ab}=c_{\ab}(\ch(t),\cch(t))>0$ and
$c_{\mt}=c_{\mt}(\ch(t),\cch(t))>0$ such that
\leqn{ambnd}{
\zeta_\Sigma \zeta_\Lambda \ab^{\Sigma\Lambda}(\xb^\lambda) \geq c_{\ab}|\zeta|^2 \AND
\omega_\mu \omega_\nu \mt^{\mu\nu}(\xb^\lambda) \geq c_{\mt}|\omega|^2
}
for all $(\xb^\lambda,\zeta_\Lambda,\omega_\mu) \in \overline{\Omega}_t\times \Rbb^3\times \Rbb^4$.
Fixing
%$\chi \in C^1(\overline{\Omega})$ and
$\Xi=(\xi,\hat{\xi}_\nu)\in C^1(\overline{\Omega},\Rbb\times\Rbb^4)$, it
then follows immediately from the bounds \eqref{ambnd},
%the fact that
%\eqn{inWc}{
%\bigl(f(\xb^\Sigma),(\phi(\xb^\lambda),J(\xb^\lambda)),\Psi(\xb^\lambda)\bigr) \in \Vc,
%}
the definitions \eqref{cbvpB.2}, \eqref{defrec.8}, \eqref{defrec.17},
 and H\"{o}lder's inequality
 that there exists constants  %$c_{\At},c_{\Ac}>0,c_{\Sc}>0$, independent of $\chi$ and $\Xi$, such that
$c_{\Ac}= c_{\Ac}(\ch(t),\cch(t))>0$, $c_{\Sc}= c_{\Sc}(\ch(t),\cch(t))>0$, independent of $\Xi$, such that
\lalin{AtAccoerc}{
%\ip{\delb{\Sigma}\chi}{\At^{\Sigma\Lambda}(\xb^0)\delb{\Lambda}\chi}_{\Omega} &\geq c_{\At} \norm{\Db{}\chi}_{L^2(\Omega)}^2, \label{AtAccoerc.1}\\
\ip{\delb{\Sigma}\Xi}{(\id+\delta\Pbb)\Ac^{\Sigma\Lambda}(\xb^0) \delb{\Lambda}\Xi}_{\Omega} &\geq c_{\Ac}\bigl(
\norm{\Db{}\xi}_{L^2(\Omega)}^2+ (1+\delta)\norm{\Db{}\hat{\xi}}_{L^2(\Omega)}^2\bigr)
\label{AtAccoerc.2}
\intertext{and}
\bigl|\ip{\delb{\Sigma}\Xi}{2\nml^{[\Sigma}\Sc^{\Lambda]}(\xb^0)\delb{\Lambda}\Xi}_{L^2(\Omega)}\bigr|
&\leq c_{\Sc}\bigl(
\norm{\Db{}\hat{\xi}}_{\Omega}^2+ \norm{\Db{}\xi}_{L^2(\Omega)}\norm{\Db{}\hat{\xi}}_{L^2(\Omega)}\bigr)
\label{AtAccoerc.3}
}
for all $\xb^0\in [0,t]$ and $\delta\geq 0$.
Applying Young's inequality to the right hand side of \eqref{AtAccoerc.3} gives
\leqn{BcoercA}{
\bigl|\ip{\delb{\Sigma}\Xi}{2\nml^{[\Sigma}\Sc^{\Lambda]}(\xb^0)\delb{\Lambda}\Xi}_{L^2(\Omega)}\bigr|
\leq c_{\Sc}\biggl(1 + \frac{c_{\Sc}}{2c_{\Ac}}\biggr)
\norm{\Db{}\hat{\xi}}_{L^2(\Omega)}^2+ \frac{c_{\Ac}}{2}\norm{\Db{}\xi}_{L^2(\Omega)}^2.
}
Setting
\eqn{deltaset}{
\delta=  \frac{c_{\Sc}}{c_{\Ac}}\biggl(1+\frac{c_{\Sc}}{2c_{\Ac}}\biggr)+\frac{1}{2},
}
the coercive estimate
\eqn{BcoercB}{
\ip{\del{\Sigma}\Xi}{\Bc^{\Sigma\Lambda}(\xb^0) \del{\Lambda}\Xi}_{\Omega} \geq \frac{c_{\Ac}}{2}
\norm{\Db{}\Xi}_{L^2(\Omega)}^2, \qquad 0\leq \xb^0 \leq t\leq T,
}
then follows directly from \eqref{cbvpB.3}, and the estimates \eqref{AtAccoerc.2} and \eqref{BcoercA}.
\end{proof}

\begin{lem} \label{coercAlem}
Suppose $\ep >0$ and let $\delta>0$ be as in Lemma \ref{coerclem}.
Then there exist constants $c^{\pm}_{B}=c^{\pm}(\ch(t),\cch(t)) > 0$, independent of $\kappa \in \Rbb$, such that
\eqn{coercAlem1}{
\ip{\delb{\Sigma}\xi}{B^{\Sigma\Lambda}(\xb^0) \delb{\Lambda}\xi}_{\Omega} \geq c^+_{B}
\norm{\Db\xi}_{L^2(\Omega)}^2-c^-_{B}\norm{\xi}_{L^2(\Omega)}^2
}
for all $\xb^0\in [0,t]$, $0\leq t\leq T$, and  $\xi \in C^1(\overline{\Omega},\Rbb^4)$.
\end{lem}
\begin{proof}
%Since $\bigl(f(\xb^\Lambda),(\phi(\xb^0,\xb^\Lambda),J(\xb^0,\xb^\Lambda)),\Psi(\xb^0,\xb^\Lambda)\bigr)
%\in \Vc$ for all $(\xb^0,\xb^\Lambda)\in \Omega_T$,
From the bounds \eqref{ambnd}, the definition  \eqref{cbvpBa.1} and the fact that $\pi^{\mu\nu}$ is non-negative,
it is clear that
there exists a constant $c_{B}=c_{B}(\ch(t),\cch(t))>0$, independent of $\ep >0$, such that
\eqn{coercAlem2}{
\zeta_\Sigma\zeta_\Lambda\ipe{\omega}{B^{\Sigma\Lambda}(\xb^\lambda)\omega} \geq c_{B}|\zeta|^2|\omega|^2
}
for all $(\xb^\lambda,\zeta_\Lambda,\omega_\mu)\in \Omega_t\times\Rbb^3\times\Rbb^4$. By definition,
 this establishes the strong ellipticity of $B^{\Sigma\Lambda}$.
From Theorem 3 in Section 6 of \cite{SimpsonSpector:1987}, we see that the proof
follows if we can verify that the BVP
\lalin{coercAlem3}{
\delb{\Sigma}\bigl(B^{\Sigma\Lambda}(\xb^0)\delb{\Lambda}\xi\bigr) &= 0 \hspace{0.5cm} \text{in $\Omega$,} \label{coercAlem3.1}\\
\nu_\Sigma B^{\Sigma\Lambda}(\xb^0)\delb{\Lambda}\xi &= 0 \hspace{0.5 cm} \text{in $\del{}\Omega$,} \label{coercAlem3.2}
}
satisfies the \emph{strong complementing condition}; see \cite[\S 4]{SimpsonSpector:1987} for a
precise definition, for each $\xb^0\in [0,t]$.

To verify that \eqref{coercAlem3.1}-\eqref{coercAlem3.2} satisfies the complementing condition, we ``freeze''
the coefficients at a point
$(\xb^0,\xb^\Lambda)\in [0,t]\times \del{}\Omega$, and consider the following BVP on the half-plane:
\lalin{coercAlem4}{
\delb{\Sigma}\bigl(B^{\Sigma\Lambda}_{\mu\nu}\delb{\Lambda}\xi^\nu\bigr) &= \alpha^2 \mt_{\mu\nu}\xi^\nu \hspace{0.5cm} \text{in $\Rbb^2\times \Rbb_{>0}$,} \label{coercAlem4.1}\\
\nu_\Sigma B^{\Sigma\Lambda}_{\mu\nu}\delb{\Lambda}\xi^\nu &= 0 \hspace{1.75 cm} \text{in $\Rbb^2$,} \label{coercAlem4.2}
}
where $\alpha \in \Rbb$, $\nu_\Sigma=-\delta^3_\Sigma$ is an outward pointing co-normal, $\nu^\Sigma = -\delta^\Sigma_3$,
\eqn{Bdndef}{
B^{\Sigma\Lambda}_{\mu\nu} = \bigl(\mt_{\mu\nu}+\ep \pi_{\mu\nu}\bigr)\At^{\Sigma\Lambda}+2 \St_{\mu\nu}{}^{[\Sigma}\nu^{\Lambda]},
}
and for notational simplicity,  we use the same notation for any of the previously
defined geometric objects and their frozen versions, e.g. we denote $\pi^{\mu\nu}(\xb^0,\xb^\Lambda)$ by $\pi^{\mu\nu}$. In the following, upper case calligraphic and Fraktur letters, e.g. $\Ac,\Bc$ and $\Af,\Bf$, will run from 1,2 and index the boundary coordinate and frame indices, respectively.

We proceed by making the following definitions:
\alin{frozenA}{
S_{ij}{}^\Lambda &= \St_{\mu\nu}{}^\Lambda \et^\mu_i \et^\nu_j, \\
m_{ij} &= \mt_{\mu\nu}\et^\mu_i\et^\mu_j,\\
\pi_{ij} &= \pi_{\mu\nu}\et^\mu_i\et^\mu_j,
\intertext{and}
\xi^i = &\thetat^i_\nu \xi^\nu.
}
From these definitions,  it is not difficult to verify that \eqref{coercAlem4.1}-\eqref{coercAlem4.2} are equivalent to
\lalin{coercAlem5}{
(m_{ij}+\ep\pi_{ij})\At^{\Sigma\Lambda}\delb{\Sigma}\delb{\Lambda} \xi^j &= \alpha^2 m_{ij}\xi^j \hspace{0.6cm} \text{in $\Rbb^2\times \Rbb_{>0}$,} \label{coercAlem5.1}\\
-(m_{ij}+\ep\pi_{ij}) \At^{3\Lambda}\delb{3}\xi^j &= S_{ij}{}^\Lambda\delb{\Lambda}\xi^j \hspace{0.5 cm} \text{in $\Rbb^2$.} \label{coercAlem5.2}
}
It then follows from \eqref{defrecLag.12},
 \eqref{defrecproj.1}-\eqref{defrecproj.2},
\eqref{defrec.13}, \eqref{conprev.1}-\eqref{conprev.3}, and \eqref{conprev.6}-\eqref{conprev.7}  that
\lalin{frozenB}{
(m_{ij}) &= \begin{pmatrix} 1 & 0  & 0\\ 0 & \fmt_{\Af\Bf} & \fmt_{3\Bf} \\
0 & \fmt_{\Af 3} & \fmt_{33} \end{pmatrix}, \label{frozenB.1}\\
(\pi_{ij}) &=  \begin{pmatrix} 1 & 0  & 0\\ 0 & \fmt_{\Af\Bf} & \fmt_{3\Bf} \\
0 & \fmt_{\Af 3} & \fmt_{33}-\begin{displaystyle}\frac{1}{\fmt_{33}}\end{displaystyle} \end{pmatrix},
\label{frozenB.2}\\
\intertext{and}
(S_{ij}{}^\Lambda \omega_\Lambda) &=  \det(\thetab)\ep_{0123}\begin{pmatrix}0 & 0 & 0 \\
0 & 0 &  \eb^\Cc_\Bf \omega_\Cc \\
0 & -\eb^\Cc_\Af \omega_\Cc &  0
\end{pmatrix}, \label{frozenB.3}
}
where
\eqn{frozenC}{
\ep_{0123} = \ep_{\mu\nu\gamma\lambda}\et^\mu_0 \et^\nu_1 \et^\gamma_2 \et^\lambda_3 = \det(\et).
}
The fact that $\nu_\mu = -\delta^3_\mu$ is an outward
pointing co-normal to $\Gamma_T$ implies via \eqref{OmegaSpanB} that the frame coefficients
$\eb^\mu_\Af$ satisfy $\eb^3_\Af=0$. Using this together with
the boundary condition \eqref{conprev.8} and the identity
\eqn{frozenD}{
\det(J) = \det(J)\det(\eb\thetab)=\det(J\eb)\det(\thetab)=\det(\et)\det(\thetab)
}
allows us to express \eqref{frozenB.3} as
\eqn{frozenE}{
(S_{ij}{}^\Lambda \omega_\Lambda) =  \ell \begin{pmatrix}0 & 0 & 0 \\
0 & 0 &  \eb^\Cc_\Bf \omega_\Cc \\
0 & -\eb^\Cc_\Af \omega_\Cc &  0
\end{pmatrix},
}
where
\leqn{frozenF}{
\ell = -\frac{\det(J)}{\fft} > 0.
}
We note also that, by making a linear change of coordinates if necessary, we can always arrange that
\leqn{frozenG}{
\ab^{\Sigma\Lambda} = \gb^{\Sigma\Lambda} = \delta^{\Sigma\Lambda} \AND \gb^{33}=\fmt^{33}=1.
}

Next, we look for bounded exponential solutions to \eqref{coercAlem5.1}-\eqref{coercAlem5.2} that are of the form
\leqn{frozenH}{
\xi^j= z^j(\xb^3)\exp(i\omega_\Ac\xb^\Ac), \qquad (\omega_\Ac)\in \Rbb^2.
}
Since $\pi_{ij}$ arises from lowering the index of the projection operator $\pi^i_j$ using
the positive definite metric $m_{ij}$, it follows that there exists an invertible matrix $U^{i}_j$
such that
\eqn{frozenK}{
U_i^k U_j^l m_{kl} = \delta_{ij} \AND  U_i^k U_j^l \pi_{kl} = \delta_{ij}-\delta_i^3\delta_j^3.
}
Letting $\Uch^i_j$ denote the inverse of $U^i_j$, we can write \eqref{frozenH} as
\eqn{frozenL}{
\xi^j= U^j_k \zt^k(\xb^3)\exp(i\omega_\Ac\xb^\Ac), \qquad \text{where} \quad \zt^k = \Uch^k_l z^l.
}
Substituting this into \eqref{coercAlem5.1}, while  noting that $\At^{\Sigma\Lambda}=\ell \delta^{\Sigma\Lambda}$
by \eqref{defrec.1}, \eqref{frozenF} and \eqref{frozenG}, we find, after a short calculation,
that $\zt^j$ satisfies the differential equation
\leqn{frozenM}{
\zt^j{}''(\xb^3)-(|\vec{\omega}|^2+\alphat^2)\zt^j(\xb^3) = 0, \qquad \xb^3 > 0,
}
where
\eqn{frozenN}{
\text{$\vec{\omega}=(\omega_\Ac)\quad$ and  $\quad\alphat \in \begin{displaystyle}\left[\frac{\alpha}{\sqrt{\lambda}}
,\frac{\alpha}{\sqrt{(1+\ep)\lambda}}
\right]\end{displaystyle}$}
}
for some $\lambda > 0$.
From standard ODE theory, $\zt(\xb^3)=(\zt^j(\xb^3))$ must be expressible as a linear combination of exponential
solutions of the form
\eqn{frozenO}{
\zt(\xb^3) = \exp(i \xb^3 \omega_3)\tilde{\Upsilon}, \qquad \omega_3\in \Cbb, \; \tilde{\Upsilon}=(\tilde{\Upsilon^j}) \in \Rbb^4.
}
Substituting this into  \eqref{frozenM} gives
\eqn{frozenP}{
\big(\omega_3^2 + |\vec{\omega}|^2+\alphat^2)\tilde{\Upsilon} = 0.
}
Assuming that $\tilde{\Upsilon} \neq 0$, we see that
\eqn{frozenQ}{
\omega_3 = \pm i \sqrt{ |\vec{\omega}|^2+\alphat^2}.
}
Of these two solutions, only
\leqn{frozenR}{
\omega_3 = i \sqrt{ |\vec{\omega}|^2+\alphat^2}
}
is compatible with $\zt(\xb^3)$ being bounded as $\xb^3\rightarrow \infty$; consequently, every bounded solution of \eqref{coercAlem5.1}  must be a linear combination of terms of the form
\leqn{frozenT}{
\xi=(\xi^j) = \exp(i \xb^3 \omega_3)\exp(i\omega_\Ac\xb^\Ac)\Upsilon
}
with $\omega_3$ given by \eqref{frozenR} and $\Upsilon := (U^j_k \tilde{\Upsilon}^k) \in \Rbb^4_{\times}$.

Substituting \eqref{frozenT} into \eqref{coercAlem5.2}, we see,
using \eqref{frozenB.1}-\eqref{frozenB.2} and \eqref{frozenG}, that
\lgath{frozenU}{
(1+\ep)\sqrt{ |\vec{\omega}|^2+\alphat^2}\Upsilon^0 = 0 \label{frozenU.1}
\intertext{and}
M\vec{\Upsilon} = 0, \label{frozenU.2}
}
where
\eqn{frozenV}{
M = \sqrt{ |\vec{\omega}|^2+\alphat^2}\begin{pmatrix}  (1+\ep)\fmt_{\Af\Bf} & (1+\ep)\fmt_{3\Bf} \\
(1+\ep)\fmt_{\Af 3} & (1+\ep)\fmt_{33}-\ep \end{pmatrix} - i  \begin{pmatrix}  0  & \eb^\Cc_\Bf \omega_\Cc \\
-\eb^\Cc_\Af \omega_\Cc &  0 \end{pmatrix}, \AND \vec{\Upsilon}=\begin{pmatrix} \Upsilon^\Bf \\ \Upsilon^3 \end{pmatrix}.
}
With the help of Lemma \ref{matrixlem}, we compute
\leqn{frozenW}{
\det(M) = (|\vec{\omega}|^2+\alphat^2)^{3/2}\begin{pmatrix}  (1+\ep)\fmt_{\Af\Bf} & (1+\ep)\fmt_{3\Bf} \\
(1+\ep)\fmt_{\Af 3} & (1+\ep)\fmt_{33}-\ep \end{pmatrix}
-(1+\ep)(|\vec{\omega}|^2+\alphat^2)^{1/2}\det(\fmt_{\Af\Bf})\check{\fmt}^{\Af\Bf}\eb_\Af^{\Cc} \eb_\Bf^{\Dc}\omega_\Cc\omega_\Dc,
}
where
\eqn{frozenX}{
(\check{\fmt}^{\Af\Bf}) := (\fmt_{\Af\Bf})^{-1}.
}
Next, we decompose $\fmt_{IJ}$ as
\eqn{frozenY}{
(\fmt_{IJ}) = \begin{pmatrix} \fmt_{\Af\Bf} & \fmt_{3\Bf} \\ \fmt_{\Af 3} & \fmt_{33} \end{pmatrix} = \begin{pmatrix} \fmt_{\Af\Bf} & \beta_\Bf\\ \beta_\Af & n^2+\beta_\Cf\beta^{\Cf}\end{pmatrix},
}
where we have defined
\eqn{frozenZ}{
\beta^{\Af} := \check{\fmt}^{\Af\Bf}\beta_\Bf.
}
Inverting $(\fmt_{IJ})$, we find that
\eqn{frozenAA}{
(\fmt^{IJ}) = \begin{pmatrix} \check{\fmt}^{\Af\Bf} + \begin{displaystyle}\frac{\beta^{\Af}\beta^{\Bf}}{n^2}\end{displaystyle} & \beta^{\Bf} \\
\beta^{\Af} & \begin{displaystyle} \frac{1}{n^2} \end{displaystyle}\end{pmatrix}.
}
But, $\fmt^{33} = 1$, and so we have that
\lgath{frozenBB}{
n=1 \notag
\intertext{and}
(\fmt^{IJ})=\begin{pmatrix} \check{\fmt}^{\Af\Bf} + \beta^{\Af}\beta^{\Bf} & \beta^{\Bf}  \\
\beta^{\Af} & \begin{displaystyle} 1 \end{displaystyle}\end{pmatrix}. \label{frozenBB.1}
}
Using Lemma \ref{matrixlemA}, we compute
\alin{frozenCC}{
\det\begin{pmatrix}  (1+\ep)\fmt_{\Af\Bf} & (1+\ep)\fmt_{3\Bf} \\
(1+\ep)\fmt_{\Af 3} & (1+\ep)\fmt_{33}-\ep \end{pmatrix} &= (1+\ep)^2\det(\fmt_{\Af\Bf})\bigl((1+\ep)\fmt_{33}-\ep - (1+\ep)\beta_\Af\beta^\Af\bigr) \\
&= (1+\ep)^2\det(\fmt_{\Af\Bf}),
}
where in deriving the last equality, we used $\fmt_{33}-\beta_\Af\beta^\Af = n^2 =1$. Substituting the above expression into \eqref{frozenW} yields
\eqn{frozenDD}{
\det(M) = (1+\ep)(|\vec{\omega}|^2+\at^2)^{1/2}\det(\fmt_{\Af\Bf})
\bigl[(1+\ep)(|\vec{\omega}|^2+\at^2)-\check{\fmt}^{\Af\Bf}\eb_\Af^{\Cc}\eb_\Bf^\Dc\omega_{\Cc}\omega_{\Dc}\bigr].
}
However,
\alin{frozenEE}{
\check{\fmt}^{\Af\Bf}\eb_\Af^{\Cc}\eb\Bf^\Dc\omega_{\Cc}\omega_{\Dc} &= \bigl(\fmt^{\Af\Bf} -\beta^{\Af}\beta^{\Bf} \bigr)\eb_\Af^{\Cc}\eb_\Bf^\Dc\omega_{\Cc}\omega_{\Dc}
&&\text{(by \eqref{frozenBB.1})} \\
&= \fmt^{IJ}\eb_I^\Cc\eb_J^\Dc \omega_{\Cc}\omega_{\Dc}-2\fmt^{3\Af} \eb_3^{\Cc}\eb_\Af^\Dc\omega_{\Cc}\omega_{\Dc}- \fmt^{33}\bigl(\eb_3^\Bc \omega_\Bc\bigr)^2 - \bigl(\beta^{\Af}\eb_{\Af}^\Bc \omega_\Bc\bigr)^2\\
&= \gb^{\Cc\Dc} \omega_{\Cc}\omega_{\Dc}-2(\beta^{\Af}\eb_\Af^\Cc \omega_\Dc)\eb_3^\Cc\omega_\Cc-\bigl(\eb_3^\Bc \omega_\Bc\bigr)^2- \bigl(\beta^{\Af}\eb_{\Af}^\Bc \omega_\Bc\bigr)^2
&&\text{(by \eqref{frozenBB.1})} \\
&=|\vec{\omega}|^2-\bigl[(\beta^{\Af}\eb_\Af^\Bc + \eb_3^\Bc)\omega_\Bc\bigr]^2 && \text{(by \eqref{frozenG}),}
}
and so, we see that
\eqn{frozenFF}{
\det(M) = (1+\ep)(|\vec{\omega}|^2+\at^2)^{1/2}\det(\fmt_{\Af\Bf})\Bigl(\ep|\vec{\omega}|^2+(1+\ep)\alphat^2+\bigl[(\beta^{\Af}\eb_\Af^\Bc + \eb_3^\Bc)\omega_\Bc\bigr]^2\Bigr).
}
Since $\ep >0$ and $\det(\fmt_{\Af\Bf}) >0$, we conclude that
\eqn{frozenGG}{
\det(M) = 0 \quad \Longleftrightarrow \quad \vec{\omega} = 0 \AND \alpha = 0,
}
 thereby establishing that all solutions $\Upsilon=(\Upsilon^0,\vec{\Upsilon})$ to \eqref{frozenU.1}-\eqref{frozenU.2} with $\vec{\omega}\neq 0$ are trivial.
This, in turn, implies that only bounded solutions of the type \eqref{frozenH} to the frozen BVP \eqref{coercAlem4.1}-\eqref{coercAlem4.2} are
those for which  the $z^j$ are constant (zero if $\alpha\neq 0$) and $\omega_\Ac = 0$. By definition, this verifies that the BVP \eqref{coercAlem3.1}-\eqref{coercAlem3.2} satisfies the
strong complementing condition for each $(\xb^\lambda)\in \Gamma_t$, and the proof is complete.
\end{proof}

\begin{lem} \label{pqlem}
Suppose $\tau \geq 0$ and let $\delta > 0$ be as in Lemma \ref{coerclem}.
Then
\eqn{pglem1d}{
Q^{\tr}=Q, \qquad Q \leq 0, \qquad \text{\em rank}\,Q = 3,
}
and there exist constants $c_Q=c_Q(\ch(t),\cch(t)),\rc=\rc(\ch,\cch(t),\tau)>0$ and $\kappa=\kappa(\ch(t),\cch(t),\tau) < 0$ such that
\eqn{pqlem1c}{
P+\tau \delb{0}Q + \rc Q \leq 0 \AND
\frac{1}{c_Q} Q\leq Q^2 \leq c_Q Q
}
in $\Gamma_t$, $0\leq t \leq T$.
\end{lem}
\begin{proof}
From \eqref{chdef.1}, \eqref{cchdef}, the choice of
initial data from Section \ref{IV}, see in particular \eqref{deteb},  and the formulas
\eqref{defrecLag.12}, \eqref{defrec.5}, and \eqref{defrec.18}-\eqref{defrec.20}, it is not difficult
to verify that there exists constants $c_2>c_1>0$, where $c_1=c_1(\ch(t),\cch(t))$ and $c_2=c_2(\ch(t),\cch(t))$,
such that
\lgath{pqlem4}{
c_{1}|\xi|^2 \leq \xi_\mu\mt^{\mu\nu}(\xb^\lambda)\xi_\nu \leq c_{2}|\xi|^2, \label{pqlem4.1} \\
%c_{1}\ipe{X}{Y}^2 \leq X^\mu\mt_{\mu\nu}(\xb) Y^\nu \leq c_{2}\ipe{X}{Y}, \label{pqlem4.2} \\
c_{1} \leq -\alpha(\xb^\lambda)\mu(\xb^\lambda) \leq c_2, \label{pqlem4.3} \\
%c_{1} \leq |\psih(\xb^\lambda)|^2 \leq c_2, \label{pqlem4.4}\\
\bigl|\mu(\xb^\lambda)\bigl(\delb{0}\alpha(\xb^\lambda)+\lambda(\xb^\lambda)\bigr)\bigr| \leq c_2, \label{pqlem4.5}
}
and\footnote{Given $A\in \Mbb{n}$, the operator norm of $A$ is defined, as usual, by $\norm{A}_{\text{op}} = \sup_{|\xi|=1} |A\xi|$.}
\leqn{pqlem4.6}{
\norm{\beta(\xb^\lambda)}_{\text{op}}+ \norm{\mt(\xb^\lambda)}_{\text{op}}+ \norm{\delb{0}q(\xb^\lambda)}_{\text{op}} \leq c_2
}
for all $(\xb^\lambda)\in \Gamma_t$, $0\leq t \leq T$, and $\xi=(\xi_\mu)
%,X^\mu,Y^\mu
\in \Rbb^4$. Fixing $\xi=(\xi_\mu)\in \Rbb^4$, we then observe that
\lalin{pqlem5}{
|\pi\xi|^2 \leq \frac{1}{c_1}\pi_\mu^\lambda\xi_\lambda\mt^{\mu\nu}\pi_\nu^\omega\xi_\omega
 = \frac{c_{2}}{c_1} \xi_\lambda \pi^{\lambda\omega}\xi_\omega = \frac{c_2}{c_1}\frac{1}{\alpha\mu}\ipe{\xi}{q\xi} \label{pqlem5.1}
}
by \eqref{pqlem4.1}, \eqref{pqlem4.3} and \eqref{defrec.20}.

Next, we estimate
\lalin{pqlem6}{
| \delb{0}q^{\mu\nu}\xi_\mu\xi_\nu| &= \biggl|\biggl(\pi^\lambda_\mu \xi_\lambda +\frac{\psih^\lambda\psih_\mu}{|\psih|_{\mt}^2} \xi_\lambda\biggr)\delb{0}q^{\mu\nu}
\biggl(\pi^\omega_\nu \xi_\omega + \frac{\psih^\omega\psih_\nu}{|\psih|_{\mt}^2} \xi_\omega\biggr)\biggr| && \text{(by \eqref{pirem1})}  \notag \\
%&= \bigl|\pi^\lambda_\mu \xi_\lambda \delb{0}q^{\mu\nu}\pi^\omega_\nu \xi_\omega
%+ 2 \xi_\lambda\pi^\lambda_\mu \delb{0}q^{\mu\nu}\psih_\nu \psih^\omega \xi_\omega
%+ \psih_\mu \delb{0}q^{\mu\nu}\psih_\nu \bigl( \psih^\omega \xi_\omega\bigr)^2 \bigr| \notag \\
&\leq \bigl|\pi^\lambda_\mu \xi_\lambda \delb{0}q^{\mu\nu}\pi^\omega_\nu \xi_\omega\bigr|
+ 2 \biggl|\xi_\lambda\pi^\lambda_\mu \delb{0}q^{\mu\nu}\frac{\psih_\nu}{|\psih|_{\mt}}\biggr|
\biggl|\frac{\psih^\omega}{|\psih|_{\mt}} \xi_\omega\biggr|
+ \biggl|\frac{\psih_\mu}{|\psih|_{\mt}} \delb{0}q^{\mu\nu}\frac{\psih_\nu}{|\psih|_{\mt}}\biggr| \biggl(
\frac{\psih^\omega}{|\psih|_{\mt}|} \xi_\omega\biggr)^2  \notag \\
& \leq \norm{\delb{0}q}_{\text{op}}|\pi\xi|^2 + 2\norm{\delb{0}q}_{\text{op}}\frac{|\psih|}{|\psih|_{\mt}^2}|\pi\xi||\psih^\omega \xi_\omega|
+ \norm{\delb{0}q}_{\text{op}}\frac{|\psih|^2}{|\psih|_{\mt}^4} (\psih^\omega \xi_\omega)^2 \notag\\
&\leq 2\norm{\delb{0}q}_{\text{op}}|\pi\xi|^2
+ 2\norm{\delb{0}q}_{\text{op}}\frac{|\psih|^2}{|\psih|_{\mt}^4} (\psih^\omega \xi_\omega)^2 \notag\\
& \leq  \frac{2c_2}{c_1}\biggl(\frac{c_2}{\mu\alpha}\ipe{\xi}{q\xi}+ \biggl(\frac{\psih^\omega}{|\psih|_{\mt}} \xi_\omega\biggr)^2 \biggr),  \label{pqlem6.1}
}
where in obtaining the last inequality, we used \eqref{pqlem4.1}
and \eqref{pqlem5.1}. Additionally, we estimate
\lalin{pqlem7}{
|2\mu\alpha \pi^{\mu\omega}\beta_\omega^\lambda \pi_\lambda^\mu \xi_\mu\xi_\nu|
&= |2 \mu \alpha \pi^{\mu}_\sigma \mt^{\sigma\omega}\beta_\omega^\lambda \pi_\lambda^\mu \xi_\mu \xi_\nu| \notag \\
&\leq 2c_2^3|\pi\xi|^2 && \text{(by \eqref{pqlem4.1}, \eqref{pqlem4.3} \& \eqref{pqlem4.6})} \notag\\
&\leq \frac{2 c_2^4}{c_1}\frac{1}{\mu\alpha}\ipe{\xi}{q\xi} && \text{(by \eqref{pqlem4.3} \& \eqref{pqlem5.1})} \label{pqlem7.1}
}
and
\lalin{pqlem8}{
\biggl(\mu(\delb{0}\alpha+\lambda)\mt^{\mu\nu} + \kappa \frac{\psih^\mu\psih^\nu}{|\psih|_{\mt}^2}\biggr)\xi_\mu \xi_\nu
&= \biggl(\mu(\delb{0}\alpha+\lambda)\pi^{\mu\nu} + \bigl[\mu(\delb{0}\alpha+\lambda)+\kappa\bigr]\bigr) \frac{\psih^\mu\psih^\nu}{|\psih|^2_{\mt}}\biggr)\xi_\mu \xi_\nu
 \notag \\
&\leq \frac{c_2}{\mu\alpha} \ipe{\xi}{q\xi}+ (c_2+\kappa)\biggl(\frac{\psih^\omega}{|\psih|_{\mt}} \xi_\omega\biggr)^2\biggr), \label{pqlem8.1}
}
where in deriving the last inequality, we have used \eqref{pqlem4.5}.
Fixing $\tau \geq 0$, it then follows from \eqref{defrec.20} and \eqref{pirem2}, and the inequalities \eqref{pqlem6.1}-\eqref{pqlem8.1}
that
\eqn{pqlem9a}{
\bigl(p^{\mu\nu}+\tau\delb{0}q^{\mu\nu} \bigr)\xi_\mu\xi_\nu \leq \biggl(\kappa+c_2+\frac{2c_2\tau}{c_1}\biggr)\biggl(
\frac{\psih^\omega}{|\psih|_{\mt}} \xi_\omega\biggr)^2 +
\biggl(c_2+\frac{2c_2^4}{c_1}+\frac{2c_2^2 \tau}{c_1}\biggr)\frac{1}{-\mu\alpha}\ipe{\xi}{\!-\!q\xi}.
}
Setting $\kappa = -c_2-\frac{2c_2\tau}{c_1}<0$ then yields, with the help of \eqref{pqlem4.3}, the inequality
\eqn{pqlem9}{
\ipe{\xi}{p\xi}+\tau\ipe{\xi}{\delb{0}q\xi} + \biggl(\frac{c_2}{c_1}+\frac{2c_2^4}{c^2_1}+\frac{2c_2^2 \tau}{c^2_1}\biggr)\ipe{\xi}{q\xi} \leq 0,
}
or equivalently, see \eqref{cbvpB.6} and \eqref{cbvpB.7},
\eqn{pqlem10}{
P + \tau \delb{0}Q + \biggl(\frac{c_2}{c_1}+\frac{2c_2^4}{c^2_1}+\frac{2c_2^2 \tau}{c^2_1}\biggr)Q \leq 0.
}

To conclude, we note that the statements
\eqn{pqlem11}{
Q^{\tr}=Q, \qquad Q \leq 0, \AND \text{rank}\,Q = 3
}
are a direct consequence of the formulas \eqref{cbvpB.6} and \eqref{defrec.20}, the fact that $\pi^\mu_\nu$ is a projection operator with a 1-dimensional kernel, and the positivity of $\delta$, $\mu$, and $-\alpha$.
Moreover, from \eqref{pqlem4.1} and the definition $(\mt_{\mu\nu})=(\mt^{\mu\nu})^{-1}$, we see
that
\eqn{pqlem12}{
\frac{1}{c_{2}}|\omega|^2\leq \omega^\mu\omega^\nu \mt_{\mu\nu} \leq \frac{1}{c_1}|\omega|^2 , \quad \forall \: \omega=(\omega^\mu)\in \Rbb^4.
}
From this, \eqref{defrecproj.2}, and the fact that $\pi^{\mu}_{\nu}$ is a projection operator, we obtain
\alin{pqlem13}{
\frac{1}{c_{2}}\omega_\mu\pi^{\mu\alpha}\delta_{\alpha\beta}\pi^{\beta\nu}\omega_\nu \leq
\omega_\mu\pi^{\mu\alpha}\mt_{\alpha\beta}\pi^{\beta\nu}\omega_\nu
= \omega_\mu\pi^{\mu\nu}\omega_\nu
\leq \frac{1}{c_{1}} \omega_\mu\pi^{\mu\alpha}\delta_{\alpha\beta}\pi^{\beta\nu}\omega_\nu \,
\quad \forall\: \omega=(\omega_\mu) \in \Rbb^4,
}
or equivalently
\eqn{pqlem14}{
c_1 \pi \leq \pi^2 \leq c_{2} \pi,
}
where $\pi=(\pi^{\mu\nu})$. We see immediately from this result and  \eqref{pqlem4.3}
that $Q$ satisfies
\eqn{pqlem15}{
c_1 Q \leq Q^2 \leq c_{2} Q.
}
\end{proof}

Next, we need to estimate the various coefficients
from the equations \eqref{cbvpA.1}-\eqref{cbvpA.4}.
We start by differentiating \eqref{cfsmoothA.5} and \eqref{cfsmoothA.6} to obtain
\leqn{aplemA0a}{
\delb{0}^{2s}P = \lc^\mu_1\delb{\mu}\theta + h_1 \AND
\delb{0}^{2s}G = \lc^\mu_2\delb{\mu}\theta + h_2
}
where
\leqn{aplemA0b}{
\theta = \delb{0}^{2s-1}(\Jf(\phi,\thetat^0,\Psi)) \AND
h_a = [\delb{0}^{2s},\lc^\mu_a]\delb{\mu}\Jc(\phi,\thetat^0) + \delb{0}^{2s-1}\dot{\kc}_a, \quad a=1,2.
}
We then define
\alin{aplemA0c}{
\alpha_1(\xb^0) & = 1+\norm{\,\dot{\!\Bc}(\xb^0)}_{E^{s}}^2 +
\norm{\delb{0}P(\xb^0)}_{E^{s,2s-2}}^2 +\norm{\dot{Q}(\xb^0)}_{\Ec^{s+\frac{1}{2}}}^2
+\norm{\theta(\xb^0)}_{H^1(\Omega)}^2+\norm{\delb{0}\theta(\xb^0)}_{H^1(\Omega)}^2+
\norm{\vec{h}(\xb^0)}_{\Ec^1}^2,
\intertext{and}
\alpha_2(\xb^0) &= 
\norm{\Hc(\xb^0)}_{E^{s-1,2s-2}}^2+ \norm{\delb{0}\Hc(\xb^0)}_{E^{s-1,2s-2}}^2+ \norm{\delb{0}^{2s}\Hc(\xb^0)}_{L^2(\Omega)}^2
+\norm{\Mcal(\xb^0)}_{E^{s}}^2+\norm{\,\,\dot{\!\!\Mcal}(\xb^0)}_{E^{s}}^2 \notag \\
&\hspace{1.0cm}+\norm{\Kc(\xb^0)}_{E^{s,2s-2}}^2
 +\norm{\delb{0}\Kc(\xb^0)}_{E^{s,2s-2}}^2
+ \norm{\vec{h}(\xb^0)}_{\Ec^1}^2
+\norm{\theta(\xb^0)}_{H^1(\Omega)}^2+\norm{\delb{0}\theta(\xb^0)}_{H^1(\Omega)}^2,
}
where $\vec{\lc}=(\lc^\mu_1,\lc^\mu_2)$, $\vec{h}=(h_1,h_2)$, and the coefficients with a dot, e.g.
$\dot{\!\Bc}$, $\dot{Q}$, etc., are defined by \eqref{cfdsmoothA.1}-\eqref{cfdsmoothA.6}, that is
via time differentiation.

\begin{lem} \label{aplemA}
Let
\eqn{aplemA1}{
\bar{\Rc} = \norm{\phi}_{E^{s+\frac{3}{2},2s-1}}^2+\norm{\thetat^0}_{E^{s+\frac{3}{2},2s-1}}^2+
\norm{\phi}_{\Ec^{s+1}}^2+\norm{\thetat^0}_{\Ec^{s+1}}^2+\norm{\Psi}_{\Ec^{s+1}}^2.
}
Then
\alin{aplemA2}{
\alpha_1+\alpha_2 \leq C(\bar{\Rc}).
}
\end{lem}
\begin{proof} We will only estimate the terms involving $\Kc$ and $h_a$ since the estimates for the remainder of
the terms follow from similar arguments. First, we estimate $\Kc$ by
\eqn{aplemA3}{
\norm{\Kc}_{E^{s,2s-2}}+\norm{\delb{0}\Kc}_{E^{s,2s-2}}
\lesssim \norm{\Kc}_{\Ec^{s+\frac{1}{2}}} \leq C(\bar{\Rc})
}
where the last inequality follows from \eqref{cfsmoothA.5} and an application of Proposition \ref{stpropEa}.
Next, by \eqref{cfdsmoothA.2} and another application of Proposition \ref{stpropEa}, we see that
\leqn{aplemA4}{
\norm{\delb{0}^{2s-1}\dot{\kc}_a}_{\Ec^1} \lesssim
\norm{\dot{\kc}_a}_{\Ec^{s+\frac{1}{2}}} \leq C(\bar{\Rc}),
}
while
\lalin{aplemA5}{
\norm{[\delb{0}^{2s},\lc_a^\mu]\delb{\mu}\Jc(\phi,\thetat^0)}_{H^1(\Omega)}
&\lesssim \norm{\dot{\lc}_a}_{E^{s+\frac{1}{2},2s-1}}\norm{\delb{} \Jc(\phi,\thetat^0)}_{E^{s+\frac{1}{2},2s-1}}
&&\text{(by Prop. \ref{stcomA})}\notag \\
& \lesssim \norm{\dot{\lc}_a}_{\Ec^{s+\frac{1}{2}}}\norm{\delb{} \Jc(\phi,\thetat^0)}_{\Ec^{s+\frac{1}{2}}} \notag \\
&\leq C(\bar{\Rc}) \label{aplemA5.1}
}
where the last inequality follows from \eqref{cfdsmoothA.1} and Proposition \ref{stpropEa}.
Noting that
\eqn{aplemA6}{
\delb{0}\bigl([\delb{0}^{2s},\lc_a^\mu]\delb{\mu}\Jc(\phi,\thetat^0)\bigr)
= [\delb{0}^{2s+1},\lc_a^{\mu}]\delb{\mu}\Jc(\phi,\thetat^0) - \dot{\lc}_a^{\mu}\delb{0}^{2s}\delb{\mu}\Jc(\phi,\thetat^0),
}
we find from \eqref{cfdsmoothA.1}, and Propositions \ref{stpropA} and \ref{stpropE} that
\leqn{aplemA7}{
\norm{\delb{0}\bigl([\delb{0}^{2s},\lc_a^\mu]\delb{\mu}\Jc(\phi,\thetat^0)\bigr)}_{L^2(\Omega)}
\lesssim \norm{\dot{\lc}_a}_{E^s}\norm{\delb{}\Jc(\phi,\thetat^0)}_{E^s} \leq C(\bar{\Rc})
}
Combining the estimates \eqref{aplemA4}, \eqref{aplemA5.1} and \eqref{aplemA7}, we see that $h_a$ satisfies
\eqn{aplemA8}{
\norm{h_a}_{\Ec^1} \leq C(\bar{\Rc}).
}
\end{proof}

\begin{lem} \label{aplemC}
Let
\eqn{aplemC1}{
%\hat{\Rc} = \norm{\phi}_{E^{s+\frac{3}{2},2s-1}(\Omega)}^2+\norm{\thetat^0}_{E^{s+1,2s-2}}^2+\norm{\Psi}_{\Ec^{s+1}}^2.
\Rh = \norm{\phi}_{E^{s+\frac{3}{2},2s-1}}^2+
\norm{\phi}_{\Ec^{s+1}}^2+\norm{\thetat^0}_{\Ec^{s+1}}^2+\norm{\Psi}_{\Ec^{s+1}}^2.
}
Then
\alin{aplemC2}{
\norm{\delb{0}(B^{0\beta}\delb{\beta}\thetat^0)}_{H^{s-\frac{1}{2}}(\Omega)}+&\norm{\delb{0}M^0}_{H^{s-\frac{1}{2}}(\Omega)}+ \norm{H}_{H^{s-\frac{1}{2}}(\Omega)} \\
\qquad +\norm{B^{\Sigma\Lambda}}_{H^{s+\frac{1}{2}}(\Omega)}& + \norm{B^{\Sigma 0}\delb{0}\thetat^0}_{H^{s+\frac{1}{2}}(\Omega)}+\norm{M^\Sigma}_{H^{s+\frac{1}{2}}(\Omega)}
+ \norm{K}_{H^{s+\frac{1}{2}}(\Omega)}
\leq %1 + C(\hat{\Rc})\hat{\Rc}.
C(\hat{\Rc}).
}
\end{lem}
\begin{proof}
The proof follows from the same arguments used to establish Lemma \ref{aplemA}.
\end{proof}
Viewing equation \eqref{cbvpA.1} as an elliptic equation for $\thetat^0_\mu$, we find that
\leqn{apellip}{
\norm{\thetat^0}_{H^{s+\frac{3}{2}}(\Omega)}^2 \leq C\bigl(\Rh\bigr).
}
by Theorem \ref{ellipthmB} and Lemma \ref{aplemC}. But then we note that
\alin{apellipA}{
\norm{\thetat^0}^2_{E^{s+\frac{3}{2},2s-1}} &= \norm{\thetat^0}^2_{H^{s+\frac{3}{2}}(\Omega)}
+\norm{\delb{0}\thetat^0}_{E^{s+1,2s-2}}^2, \\
& \leq  \norm{\thetat^0}^2_{H^{s+\frac{3}{2}}(\Omega)} +
C(\Rh)
}
where in the last inequality we used the evolution equation \eqref{cbvpA.6} and Proposition \ref{stpropB}.
Combining the above estimate with \eqref{apellip}, we see that
\leqn{apellipB}{
\norm{\thetat^0}^2_{E^{s+\frac{3}{2},2s-1}} \leq C(\Rh),
}
which in turn, implies that
\leqn{apellipC}{
\Rb \leq C(\Rh).
}
We also observe that the estimate
\leqn{apellipD}{
\norm{\phi(t)}_{E^{s+\frac{3}{2},2s-1}}^2+
\norm{\phi(t)}_{\Ec^{s+1}}^2+\norm{\thetat^0(t)}_{\Ec^{s+1}}^2
\leq \Rh(0) + \int_0^t C(\Rh(\tau))\, d\tau
}
follows directly from the evolution equations \eqref{cbvpA.5}-\eqref{cbvpA.6} and Proposition \ref{stpropF}.

Applying the energy estimates from Theorem \ref{linlocthmC} to \eqref{cbvpA.3}-\eqref{cbvpA.4}, which is possible in view
of Lemma \ref{coerclem} and Lemma \ref{pqlem} with $\tau=2s-\frac{3}{2}$,
we obtain, with the help of \eqref{apellipD} and Lemma \ref{aplemA}, the energy estimate
\lalin{PsiA}{
&\norm{\Psi(\xb^0)}_{\Ec^{s+1}}^2+\ip{\delb{0}^{2s+1}\Psi(\xb^0)}{(-Q)\delb{0}^{2s+1}\Psi(\xb^0)}_{\Omega}
\leq C(\rho_{t})\biggl(\Rc(0)+ \int_0^{\xb^0} C(\Rc(\tau)) d\tau\biggr), \quad 0\leq \xb^0 \leq t, \label{PsiA.1}
}
where $t\in (0,T]$,
\eqn{apRcdef}{
\Rc = \Rh+ \ip{\delb{0}^{2s+1}\Psi}{(-Q)\delb{0}^{2s+1}\Psi}_{\Omega},
}
and
\alin{rhoconstA}{
\rho(t) &= \norm{\Bc}_{X_{t}^{s}}+ \norm{P}_{X_{t}^{s,2s-2}}+
\norm{\delb{0}^{2s-1}P}_{L^{\infty}([0,t],L^2(\del{}\Omega))]} +\norm{Q}_{L^\infty([0,t],H^s(\Omega))}\\
&\hspace{4.0cm}+
 \norm{\dot{Q}}_{X_{t}^{s,2s-2}}+ \norm{\vec{\lc}}_{W^{1,\infty}([0,t],W^{1,\infty}(\Omega))} +
 \ch(t)+\cch(t)
}
with $\ch(t)$ and $\cch(t)$ as defined previously by \eqref{chdef.1} and \eqref{cchdef}.

Collectively, the estimates \eqref{apellipD}-\eqref{PsiA.1} imply
that
\eqn{apRcA}{
\Rc(\xb^0) \leq C(\rho(t))\biggl(\Rc(0)+\int_0^{\xb^0} C(\Rc(\tau))\, d\tau\biggr), \quad 0\leq \xb^0 \leq t \leq T,
}
and hence, that
\leqn{apRcAa}{
\norm{\Rc}_{L^\infty([0,t])} \leq C(\rho(t))\biggl(\Rc(0)+\int_0^{t}
C\Bigl(\norm{\Rc}_{L^\infty([0,\tau])}\Bigr)\, d\tau\biggr).
}
\begin{lem} \label{aplemD}
\leqn{aplemD1}{
\rho(t) \leq C(\Rc(0)) + t C\bigl(\norm{\Rc}_{L^\infty([0,t])}\bigr), \quad 0\leq t \leq T.
}
\end{lem}
\begin{proof}
We only estimate the terms involving $P$ since the rest of the terms can be bounded in a similar fashion. We begin
with the inequality:
\alin{aplemD2}{
\norm{P(\xb^0)}_{E^{s,2s-2}}^2 &\leq \norm{P(0)}_{E^{s,2s-2}}^2
+ \int_0^{\xb^0} \norm{P(\tau)}_{E^{s,2s-2}}^2 + \norm{\delb{0}P(\tau)}_{E^{s,2s-2}}^2\, d\tau
&& \text{(by Prop. \ref{stpropF})} \\
& \leq \Rc(0) + \int_0^{\xb^0} C(\bar{\Rc}(\tau))\, d\tau && \text{(by Lemma \ref{aplemA})} \\
& \leq \Rc(0) + \int_0^{\xb^0}  C(\Rc(\tau))\, d\tau && \text{(by \eqref{apellipC})}.
}
From this, it is then clear that
\eqn{aplemD3}{
\norm{P}_{X^{s,2s-2}_t}^2 \leq \Rc(0) + t C\bigl(\norm{\Rc}_{L^\infty([0,t])}\bigr), \quad 0\leq t \leq T.
}

Next, from \eqref{aplemA0a} and the argument used to derive \eqref{linlocthmA16e.1},
we have that
\lalin{aplemD4}{
&\norm{\delb{0}^{2s-1}P(\xb^0)}^2_{L^2(\del{}\Omega)} \lesssim \norm{\delb{0}^{2s-1}P(0)}^2_{L^2(\del{}\Omega)} \notag \\
&\hspace{1.0cm}+
\norm{\lc_{1}}_{L^\infty([0,t],W^{1,\infty}(\Omega))}
\int_0^{\xb^0} \norm{\delb{0}^{2s-1}P(\tau)}_{H^1(\Omega)}^2
+\norm{h_1(\tau)}_{H^1(\Omega)}^2+\norm{\theta(\tau)}^2_{H^1(\Omega)}\, d\tau \label{aplemD4.1}
}
for $0\leq \xb^0\leq t \leq T$. Noting that
\alin{aplemD5}{
\norm{\delb{0}^{2s-1}P(0)}_{L^2(\del{}\Omega)} &\lesssim \norm{\delb{0}^{2s-1}P(0)}_{H^1(\Omega)}
\lesssim \norm{\delb{0}P(0)}_{E^{s,2s-2}} && \text{(by Theorem \ref{trace})},\\
\norm{\lc_{1}}_{L^\infty([0,t],W^{1,\infty}(\Omega))} &\lesssim
\norm{(\phi,\Db\phi,\Psi)}_{L^\infty([0,t],W^{1,\infty}(\Omega))}
&& \text{(see \eqref{cfsmoothA.6})}
\intertext{and}
\norm{\theta}_{H^1(\Omega)} &\leq \norm{\Jf(\phi,\thetat^0,\Psi)}_{\Ec^{s+\frac{1}{2}}},
&&\text{(see \eqref{aplemA0b})}
}
we conclude from Lemma \ref{aplemA}, Theorem \ref{FSobolev}, Proposition \ref{stpropD},  and the estimates \eqref{apellipB} and \eqref{aplemD4.1} that
\eqn{aplemB6}{
\norm{\delb{0}^{2s-1}P}^2_{L^\infty([0,t],L^2(\del{}\Omega))} \lesssim \Rb(0)
+ t C\bigl(\norm{\Rc}_{L^\infty([0,t])}\bigr), \quad 0\leq t \leq T.
}
\end{proof}

Taken together, \eqref{apRcAa} and \eqref{aplemD1} imply that $\Rc(\xb^0)$ satisfies an estimate of the form
\leqn{apRcB}{
\norm{\Rc}_{L^\infty([0,t])} \leq c_1\Bigl(\Rc(0)+tc_2\bigl(\norm{\Rc}_{L^\infty([0,t])}\bigr)\Bigr)
}
for some non-decreasing, continuous functions $c_{i}\, :\, [0,\infty) \rightarrow [1,\infty)$, $i=1,2$.
Setting
\eqn{aprdef}{
r(t) = t - \frac{T}{1+c_2\bigl(\norm{\Rc}_{L^\infty([0,t])}\bigr)},
}
we see that $r\, : \, [0,T] \rightarrow \Rbb$ is continuous and satisfies $r(0) <0$
and $r(T)>0$. The Intermediate Value Theorem then guarantees the existence of a $T_*\in (0,T)$ such
that $r(T_*)=0$, or equivalently
\eqn{T*boundA}{
T_* = \frac{1}{1+c_2\bigl(\norm{\Rc}_{L^\infty([0,T_*])}\bigr)}.
}
Substituting this into \eqref{apRcB} yields the estimate
\eqn{apRcC}{
\norm{\Rc}_{L^\infty([0,T_*])} \leq c_1(\Rc(0)+1),
}
which, in turn, implies that
\eqn{T*boundB}{
T_* \geq \frac{T}{1+c_2\bigl(c_1(\Rc(0)+1)\bigr)}.
}
\end{proof}

\begin{rem}
From the proof of Theorem \ref{apthm}, it is clear that the same result also holds
in the non-physical dimensions $n\neq 3$ for $s$ satisfying $s>n/2+1$ and $s=k/2$ with
$k\in \Zbb$.
\end{rem}

\bigskip

\noindent \textit{Acknowledgements:} This work was partially supported by the Australian Research Council grant FT1210045.

\appendix

\sect{calculus}{Calculus inequalities}

In this appendix, we state, for the convenience of the reader, some well known calculus inequalities.
As above, $\Omega$ is a bounded, open subset of $\Rbb^n$,
$n\geq 2$, with
smooth boundary.

\subsect{Sineq}{Spatial inequalities}

The proof of the following calculus inequalities are well known and may be found, for example, in
the references \cite{AdamsFournier:2003,Friedman:1976,RunstSickel:1996,TaylorIII:1996}.
In the following, $M$ will denote either $\Omega$, or a closed $n$-manifold.

\begin{thm}{\emph{[H\"{o}lder's inequality]}} \label{Holder}
If $0< p,q,r \leq \infty$ satisfy $1/p+1/q = 1/r$, then
\eqn{HolderA}{
\norm{uv}_{L^r(M)} \leq \norm{u}_{L^p(M)}\norm{v}_{L^q(M)}
}
for all $u\in L^p(M)$ and $v\in L^q(M)$.
\end{thm}

\begin{thm}{\emph{[Integral Sobolev inequalities]}} \label{ISobolev}
Suppose $s\in \Zbb_{\geq 1}$ and $1\leq p < \infty$.
\begin{enumerate}[(i)]
\item If $sp<n$, then
\eqn{ISobolev1}{
\norm{u}_{L^q(M)} \lesssim \norm{u}_{W^{s,p}(M)}, \qquad p\leq q \leq \frac{np}{n-s p},
}
for all $u\in W^{s,p}(M)$.
%\item[(ii)] If $sp=n$, then
%\eqn{Sobolev2}{
%\norm{u}_{L^q(\Omega)} \lesssim \norm{u}_{W^{s,p}(\Omega)} \qquad p\leq q < \infty
%}
%for all $u\in W^{s,p}(\Omega)$.
\item (Morrey's inequality)
If  $sp > n$, then
\eqn{ISobolev2}{
\norm{u}_{C^{0,\mu}(M)} \lesssim \norm{u}_{W^{s,p}(M)}, \qquad 0 < \mu \leq \min\{1,s-n/p\},
}
for all $u\in W^{s,p}(M)$.
\end{enumerate}
\end{thm}

\begin{thm}{\emph{[Fractional Sobolev inequalities]}} \label{FSobolev} Suppose $s\in \Rbb$ and
$1< p < \infty$.
\begin{itemize}
\item[(i)] If $sp<n$, then
\eqn{FSobolev1}{
\norm{u}_{L^q(M)} \lesssim \norm{u}_{W^{s,p}(M)}, \qquad p\leq q \leq \frac{np}{n-s p},
}
for all $u\in W^{s,p}(M)$.
%\item[(ii)] If $sp=n$, then
%\eqn{Sobolev2}{
%\norm{u}_{L^q(\Omega)} \lesssim \norm{u}_{W^{s,p}(\Omega)} \qquad p\leq q < \infty
%}
%for all $u\in W^{s,p}(\Omega)$.
\item[(ii)] If $sp > n$, then
\eqn{FSobolev2}{
\norm{u}_{L^\infty(M)} \lesssim \norm{u}_{W^{s,p}(M)}
}
for all $u\in W^{s,p}(M)$.
\end{itemize}
\end{thm}

\begin{comment}
\begin{thm}{\emph{[Product and commutator estimates]}} \label{Iprod} $\;$
\begin{enumerate}[(i)]
\item
Suppose $1\leq p_1,p_2,q_1,q_2\leq \infty$, $s=|\alpha|\in \Zbb_{\geq 1}$, and
\eqn{Iprod.1}{
\frac{1}{p_1}+\frac{1}{p_2} = \frac{1}{q_1} + \frac{1}{q_2} = \frac{1}{r}.
}
Then
\alin{Iprod.2}{
\norm{D^\alpha(uv)}_{L^r(M)} \lesssim \norm{u}_{W^{s,p_1}(\Omega)}\norm{v}_{L^{q_1}(\Rbb^n)} + \norm{u}_{L^{p_2}(M)}\norm{v}_{W^{s,q_2}(M)} \label{clacpropB.2.1}
\intertext{and}
\norm{D^\alpha(uv)-uD^\alpha v}_{L^r(M)} \lesssim \norm{Du}_{L^{p_1}(\Omega)}\norm{v}_{W^{s-1,q_1}(M)} + \norm{Du}_{
W^{s-1,p_2}(M)}\norm{v}_{L^{q_2}(M)}
}
for all $u,v \in C^\infty(\overline{M})$.
\item[(ii)]  If $s_1,s_2,s_3\in \Zbb_{\geq 0}$, $s_1,s_2\geq s_3\geq 0$, $1\leq p \leq \infty$, and $s_1+s_2-s_3 > n/p$, then
\eqn{Iprod.3}{
\norm{uv}_{W^{s_3,p}(M)} \lesssim \norm{u}_{W^{s_1,p}(M)}\norm{v}_{W^{s_2,p}(M)}
}
for all $u\in W^{s_1,p}(M)$ and $W^{s_2,p}(M)$.
\end{enumerate}
\end{thm}
\end{comment}

\begin{thm} {\emph{[Trace theorem]}} \label{trace}
If $s>1/2$, then the trace operator
\eqn{trace1}{
H^s(\Omega) \ni u \longmapsto u|_{\del{}\Omega} \in H^{s-\frac{1}{2}}(\del{}\Omega)
}
is well-defined, continuous (i.e. bounded), and surjective.
\end{thm}

\begin{lem} {\emph{[Ehrling's lemma]}} \label{Ehrling}
Suppose $1<p<\infty$, $s_0 < s < s_1$. Then for any $\delta>0$ there exists a constant $C=C(\delta)$ such
that
\eqn{Ehrling1}{
\norm{u}_{W^{s,p}(M)} \leq \delta \norm{u}_{W^{s_1,p}(M)} + C(\delta)\norm{u}_{W^{s_0,p}(M)}
}
for all $u\in W^{s_1,p}(M)$.
\end{lem}
\begin{comment}

\begin{thm}{\emph{[Interpolation]}} \label{interp}
Suppose $\ep_0 >0$, $1\leq p \leq \infty$, $k,s\in \Zbb_{\geq 0}$ and $k\leq s$. Then there exists
a constant $K>0$ such that
\eqn{thm}{
|u|_{k,p} \leq K\bigl(\ep |u|_{s,p} + \ep^{-k/(s-k)}\norm{u}_{L^p(\Omega)}\bigr)
}
for $0<\ep \leq \ep_0$, where $|\cdot|_{k,p}$ is the seminorm defined by
\eqn{interpB}{
|u|_{k,p} = \left(\sum_{|\alpha|=k}\norm{D^\alpha u}_{L^p(\Omega)}^p\right)^{1/p}.
}
\end{thm}

\end{comment}
%------------------------- end comment -00000000000000000000000

\begin{thm}{\emph{[Integral multiplication inequality]}} \label{Iprod} $\;$
Suppose $s_1,s_2,s_3\in \Zbb_{\geq 0}$, $s_1,s_2\geq s_3$, $1\leq p \leq \infty$, and $s_1+s_2-s_3 > n/p$. Then
\eqn{Iprod.3}{
\norm{uv}_{W^{s_3,p}(M)} \lesssim \norm{u}_{W^{s_1,p}(M)}\norm{v}_{W^{s_2,p}(M)}
}
for all $u\in W^{s_1,p}(M)$ and $v\in W^{s_2,p}(M)$.
\end{thm}

\begin{thm}{\emph{[Fractional multiplication inequality]}} \label{calcpropB}
Suppose $1<p <\infty$, $s_1,s_2,s_3\in \Rbb$, $s_1+s_2 > 0$, $s_1,s_2 \geq s_3$, and
$s_1+s_2 - s_3 > n/p$. Then
\eqn{calcpropB.1}{
\norm{u v}_{W^{s_{3},p}(M)} \lesssim \norm{u}_{W^{s_1,p}(M)} \norm{v}_{W^{s_2,p}(M)}
}
for all $u \in W^{s_1,p}(M)$ and $v \in W^{s_2,p}(M)$.
\end{thm}

\subsect{STineq}{Spacetime inequalities}
In this section, we establish a number of product, commutator, and related estimates for the
spaces $X_T^{s,r}(M)$ and $\Xc_T^s(M)$. We begin with a product estimate.
\begin{prop} \label{stpropA}
Suppose $s_1=k_1/2$ and  $s_2=k_2/2$ for $k_1,k_2\in \Zbb_{\geq 0}$, $s_3\in \Rbb$, $s_1,s_2 \geq s_3$,
$s_1+s_2 - s_3 > n/2$, $r\in \Zbb$, and $0\leq r \leq 2 s_3$.
Then
\eqn{stpropA1}{
\norm{\del{t}^\ell(u(t)v(t))}_{H^{s_3-\frac{\ell}{2}}(M)}\leq
\norm{u(t)v(t)}_{E^{s_3,r}}  \lesssim \norm{u(t)}_{E^{s_1,r}} \norm{v(t)}_{E^{s_2,r}}
}
for $0\leq t \leq T$, $0\leq \ell \leq r$, and all $u\in X^{s_1,r}_T(M)$ and $v\in X^{s_2,r}_T(M)$.
\end{prop}
\begin{proof}
By the assumptions on $s_1$, $s_2$ and $s_3$, we note that
\eqn{stpropA2}{
\biggl[s_1 - \biggl(\frac{\ell-m}{2}\biggr)\biggr]+\biggl[s_2 - \biggl(\frac{m}{2}\biggr)\biggr]
-\biggl[s_3 - \biggl(\frac{\ell}{2}\biggr)\biggr] - \frac{n}{2} = s_1+s_2-s_3 -\frac{n}{2} > 0,
}
\eqn{stpropA3}{
\biggl[s_1 - \biggl(\frac{\ell-m}{2}\biggr)\biggr]
-\biggl[s_3 - \biggl(\frac{\ell}{2}\biggr)\biggr]  = s_1-s_3 + \frac{m}{2} \geq 0, \quad m\geq 0,
}
\eqn{stpropA4}{
\biggl[s_2 - \biggl(\frac{m}{2}\biggr)\biggr]
-\biggl[s_3 - \biggl(\frac{\ell}{2}\biggr)\biggr] = s_2-s_3 +\biggl(\frac{\ell-m}{2}\biggr) \geq 0, \quad m\leq \ell,
}
and
\eqn{stpropA5}{
\biggl[s_1 - \biggl(\frac{\ell-m}{2}\biggr)\biggr]+\biggl[s_2 - \biggl(\frac{m}{2}\biggr)\biggr] = s_1+s_2 - \frac{\ell}{2} > s_3 +\frac{n}{2} - \frac{\ell}{2}
\geq 0, \quad \ell \leq 2s_3+n.
}
From these observations, we see via the fractional multiplication inequality from Theorem \ref{calcpropB} that
\leqn{stpropA6}{
\norm{\del{t}^{\ell-m}u \del{t}^m v}_{H^{s_3-\frac{\ell}{2}}(M)} \lesssim
\norm{\del{t}^{\ell-m}u}_{H^{s_1-\frac{\ell-m}{2}}(M)}
\norm{\del{t}^{m}v}_{H^{s_2-\frac{m}{2}}(M)}, \quad 0\leq m \leq \ell \leq 2 s_3 + n.
}

Next, differentiating the product $uv$ $\ell$-times with respect to $t$, we obtain
\eqn{stpropA7}{
\del{t}^\ell(uv) = \sum_{m=0}^\ell \binom{\ell}{m}\del{t}^{\ell-m}u\del{t}^m v.
}
Applying the triangle inequality and \eqref{stpropA6} to this expression, we see that
\eqn{stpropA7}{
\norm{\del{t}^\ell(uv)}_{H^{s_3-\frac{\ell}{2}}(M)} \lesssim \sum_{m=0}^\ell \norm{\del{t}^{\ell-m}u}_{H^{s_1-\frac{\ell-m}{2}}(M)}
\norm{\del{t}^{m}v}_{H^{s_2-\frac{m}{2}}(M)}, \quad 0\leq \ell \leq 2 s_3 +n,
}
and hence, that
\eqn{stpropA8}{
\norm{\del{t}^\ell(uv)}_{H^{s_3-\frac{\ell}{2}}(M)} \lesssim
\norm{u}_{E^{s_1,r}}\norm{v}_{E^{s_2,r}}, \quad 0\leq \ell \leq r \leq 2 s_3 +n.
}
\end{proof}

\begin{prop} \label{stcomA}
Suppose $s_1=k_1/2$ and  $s_2=k_2/2$ for $k_1,k_2\in \Zbb_{\geq 0}$, $\ell\in \Zbb_{\geq 1}$,
$s_1+s_2 -\ell/2>0$,  $s_3\in \Rbb$, $s_1\geq s_3$,  $s_2\geq s_3-1/2$ and $s_1+s_2-s_3>n/2$. Then
\eqn{stcomA1}{
\norm{[\del{t}^\ell,u(t)]v(t)}_{H^{s_3-\frac{\ell}{2}}(M)} \lesssim \norm{\del{t}u(t)}_{E^{s_1-\frac{1}{2},\ell-1}}
\norm{v(t)}_{E^{s_2,\ell-1}}
}
for $0\leq t \leq T$ and all $u\in X^{s_1,\ell}_T(M)$ and $v\in X^{s_2,\ell-1}_T(M)$.
\end{prop}
\begin{proof}
By definition,
\eqn{stcomA2}{
[\del{t}^\ell,u]v=\del{t}^\ell(uv)-u\del{t}^\ell v = \sum_{r=0}^{\ell-1}\del{t}^{\ell-r}u\del{t}^r v.
}
Using this together with the multiplication inequality as in the proof of Proposition \ref{stpropA} above, we find
that
\eqn{stcomA3}{
\norm{[\del{t}^\ell,u]v}_{H^{s_3-\frac{\ell}{2}}(\Omega)}
\lesssim \sum_{r=0}^{\ell-1}\norm{\del{t}^{\ell-r}u}_{H^{s_1-\frac{\ell-r}{2}}(\Omega)}
\norm{\del{t}^r v}_{H^{s_2-\frac{r}{2}}(\Omega)}
}
since $s_1\geq s_3$, $s_2\geq s_3-1/2$ and $s_1+s_2-\ell/2>0$ by assumption. The proof follows since
\eqn{stcomA3}{
\sum_{r=0}^{\ell-1}\norm{\del{t}^{\ell-r}u}_{H^{s_1-\frac{\ell-r}{2}}(\Omega)}
\norm{\del{t}^r v}_{H^{s_2-\frac{r}{2}}(\Omega)}  \lesssim \norm{\del{t}u}_{E^{s_1-\frac{1}{2},\ell-1}}
\norm{v}_{E^{s_2,\ell-1}}.
}
\end{proof}

In addition to the product estimate from Proposition \ref{stpropA}, we need estimates for more general
non-linear maps. The first estimate of this type is contained in the next proposition and the proof follows from a straightforward adaptation of the proof of Proposition A.8 from \cite{AnderssonOliynyk:2014} that involves replacing the integral multiplication
inequalities used there with their fractional versions, i.e. Theorem \ref{calcpropB}, in a similar fashion to
the proof of the product estimate above, i.e. Proposition \ref{stpropA}. We omit the proof.
\begin{prop} \label{stpropB}
Suppose $s =k/2$, $s>n/2$, $r\in \Zbb$, $0\leq r \leq 2s$, $f\in C^r(\Rbb)$, and $f(0)=0$.
Then
\eqn{stpropB1}{
\norm{\del{t}^\ell f(u(t))}_{H^{s-\frac{\ell}{2}}(M)}\leq
\norm{f(u(t))}_{E^{s,r}} \leq C(\norm{u(t)}_{E^{s,r}})\norm{u(t)}_{E^{s,r}}
}
for $0\leq t \leq T$, $0\leq \ell \leq r$, and all $u\in X^{s,r}_T(M)$.
\end{prop}

Next, we consider a number of variations on the above basic estimates.
\begin{prop} \label{stpropC}
Suppose $s_1=k_1/2$ and  $s_2=k_2/2$ for $k_1,k_2\in \Zbb_{\geq 0}$, $s_3\in \Rbb$, $s_1,s_2 \geq s_3$,
and $s_1+s_2 - s_3 > n/2$.
Then
\eqn{stpropC1}{
\norm{u(t)v(t)}_{\Ec^{s_3}}  \lesssim \norm{u(t)}_{\Ec^{s_1}} \norm{v(t)}_{\Ec^{s_2}}
}
for $0\leq t \leq T$, and all $u\in \Xc^{s_1}_T(M)$ and $v\in \Xc^{s_2}_T(M)$.
\end{prop}
\begin{proof}
First, we note that
\lalin{stpropC2}{
\norm{uv}_{\Ec^{s_3}} &\lesssim \sum_{\ell=0}^{2s_3-2}\norm{\del{t}^\ell(uv)}_{H^{s_3-\frac{\ell}{2}}(M)}
+ \norm{\del{t}^{2s_3-1}(uv)}_{L^2(M)}\notag \\
&\lesssim \norm{u}_{E^{s_1,2s_3-2}}\norm{v}_{E^{s_2,2s_3-2}} +  \norm{\del{t}^{2s_3-1}(uv)}_{L^2(M)}
&& \text{(by Proposition \ref{stpropA})} \notag\\
&\lesssim \norm{u}_{\Ec^{s_1}}\norm{v}_{\Ec^{s_2}} +  \norm{\del{t}^{2s_3-1}(uv)}_{L^2(M)}. \label{stpropC2.1}
}
Next, we find from
\eqn{stpropC3}{
\del{t}^{2s_3-1} = \sum_{\ell=0}^{2s_3-1}\binom{2s_3-1}{\ell}\del{t}^{2s_3-1-\ell}u\del{t}^\ell v
}
that
\lalin{stpropC4}{
\norm{\del{t}^{2s_3-1}(uv)}_{L^2(M)}
&\lesssim \sum_{\ell=1}^{2s_3-2}\norm{\del{t}^{2s_3-1-\ell}u\del{t}^\ell v}_{L^2}
+\norm{v\del{t}^{2s_3-1}u}_{L^2(M)} + \norm{v\del{t}^{2s_3-1}u}_{L^2(M)}\notag \\
&\lesssim \norm{u}_{E^{s_1,2s_3-2}}\norm{v}_{E^{s_2,2s_3-2}} + \norm{v\del{t}^{2s_3-1}u}_{L^2(M)} + \norm{v\del{t}^{2s_3-1}u}_{L^2(M)} \notag\\
&\lesssim \norm{u}_{\Ec^{s_1}}\norm{v}_{\Ec^{s_2}} + \norm{u\del{t}^{2s_3-1}v}_{L^2(M)} + \norm{v\del{t}^{2s_3-1}u}_{L^2(M)}
\label{stpropC4.1}
}
where in deriving the second inequality we used Proposition \ref{stpropA}.

Letting
\eqn{stpropC5.1}{
\Ic = \norm{u\del{t}^{2s_3-1}v}_{L^2(M)} + \norm{v\del{t}^{2s_3-1}u}_{L^2(M)},
}
three cases follow:
\smallskip

\noindent\textit{Case 1: $s_1=s_2=s_3$}

\smallskip
\lalin{stpropC6}{
\Ic &\lesssim \norm{u}_{L^\infty(M)}\norm{\del{t}^{2s_2-1}v}_{L^2(M)} + \norm{v}_{L^\infty(M)}
\norm{\del{t}^{2s_1-1}u}_{L^2(M)}  \notag\\
&\lesssim \norm{u}_{H^{s_1}(M)}\norm{\del{t}^{2s_2-1}v}_{L^2(M)} + \norm{v}_{H^{s_2}(M)}
\norm{\del{t}^{2s_1-1}u}_{L^2(M)} && \text{(by Theorem \ref{FSobolev})} \notag\\
&\lesssim \norm{u}_{\Ec^{s_1}}\norm{v}_{\Ec^{s_2}}. \label{stpropC6.1}
}
\smallskip

\noindent\textit{Case 2: $s_1=s_3$ and $s_2>s_3$, or $s_2=s_3$ and $s_1>s_3$}

\smallskip

Here, we assume that $s_2=s_3$ and $s_1>s_3$ so that $s_1>n/2$. The other case,
$s_1=s_3$ and $s_2>s_3$, follows from switching $u$ and $v$.
\lalin{stpropC7}{
\Ic &\lesssim \norm{u}_{L^\infty(M)}\norm{\del{t}^{2s_2-1}v}_{L^2(M)} +
\norm{v\del{t}^{2s_3-1}u}_{L^2(M)} \notag \\
& \lesssim \norm{u}_{H^{s_1}(M)}\norm{\del{t}^{2s_2-1}v}_{L^2(M)} +
\norm{v}_{H^{s_2}(M)}\norm{\del{t}^{2s_3-1}u}_{H^{s_1-\frac{2s_3-1}{2}}(M)} \notag \\
& \lesssim  \norm{u}_{\Ec^{s_1}}\norm{v}_{\Ec^{s_2}} \hspace{3.0cm} \text{(since $2s_3-1\leq 2s_1-2$)},\label{stpropC7.1}
}
where in deriving this second inequality we used Theorems \ref{FSobolev} and \ref{calcpropB}.
\smallskip

\noindent\textit{Case 3: $s_1>s_3$ and $s_2 > s_3$}

\smallskip

\lalin{stpropC8}{
\Ic &\lesssim
\norm{u}_{H^{s_1}(M)}\norm{\del{t}^{2s_3-1}v}_{H^{s_2-\frac{2s_3-1}{2}}(M)}+
\norm{v}_{H^{s_2}(M)}\norm{\del{t}^{2s_3-1}u}_{H^{s_1-\frac{2s_3-1}{2}}(M)} \notag \\
&\lesssim \norm{u}_{\Ec^{s_1}}\norm{v}_{\Ec^{s_2}} \hspace{3.0cm} \text{(since $2s_3-1\leq 2s_1-2$ and
$2s_3-1\leq 2s_2-2$)}
\label{stpropC8.1}
}
where in deriving the second inequality we used Theorem \ref{calcpropB}.

Combining the inequalities \eqref{stpropC2.1}-\eqref{stpropC8.1} completes the proof.
\end{proof}

\begin{prop} \label{stpropD}
Suppose $s =k/2$, $s>n/2$,  $f\in C^{2s-1}(\Rbb)$ and $f(0)=0$.
Then
\eqn{stpropB1}{
\norm{f(u(t))}_{\Ec^{s}}\leq  C(\norm{u(t)}_{\Ec^{s}})\norm{u(t)}_{\Ec^{s}}
}
for $0\leq t \leq T$ and all $u\in \Xc^{s}_T(M)$.
\end{prop}
\begin{proof}
Again, the proof follow from an adaptation of the proof of Proposition A.8 from \cite{AnderssonOliynyk:2014} that involves replacing the integral multiplication
inequalities used there with their fractional versions, i.e. Theorem \ref{calcpropB}, in a similar fashion to
the proof of the product estimate from Proposition \ref{stpropC}.
\end{proof}

\begin{prop} \label{stpropE}
Suppose $s =k/2$, $s>n/2$,  $f\in C^{2s}(\Rbb\times \Rbb^{n+1})$ and $f(0,0)=0$. Then
\eqn{stpropE1}{
\norm{f(u(t),\del{}u(t))}_{E^{s}} \leq C\Bigl(\norm{u(t)}_{\Ec^{s+1}}\Bigr)\norm{u(t)}_{\Ec^{s+1}}
\quad \bigl(\del{}u=(\del{t}u,\del{i}u)\bigr)
}
for $t\in [0,T]$ and $u \in  \Xc^{s+1}_T(M)$.
\end{prop}
\begin{proof}
First, we note the following simple inequalities, which follow directly from the definition
of the norms $\norm{\cdot}_{E^s}$ and $\norm{\cdot}_{\Ec^{s+1}}$.
\eqn{stpropE1}{
\norm{u}_{E^s}^2 \leq \norm{u}_{\Ec^{s+1}}^2,
\quad \norm{\del{t}u}_{E^s}^2 \leq \norm{u}_{\Ec^{s+1}}^2
\AND
\norm{Du}^2_{E^s} \leq \norm{u}_{\Ec^{s+1}}^2.
}
Using these in conjunction with Proposition \ref{stpropB}, we obtain
\eqn{stpropE2}{
\norm{f(u,\del{}u)}_{E^s} \leq C\bigl(\norm{u}_{E^{s}}+\norm{\del{}u}_{E^s}\bigr)
\bigl(\norm{u}_{E^{s}}+\norm{\del{}u}_{E^s}\bigr) \leq C(\norm{u}_{\Ec^{s+1}})\norm{u}_{\Ec^{s+1}},
}
which completes the proof.
\end{proof}

\begin{prop} \label{stpropEa}
Suppose $s =k/2$, $s>n/2$,  $f\in C^{2s}(\Rbb\times \Rbb^{n+1}\times \Rbb\times \Rbb)$ and $f(0,0,0,0)=0$. Then
\alin{stpropEa1}{
&\norm{f(u(t),\del{}u(t),v(t),\del{t}v(t))}_{\Ec^{s+\frac{1}{2}}} \leq
C\Bigl(\norm{u(t)}_{E^{s+\frac{3}{2},2s-1}}+\norm{u(t)}_{\Ec^{s+1}}+\norm{v(t)}_{\Ec^{s+1}}\Bigr) \\
&\hspace{6.0cm}\times
\Bigl(\norm{u(t)}_{E^{s+\frac{3}{2},2s-1}}+\norm{u(t)}_{\Ec^{s+1}}+\norm{v(t)}_{\Ec^{s+1}}\Bigr)
}
for all $t\in [0,T]$, $u\in X_T^{s+1}(M)\cap X_T^{s+\frac{3}{2},2s-1}(M)$, and $v\in \Xc_T^{s+1}(M)$.
\end{prop}
\begin{proof}
We begin by observing the following simple inequalities that follow directly from the defintion
of the norms:
\alin{stpropEa1}{
\norm{v}_{\Ec^{s+\frac{1}{2}}}+\norm{\del{t}v}_{\Ec^{s+\frac{1}{2}}} &\lesssim \norm{v}_{\Ec^{s+1}},
\intertext{and}
\norm{u}_{\Ec^{s+\frac{1}{2}}}+\norm{\del{t}u}_{\Ec^{s+\frac{1}{2}}} +
\norm{Du}_{\Ec^{s+\frac{1}{2}}} & \lesssim \norm{u}_{E^{s+\frac{3}{2},2s-1}}+ \norm{u}_{\Ec^{s+1}}.
}
The desired inequality then follows from these estimates and an application of Proposition \ref{stpropD}:
\alin{stpropEa2}{
\norm{f(u,\del{}u,v,\del{t}v)}_{\Ec^{s+\frac{1}{2}}}
&\leq C\Bigl(\norm{u}_{\Ec^{s+\frac{1}{2}}}+ \norm{\del{}u}_{\Ec^{s+\frac{1}{2}}}+
\norm{v}_{\Ec^{s+\frac{1}{2}}}+ \norm{\del{t}v}_{\Ec^{s+\frac{1}{2}}} \Bigr) \\
&\hspace{2.0cm} \times \Bigl(\norm{u}_{\Ec^{s+\frac{1}{2}}}+ \norm{\del{}u}_{\Ec^{s+\frac{1}{2}}}+
\norm{v}_{\Ec^{s+\frac{1}{2}}}+ \norm{\del{t}v}_{\Ec^{s+\frac{1}{2}}}\Bigr), \\
&\leq C\Bigl(\norm{u}_{E^{s+\frac{3}{2},2s-1}}+ \norm{u}_{\Ec^{s+1}}+\norm{v}_{\Ec^{s+1}}\Bigr)
\Bigl(\norm{u}_{E^{s+\frac{3}{2},2s-1}}+ \norm{u}_{\Ec^{s+1}}+\norm{v}_{\Ec^{s+1}}\Bigr).
}
\end{proof}

We conclude with the following integral estimates for ODEs.

\begin{prop} \label{stpropF}
Suppose that $s=k/2$, $k\in \Zbb_{\geq 0}$, $0\leq r \leq 2s$, $u,\del{t}u\in X_T^{s,r}(M)$, $f\in X_T^{s,r}(M)$, $v,\del{t}v \in X_T^{s}(M)$, $g\in X^{s}_T(M)$, and
\alin{stpropF1}{
\del{t}u = f \AND \del{t}v = g \quad \text{in $[0,T]\times M$.}
}
Then
\alin{stpropF2}{
\norm{u(t)}_{E^{s,r}}^2 &\leq \norm{u(0)}_{E^{s,r}}^2 + \int_0^t
\norm{u(\tau)}_{E^{s,r}}^2 + \norm{f(\tau)}_{E^{s,r}}^2 \, d \tau
\intertext{and}
\norm{v(t)}_{\Ec^{s}}^2 &\leq \norm{v(0)}_{\Ec^{s}}^2 + \int_0^t
\norm{v(\tau)}_{\Ec^{s}}^2 + \norm{g(\tau)}_{\Ec^{s}}^2 \, d \tau
}
for $t\in [0,T]$.
\end{prop}
\begin{proof}
We only prove the first inequality with the second following from similar arguments. Fixing
$\ell \in \{0,1,\ldots,r\}$, we compute
\eqn{stpropF3}{
\del{t}\norm{u(t)}_{H^{s-\frac{\ell}{2}}(M)}^2 = 2\ip{u}{\del{t}u}_{H^{s-\frac{\ell}{2}}(M)}
\leq \norm{u(t)}_{H^{s-\frac{\ell}{2}}(M)}^2 + \norm{\del{t}u(t)}_{H^{s-\frac{\ell}{2}}(M)}^2
}
where in deriving the last inequality we used Cauchy-Schwartz and Young's inequalities. But, $\del{t}u=f$
by assumption, and hence we have that
\eqn{stpropF4}{
\del{t}\norm{u(t)}_{H^{s-\frac{\ell}{2}}(M)}^2 \leq
\norm{u(t)}_{H^{s-\frac{\ell}{2}}(M)}^2 + \norm{f(t)}_{H^{s-\frac{\ell}{2}}(M)}^2.
}
Integrating this expression in time and summing the result over $\ell$ from $0$ to $r$ then gives
\eqn{stpropF5}{
\norm{u(t)}_{E^{s,r}}^2 \leq \norm{u(0)}_{E^{s,r}}^2
\leq \int_0^t \norm{u(\tau)}_{E^{s,r}}^2 + \norm{f(\tau)}_{E^{s,r}}^2\, d\tau,
}
which proves the first inequality.
\end{proof} 
\sect{elliptic}{Elliptic systems}

In this appendix, we recall some well known existence and regularity results for elliptic systems
of the form
\lalin{ellipA}{
\del{i}(b^{ij}\del{j}v+ \ep d^i v+L^i) + \ep a^i\del{i} v - \lambda c v &= F\hspace{1.35cm}\text{in $\Omega$,} \label{ellipA.1}\\
\nu_i(b^{ij}\del{j}v+ \ep d^i v + L^i) &= \ep h v + G \hspace{0.4cm}\text{in $\del{}\Omega$,} \label{ellipA.2}
}
where
\begin{enumerate}[(i)]
\item as above, $\Omega$ is open, bounded in $\Rbb^n$, $n\geq 2$, with smooth boundary,
\item $\nu_i$ is the outward pointing unit co-normal to $\del{}\Omega$,
\item $L^i \in L^2(\Omega,\Rbb^N)$, $F\in L^{2}(\Omega,\Rbb^N)$ and $G\in H^{-1/2}(\del{}\Omega,\Rbb^N)$,
\item $a^i,d^i\in L^n(\Omega,\Mbb{N})$ and $h\in L^\infty(\Omega,\Mbb{N})\cap W^{1,n}(\Omega,\Mbb{N})$,
\item $c\in L^n(\Omega,\Mbb{N})$ and satisfies
\eqn{cellip1}{
c \geq \sigma > 0 \quad \text{in $\Omega$}
}
for some positive constant $\sigma$,
\item $b^{ij}\in L^\infty(\Omega,\Mbb{N})$ and satisfies
\leqn{bellip}{
(b^{ij})^{\tr} = b^{ji},
}
\item and there exists a $\kappa_1 >0$ and $\mu \geq 0$ such that
\leqn{coercellip}{
\ip{\del{i}v}{b^{ij}\del{j}v}_{L^2(\Omega)} \geq \kappa_1 \norm{v}^2_{H^1(\Omega)} -\mu \norm{v}^2_{L^2(\Omega)}
}
for all $v\in H^1(\Omega)$.
\end{enumerate}

\begin{Def} \label{eweakdef}
Under the assumptions (i)-(vii) above, $v\in H^1(\Omega,\Rbb^N)$ is called a weak solution
of \eqref{ellipA.1}-\eqref{ellipA.2} if
\lalin{eweakdef1}{
\ip{b^{ij}\del{i}v}{\del{j}\phi}_{\Omega} &+\ep\ip{d^j v}{\del{j}\phi}_{\Omega} -\ep\ip{a^i\del{i}v}{\phi}_{\Omega}
+\lambda\ip{c v}{\phi}_{\Omega} \notag  \\
&\hspace{0.5cm} -\ep\ip{(hv)|_{\del{}\Omega}}{\phi|_{\del{}\Omega}}_{\del{}\Omega}
= -\ip{F}{\phi}_{\Omega} - \ip{L^j}{\del{j}\phi}_{\Omega}
+ \ip{G}{\phi|_{\del{}\Omega}}_{\del{}\Omega} \label{eweakdef1.1}
}
for all $\phi \in H^1(\Omega,\Rbb^N)$.
\end{Def}

\begin{rem} \label{weakremAa}
That the above definition makes sense follows from repeated use of
H\"{o}lder's inequality, the Trace theorem,
the Sobolev inequalities, Ehrling's lemma and the duality relation $(H^s(\del{}\Omega))^* \cong H^{-s}(\del{}\Omega)$.
To see this, we observe from Theorems \ref{Holder} and \ref{ISobolev}
that
\leqn{weakremA1}{
\norm{uv}_{L^2(\Omega)} \leq \norm{u}_{L^n(\Omega)} \norm{v}_{L^{\frac{2n}{n-2}}(\Omega)}
\lesssim \norm{u}_{L^n(\Omega)}\norm{v}_{H^1(\Omega)},
}
which in turn, implies via the Cauchy-Schwartz inequality that
\leqn{weakremA2}{
\ip{uv}{w}_{\Omega} \lesssim  \norm{u}_{L^n(\Omega)}\norm{v}_{H^1(\Omega)}\norm{w}_{L^2(\Omega)}.
}
Using Theorems \ref{Holder}, \ref{ISobolev} and \ref{Iprod}, we observe also that
\leqn{weakremA3}{
\norm{uv}_{H^1(\Omega)} \lesssim \norm{u}_{L^\infty(\Omega)}\norm{v}_{H^1(\Omega)}
+ \norm{u}_{W^{1,n}(\Omega)}\norm{v}_{L^{\frac{2n}{n-2}}(\Omega)} \lesssim
\norm{u}_{L^\infty(\Omega)\cap W^{1,n}(\Omega)} \norm{v}_{H^1(\Omega)}.
}
This inequality, together with Theorem \ref{trace}, the Cauchy-Schwartz inequality, and
Lemma \ref{Ehrling}, implies that
\lalin{weakremA3a}{
\ip{(uv)|_{\del{}\Omega}}{w|_{\del{}\Omega}}_{\del{}\Omega}
&\lesssim \norm{uv}_{H^1(\Omega)}\norm{w}_{H^{\frac{3}{4}}(\Omega)} \notag \\
&\lesssim  \norm{u}_{L^\infty(\Omega)\cap W^{1,n}(\Omega)} \norm{v}_{H^1(\Omega)}
\bigl(\gamma \norm{w}_{H^{1}(\Omega)} + %\frac{1}{\gamma^3}
C(\gamma)\norm{w}_{L^2(\Omega)} \bigr) \label{weakremA3.1}
}
for any $\gamma > 0$. Finally,
\leqn{weakremA5}{
\ip{u}{v|_{\del{}\Omega}}_{\del{}\Omega} \lesssim \norm{u}_{H^{-\frac{1}{2}}(\del{}\Omega)}\norm{v}_{H^1(\Omega)}
}
follows by the duality relation $(H^s(\del{}\Omega))^* \cong H^{-s}(\del{}\Omega)$ and Theorem \ref{trace}.
\end{rem}

The following existence result is a slight modification of Theorem B.2 from \cite{Koch:1990}, and
is proved using similar arguments.
\begin{thm} \label{ellipthmA}
Suppose assumptions (i)-(vii) above are satisfied. Then
\begin{enumerate}[(i)]
\item there exists an $\delta^*=\delta^*\bigl(\kappa_1,\norm{a}_{L^n(\Omega)},
    \norm{d}_{L^n(\Omega)},\norm{h}_{L^\infty(\Omega)\cap W^{1,n}(\Omega)}\bigr)\geq 1$ such that
every weak solution $v$ of \eqref{ellipA.1}-\eqref{ellipA.2} with $\ep \in [0,\frac{1}{\delta^*}]$ satisfies the estimate
\eqn{ellipthmA1}{
\norm{v}_{H^1(\Omega)} \leq C\Bigl(\norm{v}_{L^2(\Omega)}
+ \norm{L}_{L^2(\Omega)}
+ \norm{F}_{L^{2}(\Omega)} + \norm{G}_{H^{-\frac{1}{2}}(\del{}\Omega)}\Bigr),
}
where
\eqn{ellipthmA2}{
C = C\bigl(\kappa_1,\mu,\lambda,\norm{a}_{L^n(\Omega)},\norm{d}_{L^n(\Omega)},\norm{h}_{L^\infty(\Omega)\cap W^{1,n}(\Omega)},\norm{c}_{L^n(\Omega)}\bigr),
}
and
\item there exists a $\lambda^*= \lambda^*(\sigma,\mu) > 0$ such that for $\lambda \geq \lambda^*$
and $\ep \in [0,\ep^*]$
there exists a unique weak solution of  \eqref{ellipA.1}-\eqref{ellipA.2}.
\end{enumerate}
\end{thm}
\begin{proof}
\textit{(i):} Given a weak solution $v$ of \eqref{ellipA.1}-\eqref{ellipA.2}, we
see, after setting $v=\phi$ in \eqref{eweakdef1.1}, and using \eqref{coercellip} and \eqref{weakremA1}-\eqref{weakremA5}, that
\alin{ellipthmA3}{
\kappa_1 &\norm{v}^2_{H^1(\Omega)} \leq K\bigl( \ep \norm{a}_{L^n(\Omega)}+
\ep \norm{d}_{L^n(\Omega)}+ \ep\norm{h}_{L^\infty(\Omega)\cap W^{1,n}(\Omega)}\bigr)
\norm{v}^2_{H^1(\Omega)} + K\biggl(\biggl[|\lambda|\norm{c}_{L^n(\Omega)}
 \\
&+\ep\norm{h}_{L^\infty(\Omega)\cap W^{1,n}(\Omega)}+\mu\biggr]\norm{v}_{L^2(\Omega)}+ \norm{F}_{L^2(\Omega)}+\norm{L}_{L^2(\Omega)}
+ \norm{G}_{H^{-\frac{1}{2}}(\del{}\Omega)}\biggr)\norm{v}_{H^1(\Omega)}
}
for some constant $K>0$ independent of $\ep > 0$. Setting
\eqn{ep*def}{
\delta^*=\max\left\{1,\frac{2K\Bigl(\norm{a}_{L^n(\Omega)}+\norm{d}_{L^n(\Omega)}+\norm{h}_{L^\infty(\Omega)\cap W^{1,n}(\Omega)}\Bigr)}{\kappa_1} \right\},
}
and choosing $\ep \in [0,\frac{1}{\delta^*}]$, it follows immediately from the above estimate that
\eqn{ellipthmA4}{
\norm{v}_{H^1(\Omega)} \leq C\bigl( \norm{v}_{L^2(\Omega)}+ \norm{F}_{L^2(\Omega)}+\norm{L}_{L^2(\Omega)}
+ \norm{G}_{H^{-\frac{1}{2}}(\del{}\Omega)}\bigr),
}
where
\eqn{ellipthmA5}{
C = C\bigl(\kappa_1,\mu,\lambda,\norm{a}_{L^n(\Omega)},\norm{d}_{L^n(\Omega)},\norm{h}_{L^\infty(\Omega)\cap W^{1,n}(\Omega)},\norm{c}_{L^n(\Omega)}\bigr).
}

\noindent\textit{(ii:)} Setting
\eqn{ellipthmA6}{
B(v,\phi) = \ip{b^{ij}\del{i}v}{\del{j}\phi}_{\Omega} +\ip{d^j v}{\del{j}\phi}_{\Omega} -\ip{a^i\del{i}v}{\phi}_{\Omega}
+\lambda\ip{c v}{\phi}_{\Omega} -\ep\ip{(fv)|_{\del{}\Omega}}{\phi|_{\del{}\Omega}}_{\del{}\Omega},
}
and
\eqn{ellipthmA7}{
\Lambda(\phi) = -\ip{F}{\phi}_{\Omega} - \ip{L^j}{\del{j}\phi}_{\Omega}
+ \ip{G}{\phi|_{\del{}\Omega}}_{\del{}\Omega},
}
it follows from the estimates \eqref{coercellip} and \eqref{weakremA1}-\eqref{weakremA5} that $B$ and $\Lambda$
define bounded forms on $H^1(\Omega)$, and, moreover, that $B$ satisfies the estimate
\alin{ellipthmA8}{
B(\phi,\phi) \geq& \bigl(\kappa_1-K\bigl[\ep\norm{d}_{L^\infty(\Omega)}+\ep\norm{a}_{L^\infty(\Omega)}+\ep\norm{h}_{L^\infty(\Omega)\cap W^{1,n}(\Omega)}\bigr]\bigr)\norm{\phi}^2_{H^1(\Omega)}\\
& +(\lambda\sigma-\mu)\norm{\phi}^2_{L^2(\Omega)}
 -\ep\norm{h}_{L^\infty(\Omega)\cap W^{1,n}(\Omega)}\norm{\phi}_{L^2(\Omega)}\norm{\phi}_{H^1(\Omega)}.
}
From this estimate and Young's inequality, i.e. $\alpha\beta$ $\leq$ $\frac{\alpha^2}{2\gamma}$ $+$ $\frac{\gamma \beta^2}{2}$ ($\alpha,\beta\geq0$, $\gamma>0$), it is clear that there exists a constant
$\lambda^* = \lambda^*(\sigma,\mu) >0$
such that
\eqn{ellipthmA10}{
B(\phi,\phi) \geq \frac{\kappa_1}{2}\norm{\phi}_{H^1(\Omega)}^2
}
whenever $\lambda \geq \lambda^*$ and $0\leq \ep \leq \frac{1}{\delta^*}$. Fixing $\lambda \in [\lambda^*,\infty)$ and $\ep \in [0,\frac{1}{\delta^*}]$,
we can apply the Lax-Milgram theorem, see
Theorem 1 from Section 6.2.1 of \cite{Evans:2010}, to
obtain the existence of a unique $v\in H^1(\Omega)$ satisfying $B(v,\phi)=\Lambda(\phi)$ for all
$\phi\in H^1(\Omega)$. By definition of $B$ and $\Lambda$, it is clear that $v$ is the desired weak solution
of \eqref{ellipA.1}-\eqref{ellipA.2}.
\end{proof}

In addition to the above existence result, we will also require the following version of
elliptic regularity, where we define
\eqn{r*def}{
r^* = \begin{cases} r-1 & \text{if $r> 1$}\\
                   0 & \text{otherwise}
\end{cases}
}
\begin{thm} \label{ellipthmB}
Suppose $r,s \in \Rbb$, $s>n/2$, $0\leq r \leq s$, $b^{ij} \in H^s(\Omega,\Mbb{N})$, $L^i\in H^{r}(\Omega,\Rbb^N)$,
$F\in H^{r^*}(\Omega,\Rbb^N)$, $G\in H^{r-\frac{1}{2}}(\del{}\Omega,\Rbb^N)$,
the $b^{ij}$ satisfy \eqref{bellip} and \eqref{coercellip}, $a^i=c=d=h=0$, and
$v$ is a weak solution of  \eqref{ellipA.1}-\eqref{ellipA.2}. Then
$v\in H^{r+1}(\Omega,\Rbb^N)$ and satisfies
\eqn{ellipthmB.1}{
\norm{v}_{H^{r+1}(\Omega)} \leq C\Bigl(\norm{v}_{H^{r}(\Omega)}+\norm{F}_{H^{r^*}(\Omega)}
+ \norm{L}_{H^{r}(\Omega)} + \norm{G}_{H^{r-1/2}(\del{}\Omega)} \Bigr),
}
where $C = C(\kappa_1,\mu,\norm{b}_{H^s(\Omega)})$.
\end{thm}

\sect{LA}{Determinant formulas}

\begin{lem} \label{matrixlem}
Suppose that $X,Y\in \Cbb^2$, $N\in \Cbb$, and $L\in \text{\em Gl}(2,\Cbb)$ satisfies $L^{\tr}=L$. Then
\eqn{matrixlem1}{
\det\begin{pmatrix} L & X+Y \\ -X^{\tr}+Y^{\tr} & N\end{pmatrix} = \det\begin{pmatrix} L & Y \\ Y^{\tr} & N \end{pmatrix}
+ \det(L)X^{\tr} L^{-1} X.
}
\end{lem}
\begin{proof}
Direct computation.
\end{proof}

\begin{lem} \label{matrixlemA}
Suppose that $X\in\Cbb^2$, $N\in \Cbb$, and $L\in \text{\em Gl}(2,\Cbb)$ satisfies $L^{\tr}=L$. Then
\eqn{matrixlemA1}{
\det\begin{pmatrix} L & X \\ X^{\tr} & N+X^{\tr}L^{-1}X\end{pmatrix} = N \det(L).
}
\end{lem}
\begin{proof}
Noting that
\alin{matrixlemA2}{
\det\begin{pmatrix} L & X \\ X^{\tr} & N+X^{\tr}L^{-1}X\end{pmatrix} &= 
\det\left(\begin{pmatrix} \id & 0 \\ 0 & -1 \end{pmatrix} \begin{pmatrix} L & X \\ -X^{\tr} & -N-X^{\tr}L^{-1}X\end{pmatrix}\right) \\
&= -\det\begin{pmatrix} L & X \\ -X^{\tr} & -N-X^{\tr}L^{-1}X\end{pmatrix},
}
the proof follows since
\eqn{matrixlemA3}{
\det\begin{pmatrix} L & X \\ -X^{\tr} & -N-X^{\tr}L^{-1}X\end{pmatrix} = -N\det(L)
}
by Lemma \ref{matrixlem}.
\end{proof}

\bibliographystyle{amsplain}
\bibliography{refs}

%\end{spacing}

\end{document}